\documentclass[12pt]{amsart}
\usepackage{amsbsy}
\usepackage{graphicx,epsfig,subfigure,psfrag}
\begin{document}

\newtheorem{theorem}{Theorem}
\newtheorem{proposition}{Proposition}
\newtheorem{lemma}{Lemma}
\newtheorem{corollary}{Corollary}
\newtheorem{definition}{Definition}
\newtheorem{remark}{Remark}
\newcommand{\tex}{\textstyle}
\numberwithin{equation}{section} \numberwithin{theorem}{section}
\numberwithin{proposition}{section} \numberwithin{lemma}{section}
\numberwithin{corollary}{section}
\numberwithin{definition}{section} \numberwithin{remark}{section}
\newcommand{\ren}{\mathbb{R}^N}
\newcommand{\re}{\mathbb{R}}
\newcommand{\n}{\nabla}
\newcommand{\iy}{\infty}
\newcommand{\pa}{\partial}
\newcommand{\fp}{\noindent}
\newcommand{\ms}{\medskip\vskip-.1cm}
\newcommand{\mpb}{\medskip}
\newcommand{\AAA}{{\bf A}}
\newcommand{\BB}{{\bf B}}
\newcommand{\CC}{{\bf C}}
\newcommand{\DD}{{\bf D}}
\newcommand{\EE}{{\bf E}}
\newcommand{\FF}{{\bf F}}
\newcommand{\GG}{{\bf G}}
\newcommand{\oo}{{\mathbf \omega}}
\newcommand{\Am}{{\bf A}_{2m}}
\renewcommand{\a}{\alpha}
\renewcommand{\b}{\beta}
\newcommand{\g}{\gamma}
\newcommand{\G}{\Gamma}
\renewcommand{\d}{\delta}
\newcommand{\D}{\Delta}
\newcommand{\e}{\varepsilon}
\newcommand{\var}{\varphi}
\renewcommand{\l}{\lambda}
\renewcommand{\o}{\omega}
\renewcommand{\O}{\Omega}
\newcommand{\s}{\sigma}
\renewcommand{\t}{\tau}
\renewcommand{\th}{\theta}
\newcommand{\z}{\zeta}
\newcommand{\wx}{\widetilde x}
\newcommand{\wt}{\widetilde t}
\newcommand{\noi}{\noindent}
\newcommand{\uu}{{\bf u}}
\newcommand{\xx}{{\bf x}}
\newcommand{\yy}{{\bf y}}
\newcommand{\zz}{{\bf z}}
\newcommand{\aaa}{{\bf a}}
\newcommand{\cc}{{\bf c}}
\newcommand{\jj}{{\bf j}}
\newcommand{\ff}{{\bf f}}
\newcommand{\ggg}{{\bf g}}
\newcommand{\UU}{{\bf U}}
\newcommand{\YY}{{\bf Y}}
\newcommand{\WW}{{\bf W}}
\newcommand{\HH}{{\bf H}}
\newcommand{\GGG}{{\bf G}}
\newcommand{\VV}{{\bf V}}
\newcommand{\ww}{{\bf w}}
\newcommand{\vv}{{\bf v}}
\newcommand{\hh}{{\bf h}}
\newcommand{\di}{{\rm div}\,}
\newcommand{\inA}{\quad \mbox{in} \quad \ren \times \re_+}
\newcommand{\inB}{\quad \mbox{in} \quad}
\newcommand{\inC}{\quad \mbox{in} \quad \re \times \re_+}
\newcommand{\inD}{\quad \mbox{in} \quad \re}
\newcommand{\forA}{\quad \mbox{for} \quad}
\newcommand{\whereA}{,\quad \mbox{where} \quad}
\newcommand{\asA}{\quad \mbox{as} \quad}
\newcommand{\andA}{\quad \mbox{and} \quad}
\newcommand{\withA}{,\quad \mbox{with} \quad}
\newcommand{\orA}{,\quad \mbox{or} \quad}
\newcommand{\ef}{\eqref}
\newcommand{\ssk}{\smallskip}
\newcommand{\LongA}{\quad \Longrightarrow \quad}
\def\com#1{\fbox{\parbox{6in}{\texttt{#1}}}}
\def\N{{\mathbb N}}
\def\A{{\cal A}}
\newcommand{\de}{\,d}
\newcommand{\eps}{\varepsilon}
\newcommand{\be}{\begin{equation}}
\newcommand{\ee}{\end{equation}}
\newcommand{\spt}{{\mbox spt}}
\newcommand{\ind}{{\mbox ind}}
\newcommand{\supp}{{\mbox supp}}
\newcommand{\dip}{\displaystyle}
\newcommand{\prt}{\partial}
\renewcommand{\theequation}{\thesection.\arabic{equation}}
\renewcommand{\baselinestretch}{1.1}

\title
{\bf On blow-up ``twistors"\\
for the Navier--Stokes equations in $\re^3$:\\ a view from
reaction-diffusion theory}



\author{
V.A.~Galaktionov}

\address{Department of Mathematical Sciences, University of Bath,
 Bath BA2 7AY, UK}
\email{vag@maths.bath.ac.uk}



  \keywords{Navier--Stokes equations in $\re^3$,
 solutions on linear subspaces, blow-up angular swirl mechanism,
 blow-up periodic patterns, solenoidal Hermite polynomials,
 multiple zero structures, Slezkin--Landau
 singular steady solution,   blow-up tornado.
 {\bf Submitted to:} Adv. Differ. Equat.}
 \date{\today}




 \begin{abstract}

Formation of blow-up singularities for the Navier--Stokes
equations (NSEs)
 $$
 \uu_t +(\uu \cdot \n)\uu=- \n p +\D  \uu, \quad \di \uu=0 \inB \re^3 \times \re_+,
 $$
 with bounded data $\uu_0$
 is discussed.
 Using natural links with
 blow-up theory for nonlinear reaction-diffusion PDEs, some possibilities
 to construct
  special self-similar  and other related solutions
  that are characterized by
 {\em blow-up  swirl} with the angular speed near the blow-up
 time (this represents simplest $\o$-limits of rescaled orbits as {\em periodic} ones)
  $$
   \tex{
  \var(t) \sim -\sigma  \,\ln(T-t) \LongA
  \dot \var(t) \sim \frac {\sigma }{T-t}
  \to  \iy \asA t \to T^-
   \quad (\sigma  \not =0).
  }
   $$
  This is done in cylindrical polar coordinates $\{r,\var,z\}$ in $\re^3$,
using the restriction of the NSEs to the linear subspace $W_2={\rm
 Span}\{1,z\}$. Similarly,
 blow-up twistors with axis precessions in the spherical geometry $\{r,\th,\var\}$
  are  introduced.

It is shown that other blow-up patterns (a ``screwing in tornado")
may correspond to a slow ``centre-stable manifold-like
  drift"
 about  Slezkin--Landau singular or other  equilibria of the
 NSEs.
Some  approaches to blow-up singularities can be  applied to 3D
Euler's equations and to well-posed Burnett equations in 7D (i.e.,
the NSEs with $\D \mapsto -\D^2$). Though  most of blow-up
 scenarios were not  justified even at a qualitative level,
 the author hopes that the proposed approaches to  families of blow-up and other
 patterns, including those with blow-up swirl,
 will give some extra insight into
 the micro-scale ``turbulent" structure  of the NSEs.

The discussion of possible types of blow-up patterns for the NSEs
is going
 in conjunction with some other classic nonlinear PDEs of mathematical physics.

\end{abstract}

\maketitle


\section{Introduction: Navier--Stokes equations and blow-up}
 \label{SInt}


\subsection
{Two faces of the open problem via blow-up formulation: first
discussion around self-similar blow-up}

 It is well-understood  (see e.g., \cite[Ch.~5]{Maj02}) that
  the fundamental open problem of fluid
 mechanics\footnote{The Millennium Prize Problem for the Clay Institute; see Fefferman
 \cite{Feff00}.}
  and  PDE theory on {\em global existence or nonexistence of
 bounded smooth solutions} of the
 {\em Navier--Stokes}\footnote{Claude Louis Marie Henri Navier, 1785-1836, and George
Gabriel Stokes, 1819-1903.} {\em equations} (the NSEs),
  \be
  \label{NS1}
  \left\{
   \begin{matrix}
 \uu_t +(\uu \cdot \n)\uu=- \n p + \D \uu, \qquad\qquad\qquad
 \ssk \\
 \di \uu=0 \inB \re^3 \times
 \re_+ \quad(\uu=(u,v,w)),
  \end{matrix}
  \right.
  \ee
 with {\em arbitrary bounded divergence-free $L^2$-data} $\uu_0$, is
 two-fold:

\ssk

 From a standard evolution  PDE and blow-up point of view, there exist  two
 possibility: the
 {\bf positive answer},
  i.e., {\em existence of a global
 bounded solution}
 (then necessarily it is unique by the existing
 local semigroup of bounded smooth solutions) is equivalent to
  \be
  \label{Pos}
  \fbox{$
  \mbox{\bf $\exists$ global bounded solution}
  = \mbox{\bf $\not \exists$  finite-time blow-up at any $T>0$}.
  $}
  \ee

  \ssk

  On the other hand, the {\bf negative
answer}, i.e., {\em nonexistence in general of a global bounded
solution}, is equivalent to the following:
 \be
  \label{Neg}
  \fbox{$
  \mbox{\bf ${\mathbf \not \exists}$ global bounded solution}
  = \mbox{\bf $\exists$ a finite-time blow-up pattern,}
  $}
  \ee
 which corresponds to bounded $L^2$-initial data.
In both cases, we mean blow-up in $L^\iy(\ren)$ of bounded smooth
solutions at the first blow-up time $t=T$.

\ssk

It seems obvious that both scenarios (\ref{Pos}) and \ef{Neg}
assume a detailed study of possible  blow-up behaviour of
solutions as $t \to T^-$, and  this cannot be avoided and
represents a clear alternative face of this  fundamental open
problem of PDE theory.
Of course, \ef{Pos} would be achieved without any blow-up study
provided that a new technique (e.g., a new conservation  law
and/or a monotonicity formulae) could be invented for the NSEs.
However, the long history of this Millennium Open Problem suggests
that this is expected to be extremely difficult. Anyway, this can
happen
since the Navier--Stokes equations can indeed inherit from their
universal physical nature (a system comprising  Newton's Second
Law, the continuity, and a basic viscosity) some extra still
hidden and unknown continuous, discrete, or other
symmetries/monotonicity/symplectic/dissipative, etc. features.

\ssk


 It is definite that, during a few  last years,
 the direction of
 the attacking this open
problem was clearly partially changed and
 a seriously increasing number of papers along the evolution
 blow-up scenarios \ef{Pos} and \ef{Neg} were published.
 In particular, a  rather complete negative answer was achieved
 supporting somehow \ef{Pos} in the following way:
  \be
  \label{Sim}
   \fbox{$
  \mbox{for (\ref{NS1}), a standard self-similar  blow-up as $t \to T^-$  is impossible}.
   $}
   \ee
 In a most general evolution setting, this is due to the work by
  Hou and Li \cite{Hou07}, based on the crucial nonexistence  result in
  Ne\v{c}as--Ru\v{z}i\v{c}ka--\v{S}ver\'ak
   \cite{Nec96}
  proved via the Maximum Principle (the MP). The ban \ef{Sim} was proven in \cite{Hou07}
  by
  solid semigroup theory that is adequate
 to classic evolution approaches to nonlinear PDEs with blow-up
 (we present  a spectral discussion  of \ef{Sim} in Section
 \ref{SLer}).
 In fact, it was proved that any convergence in $L^p(\re^3)$, with
 $p>3$,
as $t \to T^-$ to a rescaled self-similar profile means that the
solution remains bounded at  $t=T$ (so it is a {\em removable
singularity}). This approach does not exhaust all the
possibilities of, say, almost self-similar or other (Type II)
blow-up,
 but, nevertheless, is a new convincing fact
on possible singularities in the fluid  model \ef{NS1}.

\ssk

 The negative result (\ref{Sim}) completes a long remarkable  history  of the
 study of
 similarity blow-up singularities for the NSEs
 that was initiated by J.~Leray in 1933-34 \cite{Ler33, Ler34}, who actually posed a deeper
 problem on both backward and forward phenomena:
 \be
 \label{for1}
 \fbox{$
 \begin{matrix}
 \mbox{{\bf Leray's blow-up scenario:} \quad
self-similar blow-up as $t \to T^-$ ($t<T$)} \ssk\\
 \mbox{and similarity collapse of
singularities as $t \to T^+$ ($t
> T$)};
 \end{matrix}
 $}
 \ee
  see his precise statements and a discussion on these principal issues
  in Section \ref{SLer}.


\subsection{Similar open regularity problems for $N=4$ and for the well-posed Burnett
equations}

 It is important to recognize that, unlike its great sounding and well-established reputation even among non-experts,
  the NSEs Millennium problem is
 just a ``remnant" of a general fundamental difficulty of PDE theory
 in the twenty first century, which seems cannot principally be
 understood by classic modern techniques, to say nothing about a rigorous proof.
  Therefore,  most general
new concepts covering a wider PDE area are in great demand.

 To justify that, one does not need to address essentially other
 classes of nonlinear PDEs and/or systems, and
 can just slightly generalize the classic NSEs.
 First, let us mention that a similar and not less difficult (but indeed more exotic)  open
 problem on existence/nonexistence of global classical solutions
 persists  for the Navier--Stokes equations in dimension $N=4$; see
 Scheffer (1978) \cite{Sc78} and recent developments in \cite{Dong07N4}.

\ssk

 Second, as
a next neighbouring example of a similar  nature associated with
applications,
 consider
 the {\em well-posed Burnett
equations} 
 \be
  \label{NS1m}
 \uu_t +(\uu \cdot \n)\uu=- \n p - \D^2 \uu, \quad \di \uu=0
 \inA,
  \ee
where, in comparison with the standard model \ef{NS1}, the
Laplacian squared with the correct sign $-\D^2$ on the right-hand
side is posed. Note that here one deals with the linear solenoidal
bi-harmonic flow induced by the operator $D_t+\D^2$, which is
parabolic but assumes no any order-preserving and other properties
of the semigroup ${\mathrm e}^{-\D^2 t}$ (which  are somehow
naturally partially  inherited from the positive Gaussian for
$D_t-\D$ in the
NSEs \ef{NS1}). Surely, both problems, the NSEs \ef{NS1} in
$\re^4$ and the Burnett  equations \ef{NS1m} in $\re^7$ look
rather ridiculous and seem cannot have a real application. But for
PDE theory these can be considered as some key and canonical
representatives exhibiting the necessary principal difficulties
(or course, there are other examples like that, e.g.,
supercritical nonlinear Schr\"odinger or Ginzburg--Landau
equations to be used and discussed as well).

\ssk

 The more complicated model \ef{NS1m}
 appears on the basis of Grad's
method in
 Chapman--Enskog expansions for hydrodynamics. In particular,
 equations \ef{NS1m} were studied  in
 \cite{G3} as a particular case of {\em Kuramoto--Sivashinsky}-type PDEs
  (see on a derivation therein also), where it was shown
 that, for the smooth orbits, the following embedding holds:
  \be
  \label{pp1}
   L^{\frac N{3}}(\ren) \LongA L^\iy(\ren) \quad (N>6).
   \ee
   This is an analogy for \ef{NS1}, where the idea of
   a similar transition $L^3(\re^3)\Longrightarrow L^\iy(\re^3)$
  goes back to Leray (1934); see a survey in the next section.
 Fixing
 \be
 \label{N73}
 N=7 \quad (\mbox{or any \,\,$N \ge 7$}),
  \ee
 we arrive at a similar  open problem on
existence/nonexistence of global smooth solutions. For $N \le 6$
the proof is easier. In fact, this is the analogy of $N=2$ for
(\ref{NS1}),
 where existence-uniqueness in the Cauchy problem is due to Leray (1933)
 \cite{Ler33N2} (extended by Hopf in 1951 \cite{Hopf51})
 and Ladyzhenskaya (1958) for the IBVPs
 \cite{Lad58, Lad61}. Obviously, with such a huge
growth of the order of PDEs involved and the dimension $N$ of the
Euclidean space for the Laplacian and the convective term (to say
nothing about the non order-preserving properties of the
higher-order semigroup
 ${\mathrm e}^{-\D^2\, t}$), a detailed and convincing analysis of
 \ef{Pos} and \ef{Neg}, together with a complete description of blow-up
 patterns, seems entirely illusive and non-achievable.

\subsection{
On universality
 of the open $L^p \Longrightarrow L^\iy$ problem in PDE theory}

As we have just mentioned above, roughly speaking, the Millennium
Prize Problem, posed specially for the NSEs, is, in a loose sense,
``non-unique", since similar open regularity problems (or not that
lighter significance) occur for many evolution PDEs of various
types. We list a few of them, where the difficult open
mathematical aspects of global existence and/or blow-up are
associated with the following factors:

(i) supercritical Sobolev parameter range of the principal
operator (hence, standard or very enhanced embedding-interpolation
techniques fails), and, in fact, as a corollary,

(ii) multi-dimensional space $x \in \ren$, with $N \ge 3$, at
least (this leaves a lot of room for constructing various $L^\iy$
blow-up patterns via self-similarity, angular swirl, axis
precessions, linearization, matching, etc.).

\ssk

We now list those PDEs, where we give a few recent basic
references to feel the subject. (I) {Supercritical defocusing {\em
nonlinear Schr\"odinger equation} (NLSE)} (see \cite{Mer05,
Vis07})
\be
\label{Sr1}
   \tex{
-{\rm i}\, u_t= \D u - |u|^{p-1} u \withA p>p_{\rm S}(2)=
\frac{N+2}{N-2} \quad (N \ge 3);
 }
 \ee
 (II) $2m$th-order supercritical {\em semilinear heat equation with
 absorption} ($m=1$ is covered by the MP; see \cite {GW2} and \cite{ChGal2m},
  where the result in \S~4
 for $p>p_{\rm S}(2m)$ applies to small
  solutions only):
\be
\label{Sr2}
   \tex{
 u_t= - (-\D)^m u - |u|^{p-1} u \withA p>p_{\rm S}(2m)=
 \frac{N+2m}{N-2m}\quad
 (N>2m, \,\,\,m \ge 2);
 }
 \ee
 (III) The semilinear {\em supercritical wave equations} (see
 \cite{Ikeh08, Yang08}, as most recent guides)
\be
\label{Sr3}
   \tex{
 u_{tt}= \D u - |u|^{p-1} u \withA p>p_{\rm S}(2)=
 \frac{N+2}{N-2}\quad
 (N \ge 3).
 }
 \ee

One can add to those ``supercritical" PDEs some others of a
different structure such as the {\em Kuramoto--Sivashinsky
equations} for $l=1,2,...$ \cite{G3}
 \be
 \label{hSr1}
  \tex{
 u_t=-(-\D)^{2l} u + (-\D)^l u+ \frac 1{p} \, \sum_{(k)} d_k
 D_{x_k}(|u|^p), \,\,\, |{\bf d}|=1,  \,\,\,  p> p_0= 1+ \frac{2(4l-1)}{N}.
 }
 \ee
Here, $p_0$ is not the Sobolev critical exponent, though precisely
for $p>p_0$, $L^2 \not \Rightarrow L^\iy$ by blow-up scaling,
\cite[\S~5]{G3}. On the other hand, a more exotic applied models
exhibit similar fundamental difficulties such as the following
{\em nonlinear dispersion equation} (see \cite{GalNDE5, GPnde} for
references and some details)
 \be
 \label{hS2}
  \tex{
  u_t= - D_{x_1}[(-\D)^m u]- D_{x_1} (|u|^{p-1}u)
   \withA p>p_{\rm S}(2m)=
 \frac{N+2m}{N-2m}.
 }
 \ee
 In view of the  conservation properties for the
models \ef{hSr1} and \ef{hS2}, these, though being local, can be
more adequate to the nonlocal NSEs \ef{NS1}, than the others
above.

In most of the cases, the operator on the right-hand sides
satisfying for  $u \in C_0^\iy(\ren)$
 \be
 \label{Sr4}
 \tex{
 \AAA(u) = - (-\D)^m u - |u|^{p-1}u
 \LongA
 \langle
 \AAA(u),u\rangle=- \int |D^m u|^2 - \int |u|^{p+1} \le 0,
 }
 \ee
is indeed coercive and monotone in the metric of $L^2(\ren)$,
which always helps for global existence-uniqueness of sufficiently
smooth solutions of these evolution PDEs. For the NLS \ef{Sr1},
this gives a stronger conservation laws than for the focusing
equation with the ``source-like" term $=|u|^{p-1}u$.  Evidently,
replacing $\D$ in \ef{Sr1} and \ef{Sr3} by $-(-\D)^m$, $m \ge 2$
moves the supercritical range to that in \ef{Sr2}. On the other
hand, introducing quasilinear differential operators $-(-\D)^m
|u|^\s u$ with $\s>0$ moves the critical exponent to $p_{\rm
S}(2m,\s)=(\s+1) \frac {N+2m}{N-2m}$. Similar supercritical PDEs
can contain $2m$th-order $p$-Laplacian operators, such as the one
for $m=2$, with $\s>0$,
\be
 \label{Sr4NN}
 \tex{
 \AAA(u) = -
\D(|\D u|^\s \D u) - |u|^{p-1}u, \quad
 \langle
 \AAA(u),u\rangle=- \int |\D u |^{\s+2} - \int |u|^{p+1} \le 0.
 }
 \ee


 However, the
lack of embedding-interpolation techniques to get $L^\iy$-bounds,
which can be expressed as the lack of
  compact Sobolev embedding of the corresponding spaces
  for bounded domains $\O \subset \ren$ (this analogy is not
  straightforward and is used as a consistent illustration only)
 \be
 \label{Sr5}
 \tex{
 H^m(\O) \not \subset L^{p+1}(\O) \forA p
> p_{\rm S}(2m),
 }
 \ee
 actually presents the core of the problem: it is not clear how
 and when bounded solutions can attain in a finite blow-up time a
 ``singular blow-up component" in $L^\iy$.
 For the operator in \ef{Sr4NN}, a similar supercritical demand reads
\be
 \label{Sr5NN}
 \tex{
 W^2_{\s+2}(\O) \not \subset L^{p+1}(\O) \forA p
> p_{\rm S}(4,\s) = \frac{(\s+1)N+2(\s+2)}{N-2(\s+2)}, \quad
N>2(\s+2).
 }
 \ee
 In the given supercritical Sobolev
 ranges, finite mass/energy blow-up patterns for \ef{Sr1}--\ef{hS2} are unknown, as
 well as global existence of arbitrary (non-small) solutions.

It is curious that for the NSEs with the same absorption
mechanism as above,
 \be
 \label{Sr6}
  \tex{
  \uu_t +(\uu \cdot \n)\uu= - \n p + \D \uu -|\uu|^{p-1}\uu,
  \quad {\rm div}\, \uu=0 \inB \re^3 \times \re_+,
  }
  \ee
  by the same reasons and similar to \ef{Sr2},
the global existence of smooth solutions is guaranteed
\cite{Cai08} in the subcritical Sobolev range only: for
 \be
 \label{Sr7}
  \tex{
  p \le 5= \frac{N+2}{N-2}\big|_{N=3} \quad \big(\mbox{and $p \ge
  \frac 72$ by another natural reason}\big).
   }
   \ee

 We thus  claim that, even for the PDEs with local nonlinearities
 \ef{Sr1}--\ef{Sr3} (and similar higher-order others), the study
 of the admissible types of possible blow-up patterns can
 represent an important and constructive problem, with the results
 that can be key also for the non-local parabolic flows such as
 \ef{NS1}, \ef{NS1m}, \ef{Sr6}, etc. Moreover, it seems reasonable
 first to clarify the blow-up origins in some of looking similar and simpler (hopefully, yes, since
 \ef{NS1} is both nonlocal and vector-valued unlike the others)
local supercritical PDEs, and next to extend the approaches to the
non-local NSEs \ef{NS1}; though, obviously, the former ones are
not that attractive and, unfortunately, are not related to
``millennium" issues (however, many PDE experts very well
recognize how important these are for general PDE theory).

\subsection{Main synthetic  goal of the paper and on the proposed style of research}
 \label{S1.4}

For the author, who for almost thirty years dealt with
blow-up singularities in various nonlinear PDEs and mainly in
reaction-diffusion systems of different orders, the appearance of
the paper \cite{Hou07}\footnote{The first preprint of this paper,
which was available to the author, dated 24th August 2007.} was a
crucial sign. Actually, this announced that, if blow-up
singularities in the model \ef{NS1} are possible, these must be of
interesting, complicated enough, non-similarity, and non-symmetric
nature, though anyway a long story of various unsuccessful
attempts to reconstruct those blow-up singularities (if any)
suggested that this is expected to be a difficult affair. This
convinced us  to look now carefully at such a complicated model at
the Navier--Stokes equations
 (a system of four PDEs with $(x,t)
\in  \re^3 \times \re_+$!)
  from the point of view of  standard blow-up
 theory.
   Though the author,
who was dealing with  uncomparably simpler PDEs, which
  however generated
 a number of still open problems
 (further comments on this matter will be given later), naturally expected that a
definite, to say nothing about a rigorous, singularity
construction might not be convincingly done just in view of a
general complexity of the model.

\ssk

Let us specify the actual main goal of this essay.
 Of course, clearly, some part
of the author general motivation is  associated with the
Millennium (\ref{Pos})--(\ref{Neg}) Problem. On the other hand,
there exists the second half of the motivation, which the author
honestly regards as not as less important (and even more valuable
for general PDE theory).

More clearly, in his opinion,
 performing a partial
or most complete classification of singularities for \ef{NS1} and
other related nonlinear models of practical interest
 along the lines of blow-up theory becomes nowadays a
 fundamental mathematical direction.
 Overall, this is about a description of {\em possible complicated
 micro configurations that an evolution system can create}, or,
 \be
   \label{Turb1}
   \underline{\bf Goal:}\quad
  \fbox{$
  \begin{matrix}
 \mbox{for NSEs: \, describe all  {\em ``turbulent" incompressible fluid}} \\
  \mbox{\em configurations
 on possibly minimal micro-$(x,t)$-scales},
  \end{matrix}
 $}
  \ee
    {\em regardless}  whether these are developed in finite-time
 blow-up or non blow-up manners.
 Actually, (\ref{Turb1}) does not assume that the solutions  must be of a finite
 kinetic energy; it is just necessary that the behaviour of
 solutions at infinity (as $x \to \iy$) does not play any
 essential role for formation of patterns. Otherwise, this would
 mean posing various ``boundary conditions" at infinity that can
 immensely increase the variety of admissible patterns.
 Describing some approaches to  (\ref{Turb1}) is our actual {\sc main synthetic goal}.

\ssk

Concerning the style of a structural characterization of our
concepts and ideas, we are oriented to perform our research in a
such unified manner that, at least, formally,
 \be
  \label{meth1}
  \underline{\bf Style:}\quad
  \fbox{$
 \mbox{approaches cover both NSEs (\ref{NS1}) and the
Burnett equations (\ref{NS1m}).
 }
 $}
 \ee
This gives a room for possible future author's apology in the
sense that
 \be
  \label{meth2}
   \fbox{$
   \mbox{if (\ref{Pos}) were proved, the goal would be the 7D
   Burnett equations (\ref{NS1m})}
    $}
    \ee
 (or other complicated nonlinear evolution PDEs with similar
 principal difficulties of $L^\iy$-bounds).
 Of course, \ef{NS1m} will be studied  less, but
 we will indeed seriously present basics of necessary related linear and nonlinear
operator theory covering the case of  the bi-harmonic operator
$-\D^2 \uu$ or the $2m$th-order one $-(-\D)^m \uu$ for any $m \ge
2$.
 It is clear that any rigorous
justifying the conclusions for \ef{NS1m} would be incredibly
difficult or even will be never achieved, but anyway this cannot
prevent us from performing the research: {\em the necessity of a
deep mathematical study of various nonlinear PDEs {\sc without}
any hope of strict formulations of many results is already
clearly an inevitable feature of modern PDE theory associated with
several types of higher-order singular or degenerate
 equations\footnote{``The main goal of a mathematician is not
 proving a theorem, but an effective investigation of the
 problem...", A.N.~Kolmogorov, 1980s 
.}.}

 \subsection{Layout of the paper} The main steps of our analysis
 are as follows:

 \ssk

\noi$\bullet$\underline{\em  Blow-up survey (Section
$\ref{SLer}$)}. We begin with a necessary  short survey devoted to
   classic and recent blow-up results for the Navier--Stokes and
   some for
 Euler equations.

\ssk

\noi$\bullet$\underline{\em  First application of blow-up scaling:
Hermitian structure of zero sets  (Section $\ref{SHerm}$)}. This
is devoted
 to a first application of blow-up scaling to the NSEs. Namely, we
show that the local structure of multiple zeros of regular
solutions can be governed by special vector solenoidal {\em
Hermite polynomials} as eigenfunctions of the adjoint rescaled
linear operator. We claim that this study is a natural and
unavoidable step for further deeper discussion and classification
of all types of micro-scale, single-point configurations that can
be generated by the evolution system in finite time. These sets of
singular patterns are assumed to include also possible
$L^\iy$-blow-up patterns, which are essentially nonlinear and, for
revealing of those, demands complicated matching procedures, where
the Hermitian polynomial space-time structures will be  key.





\ssk

\noi$\bullet$\underline{\em  Third  blow-up scaling: singularities
in NSEs and EEs (Section $\ref{SEE2}$)}. This is another version
of blow-up scaling showing that any blow-up in the NSEs must be
supported by a bounded ``NS-entropy" solution of the EEs, which is
defined on larger
 space-time subsets.

 \ssk

\noi$\bullet$\underline{\em  Blow-up twistor mechanism   (Section
$\ref{STw}$)}. According to typical blow-up results, which are
well-known for a wide audience of mathematicians working with
nonlinear parabolic, hyperbolic, dispersion,  Boussinesq, and
other evolution PDEs, it is clear that,
regardless a pretty strong negative result in \ef{Sim}, the story
of the open problem is far away not only  from being solved but
even reasonable understood.
We intend to show some extra ways how the Navier--Stokes equations
can create complicated blow-up patterns.
 The main difficulty is indeed to detect a suitable  and adequate for fluid vortex models
 mechanism of blow-up
 swirl (rotation) leading to an essential reduction of the system
 \ef{NS1}, which can be at least discussed to clarify the way of
 construction of such blow-up patterns that have been called in Section \ref{STw}  {\em
blow-up twistors}. Namely, we show that, in cylindrical polar
coordinates  $\{r, \var,z\}$, NSEs \ef{NS1} allow a consistent
restriction to the 2D
linear subspace
 \be
 \label{lin11}
 W_2= {\rm Span}\{1,z\},
  \ee
and moreover this subspace is {\em partially invariant} under the
nonlinear operators involved (this is understood in the sense of
mappings, \cite[Ch.~3,~7]{GSVR}). Such blow-up patterns belong to
$W_2$ and lead to  involved nonlinear systems, which are very
difficult to study, but, anyway, this admits  the natural
similarity {\em logarithmic travelling wave} (TW) mechanism of
blow-up vortex swirl
 about the $z$-axis, with the angular dependence
 \be
 \label{ang1}
  \tex{
  \var(t)= -\sigma  \ln(T-t) \LongA \dot \var(\t) = \frac \sigma {T-t} \to
  \iy \asA t \to T^- \whereA \sigma  \not = 0.
 }
 \ee
 In other words, we pose an extra logTW angular dependence in
 cylindrical coordinates\footnote{Recently, the blow-up $\ln(T-t)$ factors are more boldly appear
  for the NSEs; see e.g., \cite[\S~1.2]{Chae08}.}
  \be
  \label{ang2}
   \var= \mu - \sigma  \ln(T-t) \whereA \mu \in (0,2\pi) \,\,\,
   \mbox{is the rescaled angle},
    \ee
which inserts  into the system a new ``nonlinear eigenvalue"
$\sigma \in \re$. As usual, assuming such an extra evolution
freedom in this swirling dynamical system extends the overall
possibility to get suitable blow-up patterns, possibly even in the
self-similar form\footnote{Introducing ``twistors", we mean a
specific self-similar angular blow-up mechanism \ef{ang2} and
around with further perturbations and approximate similarity
features to be described shortly. The fact that blow-up in the
NSEs could be connected with a ``tornado-type" structures is
well-recognized. {\em Q.v.} Sinai \cite[p.~730]{Sinai02}:
``... A negative answer [to existence of strong solutions of the
3D NSEs] could be connected with solutions which develop
singularities in finite time like a \underline{tornado-type
solution} where infinite vorticity appears at some particular
points in time and space." (Underlying is author's.)}. Inserting
the blow-up swirl \ef{ang2} into \ef{for1}
 revives this Leray scenario, since now, in the rescaled
 variables, we observe not a stabilization to a point (already prohibited
 by the MP), but convergence to a {\em periodic orbit}. Roughly
 speaking, this falls into the scope of a much more difficult
 dynamics corresponding to the ``Poincare--Bendixson Theorem"
 (existence of blow-up), or to ``Dulac's Negative Criterion"
 (nonexistence), as in classic ODE theory \cite{Perko}, but the current
 PDEs are infinite-dimensional and nonlocal.
But this is not the end of the story even if a nonexistence
Dulac-like result would have been proved: the rescaled orbits may
converge to various {\em quasi-periodic orbits} with arbitrarily
large number of fundamental frequencies (some details to be also
discussed).

 Meantime,
 we state this our observation as: the NSEs
  \be
  \label{NegT}
  \fbox{$
   \mbox{(\ref{NS1}) on $W_2$ admit a (self-similar) mechanism of
 blow-up swirl at a point}.
 $}
 \ee
 We are not able to study somehow
 rigorously the evolution dynamical systems that occur
 and propose
 a few ideas how various  blow-up twistors can occur via evolution close to certain invariant
 manifolds associated with the similarity flow. Though the fluid
 model \ef{NS1} naturally supports  formation and evolution of
 vortices (``{von K\'arm\'an's streets}"), it should be noted that, mathematically speaking,
 {\em taking into account the rotational torsion-like ``spiral wave" mechanism \ef{NegT}
 and other related axis and vertex precession phenomena
 includes into the nonstationary rescaled system  extra velocity and other parameters,  being  nonlinear
 eigenvalues, that, as usual, improves the overall probability to get
 a necessary pattern by matching of various local blow-up and non-singular flows.}


\ssk

  We introduce and  discuss  these examples as a warning
 showing that a reasonable simple treatment of the scenario
 \ef{Neg} cannot be expected.
 We postpone until Section \ref{SLanFin} construction and analysis
 of  a most involved blow-up pattern  with the swirl for \ef{NS1}
  in the spherical coordinates  $\{r,\th,\var\}$, where this gets much more
 complicated and seems does not admit any lower-dimensional
 reductions. Overall, these
 lead to \ef{NegT}.


\ssk

Actually, our blow-up concepts
  are rather
general and are not bounded by the framework of the classic model
\ef{NS1}.
In particular, as a key ingredient, we show possible generating
mechanisms  of
 countable families of other blow-up patterns, which can exhibit
 different properties. In other words, we propose the 
 following statement, which
 is
 well-understood in reaction-diffusion theory (see comments
 below):
 \be
 \label{nn1}
 \fbox{$
 \mbox{in general, ``self-similarity ban" (\ref{Sim}) does not prevent
 blow-up in NSEs (\ref{NS1}).}
  $}
   \ee
For  higher-order
 systems of PDEs
with \ef{Sim}, {\bf the blow-up story is about to begin}, and the
present paper pretends to be just a first step along a ``blow-up
R--D direction".

\ssk

 Of course, \ef{NS1} is a dynamical system that admits the strong {\em a
 priori} control of the $L^2$-norm of the solutions at the blow-up
 time $t=T$, to say nothing about the evident presence of the MP \cite{Nec96}
 in the stationary rescaled form of the equations,
  but, possibly, the system is complicated enough to
 get over such an obstacle in an evolution way.
 Therefore, we will present in Section \ref{S5.4} a first
  discussion of various ways how to overcame the ban \ef{Sim} keeping
 blow-up similarity rescaled variables. We show
 that a thorough study of blow-up evolution on the quasi-stationary
 manifold of {\em Slezkin--Landau's singular solutions}
 for a submerged jet\footnote{See bibliographical  comments and references in
 Section
 \ref{SLan} and on
 further extensions in Appendix C.}
 \cite{Slez34, Lan44} (see also \cite{Squ51}) could provide us  with other types  of
 patterns.


\ssk

Thus,
 the introduced blow-up twistor mechanism  shows that the Navier--Stokes equations in $\re^3$,
 which naturally support vorticity-type evolution (as a necessary feature of this basic
 fluid model), can develop  blow-up rotational  angular phenomena
 in finite time at a fixed stagnation point of the flow (Section
 \ref{STw}), which naturally leads to periodic $\o$-limits.
 It is not still clear if this twistor construction may lead to a
 truly localized blow-up swirl-like singularity.
 By the partial invariance of the twistor (it belongs to the 2D
 subspace \ef{lin11}), it has a velocity field that is
 unbounded in the $z$-direction. In our discussion, we show that
 the branching phenomena that could lead to $z$-localization of
 such a perturbed twistor, imply a very difficult evolution
 matching-like problem, which seems cannot be tackled rigorously
 still. In Section \ref{SLanFin}, in the spherical geometry, we involve
 extra precession axis and vertex mechanisms that make the problem
 more complicated but built a bridge to more realistic
 generalized quasi-periodic (or periodic) blow-up twistor
 behaviour.

\ssk

 We do not rule out  existence of
 other types of rotational singularities for \ef{NS1}, especially
 in view of the ``multiplicity curse" of blow-up scenarios;    see \ef{fk2}
 below.
 The mathematics of such singularity structures promises
 to be unbelievably difficult since assumes sorting out an infinite
(at least a countable) number of various possibilities. This
``countability curse" for blow-up asymptotics will be explained
later. In view of that, possibly, it is not an exaggeration to say
that it would be  much easier and much more pleasant to settle
\ef{Pos}, and hence would  forget about this ``awkward" multiple
blow-up stuff. This is  a warning again: in reality, blow-up
theory is well-known to be very difficult even for looking very
simple PDE models!
 However, for the NSEs,
 in author's
opinion,
 \be
 \label{aaa1}
 \fbox{$
 \mbox{both claims (\ref{Pos}) and (\ref{Neg}) assume equally
 difficult proofs by blow-up scaling.}
  $}
  \ee
In other words, to find the right answer,
one needs to pass through a sequence of very difficult steps of
analysis that in most issues are overlapping in both approaches.

\ssk

The following possible important feature of the blow-up twistors
deserves mentioning: in view of their extreme rotational nature,
they can create {\em vorticity and sometimes velocity fields that
tend to zero as $t \to T^-$} in the standard weak (integral) sense
in $L^p_{\rm loc}$. Therefore, though we observe that the
vorticity gets infinite as $t \to T^-$, the local total mass (the
integral over a small shrinking neighbourhood of the stagnation
point) becomes negligible and just disappears at the blow-up time.

\ssk

 \noi{\bf Remark:  On infinite family of
patterns with regional and global blow-up}.
 It turns  out \cite{GalJMP} that there exists an infinite countable
 family of blow-up    patterns on $W_2$ in the cylindrical coordinates.
As a compensation for a lack of proper mathematics here, in
\cite{GalJMP}, we managed
  to justify much simpler and rigorously confirmed
 construction of {\em blow-up space jets}  that exhibit {\em effective regional} or
  {\em global
 blow-up}
 in the radial variable. This example underlines another important
 feature of blow-up for \ef{NS1} (though these always are of {infinite} energy):
 \be
 \label{fk2}
  \fbox{$
 \mbox{there can exist  an infinite countable family of
  blow-up patterns}.
  $}
  \ee
In particular,  there exists {\em effective regional blow-up} of
the $z$-component of the vector field in the radial $r$-direction,
i.e., the blow-up wave does not propagate and is of a standing
type. This is a good sign, but of course the velocity field is
unbounded in the $z$-direction, so that the patterns have infinite
energy. A proper ``bending" of the $z$-axis to create a kind of
``blow-up ring" is rather suspicious, since also will exhibit
infinite kinetic energy as $t \to T^-$. A more complicated
geometry is necessary for the next refined blow-up construction.
 Another important feature  is
  \be
  \label{fk43}
  \fbox{$
 \mbox{these blow-up patterns converge  to  similarity solutions of  Euler's equations}.
  $}
  \ee
In the rescaled sense, the convergence turns out to be uniform in
$\re^3$! We do not know whether the phenomenon \ef{fk43} is
expected to be generic for other hypothetical blow-up  patterns
 of finite energy.



Similarly, in \cite{GalJMP}, for the 3D {\em Euler's
equations}\footnote{Leonhard Paul Euler, 1707-1783.} (EEs),
 \be
 \label{ee1}
 \uu_t +(\uu \cdot \n)\uu=- \n p, \quad \di \uu=0 \quad \mbox{in}
 \quad \re^3 \times \re_+,
 \ee
  the restriction to $W_2$ is shown to admit patterns with {\em single
 point blow-up} in the $r$-variable.

 \ssk

 \noi$\bullet$\underline{\em Blow-up patterns with swirl and
 convergence to EEs \ef{ee1} (Section $\ref{Snonrad}$)}.
 This is about possible extensions of the results in \cite{GalJMP}
 to the non-radial geometry with the blow-up swirl included, which
 is a training for more complicated blow-up structures to come.

\ssk

\noi$\bullet$\underline{\em Blow-up about Slezkin--Landau and
other equilibria (Section $\ref{SLan}$)}. This can lead to blow-up
patterns by a kind of ``centre-stable manifold" analysis by using
various singular or regular stationary solutions of the NSEs.
Since some results are derived on the basis of the famous
Slezkin--Landau singular steady states, we include some history of
these solutions and put translations of Slezkin's rare notes of
1934 and 1954 in Appendices A and B at the end of the paper to be
followed by a further discussion in Appendix C.

As a rather general conclusion, we state the following ``steady-TW
exercise" for the rescaled NSEs, which should be solved before
even talking about \ef{Pos} or \ef{Neg}:
 \be
 \label{Stat1}
  \fbox{$
 \begin{matrix}
  \mbox{to describe all families of equilibria in $\re^3\setminus \{0\}$ or
   $\re^3$,} \ssk\ssk\\
   \sigma  \vv_\mu +  \frac 12\, y \cdot \n \vv + \frac 12\, \vv+ (\vv \cdot \n)\vv=- \n p+ \D \vv, \quad {\rm
   div}\, \vv=0,
   \end{matrix}
 $}
 \ee
 where $\sigma  \in \re$ is the nonlinear eigenvalue introduced  as in \ef{ang2} being the blow-up angular
 velocity.
  For $\sigma =0$, \ef{Stat1} is
 Leray's classic rescaled
 stationary problem with a number of  nonexistence results,
 \cite{Nec96, Tsai98, Mil01} (though not everything is still
 known, \cite{Sv06}).

For $\sigma  \not = 0$, suitable profiles in \ef{Stat1} represent
the simplest case of {\em non self-similar blow-up with the
omega-limit set consisting of a periodic orbit}. Nevertheless,
plugging into the system a single new real eigenvalue $\sigma $
can be non-sufficient for settling the question on possible
non-trivial $\o$-limits of rescaled blow-up orbits, so we will
need to introduce further mechanisms of axis precessions leading
to quasi-periodic and other motions.

In general, even for $\sigma =0$, the stationary (nonlocal
elliptic) problem:
 \be
 \label{ST1}
  \fbox{$
  \begin{matrix}
  \mbox{describe all  singular equilibria in $\re^3 \setminus
  \{0\}$,}\ssk\\
  (\UU \cdot \n)\UU=-\n P + \D \UU, \quad {\rm div}\,\UU=0,
   \end{matrix}
   $}
   \ee
 remains still open (though some key steps have been already made, \cite{Tsai98, Sv06, Kim06, Koch07}, which
inspire some optimism that the situation is under a proper
control;
 see below). As a clue exercise to such a complexity, in Section \ref{SLan}, we demonstrate that
 complicated oscillatory sign changing singular equilibria occur
 for the simplest elliptic problem from R--D theory in the supercritical Sobolev
 range:
  \be
  \label{ST2}
   \tex{
   \D u + |u|^{p-1}u=0 \inB \ren\setminus\{0\} \whereA p \ge p_{\rm
   S}= \frac{N+2}{N-2}, \quad N \ge 3.
   }
   \ee
Non-radial singularity structures for \ef{ST2} were not addressed in the literature and
 {\em are unknown}. These can generate extremely complicated  blow-up patterns in the
  corresponding parabolic R--D flow; see \ef{Sob1} below.

  In fact, the proposed swirl behaviour
 \ef{ang2} naturally corresponds to the rotational invariance of
 some symmetric singular equilibria in \ef{ST1}. In general, for
 possible more complicated singular stationary profiles $\UU(x)$,
 new types of ``rotations" are necessary to introduce and
 understand. One of a possible scenario of blow-up in the NSEs is
 as follows:
 \be
 \label{ST3}
 \fbox{$
  \begin{matrix}
 \mbox{the orbit $\{\uu(t)\}$ evolves as $t \to T^-$ ``close" to some singular
 stationary}\\
 \mbox{ manifolds (\ref{ST1}) being  ``trapped" in their complicated structure.}
  \end{matrix}
 $}
 \ee
 Then ``swirling type" of the orbit behaviour close to $x=0$ as $t
 \to T^-$ on shrinking compact subsets will naturally and entirely
 depend on the singular steady manifold involved and can be
 extremely complicated leading to any of multi-dimensional
  quasi-periodic or even chaotic attractors.


\ssk

\noi$\bullet$\underline{\em  Twistors in spherical geometry
(Section $\ref{SLanFin}$)}.
We next perform a formal study of a blow-up twistor in the
spherical coordinates that are necessary for creating a truly
spatially localized pattern. We postpone our final conclusions
inherited from the previous
 analysis until  {Section \ref{SFin}}.
As a by-product, we again naturally arrive at the problem
 \ef{ST1},
which is difficult and open, regardless a good progress made
recently on understanding of the scaling nature and uniqueness of
famous Slezkin--Landau singular solutions of a submerged jet,
\cite{Tsai98, Sv06, MiuTs08}; see Section \ref{STw} for details.
 Thus, we again discuss the possibility that finite
kinetic energy blow-up twistors may be generated in a small
shrinking as $t \to T^-$ vicinity  of {\sc every} ``proper"
singular steady states from \ef{ST1}, which possibly is not
precisely of homogenuity $-1$ in $r=|x|$ and possessing
torsion-precession mechanism of its swirl axis. Other, more
regular ``steady states" are also of importance.



\ssk

 Thus, on the basis of the given fundamental model \ef{NS1}, we will show some
mathematical tools that are necessary to tackle general difficult
problems of blow-up and non blow-up (then a local smooth solution
becomes global and hence unique) for complicated higher-order
PDEs.




\subsection{On some reaction-diffusion and parabolic analogies to be
applied: first exercise on Type I and Type II blow-up patterns}
 \label{S1.5}

Nowadays,  blow-up PDE theory, as a self-contained subject,
embraces a wide range of various nonlinear evolution models; we
list a few monographs from the 1980s and later periods up to 2007,
\cite{Al, BebEb, GalGeom, GSVR, AMGV,   MitPoh,  Pao,   QSupl,
SGKM, Sp2}, where further extensions and references can be found.
Most of them are mainly and specially devoted to blow-up behaviour
in nonlinear partial differential equations of parabolic and
hyperbolic types (those two PDE areas
 are well-established for blow-up since the 1950s)
and contain key literature and various blow-up results achieved in
the last fifty years.
 During this long time, a huge amount of beautiful and difficult
 conclusions were obtained  resulting in
  deep and sometimes exhausting
 understanding and complete classification of
 blow-up patterns for some nonlinear parabolic and other models.

\ssk

 Therefore, for those who are constantly working in these PDE areas,
  it is not surprising that the crucial two the so-called {\em Millennium
 Problems} such as

\ssk

 \noi{\bf (I)} {\em Global existence or nonexistence of smooth solutions  of the
 NSEs in $\re^3$}, and

\ssk

\noi{\bf (II)} {\em The Poincar\'e Conjecture} (see further
comments below),

\ssk

\noi both, in their already existing (for {\bf (II)}) and possible
(for {\bf (I)}) ways of the solution, heavily rely on the study of
blow-up solutions of the nonlinear evolution equations involved.

 \ssk


\ssk

\noi\underline{\em Frank--Kamenetskii equation: self-similar (Type
I) and Type II fast blow-up patterns}.
 In
discussing further consequences related to both problems \ef{Pos}
and \ef{Neg}, we apply some general results of
 blow-up theory developed  for   various reaction-diffusion equations and systems,
 which
 turn out to be  fruitful in application to blow-up for the more
complicated PDEs such as (\ref{NS1}). In particular, it is
well-known that even for simple looking semilinear parabolic
reaction-diffusion PDEs such as
 the classic {\em
Frank-Kamenetskii equation} (1938) \cite{FrK}  developed  in
combustion theory  of solid fuels (also called the {\em solid fuel
model}; first blow-up results in  related ODE models are due to
Todes, 1933)
  \be
  \label{FK1}
u_t= \D u + {\mathrm e}^u \inA,
 \ee
 the property \ef{fk2} holds. Namely,  the first family of {\em ``linearized" blow-up
 patterns} is constructed by linearization techniques and further
 matching  of centre and stable manifolds orbits, where the latter one is
 infinite-dimensional causing the eventual countability of the
 family; see explanations below. Moreover, it is crucial that this family of blow-up
 asymptotics can
  exhaust all possible types of blow-up behaviour that is available in the model.
  Such
 delicate results for $N=1$ or 2 and for other R--D PDEs \ef{Sob1} are due to Vel\'azquez
 \cite{Vel}. Hence, the family of blow-up patterns for
 \ef{FK1} from \cite{Vel, GHPV} is {\em evolutionary
 complete} (a notion introduced in \cite{CompG}, where further
 references can be found).

 For $N \ge 3$, \ef{FK1} possesses
 non-trivial self-similar blow-up (Type I\footnote{The terms ``Type I, II" were borrowed from Hamilton
 \cite{Ham95}, where Type II is
  also called {\em slow} blow-up.}) patterns, which  are called {\em nonlinear
 eigenfunctions}, and, in addition, the set of linearized patterns
 is more involved and  can include other countable families.
 The total family of blow-up patterns
  gets more complicated, so its evolution completeness is
 unknown  representing a difficult open problem for \ef{FK1}.

 It is curious that
 that this remains open even for the radial equation \ef{FK1}
  in $\re^3$,
   \be
    \label{FK111}
   \tex{
  u_t=u_{rr} + \frac 2r \, u_r +{\mathrm e}^u \inB \re_+\times
  \re_+ \quad (u=u(r,t), \,\, r=|x|>0, \,\,\, u_r(0,t) \equiv 0).
   }
   \ee
 To illustrate this fact and our future arguments in Section \ref{S3.8},
  we briefly explain how a countable family of non self-similar
 patterns of Type II can occur.

 \ssk

 \noi\underline{\em Similarity blow-up patterns}. Performing in \ef{FK111} the standard self-similar
 blow-up scaling yields the rescaled equation, which are known to
 exist for $N \ge 3$ \cite{BebEb},
  \be
  \label{WW1}
    \begin{matrix}
   u(r,t)= - \ln(T-t) + v(y,\t), \,\,\, y= \frac r{\sqrt{T-t}},
   \,\,\,
   \t= -\ln(T-t) \ssk\ssk\\
   \LongA v_\t=\HH(v) \equiv \D_N v - \frac 12\, y v_y +{\mathrm e}^v-1.
    \end{matrix}
    \ee
 The self-similar rescaled profiles $f(y) \not \equiv 0$ are its good equilibria
 (we omit details since will  be talking a lot about similarity blow-up later on):
  \be
  \label{WW2}
  \HH(f)=0 \inB \ren, \quad f(y) \,\,\, \mbox{has at most logarithmic
  growth as $y \to \iy$}.
   \ee
 According to \ef{WW1}, each $f(y)$ defines a Type I self-similar
 blow-up (see below).

\ssk

 \noi\underline{\em Non-similarity blow-up patterns}.
 Equation \ef{WW2} in dimension $N  \ge 3$ admits a singular
 equilibrium of the form
  \be
  \label{WW3}
   \tex{
   V(y)=  \ln \frac{2(N-2)}{|y|^2} \inB \ren \setminus\{0\}.
    }
    \ee
   Then some blow-up patterns can evolve as $\t \to +\iy$ ($t \to T^-$)
   close to  this ``stationary manifold".

   \noi\underline{\em Inner Region I expansion: linearization}.
    To see the applicability of this
   idea, we perform the linearization in \ef{WW1},
    \be
    \label{WW4}
 \tex{
 v(\t)=V+ Y(\t) \LongA Y_\t= \D_N Y- \frac 12\, y Y_y + \frac
 {2(N-2)}{|y|^2} \, Y + \DD(Y),
 }
  \ee
where $\DD$ is a quadratic perturbation as $Y \to 0$ in a suitable
metric. By classic theory of self-adjoint operators \cite{BS}, the
linear operator with the inverse square potential in \ef{WW4}
 \be
 \label{WW5}
  \tex{
 \HH'(V)= \D_N  - \frac 12\, y \cdot \n + \frac
 {2(N-2)}{|y|^2} \, I
 }
  \ee
 is well posed in $H^2_{\rho^*}(\ren)$ (with the weight as in
 \ef{mm1}), i.e., has a compact resolvent and a real discrete spectrum,
  provided that the inverse square potential is not too much singular at the origin $y=0$.
  Namely, one needs that
  $2(N-2)$ is less than the constant  $c_{\rm H}$
 of embedding $H^1_0(B_1) \subset L^2_{|x|^{-1}}(B_1)$ by the
  Hardy classic inequality\footnote{It's idea  goes back to the 1920s, \cite{Hardy}, and
  in this form was already used
   by Leray in 1934 \cite{Ler34}; see Section \ref{S7.3} for extra details.}:
  \be
  \label{WW6}
   \tex{
    \int_{B_1} \frac {|u|^2}{|x|^2}\, {\mathrm d}x \le c_{\rm H}
    \int_{B_1} |\n u|^2 \, \, {\mathrm d}x \LongA
  2(N-2) \le c_{\rm H}= \big(\frac{N-2}2\big)^2, \,\,\, \mbox{i.e.,}\,\, N \ge 10.
 }
   \ee

 Thus, in integer dimensions $N \ge 11$ ($N=10$ has own peculiarities; see
  \cite{GK10}), the operator \ef{WW5} has
 a discrete spectrum  and radial eigenfunctions
 satisfying
  \be
  \label{WW7}
  \begin{matrix}
   \s(\HH'(V))= \big\{ \l_k,
   \, k=0,2,4,...\big\};
   \quad
   \psi_k^*(y) \sim b_k y^k+...\,, \,\,\, y \to +\iy, \qquad\quad
  \ssk\ssk\\
   \psi_k^*(y)\sim c_k
   y^{-\d}+...\,, \,\,\, y \to 0, \,\,\,\mbox{where} \,\,\,
   \d= \frac
    {[\, N-2-\sqrt{(N-2)(N-10)}\, ]}2
    \qquad\quad
    \end{matrix}
    \ee
 and $b_k$ and $c_k$ are some normalization constants.
 The orthonormal  set of
  eigenfunctions $\Phi^*=\{\psi_k^*\}$
   is then complete and closed in $L^2_{\rho^*}(\ren)$ in the radial setting.

    Hence, for $N \ge 11$,
  equation \ef{FK111} admits very special asymptotic patterns; see
  \cite{Dold1, GKSob} (global solutions and a survey)
  and \cite{HVsup, Fil00} (blow-up solutions and a survey).
In particular, concerning fast blow-up patterns of Type
II\footnote{In blow-up R--D theory, Type II assumes faster
non-self-similar growth.}, these are obtained by matching of the
linearized behaviour on the steady singular manifold for any
$k=2,4,...$ such that $\l_k<0$:
 \be
 \label{WW8}
  v_k(y,t) = V(y) - C_k {\mathrm e}^{\l_k \t} \psi_k^*(y)+...
 \sim F_k(y,\t) \equiv -2 \ln y -{\mathrm e}^{\l_k \t} y^{-\d}
 \asA y \to 0,
   \ee
 with a bounded flow at the origin $r=0$; see \cite{HVsup, HVsup1, Dold1}
 for first results in this direction, and further references
 above. The correct $L^\iy$-rate of blow-up behaviour is obtained by
 calculating the absolute maximum  of $F_k(y,\t)$: for $\t \gg 1$,
  \be
  \label{WW9}
   \tex{
   (F_{k})_y'=0 \,\,\, \mbox{at}\,\,\,y_k \sim {\mathrm e}^{\frac{\l_k \t}{\d}} \LongA
 \sup_y F_k(y,\t)= F_k(y_k,\t) \sim \a_k \t, \,\, \a_k=\frac
 {2|\l_k|}{\d}>1.
  }
  \ee

\noi\underline{\em Inner Region II: quasi-stationary regular
flow}. Thus,
in Region II close enough to $r=0$, one needs to solve  the
original equation \ef{FK111} with the condition as in \ef{WW9},
which suggests the following scaling therein:
 \be
 \label{W101}
  \begin{matrix}
  u(0,t)= - \a_k \ln(T-t) \LongA u(r,t)= - \a_k \ln(T-t)+
  w(\xi,s), \qquad\quad
  \ssk\\
  \mbox{where} \quad \xi= \frac r{(T-t)^{\a_k/2}} \andA s=\frac 1{\a_k-1}\,
  (T-t)^{1-\a_k} \to +\iy.\qquad\quad
  \end{matrix}
   \ee
   Substituting this into the F--K equation, after elementary
   manipulations yields the following perturbed problem for $w$:
 \be
 \label{W102}
  \tex{
  w_s= \D_\xi w + {\mathrm e}^w - \frac {\a_k}{\a_k-1}\, \frac 1s
  \,\big(1+ \frac 12\, w_\xi \xi\big) \forA s \gg 1.
  }
  \ee
Recalling that according to \ef{WW9}, we have to have that
$w(0,t)=0$, general stability theory for such blow-up singularity
problems \cite{AMGV} suggests, that since \ef{W102} is a perturbed
gradient flow, there is the stabilization to the unique bounded
stationary solution: uniformly on compact subsets in $\xi$, as $s
\to +\iy$,
 \be
 \label{W103}
  \tex{
   w(\xi,s) \to W(\xi) \whereA \D_\xi W+{\mathrm e}^W=0, \,\,
   W(0)=0 \,\, (W(\xi) \sim - 2 \ln \xi, \,\, \xi \gg 1).
   }
   \ee
Eventually, in this Region II, the current blow-up pattern behaves
as:
\be
\label{W104}
 \tex{
 u(r,t) \sim  \a_k |\ln(T-t)| + W\big(\frac
r{(T-t)^{\a_k/2}}\big) \asA t \to T^-.
 }
 \ee

Thus,  \ef{WW1} gives the asymptotic structure of such blow-up
patterns $\{u_k(r,t), \, k \ge 2\}$ with the following fast
blow-up rate:
 \be
 \label{WW10}
 \tex{
 u_k(0,t) \sim \big(1+ \frac {2|\l_k|} \d\big) |\ln(T-t)| \asA t \to T^-
 \quad \big(\l_k \sim -\frac k2, \quad k \gg 1\big).
  }
  \ee
These Type II non self-similar blow-up patterns can have
arbitrarily fast growth for $k \gg 1$ than the standard Type I
similarity divergence associated with the ODE
 \be
 \label{WW11}
  u'= {\mathrm e}^u \LongA u(t) = |\ln(T-t)| \asA t \to T^-.
   \ee

  For $N=10$, the origin $y=0$ is in the limit-point case for the
  operator \ef{WW5}, so that its deficiency indices  are
  $(2,2)$ \cite{Nai1}, and there still exists its proper self-adjoint extension
 with a discrete spectrum and compact resolvent, so that \ef{WW7} holds;
  see details of an application  in
 \cite{GK10}. This also allows to construct blow-up patterns
 with a slightly different blow-up rates.

 On the other hand, for  $N \in [3,9]$, instead of \ef{WW7}, the
 origin 0 is oscillatory: as $y \to 0$,
 \be
 \label{KK1}
  \tex{
  \psi^*(y) \sim y^{-\frac{N-2}2} [C_1\cos(b \ln y)+C_2\sin(b \ln y)]  \whereA b=
  \frac{\sqrt{(N-2)(10-N)}}2.
   }
   \ee
 This yields that any eigenfunction of any self-adjoint extension of \ef{WW5}
 (existing by classic theory \cite{Nai1}) will make the function
 \ef{WW4} to be sign-changing and oscillatory near the matching
 origin $y=0$. Indeed, for nonnegative solutions, this prohibits
 any matching at all, so that
Type II blow-up patterns are
 then nonexistent,
  which is rigorously proved in Matano--Merle
 \cite{Mat04}.
Therefore, it is plausible that for  radial nonnegative solutions
in dimension $3 \le N \le 9$, the only blow-up patterns are
exhausted by the family of linearized ones (as in Section
\ref{S3.8}) and the nonlinear self-similar solutions \ef{WW2},
though there is no proof of such a completeness. However, in the
nonradial geometry for solutions of changing sign such a
classification blow-up conclusion is far away from being even
formally justified.

In view of \ef{Sim}, one of our tricks is to
 check weather Type II non-similarity blow-up patterns can be
 constructed for the NSEs, though this should be much harder of course.



\ssk

\noi\underline{\em Quasilinear reaction-diffusion extension:
countable sets of blow-up patterns}.
 For a quasilinear extension of \ef{FK1} in 1D, with
 the so-called $p$-Laplacian operator (here $p \mapsto \s$),
  \be
  \label{FK11}
  u_t= (|u_x|^\s u_x)_x+ {\mathrm e}^u \quad (\s \ge 0),
   \ee
 a countable set of blow-up patterns was constructed in
 \cite{BuGa}. As an intrinsic feature of essentially quasilinear
 problems, it was shown that, depending on $\s>0$, first patterns
 represent {\em nonlinear eigenfunctions}, i.e., are self-similar,
 while the rest are constructed by the linearization techniques
 and matching as for \ef{FK1}. For all $\s > \s_\infty=0.60...$,
 all blow-up patterns are nonlinear, i.e., then \ef{FK11} admits
  countable families of self-similar solutions.
Moreover, it was also shown that, in view of the strong degeneracy
of the $p$-Laplacian operator in \ef{FK11}, there exist other
blow-up similarity solutions, which can compose {\em uncountable}
families.
 It is still an open problem to prove evolution
 completeness of these blow-up families.
We reflect this discussion as follows: in the family \ef{fk2},
 \be
 \label{DD1}
  \fbox{$
 \mbox{there may exist infinitely many  nonlinear or
 ``linearized" blow-up patterns.}
 $}
  \ee

\ssk

\noi\underline{\em Higher-order parabolic extensions}. Some types
of blow-up singularities are already known for higher-order
reaction-diffusion equations such as
 \be
 \label{ho1}
 u_t=-(-\D)^m u + |u|^{p-1} u \inA  \quad (m \ge 2, \,\,\, p>1)
  \ee
   \be
 \label{CH1}
  \tex{
    \mbox{and} \quad
 u_t= - \D^2 u- \D (|u|^{p-1}u) \inA \quad (p>1);
 }
 \ee
 see \cite{BGW1, EGW1, Bl4, Gal2m} and references therein
concerning other models with blow-up.
  The latter {\em limit unstable Cahn--Hilliard equation} \ef{CH1}
 has the divergence form and hence the flow preserves the total
 mass of $L^1$-solutions (in this sense this mimics a similar
 feature of \ef{NS1}).
Equation \ef{CH1} is shown to obey Leray's blow-up scenario
\ef{for1} \cite{GalJMP}.
 Note that the  questions on
 a complete description of nonlinear and linearized blow-up patterns for
 \ef{ho1} or \ef{CH1} and evolution completeness of the whole countable family
 remain open even in the one-dimensional case $N=1$ and $m=2$ (the fourth-order diffusion only).


 \subsection{Second discussion around the
 NSEs}

Thus, the negative conclusion \ef{Sim} rules out a {\em standard}
self-similar way of blow-up in the Navier--Stokes equations.
Anyway, regardless known deep ideas on interaction of ``vortices
tubes" and others (mainly coming from fluid dynamics and less
mathematically developed), it seems we cannot  still imagine how
complicated other non self-similar individual ``linearized"
blow-up patterns can be (for instance, these can be partially
invariant, i.e., ``partially nonlinear", or non-invariant at all),
to say nothing about possible interaction of a number of different
related blow-up ``vortex tubes". This is our aim  to give, on the
basis of a reaction-diffusion blow-up experience,
  a possible new insight
 into the world of imaginary singularities for \ef{NS1}.

\ssk

  In view of a more complicated nature of the Navier--Stokes
  equations that are a system of four equations in the four-dimensional
  space-time continuum\footnote{I.e.,
     a ``dynamical system $4 \times 4=(\mbox{4 dependent $\times$ 4
     independent variables})$".}
   $\re^3\times \re_+$, the
  slogan  (\ref{fk2}) should contain more surprises.  At least,
  one can expect a formally infinite number of possible types of blow-up to be checked out,
  and, to achieve the positive
  answer \ef{Pos},
   all of them should be ruled out by the assumption of  finite energy and/or others.
     It is not an exaggeration to assume that
  (\ref{NS1}) may admit a countable set of different countable (or
  even uncountable)
  scenarios  of formation of single point blow-up patterns that admit local construction (on shrinking
  subset to the blow-up point)  by
  some kind of combination of linearized and nonlinear eigenvalue
  techniques. Then the crucial fact on the finiteness of
 the energy of the globally extended patterns can be checked only
 after matching (or non-matching) of these local singularity structures with
 surrounding bundles of less blowing up or even  bounded ``tails".
  Such matching procedures (which are responsible for existence or nonexistence
 of the patterns) are also expected to be extremely
 difficult.

\ssk

   {\em Vice versa}, to achieve the negative claim
  \ef{Neg}, in order to construct a suitable blow-up pattern,
  which
  prohibits global evolution of bounded smooth finite kinetic  energy solutions, one then
  needs to sort out, possibly again, a huge (infinite and even uncountable) number of
 {\em a priori} unknown blow-up solution structures.

\ssk

 Note that, for many other models
including various nonlinear higher-order PDEs, for which there is
no any hope to control existence and/or nonexistence of blow-up
singularities by some
conservation/monotonicity/order-preserving/etc.,  mechanisms (that
are nonexistent for sure), the careful study of \ef{Pos} and
\ef{Neg} is unavoidable, is very difficult indeed, and, in many
cases, in a full generality, is not doable at all. The latter
means that this problem is not {\em analytically solvable}, i.e.,
constructive conditions that guarantee \ef{Pos} (or \ef{Neg})
cannot be derived.
  Nevertheless, this is not a manifestation of any kind of a
``PDE agnosticism", since the most stable and generic
 asymptotic structures can and must be studied and understood by any, rigorous or not, mathematical
  means.
 The illusive is a full description of all the possible singularity
 patterns and their {\em evolution completeness} (see a proper setting for this below)
 to eventually {\em guarantee} \ef{Pos} or \ef{Neg}.



\subsection{On a related blow-up parabolic area: Ricci
flows and the Poincar\'e Conjecture}

 Another amazing  geometry and PDE area, which  is not less famous nowadays,
   is the study of blow-up structure of
{\em Ricci flows} of metrics in $\re^3$.  
 Following the logic of evolution blow-up PDE analysis,
 it can be characterized that
 {\em Perel'man's
new monotonicity formulae} and {\em principals of his blow-up
surgery} at finite-time singularities of the Ricci flow
 made it possible to guarantee to get
 a proper global extension beyond all blow-ups and with
a kind of necessary ``symmetrization" at the finite extinction
time.
In particular, this allowed to prove
 the
{\em Poincar\'e Conjecture}\footnote{``{\em A closed, smooth, and
simply connected 3-manifold is homeomorphic to ${\bf S}^3$}. This
still remains the only solved Millennium Problem among other
seven.}
  see \cite{CZ06}
  for details,  history, references, and recent development.
Proposed for this kind of analysis by Richard Hamilton in 1982 the
{\em Ricci flow} for a family of Riemannian metrics
$g(t)=\{g_{ij}(t)\}$:
 \be
 \label{Ri1}
g(x,t): \quad  g_t=-2{\rm Ric}\, g, \quad \mbox{or, for
components,} \quad
 (g_{ij})_t=-2{R_{ij}},
   \ee
   where ${\rm Ric}\, g=\{R_{ij}\}$ is the Ricci curvature of $g$, represents
   a
   system of  parabolic PDEs for the components, which
   overall obeys the Maximum Principle (the MP) and other related classic
   properties of parabolic flows. The {\em scalar curvature
   equation} for the scalar curvature $R=(g^{ij}R_{ij})(x,t)$ is a semilinear
   parabolic PDE with a quadratic nonlinearity,
 \be
 \label{Ri2}
 R_t= \D R + 2 |{\rm Ric}\, g|^2,
  \ee
  where $\D \le 0$ is the corresponding Laplacian.
Moreover, in dimension 2, the scalar curvature $R(x,t)$
   then takes the form of a standard   semilinear quadratic R--D
   equation \cite[p.~227]{CZ06}
\be
  \label{SEHam1}
 \mbox{$
 R_t = \D R + R^2 \quad
 (N=2).
  $}
  \ee
This shows a  (obviously, well-known) link  between Ricci flows
singularities with standard and well-developed blow-up R--D
theory.
  Similar to \ef{FK1},
  blow-up for
  this equation has  a long history and embraces hundreds of
  publications in the 1980-90s
  in almost all leading world journals on
   nonlinear PDEs and a few  monographs mentioned above.
  For future use, let us note that a  key idea of
 its generic non-trivial blow-up behaviour with  logarithmic (the appearance  of $\ln
 (T-t)$-term is due to a centre subspace behaviour as explained in Section \ref{S3.8},
 so we avoid  further comments),
  \be
  \label{ln1}
  \sim \sqrt{|\ln(T-t)|} \asA t \to T^-,
    \ee
 deformation of similarity
blow-up structures dates back to Hocking, Stuartson, and Stuart
in 1972, \cite{HSS}.
 It was proved rigorously  twenty years later  that the actual structurally stable blow-up
 behaviour for \ef{SEHam1} as $t \to T^-$ close to the blow-up
 point $x=0$ is given by
  \be
  \label{bl11}
  \tex{
   R(x,t) = \frac 1{T-t} \, f_*\big(
   \frac{x}{\sqrt{(T-t)|\ln(T-t)|}}\big)(1+o(1))
   \whereA f_*(y)= \frac 1{1+ c_*{|y|^2}},
    }
    \ee
    and $c_*>0$ is a constant depending on the dimension $N$ only
    \big($c_*= \frac 14$ for $N=1$\big);
     see \cite{Vel} and \cite[p.~312]{SGKM} for  references
     and results. A full account of results on blow-up for
     semilinear equations such as \ef{SEHam1} can be found in the
     most recent monograph by Quittner and Suplet \cite{QSupl}.
 Formula \ef{bl11} is universal for
 dimensions $N<6$, i.e., in the subcritical Sobolev range for the
 elliptic operator in \ef{SEHam1}, with the general $|u|^{p-1}u$ term,
  \be
  \label{Sob1}
   \tex{
    u_t=\D u + |u|^{p-1}u, \quad \mbox{with}
    \quad
  p=2 < p_{\rm Sobolev}= \frac{N+2}{N-2} \quad(N>2).
   }
  \ee

As an accompanying key feature for further use here,  the
structurally stable (generic) blow-up pattern \ef{bl11} perfectly
serves as a powerful confirmation of the slogan \ef{nn1}. Indeed,
the similarity scalings for the equation \ef{SEHam1} are standard
and simple,
 \be
 \label{ss21}
  \tex{
  R(x,t)= \frac 1{T-t} \, \hat R, \quad y= \frac x{\sqrt{T-t}},
   }
   \ee
   so do not contain any of logarithmic factors as in \ef{ln1} and
\ef{bl11}, which are created in the {\em blow-up evolution} as $t
\to T^-$ for almost arbitrary suitable solutions with $L^\iy \cap
L^1$-data.

\ssk

It is indeed also worth mentioning that since the function
$f_*(y)$ in \ef{bl11} is radial, this expresses the
 phenomenon  of strong {\em symmetrization} of blow-up structures near blow-up
time that is true for almost   all admissible  non-constant
solutions of \ef{SEHam1} (classification of singularities is also
a key feature of Perel'man's proof, which the
crucial blow-up surgery is based upon, \cite[Ch.~7]{CZ06}). We
will constantly use the generic blow-up behaviour
 \ef{bl11} as a source for further speculations concerning the
 model
 \ef{NS1}.

It is remarkable that one of the key ideas of Perel'man's proof
was surgery followed by a necessary classification of possible
blow-up singularities, which can occur for the Ricci flow
\ef{Ri1}. Recall our main target \ef{Turb1} for the NSEs.
\subsection{Reaction-diffusion: a new Type II blow-up patterns for $p=p_{\rm S}$}
 \label{SRD1}

The parabolic equation \ef{Sob1} in the critical case $p=p_{\rm
S}$ for $N=3$ reveals another type of construction of Type II non
self-similar blow-up patterns. Namely, in the critical Sobolev
case, \ef{Sob1} admits the {\em Loewner--Nirenberg  stationary
solution} (1974) \cite{LN}
  \be
  \label{LN11}
   \mbox{$
  u_{\rm LN}(x)=\Big[ \frac {N(N-2)}{N(N-2) +
  |x|^2} \Big]^{\frac {N-2}2} \equiv
  \sqrt{\frac {3}{3 +
  |x|^2}}
   \forA N=3.
   $}
    \ee
This is indeed a very special and remarkable explicit solution, so
that, similar to the scenario in Section \ref{S1.5}, blow-up can
occur about the 1D manifold of such rescaled equilibria. Namely,
following \cite{Fil00}, we explain how this can happen. First, as
in \ef{WW1}, we perform the standard self-similar scaling in
\ef{Sob1}, with $p=5$ for $N=3$,
 \be
 \label{sc11}
 \begin{matrix}
 u(x,t)=(T-t)^{- \frac 1{p-1}}v(y,\t), \quad y= \frac
 x{\sqrt{T-t}}, \quad \t=-\ln(T-t), \quad \mbox{so} \ssk\ssk\\
  v_\t= \HH(v) \equiv \D v- \frac 12\, y\cdot \n v - \frac
  1{p-1}\, v+|v|^{p-1}v.
   \end{matrix}
   \ee
Second, let us assume that $v(y,\t)$ behaves for $\t \gg 1$ being
 close to the stationary manifold composed of equilibria
\ef{LN11}, i.e., for some unknown function  $\var(\t) \to +\iy$ as
$\t \to +\iy$,
 \be
 \label{sc12}
  v(y,\t) \approx \var(\t) u_{\rm
  LN}\big(\var^{\frac{p-1}2}(\t)y\big)
 \ee
on the corresponding shrinking compact subsets in the new variable
$\zeta=\var^{\frac{p-1}2}(\t)y$. It then follows that, in the case
$N=3$, on the solutions \ef{LN11} in terms of the original
rescaled variable $y$ (see computations in \cite[p.~2963]{Fil00};
our notations have been slightly changed)
 \be
 \label{sc13}
  \tex{
 |v(y,\t)|^{p-1}v(y,\t) \to \frac{4 \pi \sqrt 3}{\var(\t)} \,
 \d(y) \quad \mbox{as\,\, $\t \to +\iy$}
 }
  \ee
  in the sense of distributions.
 Therefore, equation \ef{sc11} takes asymptotically the form
  \be
  \label{tt127}
   \tex{
 v_\t= \HH(v) \equiv \D v- \frac 12\, y\cdot \n v - \frac
  1{4}\, v +\frac{ 4 \pi \sqrt 3}{ \var(\t)} \,
 \d(y) +...\forA  \t \gg 1.
 }
  \ee

Thus, we get an  eigenvalue problem for  Hermite's classic
operator \cite[p.~48]{BS}:
 \be
 \label{Her77}
  \tex{
  \BB^*=\D v- \frac 12\, y\cdot \n v, \quad \mbox{where}
  \quad \s(\BB^*)=\{\l_\b=- \frac{|\b|}2, \,\, |\b|=0,1,2,...\},
   }
  \ee
  defined
  in $L^2_{\rho^*}(\re^3)$, $\rho^*(y)={\mathrm e}^{-|y|^2/4}$,
  with the domain $H^2_{\rho^*}(\re^3)$,
 to be used
   later\footnote{The same operator occurs for the rescaled equation (\ref{SEHam1}), where
   the factor $\sqrt{|\ln(T-t)|}=\sqrt \t$ is due to the
   eigenspace behaviour with $\l_2=-1$ and the second Hermite polynomial
   $H_2(y)=
    \frac 1{\sqrt 2}\,(1- \frac {y^2}6)$.}. It follows from \ef{tt127}
 that one can try the following regular parts of such Type II blow-up
 patterns: balancing two terms in \ef{tt127} yields
  \be
  \label{Her21}
  v_\b(y,\t)={\mathrm e}^{(\l_\b-\frac 14) \t} H_\b(y)+...
 \LongA \var_\b(\t) \sim {\mathrm e}^{\frac{2|\b|+1}4\, \t}
  \forA \t \gg 1,
   \ee
   where $H_\b(y)$ are Hermite polynomials as the eigenfunctions of
   \ef{Her77}.
 Together with the scaling in \ef{sc11}, this yields a countable family of complicated
blow-up structures. To reveal the actual space-time and changing
sign structures of such Type II patterns, special matching
procedures apply. In \cite{Fil00},  this analysis has been
performed in the radial geometry, though (and this is also key for
us to recognize) still no rigorous justification of the existence
of such blow-up scenarios is available.


\subsection{Further  comments on  the ``RD-sense" of this  essay}

As was already announced, the present paper is a general view {\em
from blow-up reaction-diffusion theory to possible
``singularities" (including, e.g., multiple zeros) of the 3D
NSEs}. In other words, the main intention of the author is to
involve some RD-experience into the NSEs study.
 As is easily detectable by experts in the area, this
special issue dictated, in particular, a quite special kind of
survey and references involved in the list, which by no means
reflects the actual history and modern trends of NSEs and Euler's
equations theory (to say nothing about  R--D theory as well).
Nevertheless, the author hopes that a ``detached observer view"
can deliver some new useful accents to the area.

\ssk

 The subjects of the R--D and the NSEs do not coincide at all and seems their overlapping
is  suspicious in many places. Since the NSEs subject is huge and
extremely difficult, the author apologizes for any inconvenience
caused by his approximate constructions and/or speculations that,
from classic fluid dynamic PDE views, can be classified as
 rather non-consistent or well known somehow.
It seems understandable that,
in view of not that optimistic final conclusions (cf. \ef{PP11}),
which were predicted in advance, the author, with a completely
different background, did not have a very strong motivation to
make more clear some of the matching constructions when this
looked being possible (not often, unfortunately).
The
``semi-mathematical" (and sometimes logically incomplete) language
of presenting the speculations reflects the obvious truth: each
elementary step in this long story of singularities/blow-up for
the NSEs could take years to fix (and never known how many). So,
the author took the risk to show quite a discontinuous route   (a
``jumping frog"-style) to the end of this story, of course, with a
full understanding how risky this way could be in mathematics,
when many steps and concepts of his speculations can be attacked
by justified incinerating critics from
 some experts from the field.
Anyway,
 the author hopes that the attentive Reader will find
some ideas and concepts of blow-up from the R--D systems useful,
 even if these were not presented on a sufficiently rigorous
 basis and costs; and even will find satisfactory a rigorous
 (``almost", i.e.,  if fixed in a reasonable
 finite time) construction of a countable family of blow-up space
 jets for \ef{NS1}.


\ssk

 Indeed, the combination of all the possible tools of singularity
 analysis, leading to a success for the NSEs \ef{NS1},  will serve as a solid and
 reliable basis for further development of modern  theory of
 higher-order PDEs. In this context, blow-up scaling methods become more and more
 penetrating
 into the core of PDE theory being
 natural and unavoidable tools of evolution analysis of higher-order nonlinear equations in all its
 three  ``hipostases":
 $$
 \mbox{\bf (i) existence, \,\,(ii) uniqueness, \,\,{\rm and} \,\,
 (iii) asymptotic behaviour.}
  $$

\section{Some facts of singularity history for the Navier--Stokes
equations: Leray's blow-up solutions, singular points, spectra,
and patterns}
 \label{SLer}

We present here a short survey on global solvability, singularity
formation, and other classic facts concerning   the Navier--Stokes
equations (\ref{NS1}), and stress a special attention to some
blow-up issues,
 which will be used in what follows. A perfect and detailed overview of main
 mathematical results concerning the NSEs is available in Taylor
 \cite[Ch.~17]{Tay}, which includes several aspects to be quoted
 below without proper referencing.

\subsection{Leray--Hopf (``turbulent") solutions of finite kinetic energy}

 It is a classic matter  that the energy $L^2$-norm is natural
for (\ref{NS1}). After multiplication by $\uu$ in  the metric
$\langle \cdot,\cdot \rangle$ of $L^2(\ren)$ and integration by
parts, the convective and pressure terms vanish on smooth enough
functions $\uu(x,t)$ with sufficiently fast (say, exponential)
decay at infinity,
 \be
 \label{na1}
  \langle (\uu \cdot \n)\uu, \uu \rangle =0
  \andA -\langle \n p, \uu \rangle = \langle p, \n \cdot \uu
  \rangle=0.
   \ee
Therefore, on  such smooth solutions, we have the instantaneous
rate of dissipation of kinetic energy given by
 \be
 \label{ss1}
  \tex{
   \frac 12\, \frac{\mathrm d}{{\mathrm d}t} \, \|\uu(t)\|_2^2 = -
   \| D \uu(t)\|_2^2 \LongA \|\uu(t)\|_2 + 2 \int\limits_0^t \| D \uu(s)\|_2^2
    \, {\mathrm d}s = \|\uu_0\|_2, \,\,\, t
   \ge 0. }
    \ee
Actually, the estimate in (\ref{ss1}) with the inequality sign
$``\le "$ is the energy inequality
 for  {\em Leray--Hopf weak solutions} of (\ref{NS1}) (in 1933, Leray also called
 such weak solutions ``turbulent" \cite[p.~231, 241]{Ler34} and compared these with regular ones;
  see also Lions
 \cite[Ch.~1, \S~6]{JLi} for a discussion); {\em q.v.} e.g.,
 \cite{Qio07} and references therein.
\ef{ss1} is also a crucial identity for general turbulence theory.
E.g., the famous
  {\em Kolmogorov--Obukhov power ``K-$41$"
law} (1941) \cite{Kolm41, Obu41}  for the energy spectrum of
turbulent fluctuations for wave numbers $k$ from the so-called
{\em inertial range},
 \be
  \label{Kom}
  \tex{
 E(k)= C \e^{\frac 23} k^{-\frac 53} \whereA \e= \frac 1{\rm Re}
 \, \langle |D \uu|^2 \rangle,
 }
 \ee
 uses this rate of dissipation of kinetic energy for high Reynolds numbers
  ${\rm Re}$ ($\langle \cdot
 \rangle$ is an  invariant measure of calculating expected values).

Note another, weaker ``conservation laws", since the convective
terms in the NSEs are in divergent form ($\otimes$ is the tensor
product of vectors in $\re^3$):
 \be
 \label{Div1}
  \tex{
  (\uu \cdot \n)\uu = {\rm div}\,  (\uu\otimes \uu) \equiv 
  (u \uu)_{x_1}+(v \uu)_{x_2}+(w \uu)_{x_3}.
   }
   \ee
Hence, integrating over $\re^3$ with the necessary decay at
infinity, one can get that the total ``masses" of the velocity
components are preserved:
 \be
 \label{Mass}
  \begin{matrix}
   \frac{\mathrm d}{{\mathrm d} t} \, \int
   u(x,t)\, {\mathrm d}x,\,\,\frac{\mathrm d}{{\mathrm d} t} \,
   \int
   v(x,t)\, {\mathrm d}x, \,\, \frac{\mathrm d}{{\mathrm d} t} \, \int
   w(x,t)\, {\mathrm d}x
   =0, \quad \mbox{i.e.,}
   \ssk\ssk\ssk\ssk\\
    \int
   u(t) \equiv  \int
   u_0, \,\, \int
   v(t) \equiv  \int
   v_0, \,\, \int
   w(t) \equiv  \int
   w_0 \forA t\ge 0.
    \end{matrix}
    \ee
    Of course,
   as estimates, \ef{Mass} are  weaker than \ef{ss1},
 and  the mass ``semi-norms" by taking the absolute values of the integrals
 there are difficult to apply for solutions
of changing
 sign.
    Though, as
  sharp conservation properties, these can be valuable
 in construction (or prohibiting) special  sensitive blow-up patterns.
 Moreover, for some parabolic problems with the mass conservation (and no $L^2$-control), estimates
 such as
 \ef{Mass} can play a key role for global extension of blow-up
 solutions beyond blow-up, for $t>T$. A discussion on the
 application of Leray's blow-up scenario \ef{for1}
  to the limit
 Cahn--Hilliard equation \ef{CH1}
 is presented in  \cite{GalJMP}.

Nevertheless,
 the $L^2$-bounds as in \ef{ss1} are not sufficient to control
the $L^\iy$-non-blowing up property of solutions, and this is
the origin of extensive mathematical research in the last seventy
years or so. Recall that
 global weak solutions of (\ref{NS1}) satisfying
  $$
  u \in L^\iy(\re_+;L^2(\re^3)) \cap L^2(\re_+;H^1(\re^3))
  $$
  were already constructed by Leray \cite{Ler33, Ler34} (1933), and Hopf \cite{Hopf51}
  (1951).
We recommend recent papers \cite{Aus04, Bis07, Can04K, Chem06,
Cort07} as a guide to  various results of local and global (for
small and other data) theory of \ef{NS1}, including analyticity
results in both spatial and temporal variables; see \cite{Dong07,
Zub07} for a modern overview.



\subsection{Blow-up self-similar singularities with finite energy
are nonexistent: on Leray's scenario of backward and forward
blow-up self-similar singularities}

 It seems that the original  idea that the classic fundamental problem of the
unique solvability of (\ref{NS1}) in $\re^3$, i.e., existence of a
global smooth bounded $L^2$-solution, is associated with existence
or nonexistence of certain {\em blow-up singularities} as $t \to
T^-$, goes back to Th.~von~K\'arm\'an; see 
\cite{vonK11, vonK21, vonK34}\footnote{The author apologies for
not being able to trace out von K\'arm\'an's original work (or a
lecture?), where the idea of singularity of the velocity field
$\sim(T-t)^{-\a}$ ($\a=\frac 32$, $\frac 52$, or $\frac 25$?)
appeared first.}
 and survey
\cite{Berker}.

As usual, similarity solutions, as  a next manner to further
specify the behaviour, can be attributed to
 an invariant group of scaling transformations: if
 $\{\uu(x,t),\, p(x,t)\}$ is a solution of the NSEs \ef{NS1}, then
  \be
  \label{Symm1}
  \big\{\uu_\l(x,t)= \l \uu(\l x,\l^2 t), \, p_\l(x,t)= \l^2p(\l x,\l^2
  t)\big\} \,\,\, \mbox{is also a solution for any $\l \in \re$}.
   \ee
Setting here $\l = \frac 1{\sqrt{\mp t}}$ formally yields two
types of self-similar solutions of \ef{NS1}, the blow-up and the
global ones. Let us now focus on
 key historical aspects of such a study.

Namely, in 1933, J.~Leray \cite{Ler33, Ler34} proposed a
mathematical question to look for blow-up
in (\ref{NS1}) 
  driven by the self-similar solutions\footnote{Such a blow-up {\em backward
 continuation variable} $y$  for 1D linear parabolic equations was already systematically used by
 Strum in 1836 \cite{St36}; on his backward-forward continuation analysis, see \cite[p.~4]{GalGeom}.}
 \be
 \label{NS2}
  \tex{
  \uu(x,t)= \frac 1{\sqrt{T-t}} \, {\bf U}(y), \quad p(x,t)= \frac
  1{T-t} \, P(y) \whereA y=\frac x{\sqrt{T-t}}.
 }
 \ee
This  Leray's statement is well known and was stressed to  in many
papers \cite[p.~225]{Ler34}\footnote{Here and later, boxing and
underlying are author's.}:

\ssk

\ssk

\com{ \noi  "... la solution des \'equations de Navier dont il
s'agit est:
 $$
 (3.12) \qquad\qquad u_i(x,t)=[2\a(T-t)]^{-\frac 12} \,U_i[(2
 \a(T-t))^{-\frac12} \, x] \qquad (t<T)\qquad\qquad
 $$
 ($\l  x$ d\'esigne le point de coordonn\'ees $\l x_1$, $\l x_2$, $\l
 x_3$.)"
 }

\ssk

 \ssk

However, at the end of the same paper in {\tt Acta mathematica},
Leray returned  once more to this similarity blow-up problem and
now his question is also truly remarkable \cite[p.~245]{Ler34}
(here (3.11) is the system \ef{NS1ss} below for the similarity
profiles in \ef{NS2}):

 \ssk

 \ssk

\com{
 "{\em Remarque:}  Si le syst\`eme (3.11) poss\`ede une solution non nulle $U_i(x)$ cette
  solution permet de construire un exemple tr\`{e}s simple de solution turbulente c'est le vecteur
$U_i(x,t)$ \'egal \`{a}
 \be
 \label{Ler00}
  \tex{
  [2\a(T-t)]^{-\frac 12} \,\, U_i \big[\,[2 \a(T-t)]^{-\frac 12}\, x\, \big]
  \quad \mbox{\underline{pour $t<T$ et \`{a} o pour $t>T$;}}
  }
  \ee
  il existe une seule \'epoque d'irr\'egularit\'e: $T$."
  }

\ssk

\ssk

 Therefore, Leray in \ef{Ler00} posed {\em both} principal questions on
existence of  {\bf self-similar} solutions of {\em blow-up
backward type} for $t<T$ and of the standard {\em forward type}
for $t>T$, which are naturally supposed to ``coincide" (in which
sense? -- in general,  a difficult question of extended semigroup
theory) at the unique singularity point $t=T$. Thus, this is a
principal setting not only for a similarity way of formation of a
blow-up singularity as $t \to T^-$, but also for self-similar
continuation of the solution for $t>T$, i.e., beyond blow-up time,
when it becomes again regular and bounded. Even for simple
parabolic reaction-diffusion equations such as \ef{FK1} and
\ef{SEHam1}, though such ``peaking" blow-up self-similar solutions
are known, a theory of such an incomplete blow-up is far away from
being well-understood and contains a number of open problems; see
references and results in \cite{GV97} and in survey \cite{GVaz}
(as usual, further  recent progress already achieved in this
direction can be then traced out {\em via} {\tt MathSciNet}). For
the Cahn-Hilliard equation \ef{CH1}, this is done in \cite{GalJMP}
in the lines of Leray's scenario \ef{for1}.


\ssk

Thus,
  substituting (\ref{NS2}) into (\ref{NS1}) yields for  ${\bf U}$ and $P$
   the following
  ``stationary" system  (3.11) in \cite[p.~225]{Ler34}:
 \be
  \label{NS1ss}
   \tex{
\frac 12 \, {\bf U}+  \frac 12 \, (y \cdot \n) {\bf U}+ ({\bf U}
\cdot \n){\bf U}=- \n P + \D {\bf U}, \quad
  \di {\bf U} =0 \inB \re^3.
  }
  \ee

During last twelve years, a number of enhanced negative answers
 concerning existence of such non-trivial similarity patterns (\ref{NS2}),
  (\ref{NS1ss}) were
 obtained.
 The key ingredient  \cite{Nec96} of such nonexistence proofs is the
  Maximum Principle\footnote{
  The
  fact that $\frac 12 \,|\uu|^2+p$ satisfies  the MP for
  the stationary NSEs is well-known; see e.g., \cite{Gil74} and more references in
  \cite{Nec96}.}:
 \be
 \label{MP1}
 \tex{
 \Pi=
 \frac 12\,|\UU|^2 + \frac 12\, y \cdot \UU+P
 \quad \mbox{satisfies} \quad
 - \D \Pi + (\UU+ \frac 12\, y) \cdot \n \Pi =-|{\rm curl}\, \UU|^2 \le 0;
  }
  \ee
 see further details in \cite{Chae07,  Mil01, Tsai98},
 and the advanced and  negative nonstationary PDE answer in \cite{Hou07}.
  Let us note an
  existence result in \cite{Dong07N4} for $N=4$.

Regardless the nonexistence of the similarity blow-up \ef{NS2} and
especially in connection with,  the following Leray conclusion
deserves the attention \cite[p.~224]{Ler34}. For the function
 \be
 \label{VV1}
  \tex{
  V(t) = \sup_{x} \, |\uu(x,t)| \forA t < T:
  }
  \ee


\ssk

\com{"{\em Un premier caract\`ere des irr\'egularit\'es:}~ Si une
solution des \'equations de Navier devient irr\'eguli\`ere \`a
l'\'epoque $T$, alors $V(t)$ augmente ind\'efiniment quand $t$
tend vers $T$; et plus pr\'ecis\'ement:
 $$
 (3.9) \qquad\qquad\qquad\qquad\qquad\qquad V(t)> A \sqrt{ \frac \sigma {T-t}} \, .\,\,\,"
 \qquad\qquad\qquad \qquad\qquad
  $$
   }

\ssk


Thus, Leray was the first who derived this estimate from below for
a solution to have blow-up at $t=T$.  This estimate led Leray to
pose the question on blow-up similarity solutions \ef{NS2} (not
just the fact that this would be a simple dimensional way for the
NSEs to develop a finite-time singularity as $t \to T^-$ with an
analogous extension for $t>T$).
 Indeed, otherwise, it was shown that the integral in the
estimate on p.~223 for the successive approximations
$\{u^{(n)}_i(x,t)\}$\footnote{This Leray's integral approach can
be characterized as  a forerunner type for a number of later
famous estimates via ``control of vorticity" including that by
Beale, Kato, and Majda \cite{Bea84} for Euler's equations (on
recent extensions, see \cite{Zhou08}).}
 $$
 \tex{
 V^{(n)}(t) \le \var(t) \whereA \var(t) \ge A' \int\limits_0^t
 \frac{\var^2(t')}{\sqrt{\sigma (t-t')}} \, d t' + V(0)
 }
 $$
 will never exhibit a necessary divergence (more precisely,
  here a Dini--Osgood-type integral condition is supposed to
 occur). In other words, Leray's  estimate (3.9) above then suggests that a
 {\em fast-type blow-up}, which is not self-similar (i.e., of Type II), is the only
 possible.

 We devote a notable part of the present paper to discussing such type of fast blow-up that is
 unbounded in the rescaled variables introduced  in \ef{NS2} and \ef{ll2} below.


 Thus, nonexistence \ef{Sim} of Leray's similarity solutions (\ref{NS2})
and other related local types of self-similar blow-up
 is a definite step towards better understanding of
the singularity nature for the Navier--Stokes equations.
 This does not settle the problem of singularity formation,
 since there might be other
 ways for (\ref{NS1}) to create singularities as $t \to T^-$ rather than the purely self-similar
 scenario
  (\ref{NS2}); see the monograph \cite{Maj02} for details.
   Other concepts of such a  multiplicity  are discussed below.

\ssk

Concerning blow-up of infinite energy solutions, consider the
strain field ({\em q.v.}
\cite{Nak04}):
 \be
 \label{uuu1}
  \uu=(-\zeta(t)x+u(x,y,t), - \zeta(t)y + v(x,y,t), 2 \zeta(t)z)
  \whereA \zeta(t)=\|\o(\cdot,t)\|_\iy,
   \ee
   and the vorticity is structurally associated with famous {\em Oseen's
   vortex} (1910)
    \cite{Os10},
   \be
   \label{uuu2}
    \tex{
   \o(r,t)= \frac 5{2(T-t)} \, {\mathrm e}^{- \frac{r^2}{T-t}}
   \quad (r^2=x^2+y^2), \quad \mbox{so that}
    \quad {\rm curl} \, \uu=(0,0,w(r,t)).
  }
       \ee
The additional velocity field associated with an analogous
vorticity in the cylindrical coordinates $\{r,\var,z\}$
is \cite{Moff00} (here $\Gamma>0$ is a constant)
 \be
 \label{moff1}
 \uu=(0,v(r,t),0),\,\,\, \mbox{where} \,\,\,
 \tex{
 v(r,t)= \frac \Gamma {2 \pi r}\, \Big(1- {\mathrm e}^{-
 \frac{r^2}{4(T-t)}}\Big) \quad \big({\mathbf
 \o}=(0,0,\o(r,t))\big).
  }
  \ee



\subsection{On blow-up in  Euler's equations}
 \label{S2.3}

  A similar nonexistence of self-similar\footnote{It seems that
  first strongly physically motivated arguments in favor of blow-up phenomena for Euler's equations were due
   to Onsager \cite{Ons49} (1949).} and some other types of blow-up has been obtained
  for  the {\em Euler equations} \ef{ee1} in $\re^3$
   (derived  by L.~Euler in the middle of the eighteenth century
 \cite{Euler55}\footnote{In published
form the incompressible equations appeared only in 1761, while  a
preliminary version was presented to Berlin Academy in 1752; see
the full history in \cite{Chr07}.}; see the most recent survey by
Constantin \cite{Con07} for a modern mathematical activity
exposition around);
 see \cite{Chae07, Deng06, Hou06},
  and references therein.
Local existence of smooth solutions for \ef{ee1} has been known
from 1920s; see Lichtenstein \cite{Lich25} (and \cite{Bar07} for a
full history);  see also later results of the 1970s by Kato
\cite{Kato72}.
  Existence of global classical solutions
 is still open\footnote{\,``The blow-up problem for the Euler equations is a
  major open problem of PDE theory, of far greater physical importance than the blow-up problem for
   the Navier--Stokes equation, which of course is known to nonspecialists because it is a Clay
   Millennium Problem," \cite[p.~607]{Con07}.}. Pelz \cite{Pelz01} and Gibbon \cite{Gib07} contain
    interesting  surveys on various mathematical,  symmetry,
 and more physical ideas concerning possible
 blow-up scenarios including numerical aspects (for the latter, see  Kerr
 \cite{Kerr06} for completeness concerning recent discussions on numerical blow-up issues).

Infinite energy solutions of \ef{ee1} do
 blow-up, for instance, according to the following separable solution in the cylindrical coordinates
 \cite{Gib03}:
 \be
 \label{os11}
  \tex{
  \uu=(u(r,t),0,z \g(r,t)) \whereA
\g(r,t)= - \frac{{\mathrm e}^{- r^2}}{T-t}, \quad u(r,t)= \frac
1{2r} \big( \frac{1-{\mathrm e}^{- r^2}}{T-t}\big).
 }
 \ee
Blow-up of more general solutions of this form
$\uu(x,y,z)=(u(x,y,t),v(x,y,t),z \, \g(x,y,t))$ was studied in
detail in \cite{Con00}, where, in particular, non-power blow-up
rate was observed,
 $$
 \tex{
 \g \sim  \frac 1{(T-t)|\ln(T-t)|} \asA t \to T^-.
  }
  $$
 Notice also that,
in a bounded (interior or exterior) domain in rescaled variables,
Euler's equations (\ref{ee1})
  admit non-trivial similarity solutions \cite{He04, He07}
  \big(cf. (\ref{NS2}) for $\a=\b=\frac 12$\big)
 \be
 \label{NS2E}
  \tex{
  \uu(x,t)= \frac 1{(T-t)^\a} \, {\bf U}(y), \,\,\, y=\frac
  x{({T-t})^\b} \,\,\, \mbox{in the range} \,\,\,
  \b \in\big[ \frac 25,1\big] \andA \a+\b=1.
  }
\ee
  In the original spatial $x$-variables, this self-similar blow-up is
 supported by boundary conditions for a domain that shrinks into a
 point as $t \to 1^-$. On the other hand,
 conditions (looking still rather non-constructive) of pointwise
 (in $L^\iy$) blow-up for (\ref{ee1}) were introduced in
 \cite{Chae07N}.

We believe that some ideas of our blow-up swirl analysis can be
also applied to Euler's equations \ef{ee1}, and this is reflected
in \cite{GalJMP}.
  Moreover, in \cite{GalJMP},
  we perform
the asymptotic construction of
{\em blow-up patterns of the NSEs $\ef{NS1}$ converging as $t \to
T^-$ to similarity
 solutions of the EEs $(\ref{ee1})$}.
 However, in general, \ef{ee1} is a special
subject that, in several circumstances, demands different
approaches for constructing  families of blow-up patterns.


 \subsection{Singular blow-up set has zero measure}

  There
exists another classic direction
 of the singularity theory for the Navier--Stokes equations that
 was originated by Leray himself \cite{Ler34} (see also details in \cite{Esc03})
 and  in Caffarelli, Kohn, and Nirenberg
 \cite{Caff82}.
 It was shown that the one-dimensional Hausdorff measure of
 the singular (blow-up) points in a time-space cylinder is equal
 to zero and these are contained in a space-time set of the Hausdorff dimension $ \le \frac 12$.
   We refer to \cite{Gust07, Nes07, Ser07} for further development
 and references. In particular, at a given moment $t_0$, the admitted number of singular points
can be  finite; see \cite{Neus99, Ser01} and presented references.

  Incidentally, among other results including Leray's one in \cite{Ler34},
 a refined criterion   is obtained in \cite{Ser07}, saying
 that, if $T=1$ is the first singular (blow-up) moment for a solution
 $\uu(x,t)$ of (\ref{NS1}), then
  \be
  \label{ll1}
   \tex{
   \lim
   \limits
   _{t \to 1^-} \frac 1{1-t} \,  \int
   \limits
   _t^1 \int
   \limits
   _{\re^3} |\uu(x,t)|^3 \,
   {\mathrm d}x \, {\mathrm d}t= + \infty.
 }
  \ee
  The condition (\ref{ll1}) is consistent with
Leray--Prodi--Serrin--Ladyzhenskaya regularity  $L^{p,q}$ criteria
 and other more recent
  researches;  see
  key references, history, details, and results concerning this huge
  existence-regularity-blow-up
  business around
   the NSEs in
   \cite{Esc03,  Gala07, Lem07, Maj02, Qio07, Ser07,
    Str06,
   Yuan07}; \cite{Fan08} represents a modern panorama of such studies,
   which also commented that Ohyama's result in 1960 \cite{Oh60} was obtained
   before Serrin's one in 1962 \cite{Ser62}.
 (\ref{ll1}) is also associated with   Kato's class of unique
mild
  solutions (in $\ren$),
 \cite{Kato84};  see details and key references in \cite{Can04K, Gala07, Way05}.



\subsection{Leray rescaled variables}

  We perform the {\em nonstationary} scaling as in (\ref{NS2}), $T=1$,
   \be
   \label{ll2}
   \tex{
   \uu(x,t)=\frac 1{\sqrt{1-t}} \, \hat \uu(y,\t), \quad y = \frac x{\sqrt{1-t}}, \quad \t=
   -\ln(1-t) \to +\iy \,\,\, \mbox{as} \,\,\, t \to 1^-.
 }
 \ee
 This yields
 the
rescaled equations for $\hat \uu=(\hat u^1, \hat u^2, \hat u^3)^T$
and $P$,
 \be
  \label{S1}
   \tex{
    \hat \uu_\t= \D \hat \uu -
 \frac 12 \, (y \cdot \n) \hat \uu -\frac 12 \, \hat \uu-  (\hat \uu
\cdot \n)\hat \uu- \n P , \quad
  \di \hat \uu =0 \inB \re^3 \times \re_+.
  }
  \ee

 In particular, after scaling, (\ref{ll1}) takes the form
\be
  \label{ll3}
   \tex{
   \lim\limits_{\t \to +\infty} \,  {\mathrm e}^{\t}
    \int
    _{\t}^{+\iy} {\mathrm e}^{-s} \Big(
     \int
     _{\re^3} |\hat \uu(y,s)|^3 \,
   {\mathrm d}y \Big) \, {\mathrm d}s= + \infty,
 }
  \ee
so that, if  $t=1$ is singular, then  the solution of the rescaled
equations \ef{S1}
 must diverge (blow-up) as $\t \to +\infty$ in $L^3(\re^3)$.

\ssk

As a standard next step,  we exclude the pressure from the
equations \ef{S1},
 \be
  \label{ww1}
    \begin{matrix}
    \hat \uu_\t=\HH(\hat \uu) \equiv  \big(\BB^*- \frac 12 \, I \big)\hat \uu
 -  \mathbb{P}\,(\hat \uu
\cdot \n)\hat \uu
   \inB \re^3 \times \re_+, \ssk\ssk\\
   \mbox{where} \quad \mathbb{P} \vv=\vv - \n \D^{-1}(\n \cdot
   \vv) \quad(\| \mathbb{P} \|=1)\qquad\qquad\quad
  \end{matrix}
  \ee
  is the
  Leray--Hopf projector of $(L^2(\re^3))^3$ onto  the subspace
  $\{{\bf w} \in (L^2)^3: \,\, {\rm div}\, {\bf w}=0\}$ of solenoidal vector
   fields\footnote{This emphasizes the unpleasant fact that the NSEs are
    a {\em nonlocal} parabolic problem, so that a somehow full use of
 order-preserving  properties of the semigroup is illusive; though some ``remnants" of the
  Maximum Principle (cf. \ef{MP1}) for such flows may remain and actually appear from time to
  time.}. Another representation is
   $
   \mathbb{P} \vv=(v_1-R_1 \s,v_2-R_2 \s,v_3-R_3 \s)^T,
   $
   where $R_j$ are the Riesz transforms, with symbols
   $\xi_j/|\xi|$, and $\s=R_1 v_1+R_2 v_2+R_3 v_3$.
  We then first apply $\mathbb{P}$
 to the original velocity equation in \ef{NS1} and next
 use the blow-up rescaling \ef{ll2}. Using the fundamental
 solution of $\D$ in $\ren$, $N \ge 3$ ($\s_N$ is the surface area of the unit ball $B_1 \subset \ren$)
  \be
  \label{FF55}
   \tex{
   b_N(y)= - \frac 1{(N-2)\s_N} \, \frac 1{|y|^{N-2}} \whereA \s_N= \frac { 2
   \pi^{ N/2}}{\Gamma( N/2)},
    }
    \ee
the operator in \ef{ww1} is written in the form of Leray's
formulation
\cite[p.~32]{Maj02} 
 \be
 \label{HHH21}
  \begin{matrix}
\HH(\hat \uu) \equiv  \big(\BB^*- \frac 12 \, I \big)\hat \uu -
(\hat \uu \cdot \n )\hat \uu + C_3 \int\limits_{\re^3} \frac
{y-z}{|y-z|^3}\,\, {\rm tr} (\n \hat \uu (z,\t))^2\, {\mathrm d}z,
\ssk \\
 \mbox{where} \quad
 {\rm tr} (\n \hat
\uu (z,\t))^2= \sum_{(i,j)} \,  \hat u_{z_j}^i \hat u_{z_i}^j
\andA C_N= \frac 1{\s_N}>0.\quad
 \end{matrix}
 \ee


\subsection{Hermitian  spectral theory of the blowing up rescaled operator:
 point spectrum and  solenoidal polynomial
eigenfunctions}
 \label{SPer}

In fact, the rescaled equations \ef{ww1} (or the original ones
\ef{S1}) are truly remarkable.
 Writing first linear operators on the right-hand side of
(\ref{ww1}) in a divergent form,
 \be
 \label{mm1}
 \tex{
   \tilde \BB^* \hat \uu \equiv \big(\BB^*- \frac 12 \, I \big)\hat \uu =
 \D \hat \uu -
 \frac 12 \, (y \cdot \n) \hat \uu -\frac 12 \, \hat \uu
 \equiv \frac 1 {\rho^*}\, \n \cdot (\rho^* \n \hat \uu) - \frac 12 \,
 \hat \uu,
 }
 \ee
 where the weight is
$\rho^*(y)={\mathrm e}^{-\frac{|y|^2}4}$,
 we observe that the actual rescaled evolution is now restricted to the
 weighted $L^2$-space $L^2_{\rho^*}(\re^3)$, with the exponential
  weight $\rho^*(y)$.
 Here, $\tilde \BB^*= \BB^* - \frac 12 \, I$ is a shifted adjoint Hermite operator with
 the point spectrum \cite[p.~48]{BS}
  \be
  \label{bbb1}
   \tex{
   \s(\tilde \BB^*)= \big\{ \l_k= - \frac{k}2- \frac 12, \quad
   k=|\b|=0,1,2,...\big\} \quad (\mbox{$\b$ is a multiindex}),
    }
    \ee
    where each $\l_k$ has the multiplicity $\frac{(k+1)(k+2)}2$
     for $N=3$, or the binomial number $C_{N+k-1}^k$.
     The corresponding complete and closed set of eigenfunctions
     $\Phi^*=\{\psi_\b^*(y)\}$
is composed from separable Hermite polynomials.  Similar spectral
and eigenfunction properties  can be also attributed to {\em
 blow-up problems} for $2m$th-order parabolic PDEs such as \ef{ho1}; see e.g., \cite{Eg4}, where a
 correspondence between the ``blow-up" space $L^2_{\rho^*}$ and  $L^2_\rho$ for
 the global evolution, with $\rho^*=\frac 1 \rho$, is more clearly
 explained. The bi-orthonormality  holds:
  \be
  \label{bn1}
  \langle \psi_\b^*, \psi_\g \rangle=\d_{\b \g} \quad \mbox{for
  any} \quad \b, \, \g.
   \ee
Note another important for us property of Hermite polynomials:
 \be
 \label{bn100}
 \forall \, \psi_\b^*, \quad \mbox{any derivative} \,\, D^\g  \psi_\b^* \quad \mbox{is also
 an eigenfunction with $k=|\b|-|\g| \ge 0$}.
  \ee
 Recall that \cite{BS}
  \be
  \label{bn10}
 \mbox{polynomial set} \,\,\, \Phi^* \,\,\, \mbox{is complete and closed in \,
  $L^2_{\rho^*}(\re^3)$}.
 \ee

We need to consider eigenfunction expansions in the solenoidal
restriction
 \be
 \label{bn11}
\hat L^2_{\rho^*}(\re^3)= L^2_{\rho^*}(\re^3)^3\cap\{ \di \vv=0\}.
 \ee
Indeed, among the polynomials $\Phi^*=\{\psi_\b^*\}$ there are
many that well-suit the solenoidal fields. Namely, introducing the
eigenspaces
 $$
 \Phi_k^*= {\rm Span}\,\{\psi_\b^*, \,\, |\b|=k\}, \quad k \ge 1,
  $$
in view of \ef{bn100} $\di$ plays a role of a ``shift operator" in
the sense that
 \be
 \label{sh1}
 \di: \Phi^{*3}_k \to \Phi^*_{k-1}.
  \ee
We next define  the corresponding solenoidal eigenspaces as
follows (see also Section \ref{S2.5}):
 \be
 \label{Sol1}
  \tex{
  {\mathcal S}_k^*= \{\vv^*=(v_1^*,v_2^*,v_3^*): \quad {\rm div} \, \vv^*=0,
  \,\,v_i^* \in \Phi_k^*\}\whereA {\rm dim}\,\, {\mathcal S}_k^*=k(k+2);
  }
  \ee
see \cite{Gal02, Gal02A} and further references therein. Actually,
\cite{Gal02} deals with global asymptotics as $t \to +\iy$, where
the adjoint operator $\BB$ in \ef{DDD1} occurs. Since $\BB$ is
self-adjoint in $L^2_\rho(\re^3)$, several results from
\cite[Append.~A]{Gal02A} are applied to $\BB^*$ (cf. Section
\ref{S2.5}). For a full collection,
 see \cite{Brand07} for further asymptotic expansions and extensions of these ideas.

In particular,  those solenoidal Hermite polynomial
eigenfunctions of $\BB^*$ can be chosen as follows
\cite[p.~2166-69]{Gal02A}
(the choice is obviously not unique; normalization constants are
omitted):
 \be
 \label{HP1}
 \begin{matrix}
   \underline{\l_1=- \frac 12:} \quad
  \vv_{11}^*=\left[\begin{matrix}0\\-y_3\\y_2 \end{matrix} \right], \quad
\vv_{12}^*=\left[\begin{matrix}y_3\\0\\-y_1 \end{matrix} \right],
\quad \vv_{13}^*=\left[\begin{matrix}-y_2\\y_1\\0 \end{matrix}
\right]\,\,({\rm dim}\, {\mathcal S}_1^*=3); \quad\ssk\ssk
\\
   \underline{\l_2=- 1:} \,\,\,
  \vv_{21}^*=\left[\begin{matrix}4-y_2^2-y_3^2\\y_1y_2\\y_1y_3 \end{matrix} \right],
  \,\,
\vv_{22}^*=\left[\begin{matrix}y_1y_2\\4-y_1^2-y_3^2\\y_2y_3
\end{matrix} \right], \,\,
\vv_{23}^*=\left[\begin{matrix}y_1y_3\\y_2y_3\\4-y_1^2-y_2^2
\end{matrix} \right],
 \ssk\ssk\\
  \vv_{24}^*=-\left[\begin{matrix}0\\-y_1y_3\\y_1y_2 \end{matrix} \right], \quad
\vv_{25}^*=-\left[\begin{matrix}y_2y_3\\0\\-y_2y_1 \end{matrix}
\right], \qquad\qquad\qquad  \quad 
 \ssk\ssk\\
  \vv_{26}^*=\left[\begin{matrix}-y_2 y_3\\y_2y_3\\y_1^2-y_2^2 \end{matrix} \right],
  \,\,
\vv_{27}^*=\left[\begin{matrix}y_1y_2\\y_3^2-y_1^2\\-y_2y_3
\end{matrix} \right], \,\,
\vv_{28}^*=\left[\begin{matrix}y_2^2-y_3^2\\-y_1y_2\\y_1y_3
\end{matrix} \right]\,\,\,({\rm dim}\, {\mathcal S}_2^*=8),  \quad \mbox{etc.}
  \end{matrix}
 \ee

We need the following final conclusion.
By \ef{bn10}, the set of vectors $\Phi^{*3}$ is complete and
closed in\footnote{Note a standard result of functional analysis:
polynomials are  complete in any weighted $L^p$-space with an
exponentially decaying weight; see the analyticity argument in
Kolmogorov--Fomin \cite[p.~431]{KolF}.} $L^2_{\rho^*}(\re^3)^3$,
so that
 \be
 \label{bn12}
  \tex{
 \forall \, \vv \in L^2_{\rho^*}(\re^3)^3 \LongA \vv= \sum_{(\b)} c_\b
 \vv^*_\b, \quad \vv_\b^* \in \Phi^{*3}_k, \,\,\, k= |\b| \ge 1.
 }
  \ee
   It then follows
  from \ef{bn1}--\ef{sh1} that
   \be
  \label{bn10N}
 \mbox{polynomial set} \,\,\, \hat \Phi^*= \Phi^{*3}\cap\{\di \vv=0\} \,\,\,
 \mbox{is complete and closed in \,
  $\hat L^2_{\rho^*}(\re^3)$}.
 \ee
In what follows, we always assume that we  deal with ``solenoidal"
asymptotics involving eigenfunctions as in \ef{Sol1}.

\ssk

For Burnett equations \ef{NS1m}, as we have promised to go with in
parallel, the blow-up rescaling and elements of linear solenoidal
spectral theory are found in Section \ref{S3.2}.

\subsection{On a countable set of
quasi-periodic singularities: first formalities}

 Thus, according to the criterion (\ref{ll3}), the moment $T=1$ is not a singular (and hence
 regular) point,
   if the corresponding locally smooth solution of (\ref{ww1})
 does
  not blow-up as $\t \to +\infty$ in a suitable functional
 setting. Thus, the problem of global existence and uniqueness
 of a smooth solutions of the NSEs in $\re^3$
 reduces to nonexistence of blow-up in infinite time for the
 rescaled system (\ref{S1}) or \ef{ww1}. In such a framework, this problem falls into the scope of
  blow-up/non-blow-up theory for nonlinear evolution PDEs.

\ssk

Let us first discuss a simple corollary that follows from the
above spectral properties of \ef{mm1}.
 Since by assumption $T=1$ is the first blow-up point of $\uu(x,t)$, we study
 solutions of \ef{ww1} that are globally defined in $\re^3 \times
 \re_+$, i.e., do not blow-up in  finite $\t$.
 Moreover, the scaling \ef{ll2} implies that we are looking for
 orbits with very sharp $L^2$-divergence,
  \be
  \label{l22}
   \tex{
   \|\hat \uu(\t)\|_2^2 = c_1 \,{\mathrm e}^{\frac 12\, \t}(1+o(1))
   \asA \t \to + \infty \whereA c_1=\|u(\cdot,1)\|_2^2>0.
 }
 \ee
 Consider the energy identity for smooth rescaled solutions
  \be
  \label{l23}
   \tex{
   \frac 12 \, \frac{\mathrm d}{{\mathrm d}\t} \|\hat \uu(\t)\|_2^2 \,
   {\mathrm d}\t= -\| D\hat \uu(\t)\|_2^2 + \frac 14 \, \|\hat
   \uu(\t)\|_2^2.
 }
 \ee
Solving it together with \ef{l22} yields the following control
 of the gradient $D \hat
 \uu(\t)$:
 \be
  \label{l22D}
   \tex{
   \|\uu(\t)\|_2^2= \big(c_1+ 2\int\limits_\t^\iy  {\mathrm e}^{-\frac s 2}\|D \hat \uu(s)\|_2^2\,
   {\mathrm d} s\big){\mathrm e}^{\frac \t 2} \LongA
   \int\limits^\iy {\mathrm e}^{-\frac s 2}\|D \hat \uu(s)\|_2^2\,
   {\mathrm d} s < \iy.
 }
 \ee
In particular,  for any $\e>0$, the following measure is always
finite and satisfies:
 \be
 \label{mes1}
  \tex{
  {\rm meas} \, \big\{s \gg 1: \quad \| D \hat \uu(s)\|_2^2 \ge \e
  {\mathrm e}^{\frac \t 2} \big\} = o\big( \frac 1 \e\big) \asA \e
  \to 0^+.
  }
  \ee
 Indeed, in comparison with \ef{l22}, this shows certain
 ``degeneracy" for $\t \gg 1$ of the $L^2$-norm of the gradient $\D \hat \uu(\t)$
 relative to that of $\hat \uu(\t)$. At least, this shows that the gradient cannot have
 a uniform exponential growth as in \ef{l22} and should  be much slower.
 In general, this does not imply
 essential consequences, since, according to \ef{mm1}
 \be
 \label{mes21}
  \begin{matrix}
  \mbox{the actual blow-up evolution of $\hat \uu(\t)$ as $\t \to +\iy$}
\ssk\\
   \mbox{is
  restricted to the space $L^2_{\rho^*}(\ren)$, rather than $L^2(\ren)$.}
   \end{matrix}
   \ee
Therefore, in particular, the exponential $L^2$-divergence of the
orbit \ef{l22} on expanding spatial subsets as $y \to \iy$ can be
``invisible" in the  metric  $L^2_{\rho^*}(\ren)$ with the
exponentially decaying weight $\rho^*(y)={\mathrm
e}^{-{|y|^2}/4}$.

 Thus, this is the class of global orbits of interest,  and we are looking for the structure of
  these $\o$-limit
 sets, which, for such smooth orbits, are  defined via local uniform
 convergence.

As usual in dynamical system theory, one first discusses the case
when the orbit $\{\hat \uu(\t)\}$ approaches the simplest
invariant manifold being a point.
 Thus, assume that
  \be
  \label{ww2}
  \hat \uu(\t) \to \bar  \uu \asA \t \to + \infty,
   \ee
 where, for future use, we suppose convergence in $L^q(\re^3)$, with $q >3$
  (note that a standard topology is expected to be that of $L^2_{\rho^*}(\re^3)$
  or the corresponding Sobolev one $H^2_{\rho^*}(\re^3)$).
  Then $\bar \uu$ is necessarily a
 self-similar profile, so that by global and local nonexistence results
  \cite{Chae07,  Mil01, Nec96, Tsai98}
  $$
   \bar  \uu=0.
   $$

Note that the spectrum \ef{bbb1} of the linearized operator
implies that, with a proper control of the quadratically small
convection term (this is easy in
$L^2_{\rho^*}$ with the exponentially decaying weight), the
non-stationary small solutions $\hat \uu(\t)$ must satisfy
 \be
 \label{ww3}
 |\hat \uu(\t)| \sim O \big({\mathrm e}^{-\frac 12 \t}\big) \asA \t \to +
 \infty,
  \ee
  where $\l_0=-\frac 12$ is precisely the spectral gap in
  \ef{bbb1}. More precisely, as a new application of spectral theory, we have the following:

\begin{proposition}
 \label{Pr.Gl}
 Assume that, along a subsequence $\{\t_k\} \to +\iy$,
  \be
  \label{mm5}
  \hat \uu(y,\t_k) \to 0 \quad \mbox{uniformly in} \quad L^\iy \cap L^2_{\rho^*}.
    \ee
 Then $t=T$ is not a blow-up time for $\uu(x,t)$ (in other words, the singularity  is  removable).
  \end{proposition}

  \noi{\em Proof.} Consider the sequence of solutions
  $\{\hat \uu_k(y,s)=\hat \uu(y,\t_k+s)\}$ with vanishing initial data in $L^\iy$
  according to (\ref{mm5}).
  Using  and well-developed
  spectral properties  of the linearized operator $\BB^*$ in
  (\ref{mm1})
   defined in  $L^2_{\rho^*}(\ren)$,
  with   generalized
   Hermite polynomials as  a complete and closed set of eigenfunctions.
   Therefore, according to classic asymptotic parabolic theory
   (see e.g., \cite{Lun}), we conclude that for any sufficiently
   large $k$,
    \be
    \label{mm7}
     \tex{
    \hat \uu_k(y,s) \sim O\big({\mathrm e}^{-  \frac s2}\big), \,\,\, s \gg 1
     \LongA \hat \uu(y,\t) \sim O\big({\mathrm e}^{- \frac \t 2}\big), \,\,\, \t \gg
     1.
     }
     \ee
Overall, taking into account Leray's scaling (\ref{ll2}), this
yields ($T=1$):
 \be
 \label{ww4}  
 \tex{
 \uu(x,t) \sim (T-t)^{- \frac 12} O\big({\mathrm e}^{- \frac \t 2}\big)= O(1) \asA
 t \to T^-,
 }
 \ee
 so that $\uu(x,t)$ is uniformly bounded at $t=T$. $\qed$

 \ssk


 Thus, equation \ef{ww1} does not allow  stabilization to an
 equilibrium, since this corresponds to the no-blow-up case.
 Further, as usual in textbooks on dynamical systems, the next candidate
for being the corresponding  invariant manifold is a {\em periodic
orbit} of finite period $T_*>0$\footnote{The role of $\t$-periodic
motion for rescaled blow-up solutions of the Euler equations
\ef{ee1} was pointed out in \cite{Pom05} (however, the
nonexistence conclusion on p.~218 therein looks rather
suspicious).}. In \cite[p.~1218]{Ple01},
 this conjecture was connected with the study of
the {\em complex Ginzburg--Landau equation}
 \be
 \label{GL1}
  {\rm i} \,u_t+ (1-{\rm i} \,\e) \D u + (1+ {\rm i} \, \d)|u|^{2
  \s} u= f \inB \ren \times (0,T),
  \ee
 where
 $\e>0$, $\d \ge 0$ and
   $ \s> \frac  2N$.
  It was pointed out in \cite{Ple01} that
the model \ef{GL1}
 exhibits the same
scaling, similar energy control, and local semigroup
 theory\footnote{However, unlike \ef{NS1}, it was shown in
\cite{Ple01} and in other papers cited therein that \ef{GL1}
admits a lot of (countable families of?) blow-up similarity
profiles, and this makes it more analogous to the Cahn--Hilliard
model \ef{CH1}; see \cite{EGW1}, where a countable blow-up family
was detected even for $N=1$, and  \cite{GalJMP}, where Leray's
scenario was shown to apply.}; see also \cite[\S~4]{Maz05} for
some related estimates of such a periodic behaviour of an unknown
structure. It was pointed out in \cite{Yan99} (see also
\cite{Liu03}) that the CGLE such as \ef{GL1} can be a good PDE
system modelling regularity and other questions regarding the NSEs
\ef{NS1}. In particular, global existence of weak solutions can be
obtained analogously to Leray's proof. Incidentally, concerning
the singular set of blow-up points for \ef{GL1} of zero measure
\cite{Yan99}, it is known that this set is restricted to a bounded
domain, \cite{Liu03}.

Actually, this would mean proving for \ef{ww1} the {\em
Poincar\'e--Bendixson theorem} saying, essentially, that if a
rescaled  orbit
 of \ef{ww1} satisfying \ef{l23}, \ef{l22D}
does not stabilize to an equilibrium (and actually \ef{ww4} makes
this impossible), then it convergence to a simple closed
curve\footnote{Is there any hope that a nonstationary version of
the elliptic differential inequality \ef{MP1} applied to the
nonlocal parabolic PDE \ef{ww1} can rule out  at least some of
special quasi-periodic oscillations about the unique equilibrium
$0$? (seems, not sufficient, and no hope).}, i.e.,
 \be
 \label{ww5}
  \fbox{$
  \mbox{$\o(\hat \uu_0)$ consists of a $T_*$-periodic orbit $\Gamma_1$.}
 $}
 \ee
According to \ef{mes21}, the natural metric of convergence as $\t
\to +\iy$  in \ef{ww5} is then assumed to be that of
$L^2_{\rho^*}$, while the $L^2$-divergence \ef{l22} occurs on
subsets with $|y| \gg 1$, which does not affect the convergence. A
more clear discussion of such blow-up periodic orbits will be
postponed until Section \ref{S.7.3AH}, where a proper spectral
theory is under scrutiny.

 Of course, \ef{ww5} would be the best and very pleasant
 case.
Indeed, for general DSs of such complexity, \ef{ww5} is a very
difficult open  problem.
Hence,
 the first main principal difficulty is how to
 predict a possible structure of such blow-up
 ``periodic
orbits". Actually, this is one of our main  goals, and {\em the
blow-up angular swirl mechanism  \ef{ang2} in Section $\ref{STw}$
is a natural argument in support of the periodic motion \ef{ww5}.}

Continuing using the logic of standard dynamical system theory, we
next
 would have to assume that a periodic $\o$-limit set would  have been
 ruled out (i.e., being
nonexistent) for a given rescaled orbit $\{\hat \uu(\t)\}$. Then
we should, e.g.,  conjecture a countable set of other
possibilities of evolution geometrically related to  {\em tori} in
$\re^{n+1}$: for $n=2,3,4,...\, $,
 \be
 \label{cou1}
  \fbox{$
  \mbox{$\o(\hat \uu_0)$ consists of a quasi-periodic orbit
  $\Gamma_n$ driven by $n$ fund. frequencies,}
   $}
   \ee
 where  $n=1$ leads to \ef{ww5}. To get a quasi-periodic motion for $n=2$ on
 the invariant $\o$-limit set, we introduce in Section
 \ref{SLanFin} the idea of precessions of the swirl axis, which
 leads to extremely difficult and open mathematics. For $n \ge 3$,
 such a clear visual geometric interpretation of the scenarios
 \ef{cou1} is not that easy or straightforward. A
 spectral background for bifurcation of such patterns is
 difficult and  obscure (cf. Section \ref{S.7.3AH}).

  Eventually, under the assumption of
nonexistence of all of those types of  $\o$-limits in \ef{cou1},
passing to the limit $n \to \infty$ would then have led to a kind
of a {\em strange attractor} for the dynamical system \ef{ww1},
which can be an extremely complicated invariant manifold to be
proved to exist\footnote{The proof of existence of a robust
strange attractor for the $\fbox{$\mbox{E.~Lorenz}$}$ dynamical
system in $\re^3$ \cite{Lor63} (``$3 \times 1$", i.e.,
uncomparably easier), describing thermal fluid convection with
some relations to the NSEs, proposed in 1963 was declared by Smale
as one of the several challenging problems for the twenty-first
century (1998), which eventually took nearly four decades  to
complete; see \cite{Tuck99, Tuck02} and also \cite{Mis01}.}.
Recall that, for proving nonexistence of $L^\iy$-singularities for
\ef{NS1}, all these infinite number of possibilities have to be
ruled out (via a new energy/monotonicity/spectral, etc. control).


\subsection{Global similarity solutions defined for all $t>0$  do
exist.
 Solenoidal linearized patterns}
 \label{S2.5}

In contrast to blow-up, as usual  for typical parabolic
reaction-diffusion-absorption problems,
 global similarity solutions of \ef{NS1} without blow-up, i.e.,
 \ef{NS2} with $T-t \mapsto t>0$ occur more frequently and correspond to the following
 scaling and rescaled equations (the invariant of the scaling group involved, $y=x/\sqrt t$,
  is sometimes called the
  {\em Bolzman substitution}\footnote{Similarity
solutions were used by Weierstrass around 1870, and by Bolzman
around 1890; this rescaled variable $y$ in parabolic PDEs was
widely used by Sturm in 1836 \cite{St36} (and possibly even
before?)}):
   \be
   \label{ll2G}
    \begin{matrix}
   \uu(x,t)=\frac 1{\sqrt t} \, \hat \uu(y,\t), \quad y = \frac x{\sqrt t}, \quad \t=
   \ln t \to +\iy \,\,\, \mbox{as} \,\,\, t \to +\iy,
 \ssk\ssk\\
   \mbox{where} \quad
    \hat \uu_\t= \D \hat \uu +
 \frac 12 \, (y \cdot \n) \hat \uu +\frac 12 \, \hat \uu-  \mathbb{P}\,(\hat \uu
\cdot \n)\hat \uu
 \inB \re^3 \times \re_+.
 \end{matrix}
  \ee
Here we obtain the linear operator
 \be
 \label{DDD1}
  \tex{
 \tilde \BB = \BB -I,
  }
  \ee
  where $\BB$ is adjoint to $\BB^*$ in the metric of the dual space $L^2(\re^3)$.
   The spectrum of $\tilde \BB$, which
  has a self-adjoint (Friedrichs)  extension in $L^2_\rho$, where $\rho=
  \frac 1{\rho^*}$ \cite{BS}, is
  \be
  \label{Sp111}
   \tex{
    \s(\tilde \BB)= \{ \l_k= - \frac {k}2-1, \quad
    k=|\b|=0,1,2,...\}.
     }
     \ee

Non-trivial self-similar profiles, i.e., stationary solutions of
\ef{ll2G}, do exist and describe asymptotics as $t \to +\infty$ of
various sufficiently small solutions; see \cite{Can96, Can04K,
Brand08, Gal02, Gal02A, Giga89, Gru06, Planch98}. In
\cite{Planch98},  a simple criterion of asymptotic similarity form
for solutions was obtained. Note that Slezkin--Landau singular
stationary solutions \ef{Lan1} are self-similar.

In addition, by \ef{Sp111}, 0 is asymptotically stable, and this
makes it possible to construct fast decaying solutions on each 1D
stable manifolds with the asymptotic behaviour\footnote{We present
here only the first term of expansion; as usual in dynamical
system theory, other terms in the case of ``resonance" can contain
$\ln t$-factors ({\em q.v.} \cite{Ang88} for a typical PDE
application); this phenomenon was shown to exist
 for the NSEs in $\re^2$ \cite[p.~236]{Gal02A}.}
 \be
 \label{mn1}
  \tex{
 \uu_\b(x,t) \sim  \, t^{\l_k-\frac 12} \, \vv_\b \big( \frac x{\sqrt
 t}\big)+... \asA t \to \iy, \,\,\,\mbox{where} \,\,\, \vv_\b= {\vv_\b^*}F \in
 {\mathcal S}_k
   }
   \ee
 are solenoidal eigenfunctions of $\BB$ defined in Section \ref{SPer}. Namely, taking
    \be
    \label{kl1}
    \tex{
    \vv=(v_1,v_2,v_3)^T \in {\mathcal S}_k, \,\,\, v_i \in {\Phi}_k= \big\{\psi_\b=
    \frac{(-1)^{|\b|}}{\sqrt{\b !}}\, D^\b F(y), \, |\b|=k \big\},
 }
  \ee
  where $F$ stands for the rescaled Gaussian (see \ef{G1} with $N=3$),
we have that
 \be
 \label{kl2}
  \tex{
 \di \vv = (v_1)_{y_1} + (v_2)_{y_2}+ (v_3)_{y_3}= \di (\vv^* F)
 \equiv (\di \vv^*) F - \frac 12\, y \cdot \vv^* \, F.
  }
  \ee
 This establishes  a one-to-one correspondence between
 solenoidal eigenfunction classes ${\mathcal S}_k^*$ in \ef{Sol1}
 for $\BB^*$ and ${\mathcal S}_k$ in \ef{mn1}
 for $\BB$; see \ef{HP1}
 for the first  eigenfunctions
 $\vv_\b= \vv^*_\b F$. Therefore,
  ${\rm dim}\,\, {\mathcal
   S}_k=k(k+2)$, etc.;
see details and rather involved  proofs of the asymptotics
\ef{mn1} for $k=1$ and 2 in \cite{Gal02}. We will deal with
patterns such as \ef{mn1} later on.

 It
follows from \ef{HP1} that there exists the corresponding
eigenfunctions of $\tilde \BB$ with $\l_1=-\frac 32$ given by
(\ref{Sp111}). Hence, by scaling (\ref{ll2G}), for instance, there
exists the asymptotic pattern for $t \gg 1$ (the rate $O(t^{-2})$
in $L^\iy$ is thus sharp for $\uu_0\in \hat L^2_{\rho}(\re^3_+)$)
 \be
 \label{HP5}
 \tex{
 \uu_1(x,t) \sim \frac 1{t^2} \, \vv_1\big(\frac x{\sqrt t}\big) \whereA
\vv_1(y)= \frac 1{4\sqrt 2}\, (y_3,y_3,0)^T\, {\mathrm
e}^{-\frac{|y|^2}4}.
 }
  \ee
In the half space $x \in \re^3 \cap\{x_3>0\}$ with no slip
boundary condition $\uu|_{x_3=0}=0$, the decay rate of solutions
is also of interest (clearly, some of the polynomials in
(\ref{HP1}) are good for that); see \cite{Choe08} for recent
developments.

Note again that the calculus of solenoidal eigenfunction classes
look like being specially designed to suit the divergence-free
flows not only for both types of scalings, blow-up \ef{ll2} and
the global one \ef {ll2G}, but also for the Burnett equations
\ef{NS1m}, where  classes ${\mathcal S}_k$ can be defined for any
necessary $2m$th-order linear operators \cite{Eg4}. The adjoint
ones ${\mathcal S}_k^*$ then  are also composed from solenoidal
generalized Hermite polynomials only.





\section{First application of Hermitian spectral theory: Sturmian local
structure of zero sets of bounded solutions and unique
continuation}
 \label{SHerm}

\subsection{Nodal sets for the Stokes problem and NSEs}

Here we perform a first step towards the classification problem
\ef{Turb1}. Namely, we assume that at the point $(x,t)=(0,1)$ the
solution $\hat \uu(y,t)$ is uniformly bounded and is such that the
eigenfunction expansion of the corresponding rescaled function
satisfying \ef{ww2},
 \be
 \label{ex1}
 \tex{
  \hat \uu(y,\t)= \sum_{(\b)} \cc_\b(\t) \vv_\b^*(y) \whereA
\cc_\b \vv_{\b}^* = (c_1 v_{\b 1}^*,c_2 v_{\b 2}^*,c_3 v_{\b
3}^*)^T \in \hat L^2_{\rho^*}(\re^3),
  }
  \ee
 converges in $\hat L^2_{\rho^*}(\re^3)$, and moreover, uniformly
 on compact subsets. These convergence questions of polynomial series
  are  standard; see  \cite{Eg4,
 2mSturm}, where further references and  details are given.
 Then the expansion coefficients satisfy the following dynamical
 system:
  \be
  \label{ex2}
   \left\{
   \begin{matrix}
 \dot \cc_\b= \big(\l_\b-\frac 12 \big)\cc_\b
 + \sum_{(\a,\g)} d_{\a\g\b} \cc_\a \cc_\g
  \quad \mbox{for any} \,\,\, |\b|
\ge 0,
 \ssk\ssk\\
 \mbox{where} \quad
 d_{\a\g\b}=  -  \langle\mathbb{P}\,(\hat \vv_\a^* \cdot
\n)\hat \vv_\g^*, \vv_\b \rangle \quad \mbox{for all}\quad \a,\,
\g.\quad
 \end{matrix}
 \right.
 \ee
It is natural to assume that the quadratic sum on the right-hand
side converges for the given smooth rescaled solution $\hat
\uu(y,\t)$. Recall that, according to scaling \ef{ll2}, we deal
with bounded and uniformly exponentially small functions
satisfying
 \be
 \label{ssjj1}
 |\hat \uu(y,\t)| \le C \,{\mathrm e}^{- \frac \t 2} \quad \mbox{in} \quad \re^3 \times
 \re_+.
  \ee
The system \ef{ex2} is difficult for a general study, and, of
course, it contains the answer to the existence/nonexistence
problem, provided that $(0,1)$ is a singular point of the solution
$\uu(x,t)$. For regular points, it can provide us with a typical
classification of nodal sets of solutions. This kind of study was
first performed by Sturm in 1836 for linear 1D parabolic equations
\cite{St36}; see historical and other details in
\cite[Ch.~1]{2mSturm}.

\ssk

Thus, following these lines, we clarify local zero sets of
solutions of the NSEs at regular points.
  Assume that
 \be
 \label{ex4}
\uu(0,1)={\bf 0},
 \ee
 which can be always achieved by constant shifting $\uu(x,t)
 \mapsto \uu(x,t)-\uu(0,1)$.
 In this connection,
 recall that the first eigenfunctions with $\l_\b=0$
 \be
 \label{ex3}
 \vv_0^*(y)  \sim (1,1,1)^T, \,\,(1,1,0)^T,\,\,(1,0,1)^T,...\, ,
 \ee
  are the only ones that have empty nodal sets of some of  its components.
 Then, bearing in mind the blow-up scaling term $(1-t)^{-\frac 12} \equiv {\mathrm e}^{\frac \t 2}$
 in \ef{ll2}, we have to assume that
  \be
  \label{ex5}
\cc_0(\t)={\bf 0}\,\,\,\mbox{or}\,\,\,   \cc_0(\t) \to {\bf 0}
\asA \t \to +\iy \,\,\, \mbox{exponentially faster than ${\mathrm
e}^{-\frac \t 2}$}.
    \ee

\ssk

\noi\underline{\em Polynomial nodal sets for the Stokes problem}.
A first clue to a correct understanding of the DS \ef{ex2} is
given by the Stokes problem, i.e., without the nonlinear
convection term,
 \be
 \label{SP1}
  \uu_t=-\n p + \D \uu, \quad {\rm div} \, \uu=0.
   \ee
   Then \ef{ex2} becomes linear diagonal and is easily solved:
 \be
 \label{FF2}
  \tex{
  \dot \cc_\b= \big(\l_\b-\frac 12 \big)\cc_\b
  \LongA
  \cc_\b(\t)=\cc_\b(0) {\mathrm e}^{- \frac {(1+|\b|) \t}2} \quad
  \mbox{for any}
  \quad |\b|\ge 0.
  }
  \ee
 Therefore, according to \ef{ex1},
  all possible multiple zero asymptotics for the Stokes problem (its local
``micro-scale turbulence") is described by finite solenoidal
Hermite polynomials, and the zero sets of rescaled velocity
components also asymptotically, as $\t \to +\iy$ (i.e., $t \to
1^-$) obey the nodal Hermite structures.

\ssk

\noi\underline{\em NSEs}.
 Consider the full nonlinear dynamical system \ef{ex2}, which on
 integration is
 \be
 \label{FF4}
  \tex{
   \cc_\b(\t)= \cc_\b(0) {\mathrm e}^{- \frac {(1+|\b|) \t}2}-{\mathrm e}^{- \frac {(1+|\b|) \t}2}
   \int\limits_0^\t \sum_{(\a,\g)} d_{\a\g\b} (\cc_\a \cc_\g)(s) {\mathrm e}^{\frac {(1+|\b|) s}2}
   \, {\mathrm d}s.
    }
    \ee
It follows that the nonlinear quadratic terms in \ef{FF4}, under
certain assumptions, can affect the rate of decay of solutions
near the multiple zero. As usual in calculus, this indeterminacy
can be tackled by L'Hospital rule.

\ssk

Since we are mainly interested in the study of nodal structures of
solutions by using the eigenfunction expansion \ef{ex1}, we
naturally need to assume that it is possible to choose the leading
decaying term  (or a  linear combination of terms) in this sum  as
$\t \to +\iy$. Then obviously these leading terms will
asymptotically describe the Hermitian polynomial structure of
nodal sets as $t \to 1^-$.
 For PDEs with local nonlinearities, this is done in a
standard manner as in \cite[\S~4]{2mSturm}; in the nonlocal case,
this seems can cause technical difficulties.
 However, the DS \ef{ex2} looks (but illusionary) as being obtained
 from a problem with local nonlinearities. In other words, the
 nonlocal nature of the NSEs is hidden in \ef{ex2} in the
 structure of the quadratic sum coefficients $\{d_{\a,\g,\b}\}$,
 and this do not affect the nodal set behaviour for some class of
 multiple zeros.
We will check this as follows:

\ssk

We consider a ``resonance class" of multiple zeros. Namely, let us
assume there exist a multiindex
 subset ${\mathcal B}$ and a function ${\bf h}(\t) \to {\bf 0}$ such that
 \be
 \label{vv1}
  \begin{matrix}
  \cc_\b(\t) \sim {\bf h}(\t) \asA \t \to +\iy
  \quad \mbox{for any $\b \in {\mathcal B}$}, \quad\ssk\ssk\\
|\cc_\b(\t)| \ll |{\bf h}(\t)| \asA \t \to +\iy
  \quad \mbox{for any $\b \not \in {\mathcal B}$}.
 \end{matrix}
 \ee
 In other words, only the coefficients $\{\cc_\b(\t), \, \b \in  {\mathcal
 B}\}$ are assumed to define the nodal set via \ef{ex1}, and other terms
 are negligible as $\t \to +\iy$.
 Under the natural assumption of a strong enough convergence of the quadratic sums
 in \ef{ex2} (this can be expected not to be the case for singular
 blow-up points only), taking the ODEs from \ef{ex2} for each $\b
 \in {\mathcal B}$ yields, for $\t \gg 1$,
  \be
  \label{hh1}
   \tex{
 \dot \cc_\b= \big(\l_\b-\frac 12 \big)\cc_\b + o(\cc_\b)
 \whereA \cc_\b(\t) \sim {\bf h}(\t).
 }
 \ee
 Hence, the asymptotic balancing of these equations
must assume that as $\t \to + \iy$
 \be
  \label{hh2}
   \tex{
  \dot{\bf h} \sim \big(\l_\b-\frac 12 \big){\bf h}
 \LongA \cc_\b(\t) \sim {\bf h}(\t) \sim {\mathrm e}^{(-\frac k2-\frac
 12)\t} \andA |\b|=k,
  }
  \ee
 where we may omit lower-order multipliers. Thus, there exists a $k \ge 1$ such that
  $|\b|=k$ for
 any $\b \in {\mathcal B}$.
 One can see that for such ``resonance" multiple zeros, the
 nonlocal quadratic term in \ef{ex2} is not important.
 Thus, in the resonance zero class  prescribed by \ef{vv1},
as $\t \to +\iy$, on compact subsets in $y$, similar to Stokes'
problem,
 \be
 \label{ex6}
  \fbox{$
 \mbox{the nodal set of $\hat \uu(y,\t)$ is governed by some
 solenoidal {\bf Hermite} polynomials.}
  $}
  \ee

\ssk

Note that the conclusion that, locally, for any zero of finite
order at $(0,1)$,
  \be
  \label{ex111}
  \mbox{nodal sets of $ \uu(x,t)$ are governed by finite-degree polynomials}
   \ee
 is trivially true for any sufficiently smooth solution.
 Indeed, this follows from the Taylor expansion of such solutions
  \be
  \label{Tay44}
   \tex{
   \uu(x,t)= \sum_{(\mu,\nu)}  \frac {(-1)^{\nu}}{\mu! \, \nu!}
   \big(D^{\mu,\nu}_{x,t}\uu\big)(0,1)\, x^\mu \, (1-t)^\nu + {\mathbf R}(x,t),
 }
  \ee
  where ${\mathbf R}$ stands for a higher-order remainder.
  Translating \ef{Tay44} via \ef{ll2} into the expansion for $\hat
  \uu(y,\t)$ yields some polynomial structure, so \ef{ex111} is
  obviously true.
 Thus, the principal feature of \ef{ex6} is that the Hermite
 polynomials count only therein.

\ssk

Obviously, for the nonlocal problem \ef{ww1}, there exist other
non-resonance zeros. Indeed, let $(0,1)$ be a zero of $\uu(x,t)$
of a finite order $m \ge 1$, i.e., as $x \to 0$,
 \be
 \label{mm111}
  \uu(x,1) \sim x^\s, \quad \mbox{with} \quad |\s|=m.
   \ee
We now use the following expansion:
 \be
 \label{mm2}
  \tex{
   \uu(x,t)= \uu(x,1)- \uu_t(x,1)(1-t) + \frac 1{2!}\,
   \uu_{tt}(x,1)(1-t)^2+...\,,
   }
   \ee
   where, by \ef{ww1}, all the time-derivatives
   $D_t^\mu \uu(x,0)$ can be calculated:
    \be
    \label{mm3}
 \uu_t(x,1)= \D \uu(x,1)+ (\mathbb{P}(\uu \cdot \n)\uu)(x,1)
 \sim x^{\s-2} + (\mathbb{P}(\uu \cdot \n)\uu)(x,1),
  \ee
  with a natural meaning of $\D x^\s \sim x^{\s-2}$.
 If the nonlocal term is negligible here and for other
 time-derivatives, i.e.,
  $$
 \uu_t(x,1) \sim x^{\s-2}, \quad  \uu_{tt}(x,1) \sim x^{\s-4}, ...
 \, ,
  $$
  then according to \ef{mm2} this leads to a Hermitian structure of
  nodal sets. In fact, this corresponds to
the  pioneering zero-set calculus performed  by Sturm in 1836; see
his
 original computations in \cite[p.~3]{GalGeom}.

\ssk

In general,
 the nonlocal term in \ef{mm3} is not specified by a local
 structure of the zero under consideration, so, obviously,  can essentially
 affect the zero evolution. For instance, as a hint, we can have the
 following zero:
  \be
  \label{mm4}
   \tex{
   \uu_t(0,1)={\bf C} \not = 0\,\, \Longrightarrow\,\,
    \uu(x,t) \sim x^\s + (1-t) \sim {\mathrm
    e}^{-\t}(1+z^\s), \,\, z= \frac x{(1-t)^{1/m}},
     }
     \ee
 so this nodal set is governed by the rescaled
 variable $z$, which us different from the standard similarity one
 $y$ in \ef{ll2}. Of course, due to the nonlocality of the
 equation, many other types of zeros can be described. Actually,
 such non-resonance zeros can be governed by sufficiently
 arbitrary polynomials as the general expansion \ef{Tay44}
 suggests.


\ssk

  Finally, the proof that zeros of infinite order are
not possible (and, as usual in such Carleman and Agmon-type
uniqueness results, this occurs for $\uu \equiv 0$ only) is a
difficult technical problem; see an example in
\cite[\S~6.2]{2mSturm}. For analytic in $y$ solutions of the NSEs
 (see  references and results in \cite{Dong07, Zub07}),
 this problem is nonexistent, and then in \ef{ex6} the degree of
the solenoidal vector Hermite polynomials is always finite, though
can be arbitrarily large.

\ssk

Note another straightforward consequence of this analysis that
this gives the following conventional {\em unique continuations
result}: {\em let \ef{ex4} hold,  $(0,1)$ be a
resonance zero\footnote{Indeed, this is hard to check; for PDEs
with local nonlinearities, this assumption is not needed, so that
such a unique continuation theorem makes  full sense,
\cite{2mSturm}.}, and
at least one component of the nodal set of $\hat \uu(y,\t)$ does
not obey \ef{ex6}. Then}
 \be
 \label{ex7}
  \uu \equiv {\bf 0} \quad \mbox{\em everywhere}.
   \ee
Of course, this is not that surprising since the result is just
included in the existing and properly converging eigenfunction
expansion \ef{ex1} under the assumption \ef{vv1}.

For elliptic equations $P(x,D)u=0$, this has the natural
counterpart on {\em strong unique continuation property} saying
that nontrivial solutions cannot have zeros of infinite order; a
result first proved by Carleman in 1939 for $P=-\D + V$, $V \in
L^\iy_{\rm loc}$, in $\re^2$ \cite{Carl39}; see \cite{Dos05,
Tao08} for further references and modern extensions.


\ssk

Thus, this is the first application of solenoidal Hermitian
polynomial vector fields for regular solutions of the NSEs.
 We expect that, due to the DS \ef{ex2}, some ``traces" of such an analysis and Hermite polynomials should be seen
in the fully nonlinear study of $\hat \uu(y,\t)$ at the singular
blow-up  point $(0,1)$, where, instead of \ef{ex4}, we have to
assume that, in the sense of $\limsup_{x,t}$,
 \be
 \label{ex8}
 |\uu(0,1)|=+\iy.
  \ee

\subsection{Burnett equations}
 \label{S3.2}

For \ef{NS1m}, the blow-up scaling \ef{ll2} is replaced by
 \be
 \label{yyy}
  \tex{
   \uu(x,t)= (1-t)^{-\frac 34} \hat \uu(y,\t), \quad
   y= \frac x{(1-t)^{1/4}},
   }
   \ee
   so that
 the rescaled system \ef{ww1} takes a similar form
\be
  \label{ww1yy}
   \tex{
    \hat \uu_\t=\HH(\hat \uu) \equiv  \big(\BB^*- \frac 34 \, I \big)\hat \uu
 -  \mathbb{P}\,(\hat \uu
\cdot \n)\hat \uu
   \inB \re^3 \times \re_+.
 }
  \ee
 The spectral
 theory of the given here adjoint operator
  \be
  \label{NBN1}
   \tex{
   \BB^*=- \D^2 - \frac 14\, y \cdot \n \whereA
   \s(\BB^*)=\big\{\l_\b=- \frac{|\b|}4, \, |\b| =0,1,2,...\big\}
 }
    \ee
    with eigenfunctions being generalized Hermite polynomials
is available in \cite{Eg4}; a solenoidal extension in the same
lines is needed. Therefore, under the same assumptions, the
polynomial structure of nodal sets is guaranteed for the
corresponding Stokes-like and  Burnett equations (and for an
arbitrary $2m$th-order of the viscosity $-(-\D)^m \uu$ therein).

\section{Second application of blow-up scaling: on convergence to the
EEs}
 \label{SEE2}

We now present the results of a general application
 of another related  blow-up scaling  for establishing a connection
 between blow-up in the NSEs and singularities in the EEs.
  We use the blow-up scaling
 in the form of  \cite[\S~2]{GW2},
 which received other applications and extensions in
  \cite{GalEng, ChGal2m, G3}, etc.
   As usual, such a rescaling near blow-up time leads to
  {\em ancient solutions}
   in Hamilton's notation \cite{Ham95}, which has been  a typical technique of
   R--D theory; see various form of its application in
 \cite{SGKM, AMGV} and others. In \cite{G3}, this scaling
 technique applies to the NSEs and Burnett equations in $\ren$
 to present a simple treatment  of the
 corresponding Leray--Prodi--Serrin--Ladyzhenskaya regularity  $L^{p,q}$
 criteria and other estimates. Ancient solutions of the 3D NSEs
 allowed recently to get new non-blow-up results for axi-symmetric
 flows; see \cite{Koch07} and also Section \ref{S3.8}.

We begin with a definition, which settles an ``evolution" concept
of entropy solutions for the EEs. In blow-up R--D theory, such
concepts of ``entropy-viscosity" come from extended semigroup
theory, where proper (blow-up or singular) solutions are only
those, which can be obtained by regular approximations. For
problems with the MP, such extended semigroup theory leads to the
{\em unique continuation} (partially or completely unbounded) of
any blow-up solutions beyond blow-up time $t>T$,
\cite[\S~6.2]{GalGeom}.

The questions of the vanishing viscosity limits in the NSEs to get
the EEs are classical in fluid mechanics and lead in general to a
number of fundamental open problems; see a clear statement  of
such questions in \cite[p.~422]{Bar07}, where necessary topologies
of convergence are prescribed. For general (not necessarily
bounded) solutions, the estimates \ef{ss1} suggest the weak-$*$
topology of $L^\iy(\re_-,L^2(\re^3))$, while for bounded solutions
this can be  improved (the reason for the uniform boundedness is
associated with the blow-up scaling to be applied). However, we
are not obliged to deal with a specific convergence, especially
since the necessary minimal topology for in what follows is still
unknown.

\begin{definition}
 \label{DefEE}
A function  (``distribution") $\uu(\cdot,t)$
is said to be a bounded NS-entropy (i.e., Navier--Stokes-entropy)
solution of the EEs $(\ref{ee1})$ if it is obtained as
    the limit
 \be
 \label{wl1}
  \uu_k \to \uu \asA k \to \iy \quad \mbox{weak-$*$ in}
  \quad L^\iy(\re^3 \times \re_-)
  \ee
  of a sequence of bounded classical solutions $\{\uu_k(x,t)\}$,
 \be
 \label{wl2}
 |\uu_k(x,t)| \le 1 \andA \|\uu_k(\cdot,t)\|_2 \le C \forA k \ge 0,
  \ee
   of the NSEs
with vanishing viscosity coefficients,
 \be
 \label{wl3}
 \uu_k: \quad \uu_t + {\mathbb P}(\uu \cdot \n)\uu= \d_k \D \uu
 \whereA \d_k \to 0^+ \asA k \to \iy.
 \ee
\end{definition}

The given in \ef{wl1} weak-$*$ topology of convergence
 follows from the fact that, by \ef{wl2}, the sequence $\{\uu_k\}$ is bounded in
 $L^\iy(\re^3 \times \re_-)$.
 Note that we
do not specify in which sense then $\uu$ satisfies the EEs,
 \be
 \label{wl4}
\quad \uu_t + {\mathbb P}(\uu \cdot \n)\uu=0,
  \ee
 and  we even do not demand  $\uu$ to be any kind of weak solution of
 \ef{wl4}.
  Note that this is not
 always the case even for {\em nonlinear dispersion equations} (NDEs) such as
 \cite{GalNDE5, GPnde}
  \be
  \label{wl5}
  u_t=(u u_x)_{xx}, \quad u_t=(u u_{xx})_x, \quad u_t= u u_{xxx}, \andA u_t=-(u u_x)_{xxxx},
   \quad \mbox{etc.}
   \ee
 On the other hand, the regularity concepts are not always an
 option in singularity theory: it is known that sometimes even
 analytic extensions of  solutions are not a proper one, i.e.,
  proper (extremal) solutions have finite regularity; see
 an example in \cite[p.~140]{GalGeom}.

Thus, according to \ef{wl1}, the concept of the NS-entropy for the
EEs includes the (natural, indeed) way of regular approximations
of its solutions. In this sense, this is quite similar to
conservation laws theory, e.g., for the 1D Euler equation
 \be
 \label{wl6}
 u_t + u u_x=0 \inB \re \times \re_+,
  \ee
with $L^1$-data $u_0$. The parabolic approximation is as follows:
 \be
 \label{wl7}
  \tex{
  u= \lim_{\e \to 0^+} u_\e \whereA u_\e: \quad  u_t + u u_x= \e
  u_{xx}, \,\,\, u(x,0)=u_0(x),
   }
   \ee
with the convergence in $L^1$. Of course, classic entropy concepts
apply directly to \ef{wl6} (say, in the sense of Oleinik and
Kruzhkov), but the entropy regularization description \ef{wl7} is
completely self-consistent; see \cite{Sm} for details. Note that
conservation laws such as \ef{wl6}
 are natural zero-level models for the 3D EEs, where  both features
of the divergence and the $L^2$-control are available. Then
 Tartar--Murat's compensation compactness approach could be also
 natural, but,
unfortunately, seem not applicable to singularly perturbed EEs
\ef{wl3}.

Further ``distributional" properties of NS-entropy solutions of
the EEs are unknown, and these compose the core of the problem;
see below.

\subsection{The NSEs}

Thus, we assume that  there exist sequences $\{t_k\}
 \to T^- \le \iy$, $\{x_k\} \subset \ren$, and $\{C_k\}\to +\iy$ such that the
 solution $\uu(x,t)$ of \ef{NS1} becomes unbounded:
 \be
 \label{seq11}
 \sup_{\ren \times [0,t_k]} \, |\uu(x,t)| = |\uu(x_k,t_k)| = C_k \to + \infty
 \asA k \to \iy.
 \ee
 As in  \cite[\S~2]{GW2},
 we then perform the change
 \be
 \label{rvarl}
\uu_k(x,t) \equiv \uu(x_k +x, t_k+t) = C_k \ww_k(y ,s) \whereA x =
a_k y, \quad t = a_k^{2}s,
 \ee
where the sequence $\{a_k\}$ is such that the $L^2$-norm is
preserved after rescaling, i.e.,
 \be
\label{uqq} \|\uu_k(t)\|_2= \|\ww_k(s)\|_2 \quad \Longrightarrow
\quad a_k = C_k^{-\frac 23} \to 0.
 \ee
 Taking the NSEs in
 the nonlocal from as in \ef{ww1},
we then obtain the following rescaled equations for
$\ww=\ww_k(y,s)$:
  \be
  \label{NS1n}
 \ww_s + \d_k \,\mathbb{P} (\ww \cdot \n)\ww=
 \D \ww, \quad
   \mbox{where} \quad \d_k=C_k^{\frac 13}.
 \ee
  Next, after time shifting, $s \mapsto s-s_0$,
   with a fixed arbitrarily large $s_0 > 0$,
 the solutions and data satisfy the uniform bounds: for all $k \gg 1$
  \be
  \label{rr4}
  | \ww_k(s)| \le 1 \andA \|\ww_k(s)\|_2 \le C \quad \mbox{for all} \quad s \in [-s_0,0].
   \ee
 The principle (and obvious, otherwise global existence would
 trivially follow) fact is that
  \be
  \label{rrt1}
  \d_k= C_k^{\frac 13} \to +\iy \asA k \to \iy \quad (\mbox{actually meaning that}
  \quad L^2 \not \Rightarrow L^\iy),
   \ee
 so that \ef{NS1n} is a singularly perturbed problem to be
 analyzed as follows:

We divide the equation \ef{NS1n} by $\d_k$ and introduce the new
time $\bar s= \d_k s$ to get for $\ww_k=\ww_k(\bar s)$ the
equation
 \be
  \label{NS1nN}
   \tex{
 \ww_{\bar s} + {\mathbb P}(\ww \cdot \n)\ww=
 \frac 1{\d_k}\, \D \ww.
 }
    \ee
 In view of $(\ref{rr4})$, the sequence is converging   weak-$*$
  in $L^\iy(\re^3 \times \re_-)$:
 \be
 \label{ww1ZZ}
 \ww_k(\bar s) \rightharpoonup \bar \ww(\bar s) \asA k \to \iy
 \ee
 (the convergence also takes place in  better
 topologies).
 We then need to pass to the limit in the rescaled NSEs
 \ef{NS1nN}, where the right-hand side vanishes in the weak sense
 in view of the sufficient regularity \ef{rr4}.
 Concerning the quadratic convection term,
  in view of its divergence  \ef{Div1},
 one needs extra assumptions to get in the limit  a weak formulation of the
 resulting EEs (locally, for sufficiently smooth solutions, this
 is well-known due to pioneering results by Kato (1983), Temam and Wang,
 etc.,
 see survey \cite[\S~4]{Bar07},  \cite{Xia06}, and references therein).
 Using Definition \ref{DefEE},
as a standard conclusion  from \ef{NS1nN}, we obtain:

\begin{proposition}
 \label{Pr.EE}
 Under the above hypotheses, the following holds:

 \noi {\rm (i)} $\bar \ww$ in $(\ref{ww1ZZ})$ is a bounded ancient NS-entropy
 solution of the EEs
  \be
 \label{ww21}
  \tex{
 \bar \ww_{\bar s} + {\mathbb P}(\bar \ww \cdot \n) \bar \ww=0
 \forA \bar s\le 0,
  \quad  \sup_{y \in \re^3}|\bar \ww(y,0)|=1;
  }
  \ee

\noi {\rm (ii)} $\bar \ww$ is not a classic regular solution of
the EEs in $\re^3 \times \re_-$.
 \end{proposition}

Concerning (ii), let us mention that the smooth EEs \ef{ww21} are
a gradient (and symplectic) system with the positive definite
Lyapunov function
 \be
 \label{wl61}
  \tex{
  \frac{\mathrm d}{{\mathrm d}\bar s}\, \int  |\bar \ww|^2 =0 \quad (\le
  0),
   }
   \ee
   which is monotone on evolution orbits. Therefore, the
   omega-limit set in the topology of $C_{\rm loc}(\re^3)$ of smooth orbits consists of equilibria, on
   which $V(\bar \ww)= \int|\bar \ww|^2$ vanishes locally on compact subsets, so $\bar \ww=0$.
In other words, we use the fact that smooth $L^2$-solutions of EEs
decays to zero in time uniformly on compact subsets. Therefore,
such solutions cannot satisfy the last normalization condition in
\ef{ww21}. Actually, we do not need a deeper discussion here about
these rather obscure aspects of EEs theory, since, moreover,
writing \ef{NS1nN} without an extra time-rescaling, i.e., keeping
the variable $s$,
 \be
 \label{wl8}
   \tex{
\frac 1{\d_k}\, \ww_{ s} + {\mathbb P}(\ww \cdot \n)\ww=
 \frac 1{\d_k}\, \D \ww.
 }
    \ee
yields as $k \to \iy$ that $\bar \ww$ is a weak {\em stationary}
bounded NS-entropy solution of the EEs. Finally, one can replace
\ef{Pos} by the following formal:

  \begin{corollary}
 \label{CorEE}
If the only bounded stationary NS-entropy solution of the EEs with
the additional to $\ef{wl1}$ uniform convergence on compact
subsets
 is trivial, then $\ef{Pos}$
holds.
 \end{corollary}

 The convergence \ef{ww1ZZ} embraces {\em larger compact subsets} in
$x$  than the self-similar one according to \ef{NS2}:
 \be
 \label{xkk1}
  |x-x_k| =O\big( C_k^{- \frac 23}\big), \quad \mbox{while}
  \quad
  |x-x_k|_{(\ref{NS2})} = O( C_k^{-1}).
   \ee
Admissible solutions of the problem \ef{ww21} are  generally
unknown. Thus, \ef{ww1ZZ} shows a general relation between blow-up
in the NSEs and NS-entropy solutions of the EEs: {\em special
bounded singular solutions of the EEs must create and support
blow-up for the NSEs, which eventually must evolve as $t \to T^-$
on compact subsets that are smaller than those in \ef{xkk1} for
the EEs}. It general,
 the ``almost trivial" case (cf. Corollary \ref{CorEE})
 \be
 \label{tr1}
\bar \ww(y) = 0 \inB \re^3 \setminus\{0\}
 \ee
 cannot be excluded.
 Obviously, for \ef{tr1} to be treated seriously, one needs an extra
stronger ``$N$-topology" of convergence
 in \ef{wl1}, which should be adequate to the both NSEs and EEs
 and is unknown.
  On the other hand, in \cite{GalJMP},
 we present almost ``explicit" blow-up infinite
energy patterns, where nontrivial limits like \ef{ww1ZZ} take
place. Some principles of formation of Type II blow-up patterns
are discussed in Section \ref{S7.5}.

\subsection{Burnett equations}
This is similar, consult \cite{G3} for final estimates.

\section{Construction of blow-up twistors on a 2D linear subspace in $\re^3$}
 \label{STw}

In this section, we demonstrate  the fact coming from
reaction-diffusion theory that the Navier--Stokes equations
 can admit  blow-up behaviour generating a blow-up swirl \ef{ang1} at a stagnation point
  with accelerating and eventually
 infinite angular speed as $t \to T^-$.
 Such a blow-up behaviour  is  associated with the so-called
 {\em logarithmic travelling waves} in the tangential angular direction,
  which are group-invariant solutions that occur for
  some nonlinear diffusion-combustion PDEs.
 We explain such a behaviour in the cylindrical geometry, where the
  axis of rotation is fixed to be extended to more realistic
  spherical geometry in Section \ref{SLanFin}.
To this end, we will perform a partially invariant construction
(understood not in a standard way of invariance under a group of
transformations).


\ssk

Even on the 2D subspace $W_2$ given in \ef{lin11}, the  analysis
of existence of such blow-up twistors is very difficult and leads
to involved nonlinear systems that we are not able to study
rigorously\footnote{This blow-up swirl behaviour is naturally
attached to the Poincar\'e--Bendixson-type conclusion \ef{ww5}.}.
In view of this, we will show another hypothetical  way of
constructing
other related non-similarity patterns of unknown rotational-like
singularity nature in conjunction with
the multiplicity claim \ef{fk2}.

\ssk

We begin by noting that the Navier--Stokes equations \ef{NS1}
inheriting for  given data
 a strong spatial symmetry (for instance, {\em cylindrically axisymmetric} irrotational flows),
actually  reduce their ``effective" dimension by one (a very
formal and rough issue) and then
 get into the global existence and uniqueness case of
classical solutions for $N=2$. This idea goes back to
Ukhovskii--Yudovich \cite{Ukh68} and Ladyzhenskaya \cite{Lad68A}
concerning axisymmetric geometry; see  \cite{SerZaj07, Zaj07} for
proper detailed definitions of such global smooth flows and recent
results and \cite{Mah90} for helical symmetries.
Formally speaking, the present patterns with some features of
axisymmetry could be
   not
 complicated enough to generate  singularities with finite kinetic energy, so we
 will need further mental efforts to remove such a restriction to
 blow-up; see Section \ref{SLanFin}.
  On the other hand,
   global regularity of general axisymmetric
 flows is also not known\footnote{Though,
  it is indeed an advantage to know that,
  by \cite{Caff82}, singularities may occur
  at the $z$-axis only.},
   so that blow-up patterns
 may be possible even in this simplified geometry of the linear dependence on $z$, but, probably,
  after
 extra precession-type or other necessary manipulations with this axis of rotation; see below.


 \subsection{Navier--Stokes equations in cylindrical coordinates}

We introduce the standard cylindrical coordinates  $(r,\var,z)$ in
$\re^3$, where
 \be
 \label{c1}
 r^2=x^2+y^2 \andA
  {\bf e}_r=(\sin \var,\cos \var,0), \,\,
 {\bf e}_\var=(-\cos \var,\sin \var,0), \,\,
 {\bf e}_z=(0,0,1).
  \ee
Denote the corresponding velocity field and the pressure as
follows:
 \be
 \label{c2}
 \uu=(u_r,v_\var,w_z)=(U,V,W) \equiv U{\bf e}_r+V{\bf e}_\var
 + W{\bf e}_z
  \quad \mbox{and} \quad P,
  \ee
  where $V$ stands for the {\em swirl component} of velocity.
Then equations (\ref{NS1}) take the form
 \be
 \label{c3}
 \left\{
 \begin{matrix}
 \frac{\mathrm d}{{\mathrm d}t}\, U - \frac 1r \, V^2= -  P_r + \D U - \frac 2{r^2} \,
 V_\var - \frac 1{r^2}\, U, \quad\ssk\ssk\\
 \frac{\mathrm d}{{\mathrm d}t}\, V + \frac 1r \, UV=- \frac 1r \, P_\var + \D V + \frac
 2{r^2}\,
 U_\var - \frac 1{r^2} \, V, \ssk\ssk\\
 \frac{\mathrm d}{{\mathrm d}t}\, W=- P_z+ \D W,\qquad\qquad \qquad \qquad\,\, \qquad \ssk\ssk\\
 U_r+ \frac 1r\, U+ \frac 1r \, V_\var + W_z=0. \qquad \qquad\quad
 \qquad
 \end{matrix}
 \right.
 \ee
 Here, the full time-derivative $\frac{\mathrm d}{{\mathrm d}t}$ and the Laplacian $\D$ are
 given by
 \be
 \label{c4}
  \begin{matrix}
\frac{\mathrm d}{{\mathrm d}t}= D_t + U D_r+ \frac 1r \, V D_\var
+ W D_z, \qquad\quad\quad\ssk\ssk\\ \D= \D_2 + D_{zz} \equiv
D^2_{rr} + \frac 1r \, D_r + \frac 1{r^2} \, D^2_\var + D_{zz},
 \end{matrix}
 \ee
  $D_{(\cdot)}$ being  partial derivatives, and $\D_2$  the Laplacian
 in the polar variables $(r,\var)$.

\subsection{On axi-symmetric flows: towards global existence}
 \label{SAxi}

 A flow (or a vector field) is called {\em  axi-symmetric} if it is invariant under
 rotations in $\var$.
  The principal advantage of such flows is that the $V$-equation in \ef{c3} takes the form
  of a standard parabolic one,
   \be
   \label{cc4}
    \tex{
    V_t + UV_r +W V_z = \D_r V  - \big(\frac 1r \, U + \frac 1{r^2}\big)V,
     }
     \ee
 so that some of  the MP, Harnack's inequalities, Nash--Moser iterations, De Giorgi--Nash estimates,
  and other tools of classic parabolic theory can be expected to take a part.
   For instance, the standard change for the 3D Laplacian gives the divergent equation:
  \be
  \label{cc5}
   \tex{
    r V= \hat V \LongA \hat V_t +  U\hat V_r+W \hat V_z +\frac 2 r \, \hat V_r=\D_r \hat V ,
    }
    \ee
 and hence by  comparison, via the control on
 the parabolic boundary,
  \be
  \label{cc6}
   \tex{
   |\hat V(r,z,t)|\le C \LongA  |V(r,z,t)| \le \frac C r.
    }
  \ee
  A more involved use of the MP yields some important results; e.g.,
  to prohibit Type I blow-up \cite{Koch07, Chen07}:
   $$
    \tex{
   |\uu(x,t)| \le \frac C{\sqrt{T-t}} \LongA \mbox{no blow-up at
   $t=T^-$}.
   }
   $$
Note again that  axi-symmetric flows {\em with no swirl}, i.e.,
 with $V \equiv 0$, are regular; see
Ukhovskii--Yudovich \cite{Ukh68}, Ladyzhenskaya \cite{Lad68A}, and
more recent new results in \cite{SerZaj07, Zaj07, Koch07, Chen07}
(the regularity remains under the presence of helical symmetries
\cite{Mah90}).

\ssk

Thus, for axi-symmetric settings, the $V$-equation is simple,
\ef{cc4}, while the third $W$-equation in \ef{c3} gives the
pressure,
 \be
 \label{QQ1}
  \tex{
  P= \int_z^\iy \big(\frac{\mathrm d}{{\mathrm d}t}\,W- \D W\big)
  \quad \big(W\big|_{z=\iy}=0\big).
   }
   \ee
Finally, we calculate $W$ from the last ${\rm div}$-equation,
 \be
 \label{QQ2}
  \tex{
  W= \int_z ^\iy \frac 1r\,(r U)_r.
   }
   \ee
Plugging all this into the $U$-equation yields
 \be
 \label{QQ3}
 \begin{matrix}
 U_t+U U_r+ W U_z= \D U + \frac 1r \, V^2 - \frac 1{r^2}\, U\qquad
  \ssk\ssk\\
 + \int\limits_z^\iy \int\limits_z^\iy \big\{\big( \D  \frac 1r\,(r U)_r\big)_r- \big(\frac
 1r\,(r U_t)_r\big)_r-\big[ U \big( \frac 1r\,(r
 U)_r\big)_r\big]_r-
 \big[\frac 1r\,(r U)_r \frac 1r\,(r U)_{rz}\big]_r \big\}.\qquad
  \end{matrix}
  \ee
 As we have seen in \ef{ww1}, the first two integral linear terms
 on the right-hand side
 in \ef{QQ3} can be incorporated into the main derivatives $\D U$
 and $U_t$, since these have the necessary good signs, so that the
 equation reads
 \be
 \label{QQ4}
 \begin{matrix}
 (I+{\mathcal L}_1)U_t+U U_r+ W U_z= (\D +{\mathcal L}_2) U + \frac 1r \, V^2 - \frac 1{r^2}\, U\qquad
  \ssk\\
 - \int\limits_z^\iy \int\limits_z^\iy \big\{\big[ U \big( \frac 1r\,(r
 U)_r\big)_r\big]_r +
 \big[\frac 1r\,(r U)_r \frac 1r\,(r U)_{rz}\big]_r \big\},\qquad
  \end{matrix}
  \ee
 where ${\mathcal L}_1>0$ and ${\mathcal L}_2<0$ (in the metric of $L^2$)
  are pseudo-differential operators
 with easily computed symbols (these are directly
 related to the projector ${\mathbb P}$)
  \be
  \label{QQ5}
   \tex{
   \quad {\mathcal L}_1=- D_r\big(\frac
 1r\,D_r(r I)\big) \equiv -\D_r + \frac 1{r^2}\,I, \quad
 {\mathcal L}_2=
 \int_z^\iy \int_z^\iy D_r \big(\D \frac 1r\,D_r(r
I)\big).
 }
  \ee

Let us present some typical RD-like speculations concerning
blow-up. Assume that, in the $(U,V)$-system \ef{QQ4}, \ef{cc4}, a
blow-up occurs at the origin $r=0$, $z=0$ as $t \to T^-$. Then, in
view of the $L^2$-boundedness, this must be single-point blow-up.
Then, the linear terms in the system are not essential as $t \to
T^-$ in comparison with a number of quadratic ones. This also
related to the linear nonlocal terms, where integration over the
necessarily rescaled spatial variables diminishes their rates of
divergence. Therefore, the quadratic pointwise and integral terms
only in the above $(U,V)$-system can be responsible for blow-up.
 We then need to consider two cases:

 \ssk

 \noi\underline{\em Case I: nonlocal quadratic terms are negligible in $\ef{QQ4}$}.
 Then, in the limit $t \to T^-$, the $U$-equation becomes
 asymptotically pointwise in main terms, and the system reads:
\be
 \label{QQ6}
   \left\{
 \begin{matrix}
 U_t+U U_r+ W U_z = \D U + \frac 1r \, V^2,
  \ssk\ssk\\
 V_t+UV_r+WV_z= \D V- \frac 1r\, U V.
  \end{matrix}
   \right.
  \ee
It then follows by a standard comparison (barrier) arguments that
 the growth of all the localized solutions are controlled by the ODE
 system for supersolutions,
 \be
 \label{ccbb6}
  \left\{
  \begin{matrix}
 U_t= \frac 1r \, V^2,\,\,\,\,\ssk\\
 V_t= -\frac 1r\, UV
  \end{matrix}
   \right.
   \tex{
  \LongA \frac {{\mathrm d}U}{{\mathrm d}V}=- \frac V U  \LongA U^2+V^2 = C,
   }
   \ee
which thus assumes no blow-up at all.

\ssk

\noi\underline{\em Case II: $V$ is not essential for blow-up}.
Then, performing a preliminary passage to the limit as $\t=\t_k+s
\to \iy$ ($\t=-\ln(T-t)$), we arrive at the system with $V \equiv
0$, which is known to admit no blow-up.

\ssk

\noi\underline{\em Case III: quadratic nonlocal terms reinforce
$\frac 1{r}\, V^2$}. Thus, this is the only possible case.
However, a careful analysis of those two integral terms in
\ef{QQ4} shows that both have the wrong sign to do that.
 Without pretending to any rigorous conclusions, we
speculate about the sign of the integrals in \ef{QQ8}  using
 standard clues, which look naive, but very often give correct
 answers in many RD-type problems.

 Indeed,
integrating over  a small neighbourhood of the maximum point in
$U$ (say, $x=0$), where single-point blow-up occurs, these terms
have the signs of the following integrals (others have a similar
nature and can be checked out analogously):
 \be
 \label{QQ8}
  \tex{
   \sim - \int_z^\iy \int_z^\iy U_r \D_r U
   \andA
\sim -\int_z^\iy \int_z^\iy  \D_r U  \frac 1r \, U_z.
 }
 \ee
Thus, we assume that the main sign-dominant part in these
integrals is delivered by integration over sufficiently small
neighbourhood of the maximum point. Then, since
 there $U_r \le 0$ (we think that the internal area, where $U_r$ can be positive, is not dominated
 in the integral),
$\D_r U \le 0$, and $U_z \le 0$, these integrals cannot be
essentially positive on shrinking compact subsets to create an
additional new type of ``nonlocal" blow-up. In other words, the
quadratic integral operators have the tendencies to be
non-positive on such typical axi-symmetric flows and actually
  assure their extra stabilization.

 A difficult scrutinized analysis is necessary to fix such a non-blow-up
conclusion. However,  this is  a good sign for us, since we see
that, for blow-up, a special $\var$-acceleration as $t \to T^-$ is
crucially needed to create extra {\em positive} quadratic integral
terms to generate blow-up, where  the TW dependence \ef{ang2} to
be introduced is the simplest opportunity.

\subsection{Positive  integral quadratic $\var$-dependent
operators: a rout to blow-up}
 \label{SNAxi}

We briefly review such an opportunity. Thus,
\ef{QQ2} reads (we box $D_\varphi$-dependent operators)
\be
\label{QQ2N}
  \tex{
  W= \int_z ^\iy \big[\, \frac 1r\,(r U)_r +
 \fbox{$
  \frac 1r\, V_\var
 $}\,\,
  \big].
   }
   \ee
 This adds into the right-hand side of the $U$-equation
 \ef{QQ4} the following extra operators:
  \be
  \label{QQ4N}
   \begin{matrix}
   U_t+U U_r+ \frac 1r \, V U_\var+ W U_z= \D U + \frac 1r \,
   V^2
   -...
-  \int_z^\iy \int_z^\iy \big\{
 \fbox{$
 U\big(\frac 1r\,
V_\var\big)_r
 $}
  \qquad\qquad \ssk\\
+\, \fbox{$
   \frac 1r\, V\big(\frac
 1r\, V_\var\big)_\var
 + \big[ \frac 1r\,(r U)_r+ \frac 1r\, V_\var\big] \,
 \big[ \frac 1r\,(r U)_r+ \frac 1r,\ V_\var\big]_z
 \big\}_r.
  $} \qquad\qquad
  \end{matrix}
 \ee
 The $V$-equation, due to the pressure term $-\frac 1r\, P_\var$,
  also gets a new  complicated operator,
 \be
 \label{VVV1}
  \tex{
  V_t+UV_r+ \frac 1r \,V V_\var+ WV_z= \D V- \frac 1r\, U V
  \fbox{$
  -
   \frac 1r \, \int_z^\iy( \frac{\mathrm d}{{\mathrm d}t}\,
   W- \D W)_\var.
    $}
   }
   \ee
Therefore, plugging \ef{QQ2N} yields several $\var$-dependent
quadratic operators listed below:
 \be
 \label{c3N}
  \begin{matrix}
   V_t+UV_r+ \frac 1r \,V V_\var+ WV_z= \D V- \frac 1r\, U V-
  ... \ssk \\
 \fbox{$
  - \frac 1r\,\int_z^\iy \int_z^\iy
  \big\{
  U \big(\frac 1r\,(r U)_r + \frac 1r\,V_\var\big)_r
   + \frac 1r\, V\big[ \frac 1r\,(r U)_r+ \frac 1r\,
   V_\var\big]_\var
 $}
   \ssk \\
     \fbox{$
   + \big[ \frac 1r\,(r U)_r+ \frac 1r\,
   V_\var]\,\big[\frac 1r\,(r U)_r+ \frac 1r\, V_\var
   \big]_z\,\big\}_\var.
   $}
   \end{matrix}
   \ee

Thus, the control of those new swirl-type operators in \ef{QQ4N}
and \ef{c3N} becomes much more difficult and this settles the core
difficulty of the blow-up/non-blow-up problem for the NSEs.
 Now the $\var$-derivatives are essentially involved into
 evolution, including even the third-order one in the integral term
 containing $ \frac 1{r^2}\, VV_{\var\var\var}$.
 In particular, the simplest way to introduce the $\var$-dependence
 is the
  proposed  blow-up swirl mechanism in \ef{ang2}, which actually
  has a self-similar form and corresponds to periodic omega-limits.
Thus, it seems plausible that some of the integrals in \ef{c3N}
can be divergent enough to generate
 blow-up, unless all of them have the opposite negative sign, or
 their mutual interaction prevents singularities.
 Recall that each integration in $z$,
 according to the similarity variable $y$ in \ef{ll2} yields the
 multiplier $\sqrt{T-t}$, which is not powerful enough to
 compensate other singular scalings  in $(U,V,r)$ there.
 Note again that, since $U \to +\iy$, the
 regularizing anti-blow-up term $- \frac 1r \, UV$ in the
 $V$-equation must be defeated by some of the divergent nonlocal
 terms with a $\var$-swirl. As in R--D blow-up theory, all such similarity of
  approximate self-similar balances eventually lead to
 complicated open problems for nonlocal systems, which later on we will attack and explain by using
 other, more local, blow-up approaches.

On the other hand,
 for purposes of proving global solvability, the above
 $(U,V)$-system \ef{QQ4N}, \ef{c3N} (it is a simplified one) looks most reasonable. Then a careful step by step
 analysis of nonlinear mutual interaction of all the operators in \ef{c3N}
is unavoidable.

\subsection{Generalized von K\'arm\'an solutions on a partially invariant linear
subspace: first term of general expansion}
 \label{S4.2}

We look for velocity field and pressure of the form
 \be
 \label{c5}
 U=U(r,\var,t), \quad V=V(r,\var,t),
 \quad W = z \tilde W(r,\var,t) , \andA P=P(r,\var,t).
  \ee
This dependence models a class of dynamical {\em stretched} 3D
vortex flows including {\em Burgers' vortices} (1948)
\cite{Burg48}. The structure of such solutions also corresponds to
  the earlier (1921) classic {\em von K\'arm\'an
swirling flow solutions}
 of \ef{NS1}
exhibiting typical linear dependencies  on two independent spatial
variables
$x$ and $y$\,
  \cite{vonK21, vonK61},
 \be
  \label{VK1}
  \left\{
   \begin{matrix}
  u=f'(z)x-g(z)y,\quad\,\, \\
   v=f'(z)x+g(z)y,\quad\,\, \\
    w=-2 f(z),\qquad \qquad\,\, \\
   p=-2[f'(z)+f^2(z)],
   \end{matrix}
   \right.
  \ee
  where functions $f$ and $g$ satisfy a system of two nonlinear
  ODEs,
 \be
 \label{VK2}
 \left\{
 \begin{matrix}
 f'''+2ff''-(f')^2+g^2=0, \\
 g''+2 f g'-2f'g=0. \qquad\,\,
 \end{matrix}
  \right.
  \ee
 These solutions
have been  applied to various problems of fluid dynamics; see
Berker \cite{Berker}. Note that, unlike von K\'arm\'an similarity
solutions \ef{VK1} that induce  ODEs \ef{VK2}, the reduction
\ef{c5} leads to a more complicated system of nonlinear PDEs,
which anyway is simpler than the original model. Stationary
solutions of the Euler equations \ef{ee1} with a linear
$z$-dependence of $w$ as in \ef{c5} were also used by Oseen (1927)
\cite{Os27} (cf. the structure \ef{os11} for Oseen's vortex); see
more details in \cite{Gib03}.

Actually, (\ref{c5}) means that the vector field $\uu$ and $p$
belong to a 2D linear subspace,
 \be
 \label{c6}
 \uu, \,\, p \in W_2={\rm Span}\, \{1,z\}.
  \ee
  Then some of the equations \ef{c7} below (projections of \ef{c3} onto $W_2$)
   express the {\em partial invariance}
  of $W_2$ with respect to the nonlinear operators in \ef{c3};
  see \cite[Ch.~7]{GSVR} for further details and examples.
 Here the {\em invariance}
 of an $l$-dimensional linear  subspace $W_l$ under a given
   nonlinear operator $\AAA$ is understood in the usual  mapping sense,
   $$
   \AAA(W_l) \subseteq W_l.
   $$
   Partial invariance means that only a subset of $W_l$ satisfies
   this inclusion.
Similar 2D restrictions on $W_2={\rm Span}\, \{1,y\}$ exist for
the Navier--Stokes equations in dimension $N=2$; see
\cite[p.~34]{GSVR} and references therein. The corresponding
blow-up solutions are described in \cite[Ch.~8]{AMGV}.

It is key to note that the third velocity component $W$ in \ef{c5}
becomes the first term in the asymptotic expansion of general
solutions, for which
 \be
 \label{em1}
 \tex{
 W= z \tilde W + \sum_{(k \ge 2)} z^k \tilde W_k.
  }
\ee
   Then, as usual in asymptotic expansion theory \cite{Il92}, the
first term is governed by the most nonlinear system (to be
studied), while the rest of the coefficients solve ``linearized"
systems (that, in view of equations for other components, can be
also rather difficult, at least first ones). In general, with a
great luck, analytic expansions such as \ef{em1} might lead to
solutions that are spatially localized in $z$; see further
discussions below.

\ssk

Thus, substituting (\ref{c5}) into (\ref{c3}) yields
 \be
 \label{c7}
  \left\{
   \begin{matrix}
    U_t+ UU_r+ \frac 1r \, V U_\var - \frac 1{r} \, V^2=-P_r+
    \D_2 U - \frac 2{r^2}\, V_\var - \frac 1{r^2} \, U, \quad\ssk\ssk\\
    V_t + UV_r+ \frac 1r\, VV_\var + \frac 1r\, UV =- \frac 1r\,
    P_\var+\D_2 V + \frac 2{r^2}\, U_\var - \frac 1{r^2} \, V,
    \ssk\ssk \\
    \tilde W_t+ U \tilde W_r+ \frac 1r \, V \tilde W_\var + \tilde W^2= \D_2 \tilde W,\quad\qquad\qquad
    \qquad\quad\qquad
    \ssk\ssk\\
    U_r+ \frac 1r \, U+ \frac 1r \, V_\var +\tilde  W=0. \qquad \qquad
    \qquad\qquad\qquad\qquad\qquad
    \end{matrix}
  \right.
   \ee
The last equation in (\ref{c7}) easily determines $U$ in terms of
$W$ and $V$ (cf. \ef{QQ2N}),
 \be
 \label{c72}
  \tex{
   U= - \frac 1r \, \int\limits_0^r z \tilde W(z,\var,t) \, {\mathrm d}z
- \frac 1r \, \int\limits_0^r  V_\var(z,\var,t) \, {\mathrm d}z.
 }
 \ee

 As customary, the first equation in (\ref{c7}) is a pressure one
 that {\em a posteriori} is going to define the pressure.
 Excluding $P$ from first two equations yields the following PDE:
  \be
  \label{c71}
   \begin{matrix}
\big[U_t+ UU_r+ \frac 1r \, V U_\var - \frac 1{r} \, V^2 -\D_2 U +
\frac 2{r^2}\, V_\var + \frac 1{r^2} \, U\big]_\var\ssk\ssk\\
 =\big[r\big(V_t + UV_r+ \frac 1r\, VV_\var + \frac 1r\, UV-
\D_2 V - \frac 2{r^2}\, U_\var + \frac 1{r^2} \, V)\big]_r.
 \end{matrix}
 \ee
  It is not difficult to see that this awkward  equation is just a
 pseudo-parabolic PDE that also can be reduced to a more standard semilinear
 nonlocal form. Therefore, \ef{c71} emphasizes the fact that the
 NSEs \ef{NS1} admit a unique local semigroup of smooth solutions.

Thus, taking the third equation,
 \be
 \label{c72n}
  \tex{
 \tilde W_t+ U \tilde W_r+ \frac 1r \, V \tilde W_\var + \tilde W^2= \D_2 \tilde W \quad (\mbox{with
  the
  constraint
  (\ref{c72}))},
   }
   \ee
 we arrive at the system (\ref{c71}), (\ref{c72n}) for two unknowns
 $\tilde W$ and $V$. For convenience, we now return to the original
 system (\ref{c7}) for performing further blow-up scaling.

\subsection{Blow-up twistor variables and rescaled equations}

We introduce next the following blow-up rescaled independent
variables, where, for convenience, the blow-up time is reduced to
$T=0$:
 \be
 \label{c8}
  \tex{
 U =\frac 1{\sqrt{-t}} \, u, \quad  V =\frac 1{\sqrt{-t}} \, v, \quad \tilde
 W =\frac 1{(-t)} \, w, \andA
  P =\frac 1{(-t)} \, p.
   }
  \ee
 The rescaled dependent variables are given by:
  \be
  \label{c9}
   \tex{
    y= \frac r{\sqrt{-t}}, \quad \var=\mu - \sigma  \,\ln(-t),
    \andA \t=- \ln(-t),}
     \ee
where $\sigma  \not = 0$ is an unknown parameter, which  cannot be
scaled out and plays a role of a {\em nonlinear eigenvalue} for
the future ``stationary" problem.
 As usual in blow-up problems, the new time variable
 gets infinite,
 \be
 \label{c10}
 \t=- \ln(-t) \to + \infty \asA t \to 0^-.
  \ee
It is crucial that, according to (\ref{c9}), on compact intervals
in the rescaled angle $\mu$,
 \be
 \label{c11}
  \var= \mu-\sigma  \, \ln(-t) \equiv \mu + \sigma  \, \t \to \infty \asA \t \to
  +\infty \quad (\mbox{i.e., $t \to 0^-$}).
   \ee
 Therefore, in the rescaled angle variable, the dependence $\var=
 \mu+\sigma  \, \t$ reflects a {\em travelling wave} angular
 behaviour, where $\sigma $ is the {\em wave angular speed}.
The angular dependence in \ef{c11} corresponds to infinite
acceleration of rotation of all the velocity components about the
$z$-axis,
 and therefore we call special pattern  solutions (\ref{c8}), associated with scaling
 (\ref{c9}), {\em blow-up twistors}.

\ssk

\noi{\bf Remark 1: on non blow-up flows with spiral symmetry.}
Such solutions
  of the NSEs \ef{NS1} in the cylindrical coordinates with the usual TW-type angular
dependence
 \be
 \label{TW11}
 \varphi= \mu + \s t
  \ee
 are well known and  were studied by Bytev \cite{Byt72} in
 1972; see also the results of the group classification of the
 NSEs in \cite{Ames}.
  Unlike \ef{c8},
 the standard invariant of the group of translations associated
 with \ef{TW11} does not allow blow-up. Indeed, logTWs and simply
 TWs are different group-invariant solutions, and the former ones
  assume extra invariance relative
  to a group of scalings.

\ssk

\noi{\bf Remark 2: on tornado-type blow-up in complex NSEs.}
Tornado-type blow-up behaviour of the complex version of the NSEs,
for which blow-up was established in Li--Sinai \cite{Li08}, was
confirmed in \cite{Arn08} by numerical methods (though the
achieved numerical evidence concerning the structure of such
singularities still looks rather unsufficient).

\ssk

In view of \ef{Sim}, let us now check the nature of the scalings
in \ef{c8}, \ef{c9}.
 According to full similarity rescaling in
 \ef{NS2}, the $z$-variable is
  \be
  \label{zz1}
   \tex{
   \zeta = \frac z{\sqrt{-t}} \LongA  z \tilde W =
  \frac 1 {\sqrt{-t}} \, \zeta \, \hat w,
   }
   \ee
 and hence all these scaling factors are self-similar. However,
 \ef{Sim}
 is not a warning for us, since logarithmic angular TWs
do not belong to the framework of standard similarity patterns (in
fact,  these are related to {\em periodic orbits}).
 In addition, we are going to use
{\em nonstationary evolution} governed by the PDEs \ef{c7N}
including patterns without  a similarity stabilization.
 On the other hand,  for axisymmetric flows, it is known that
the blow-up must be of Type II, i.e., with the blow-up rate faster
than self-similar \cite{Koch07, Chen07}.



\ssk

Substituting (\ref{c8}), (\ref{c9}) into the general system
(\ref{c7}) yields the following nonstationary rescaled equations:
\be
 \label{c7N}
  \left\{
   \begin{matrix}
    u_\t + \frac 12\, y u_y +\frac 12\, u +\sigma  u_\mu + uu_y+ \frac 1 y \, v u_\mu - \frac 1{y^2} \, v^2=-p_y+
    \D_2 u - \frac 2{y}\, v_\mu - \frac 1{y^2} \, u,\,\, \ssk\ssk\\
    v_\t + \frac 12\, y v_y +\frac 12\, v +\sigma  v_\mu + u v_y+ \frac 1y\, v v_\mu + \frac 1y\, u v =- \frac 1y\,
    p_\mu+\D_2 v + \frac 2{y^2}\, u_\mu - \frac 1{y^2} \, v,
    \ssk\ssk \\
    w_\t+ \frac 12\, y w_y + w +\sigma  w_\mu+ u w_y+ \frac 1y \, v w_\mu + w^2= \D_2 w,\quad\qquad\qquad
    \qquad\quad\quad\,\,\,
    \ssk\ssk\\
    u_y+ \frac 1y \, u+ \frac 1y \, v_\mu + w=0.\qquad\qquad
    \qquad\quad\,\qquad
    \qquad\,
    \qquad\qquad\qquad\qquad\qquad
    \end{matrix}
  \right.
   \ee
 The pressure and the constraint (\ref{c72}) interpretation remain
 the same as for (\ref{c7}), so (\ref{c7N}) is a nonstationary nonlocal  system of two
 equations for unknowns $w$ and $v$. Recall that
  the 2-Laplacian $\D_2$ now takes the full form
   \be
   \label{Lap2}
    \tex{
 \D_2= D_{y}^2+ \frac 1 y \, D_y + \frac 1{y^2} D_\mu^2.
 }
  \ee



 \subsection{Remark: log-TWs in reaction-diffusion systems}

 Using blow-up {logarithmic travelling waves} comes as a fruitful idea
 from reaction-diffusion  theory; see examples in \cite[pp.~105, 308,
 411]{SGKM}. For instance,
 the classic quadratic {\em porous
 medium  equations with source and convection} in 1D (a canonical combustion
 problem with regional blow-up, \cite[Ch.~4]{SGKM})
 $$
  u_t=(u u_x)_x + uu_x+ u^2
  $$
admits the following blow-up logarithmic travelling waves:
 $$
  \tex{
 u(x,t)= \frac 1{(-t)} \,f(y), \quad y=x-  \sigma  \, \ln(-t)
 \,\, \LongA \,\,
 f+ \sigma  \,  f'=(f f')'+ f f' +f^2.
 }
  $$

The scaling group-invariant nature of such logTWs seems was first
obtained by Ovsiannikov in 1959 \cite{Ov59}, who performed a full
group classification of the nonlinear heat equation
 $$
 u_t=(k(u)u_x)_x,
  $$
  for arbitrary functions $k(u)$.
  In particular, such invariant solutions appear for the porous
  medium and fast diffusion equations for $k(u)=u^n$, $ n \ne 0$:
   $$
   \tex{
   u_t=(u^n u_x)_x \LongA \exists \,\,\, u(x,t)=t^{-\frac 1n} f(x+ \s \ln
   t)\whereA
    - \frac 1n \, f + \s f'=(f^n f')'.
    }
    $$
 Blow-up angular dependence   such
as in \ef{c11} was  studied later on in \cite{Bak88}, where the
corresponding similarity solutions for the reaction-diffusion
equation with source
 \be
 \label{an1}
 u_t=\n \cdot (u^\s \n u)+ u^\b
 \inB \re^2 \times (0,T) \quad (\b>1,\,\,\,\s>0),
  \ee
were indicated by reducing the PDE to a quasilinear elliptic
problem (it seems, there is no still a rigorous proof of existence
of such patterns). For parabolic models such as \ef{an1}, that are
 order-preserving
via the MP and
 do
not have a natural vorticity mechanism, such ``spiral waves" as $t
\to T^-$ must be generated by large enough initial data
specially ``rotationally"  distributed in $\re^2$. For the fluid
model \ef{NS1} with typical vorticity features, this can  be
different.


 \subsection{On some auxiliary properties of twistor structures}

As usual,
 these are stationary solutions of (\ref{c7N}), which are
 independent of the time-variable $\t$, i.e., in the original
 variables,
 \be
 \label{c81}
  \begin{matrix}
 \hat U(r,\var,t) = \frac 1{\sqrt{-t}}\, \hat u(y,\mu), \quad \hat V(r,\var,t) = \frac 1{\sqrt{-t}}\,
 \hat v(y,\mu),
 \ssk\ssk\\
  \quad  \hat W(r,\var,t) =\frac 1{(-t)} \,\hat w(y,\mu), \quad
  \hat P(r, \var,t) =\frac 1{(-t)} \, \hat p(y, \mu).\quad
   \end{matrix}
  \ee
 The rescaled profiles $\hat u$, $\hat v$, $\hat w$, and $\hat p$
 solve the corresponding stationary system,
\be
 \label{c7Ns}
  \left\{
   \begin{matrix}
     \frac 12\, y \hat u_y +\frac 12\, \hat u +\sigma  \hat u_\mu + \hat u \hat u_y+ \frac 1y \, \hat v \hat u_\mu
      - \frac 1{y} \, \hat v^2=-\hat p_y+
    \D_2 \hat u - \frac 2{y^2}\, \hat v_\mu - \frac 1{y^2} \, \hat u,\,\,\, \ssk\ssk\\
     \frac 12\, y \hat v_y +\frac 12\,\hat  v +\sigma  \hat v_\mu + \hat u \hat v_y+ \frac 1y\, \hat v \hat v_\mu +
     \frac 1y\, \hat u \hat v =- \frac 1y\,
    \hat p_\mu+\D_2 \hat v + \frac 2{y^2}\, \hat u_\mu - \frac 1{y^2} \, \hat v,
    \ssk\ssk \\
     \frac 12\, y \hat w_y + \hat w +\sigma  \hat w_\mu+ \hat u \hat w_y+ \frac 1y \, \hat v \hat w_\mu + \hat w^2=
     \D_2 \hat w,\quad\qquad\qquad\,\,\,
    \qquad\quad\quad\,\,\,
    \ssk\ssk\\
    \hat u_y+ \frac 1y \, \hat u+ \frac 1y \, \hat v_\mu + \hat w=0.\qquad\qquad
    \qquad\,
    \qquad\qquad\qquad\qquad\qquad\qquad\quad\,
    \end{matrix}
  \right.
   \ee

\ssk

We now present the full system for $\hat w$, $\hat v$.
  First, we perform the  reflection,
 \be
 \label{q0}
 w \mapsto -w.
 \ee
 Second,
 we have from the last equation in \ef{c7Ns}
 \be
 \label{q1}
 \tex{
 \hat u= \frac 1y \, \int\limits_0^y z \hat w\,{\mathrm d}z - \frac
 1y\,
 \int\limits_0^y \hat v_\mu \, {\mathrm d}z.
  }
  \ee
  Then the $\hat w$-equation reads
 \be
 \label{q2}
   \begin{matrix}
 \D_2 \hat w - \frac 12\, y \hat w_y -\hat w + \hat w^2-\sigma  \hat w_\mu
 \ssk\ssk\\
 - \frac 1y
 \, \big(\int\limits_0^y z \hat w\,{\mathrm d}z \big) \hat w_y +
 \frac 1y \, \big(\int\limits_0^y \hat v_\mu \, {\mathrm d}z \big)
\hat w_y- \frac 1y \, \hat v \hat w_\mu=0.
 \end{matrix}
  \ee
Thirdly, the $(\hat u, \hat v)$-equation (with  (\ref{q1})) takes
the form
 \be
 \label{q3}
  \begin{matrix}
 \big[ \frac 12\, y \hat u_y + \frac 12\,\hat u +\sigma  \hat u_\mu + \hat u \hat u_y +\frac 1y \,  \hat v
 \hat u_\mu - \frac 1{y} \, \hat v^2 - \D_2 \hat u + \frac
 2{y^2} \, \hat v_\mu + \frac 1{y^2} \, \hat u\big]_\mu \qquad\qquad \ssk\ssk\ssk\\
=
  \big[y\big( \frac 12\, y \hat v_y + \frac 12\,\hat v +\sigma  \hat v_\mu + \hat u \hat v_y + \frac 1y \,  \hat
  v
 \hat v_\mu + \frac 1{y} \, \hat u \hat v - \D_2 \hat v - \frac
 2{y^2} \, \hat u_\mu + \frac 1{y^2} \, \hat v\big)\big]_y.
 \qquad\qquad
 \end{matrix}
  \ee

  \ssk

Thus, we arrive at the system of two PDEs (\ref{q2}), (\ref{q3}),
with the nonlocal constraint (\ref{q1}), which, as we have seen,
actually means the presence of an extra first-order PDE.

\ssk

This system with the real parameter $\sigma  \not = 0$ is  indeed
complicated. We note that the Laplacian $\D_2$ in \ef{Lap2}
includes second-order derivatives in both radial $y$ and the
angular $\mu$ variables. The latter dependence is crucial for such
similarity twistors. We do not intent and do not plan to study
this system somehow rigorously, though some results, in
particular, about asymptotic distributions as $y \to +\iy$ can be
obtained and will be of importance in what follows.

\ssk

Obviously, in view of \ef{c5}, these similarity twistors as
special kind of blow-up vortices are not spatially localized in
the $z$-direction, so they have infinite energy always, and, for
their existence, a whole $\re^3$-space should be taken into
account (we skip at this moment their expansion meaning \ef{em1}).
But as typical in reaction-diffusion equations such as \ef{FK1},
\ef{FK11}, \ef{ho1}, and many others, such similarity or {\em
approximate} similarity blow-up structures can appear from local
finite and well-spatially-localized data. Here there appear hard
questions of their evolution stability (in which intermediate
sense?) or behaviour close to centre and/or stable manifolds
associated with their partial similarity space-time geometry.
These questions are addressed to the full non-stationary system
\ef{c7N}.

\ssk

As we have mentioned, in the case of nonexistence of such
``stationary" swirl patterns,
some time-dependent perturbations of these structures can lead to
various non self-similar blow-up twistors that evolve close to
invariant manifolds associated with the scaling \ef{c81}.
 Actually, precisely this {\em always} happens to the semilinear scalar
 curvature equation \ef{SEHam1} for $N<6$ and for the
 Frank--Kamenetskii equation \ef{FK1} for $N=1$ or 2.
Unfortunately,
 this will lead to much more difficult mathematics than for a simpler purely stationary
 similarity one (that does not look easy at all, of course).

\subsection{On extension beyond blow-up for $t>T$}

This is connected with Leray's blow-up scenario \ef{for1}. We do
not discuss such difficult issues here in any detail and just
mention that, as typical, the extension after blow-up is a
secondary question relative to the blow-up ones on the behaviour
as $t \to T^-$, which still remains mysterious.

\ssk

 Let us also note that,  concerning blow-up twistors with the
angular variable \ef{c11}, a natural way of extensions of possible
blow-up solutions is to use the  forward variable
 \be
 \label{fw1}
 \var= \mu + \hat \sigma  \, \ln t \forA t>0
  \ee
  and $-t \mapsto t$ in \ef{c8},
  where $\hat \sigma  \not = \sigma $,  in general
  (but ${\rm sign}\, \s= {\rm sign}\, \hat \s$ to keep the direction of rotation), is an extra  parameter of
  matching at $t=T=0$.
  Principles of proper matching of a blow-up flow
for $t<0$ with a more regular one for $t>0$ for \ef{NS1} are not
clear.
 For simpler nonlinear PDEs
such similarity matching  is possible; see \cite{GalJMP} for
\ef{CH1}.



 \subsection{Discussion: on auxiliary properties of blow-up
 twistors}
  \label{S4.8}

Let us look for some auxiliary properties of the ``stationary"
profiles $\hat u$, $\hat v$,  $\hat w$ as solutions of
\ef{q1}--\ef{q3}. To reveal  a first key feature,  assume for a
moment that these  profiles are
 of changing sign, so that
rescaling (\ref{c8}) would mean {\em non-uniform blow-up}, where
components diverge as $t \to 0^-$ to $\pm \iy$ on some subsets.

We next need to point out another curious possible feature of such
solutions that is of importance  for general understanding of the
nature of such singularities. Let us  demand that the similarity
profiles have {\em angular periodic} behaviour with zero  mean as
$y \to \iy$. For instance, assume for the component $\hat
w(y,\mu)$  that
 \be
 \label{q4}
  \tex{
  \hat w(y,\mu)= \frac {C(\mu)}{y^2}+...
  \asA y \to +\iy
  \whereA \int\limits_0^{2\pi} C(\mu)\, {\mathrm d} \mu=0.
  }
  \ee
Then, in view of the behaviour of the rescaled variables in
(\ref{c9}) near blow-up time, under some extra hypotheses, this
can imply, by the blowing up rotational behaviour, that
 \be
 \label{q5}
 W(r,\var,t)\rightharpoonup 0 \asA t \to 0^-,
  \ee
 in the weak sense on the corresponding subsets in the rescaled
 variables.
  Indeed, in view of the zero mean
 in (\ref{q4}), the equality (\ref{q5}) is then a manifestation of
 Riemann's Lemma from Fourier transform theory.
 As more usual and typical, under weaker assumptions, such
 solutions locally represent an oscillating ``building block" (see \cite{Brand04V})
  for $x \approx 0$, $t \to 0^-$
  in
 terms of the corresponding vorticity $\o =\n \times \uu$, where
 $\uu$ is then defined
 by the {\em Biot--Savart law},
  \be
  \label{vort1}
  \tex{
   \uu(x,t)= - \frac 1{4\pi} \, \int\limits_{\re^3} \frac{x-y}{|x-y|^3} \times
   \o(y,t) \, {\mathrm d}y  \inB \re^3 \times \re_+.
 }
 \ee

 The property
 \ef{q5} (and also relative to the vorticity, which is easier) then represents special type
  of singularities that are almost ``invisible" (in fact,
 efficiently ``nonexistent")
 close to the blow-up time in the natural integral (weak) sense.
 Recall that this could be a key feature since all the differential equations
 of fluid mechanics are derived from  kinetic equations
 with {\em integral} operators in collision-like terms by approximating typical
 integral kernels involved by kernels with pointwise supports
  (Grad's method in Chapman--Enskog expansions).
 Therefore, the integral, average meaning
 of coherent macro-structures and micro-singularities are of practical importance only.

Another principal question is as follows: can, in the
non-stationary setting, the ``zero mean oscillation property" such
as \ef{q5} in a neighbourhood of the origin affect and diminish
the asymptotic behaviour of the corresponding solutions as $y \to
\iy$ to get them into finite energy class? It seems that this is
not that essential for the present cylindrical geometry, but can
be key for the spherical one with the unique single point at the
origin; see Section \ref{SLanFin}.



\subsection{Discussion: on asymptotic stability of twistors}
 \label{S3.6CH}

We begin with the following first observation concerning
evolutionary significance of similarity scaling:

  Assume, not taken into account possible energy characteristics, that
 a suitable non-trivial
solution of the stationary system (\ref{q2}), (\ref{q3}) exists.
Then the next and a  more difficult step is to study its
asymptotic stability in the framework of the full non-stationary
system (\ref{c7N}).
 This stabilization phenomenon then becomes a difficult
open problem. Observe that there is no  chance to prove that
(\ref{c7N}) is a gradient dynamical system in any admissible
metric. Recall that the negative result   in \cite[Th.~1,
\S~2]{Hou07} establishing the self-similar ban
 is proved {\em under} the {\em a priori} assumption on
 convergence to a fixed rescaled profile (then necessarily
 this profile must be a stationary solution and hence it is zero by
 \cite{Nec96}); and finally, the exponential decay via spectral and semigroup
 characteristics of the rescaled infinitesimal generator $\BB^*$ in $L^2_{\rho^*}(\re^3)$
implies that the solution $\uu(x,t)$ is bounded at $t=0^-$.

In this connection, it is worth mentioning that, as
 known from reaction-diffusion theory,
regardless bad energy and other global properties of twistor
profiles, good rescaled solutions can stabilize to them in a local
topology, e.g., in view of interior regularity,
\be
\label{un1}
 \mbox{uniformly on compact subsets}
 \ee
  in the rescale
variables. Of course, uniform stabilization in $\re^3$ is then
impossible in view of the energy discrepancy.

\subsection{Discussion: on linearized blow-up patterns about  a constant equilibrium, (I)}
 \label{S3.7}

Assume  next that either (i) any suitable  non-trivial steady
profiles satisfying (\ref{q2}), (\ref{q3}) are nonexistent, or
(ii) there exists a non-trivial solution (both possibilities are
suitable for us). What kind of other behaviour do we then expect
of the nonstationary system (\ref{c7N})? In other words, is there
any hope to get a kind of entirely non-self-similar blow-up
twistor behaviour?

The idea of such ``linearized construction" is as follows. We
first need to fix a family of simpler local solutions or ``almost"
(say, slightly perturbed) solutions of the dynamical system under
consideration. As the next step, we linearize the flow about this
family and  use the  orbits on the corresponding stable or centre
manifolds to match them with the surrounding orbits of the
necessary behaviour and regularity. These matching procedures,
though  being rigorously justifies for some simpler parabolic PDEs
such as \ef{FK1} or \ef{SEHam1} (we have presented the
references), remain open for many other quasilinear and
higher-order parabolic reaction-diffusion models, and, surely,
this will be the case for \ef{NS1}.

We show how then a countable family of the so-called {\em
linearized blow-up patterns} (unlike the above similarity patterns
that are {nonlinear eigenfunctions}) can be constructed.
 Consider the nonstationary $w$-equation from (\ref{c7N}) bearing
 in mind the reflection (\ref{q0}),
 \be
 \label{w1}
  \tex{
   w_\t= \D_2 w - \frac 12\, y w_y -   w + w^2 - u w_y - \sigma   w_\mu  - \frac 1y \, v
   w_\mu.
   }
   \ee
We see that (\ref{w1}) admits the constant equilibrium
 \be
 \label{w2}
  \tex{
 w_*= 1, \quad \mbox{i.e., \,\, $\tilde W(r,\var,t) \equiv
 \frac 1{(-t)}$},
 }
  \ee
  which corresponds to {\em uniform global blow-up} as $t \to 0^-$ in the
  whole space. This is the easiest  exact solution, which we are going to
  linearize the flow about.
It is  known from reaction-diffusion theory (the proof is obvious
for systems with the MP) that in order to blow-up at the fixed
time $t=0$, the rescaled solution $w(y,\t)$ should always be
``sufficiently" close  (in a certain metric) to the constant
blow-up profile (\ref{w2}), since otherwise, essentially smaller
solutions must exhibit larger blow-up times, so these become
exponentially small as $\t \to +\infty$ after rescaling
(\ref{c8}); see a spectral justification later on.

Thus, we perform the standard linearization about (\ref{w2}) by
setting
 \be
 \label{w3}
  \tex{
 w= 1+Y
 }
 \ee
 to get the  ``linearized" equations (we have taken
 into account  the $u$-representation (\ref{q1})):
 \be
 \label{w4}
  \begin{matrix}
  Y_\t= \D_2 Y - y Y_y + Y - \sigma  Y_\mu - \frac 1y  \, v Y_\mu
  + \frac 1y\, \big(\int_0^y v_\mu\big) Y_y +
  \DD(Y), \ssk\\
  \mbox{where} \quad
 \DD(Y)= Y^2 - \frac 1y \, \big(\int_0^y zY\big) Y_y,
 \qquad\qquad\qquad\qquad\qquad
   \end{matrix}
   \ee
  is a quadratic perturbation as $Y \to 0$. Consider the radial part of the
   linear operator\footnote{The non-radial eigenvalue problem with $\sigma  \not = 0$
   is studied in Section \ref{Snonrad}.}
    in (\ref{w4}) excluding at this moment  its tangential
  angular first- and second-order operators including that in $\D_2$,
 \be
 \label{w6}
  \tex{
 \tilde \BB^*= \D_2 - y D_y +  I.
  }
  \ee
  Changing the radial variable yields
   \be
   \label{w7}
    \tex{
   y \mapsto  \frac y{\sqrt 2} \LongA \tilde \BB^*= 2\big( \D_2 - \frac 12 \, y
   D_y
   + \frac 12 \, I\big),
   }
   \ee
where we observe the adjoint Hermite operator (cf. \ef{Her77})
 \be
 \label{w8}
 \tex{
\BB^*=\D_2 - \frac 12 \, y D_y \withA \BB= \D_2 + \frac 12 \, y
D_y +I \inB L^2(\re^2).}
 \ee
On the other hand, $\BB$ is self-adjoint in the weighted space
$L^2_\rho(\re^2)$, with the weight
 \be
 \label{w9}
  \rho(y)={\mathrm e}^{\frac{y^2}4},
   \ee
so that this falls into the scope of classic theory of linear
self-adjoint operators, \cite{BS}.

In particular, the adjoint operator $\BB^*$ in the adjoint space
$L^2_{\rho^*}(\re^2)$, where
 \be
 \label{w91}
  \tex{
  \rho^*(y) =\frac 1{\rho(y)}={\mathrm e}^{-\frac{y^2}4},
   }
   \ee
has  the discrete spectrum (here we take into account the
restricted radial symmetry)
 \be
 \label{w92}
  \tex{
  \s(\BB^*)= \{\l_{2k}= -  k, \quad k=0,1,2,...\},
 }
 \ee
 and the eigenfunctions are normalized Hermite polynomials
 \cite[p.~48]{BS}
  \be
\label{eigen}
 \tex{
 \psi_{2k}^*(y) =
  c_{2k} H_{2k}(y), \quad c_{2k} =
  \frac{2^{2k}} {\sqrt{(2k) !}}, \quad k=0,1,2,... \, .
   }
\end{equation}
As was already pointed out, in the general non-radial setting in
$\ren$, all these polynomials are obtained by differentiating the
rescaled Gaussian $F(y)$ of the fundamental solution
 of the heat operator $D_t-\D_N$,
 \be
  \label{G1}
   \tex{
 F(y)=\frac 1{(4 \pi)^{N/2}}\, {\mathrm e}^{- \frac{|y|^2}4} \quad (N=2),
  }
  \ee
 so that the following generating formula holds:
  \be
  \label{G33}
   \tex{
\psi_\b^*(y)= \frac 1{F(y)} \,  c_\b D^\b F(y) \equiv c_\b H_\b(y)
\whereA c_\b= \frac {2^{|\b|}}{ \sqrt{\b !}}
 }
 \ee
 and $\b=(\b_1,...,\b_N)$, with  $|\b|=\b_1+...+\b_N$, is a multiindex.
 This set of polynomials include all the
 angular-dependent eigenfunctions, so these are complete and
 closed in the whole weighted space $L^2_{\rho^*}(\re^2)$.
 Note that, unlike Section \ref{SPer}, we do not need
 solenoidal test for eigenfunctions involved.


It then follows from (\ref{w7}) that
 \be
 \label{G2}
\tex{
  \s(\tilde \BB^*)= \big\{ \l_{2k}= 1-2k, \quad k=0,1,2,...\big\}.
 }
 \ee
 The first unstable mode with $k=0$ and $\l_0= 1$ corresponds to
 the unavoidable unstability of the scaling (\ref{c81}) with
 respect to small perturbation of the blow-up time, when we
 replace $0 \mapsto T \approx 0$. On the manifold of solutions with
 the same blow-up time, this unstable mode plays no role and is
 excluded.

Thus, the first actual mode takes place for $k=1$, with
 \be
 \label{G3}
  \tex{
 \l_2=-1 \andA \psi^*_2(y) = \hat c_2 \big(1- \frac 14\, y^2).
 }
 \ee
 This corresponds to patterns with the following behaviour for
 large $\t \gg 1$ in the {\em inner region} that is characterized
 by arbitrarily large compact subsets in $y$:
  \be
  \label{G4}
  w_2(\t)\approx 1+  {\mathrm e}^{-\t} C \psi^*_2(y)+... \quad
  (C>0),
   \ee
 where the tangential operator in \ef{w4} should be now also taken
 into account. This makes the spectral theory more involved (other
 equations are also included) and can lead to complicated
 computations. For the purely radial case, analogous computations will be
 performed in the next section.

 Meantime, we deal with the conventional expansion \ef{G4}, which
 should be matched with the outer behaviour for $y \gg 1$.
 The corresponding rescaled variable in this {\em outer region}
 is also seen from \ef{G4} by using the quadratic behaviour of $H_2(y)$ as
 $y \to \iy$. Indeed, this gives the first term of the stationary
 expansion in the new rescaled variable:
  \be
  \label{G5}
   \tex{
 w_2(\t) \sim 1 - \hat c_2 C \,{\mathrm e}^{-\t} \, y^2+...
 \sim 1 -  \hat c_2 C \zeta^2+... \,\,\, \mbox{for} \,\,\,y \gg 1  \whereA
 \zeta={\mathrm e}^{- \frac 12 \, \t} y.
 }
  \ee

In a similar manner, we find other patterns with the inner
behaviour governed by other stable 1D eigenspaces of
$\tilde\BB^*$, where
\be
  \label{G6}
  w_{2k}(\t)\approx 1+  {\mathrm e}^{ \l_{2k}\t} C
  \psi^*_{2k}(y)+...\,
  \whereA k=1,2,3,... \, ,
   \ee
   and $\psi_{2k}^*(y)$ are other higher-degree non-monotone Hermite
   polynomials. In general, this can lead to
 a countable set of different blow-up patterns; see details in
 \cite{GalJMP}.

For the non-radial case, where all the tangential operators should
be taken into account, the matching procedure with the {\em outer
region} can be very complicated and cannot be justified rigorously
for the full model. This is done in \cite{GalJMP}
  for simpler blow-up axisymmetric  jets that are
admitted by the Navier--Stokes equations in $\re^3$.

In general, on matching, we arrive at a special countable family
of blow-up patterns that are characterized by a stronger and wider
propagation than the (existing or not) self-similarity scaling
\ef{c8} suggests. For instance, for the scalar curvature equation
\ef{SEHam1}, the first generic blow-up pattern \ef{bl11}, which is
characterized by a slow drift on a centre manifold corresponding
to the similarity rescaled variables \ef{ss21}, shows a faster
propagation in the $x$-direction with the extra logarithmic
 factor (here $T=0$)
 $$
 \sim \sqrt \t= \sqrt{|\ln(-t)|} \to + \infty \asA t \to 0^-.
  $$
 Such an expanding blow-up twistor can be not that adequate to the
 nature of single-point blow-up for the Navier--Stokes equations
 (at least, what we could expect).

 Therefore, we begin discussing the second, probably, more
 realistic situation.

\subsection{Discussion: on linearized patterns about singular
equilibria, (II)}
 \label{S3.8}

We now will use other equilibria of the equations (\ref{NS1}),
which are not constant and hence more and better localized about
the stagnation point $x=0$.

\ssk

\noi\underline{\em Singular radial equilibria}.   The first
candidate is {\em singular stationary solutions} that exhibit
strong singularities at $x=0$, and are possibly not any weak, very
weak,  or mild solutions at all. This is not that important for
us, since we are going to use them just for linearizing  and next
remove the singularity by matching with a regular bundle at the
singularity point.

\ssk

As an easy illustration of the type of calculus to be performed
later on, let us look
 for simple homogenuity $-1$ separable
stationary solutions of the original system \ef{c7},
 \be
 \label{z1}
  \tex{
U(r,\var)= \frac{A(\var)}r, \quad V(r,\var)= \frac{B(\var)}r,
\quad W(r,\var)= \frac{C(\var)}{r^2}, \quad P(r,\var)=
\frac{D(\var)}{r^2}.
  }
  \ee
 These separable structures are invariant under the scaling group in
\ef{c8}, so that the same singular patterns in the variables
\ef{c9} can be used for the study of twistor behaviour in the
rescaled system \ef{c7N}.

Substituting \ef{z1} into \ef{c7} yields the following system on
the unit circle ${\mathbb S}^1$ in $\re^2$:
   \be
 \label{z2}
  \left\{
   \begin{matrix}
  A''- BA' -2B'+A^2+B^2 +2D=0, \\
  B''-B B'+2A'-D'=0, \qquad \qquad\quad \,\\
  C''-BC'+4C +2AC -C^2=0, \qquad \\
  B'+C=0. \qquad\qquad\qquad\qquad\qquad\quad\,\,
 \end{matrix}
  \right.
  \ee
The second equation is integrated once that, on using the last
one,  gives $D$ as a quadratic function of other unknowns,
   \be
 \label{z3}
  \tex{
 D=2A- \frac 12\, B^2 -C + M_0 \quad (M_0={\rm const.})
  }
    \ee
    On substitution into \ef{z2}, we obtain an easier system
  \be
 \label{z4}
  \left\{
   \begin{matrix}
  A''- BA' +A^2+4A+2M_0=0,\,\, \\
  C''-BC'+4C +2AC -C^2=0,  \\
  B'+C=0.
  \qquad\qquad\qquad\qquad\quad\,\,
 \end{matrix}
  \right.
  \ee
 Substituting $C=-B'$ from the last equation yields two equations
 for functions $A$ and $B$,
 \be
 \label{z55}
  \left\{
   \begin{matrix}
  A''- BA'+A^2+4A+2M_0=0,\qquad \\
  B'''-BB''+ (B')^2+ 2(A+2)B'=0,
 \end{matrix}
  \right.
  \ee
which still remain a difficult fifth-order ODE on ${\mathbb S}^1$
with an arbitrary parameter $M_0 \in \re$, which can play a role
of a nonlinear eigenvalue. In particular, it is easy to check that
there are no explicit trigonometric solutions, where
 $$
 A,B \in W_3= {\rm Span}\, \{1,\cos l\var, \, \cos l \var\} \quad (l
 \in
 {\mathbb N}).
 $$
 Finally, we note that there exists simpler ``irrotational" (no swirl)
 solutions \ef{z1}:
  \be
  \label{ind1}
  C=0 \andA A^2+B^2+2D=0,
   \ee
 about which first linearized analysis can be naturally began
 with.

  It is worth mentioning here that \v{S}ver\'ak \cite{Sv06}  proved that
 all  singular stationary homogenuity $-1$ equilibria of the NSEs
  in $\re^3 \setminus \{0\}$
 are equivalent, up to orthogonal transformations, to
 Slezkin--Landau's solutions; a stronger result was obtained in \cite{MiuTs08}, see Section \ref{SLan}.
 Since singularities of \ef{z1} are concentrated at the $z$-axis
 $\{r=0\}$, we cannot use the result of \cite{Sv06}, though it does
 give a hope to avoid
 to scrutinize the general system (\ref{z4}),
 which looks not that easy at all.
In addition, there is a strong nonexistence result in
 \cite[\S~5]{Koch07}
 for axi-symmetric
{\em ancient} solutions\footnote{I.e., defined for all $t<0$, \cite{Ham95}. By
scaling \ef{Symm1}, any $L^\iy$ blow-up solution $\uu(x,t)$
generates a non-trivial ancient uniformly bounded one $\vv(x,t)=
\frac 1{C_k}\, \uu(x_k+ \frac x{C_k},t_k+\frac t{C_k^2})$, where
$C_k=\sup_x|\uu(x,t_k)|=|\uu(x_k,t_k)| \to + \iy$ (on the
applications, see \cite{Koch07}). This reflects
typical scaling tools
 of R--D theory \cite{AMGV}; cf. Section \ref{SEE2} and   applications to
 solvability and bounds for nonlinear parabolic equations in \cite[\S~2]{GW2}.}
  of the NSEs:
  \be
  \label{NNmm1}
  \tex{
  |\uu(x,t)| \le \frac C{\sqrt{x_1^2+x_2^2}}
  \inB \re^3 \times (-\iy,0)
   \LongA \uu \equiv 0.
   }
   \ee
   It seems
    this does not directly apply to equilibria \ef{z1} and
    similar  others, since these are singular at $r=0$ (i.e., are not
    bounded ancient solutions).
   Anyway, it seems that some singular equilibria can be eventually ruled out
    but not all of
   them.

\ssk

 Thus, by $\UU$ we denote a certain singular equilibrium of the
 NSEs,  not necessarily given by the homogeneous formulae \ef{z1}.
 In general, the description of all the possible singular
 stationary orbits at $x=0$ (or at $r=0$) leads to very difficult study of
 ill-posed elliptic evolution equations as explained in Section
 \ref{S5.1}.

\ssk

\noi\underline{\em Construction of blow-up  patterns: in need of
non self-adjoint theory}.
 Thus, we are not
going to study the system \ef{z4} in detail here, and will
concentrate on the principles of construction of blow-up twistors
using linearization about the manifold corresponding to such
singularities.
Denoting by $\uu$ the 3-vector of the rescaled variables
 $\uu=(u,v,w)$, we write \ef{c7N} as a dynamical system for the
 rescaled variable
 (cf. \ef{ww1})
   \be
 \label{z5}
  \tex{
  \hat \uu'= \HH(\hat \uu) \forA \t>0,
  }
  \ee
 where, as in \ef{ww1} or \ef{HHH21}, the pressure variable is excluded
 by projecting the solution space onto the kernel of the gradient operator.
Thus, we have  fixed $\UU(y,\mu) \not \equiv 0$ as a solenoidal
stationary singular solutions of the rescaled equation \ef{z5}.

We next perform the linearization by setting
   \be
 \label{z6}
  \tex{
\hat \uu(\t)= \UU + \YY(\t)
  }
  \ee
 to get the linear equation for $\YY$,
   \be
 \label{z7}
  \tex{
  \dot \YY= \HH'(\UU) \YY + \DD(\YY) \forA \t>0,
  }
  \ee
 where $\DD$ denotes a nonlinear  perturbation, which is quadratic as $\YY \to {\bf 0}$.

It follows from the general steady structure such as  \ef{z1} that
the first step of such a construction is to find the point
spectrum of the linear integral (pseudo-differential)  matrix
operator $\BB^*=\HH'(\UU)$ for vector-valued functions with
periodic in $\mu$ (and possibly singular at the origin)
coefficients. Thus, we define the {\em Inner Region I} as  a
family of special space-time  subsets, where the linearized
equation holds asymptotically:
 \be
 \label{As12}
  \tex{
  \mbox{Inner Region I:} \quad \dot \YY = \HH'(\UU) \YY + ...
  \forA\t \gg 1.
  }
  \ee
 Recall that, by construction, $\HH'(\UU)$ is assumed to
act in a solenoidal vector field.

 A number of accompanying questions arise. E.g.,
 one can pose the following:

\ssk

\noi\underline{\em Question {\rm (i)}:} $\BB^*=\HH'(\UU)$ is not
self-adjoint in no weighted $L^2$-spaces;

\ssk

 \noi\underline{\em  Question {\rm (ii)}:} what is the domain of $\BB^*$ (hopefully a kind of
Sobolev space $H^2_{\rho^*}(\re^3)$)?

\ssk

 \noi\underline{\em  Question {\rm (iii)}:} what is the space and domain of the adjoint operator
$\BB$?

\ssk

 \noi\underline{\em  Question {\rm (iv)}:} it is not clear why $\BB^*$ and $\BB$ could have
 enough or at least some
real (or complex) eigenvalues, compact resolvent (by compact
embedding $H^2_{\rho^*} \subset  L^2_{\rho^*}$? which conditions
at the singular origin?), and bi-orthogonality property of bases
(if any);

\ssk

 \noi\underline{\em  Question {\rm (v)}:} why the eigenfunctions of $\BB^*$ can be at least approximately
  associated with
some structures of finite polynomials (possibly separable in the
angular direction)?

\ssk

\noi\underline{\em  Question {\rm (vi)} on solenoidal
eigenfunctions and classes}. Though this is a natural part of the
eigenvalue problem for the pair $\{\BB^*,\, \BB\}$, let us recall
the key:
assuming the $-1$ homogenuity (and if not?) of the equilibrium
such as \ef{z1} employed,
 we are obliged to perform a ``solinoidalization"
  to get a sufficient amount of (possibly even some  {\em
 polynomial} where again the angle separation is assumed?)
 eigenfunctions $\psi_\b^*$ ,
 which will generated eigenspaces of {\em solenoidal} fields as in
 \ef{Sol1}; ...\,, etc.



\ssk

For the singular S--L solutions \ef{Lan1}, this will be continued
in
  Section \ref{SLan}, where we also discuss in greater detail the questions of matching of various
   asymptotic inner and outer regions.

\ssk

 Nevertheless, the parabolic experience of doing linearized
theory for $2m$th-order PDEs such as \ef{ho1} shows that the
corresponding linear operators $\BB$ and $\BB^*$ are also non
self-adjoint (the only self-adjoint case is for $m=1$; cf.
\ef{mm1}), but admit real point spectrum only, \cite{Eg4}.
However, the matrix spectrum problem for the pseudo-differential
operator $\BB^*$ in \ef{z7} and for the corresponding formally (in
the topology of $L^2$) adjoint operator $\BB$ are much more
difficult and remain open;
 see Section \ref{SLan} for  further details and related important comments in \cite{Sv06}.

Thus, as for an illustration,
 assume that $\BB^*$ has a real point
spectrum
   \be
 \label{z8}
  \tex{
  \s(\BB^*})=\{\l_k, \quad k=|\b|=0,1,2,...\} \quad (\mbox{$\b$ is a
  multiindex),
  }
  \ee
  where the first positive eigenvalue,
   \be
 \label{z9}
  \tex{
   \l_0> 0 \quad \big(=  1, \,\,\, \mbox{it seems}\big),
  }
  \ee
 reflects the natural unstability with respect to the change of
 the blow-up time. We also assume that, as a standard fact from
 linear operator theory (that also needs a difficult approving),
 the eigenfunctions $\{\vv_\b\}$ are assumed to be
 bi-orthonormal to the corresponding adjoint basis
 $\{\vv_\b^*\}$, and both bases are complete, closed, and form
 Riesz-type bases in the corresponding solenoidal  spaces.

We are interested in eigenvalues $\l_k \le 0$. The most
interesting case occurs when $\l_k=0$, since it corresponds to a
behaviour close to the centre subspace of $\BB^*$. Note that
proving existence of a centre or stable invariant manifolds for
the full nonlinear problem \ef{z5} is a very difficult problem. As
a formal hint, assuming that the kernel of $\BB^*$ is
one-dimensional spanned by an eigenfunctions $\vv_\b^*$ (for any
dimension, the analysis is similar with a system to occur), we are
looking for solutions in the form\footnote{Recall
\ef{Sol1} for such eigenfunction expansions.}
 \be
 \label{zzz11}
 \tex{
 \hat  \YY_0(\t) \sim \cc_\b(\t)  \vv_\b^* + \hat \ww(\t) \forA \t \gg 1
 \quad(\vv_\b \in {\mathcal S}_k^*, \,\, \cc \vv^*=(c_1 v^*_1,c_2 v^*_2, c_3 v^*_3)),
 }
 \ee
 where the remainder $\hat \ww(\t)$ is supposed to be orthogonal to ${\rm ker} \,
 \BB^*$ in the dual metric between the spaces for $\BB$ and
 $\BB^*$ denoted by $\langle \cdot,\cdot \rangle$. Substituting
 into \ef{z7} and multiplying by $\vv_\b$ yields the following
 asymptotic ODE system:
 \be
 \label{zzz1}
 \tex{
 \dot \cc_\b= \langle \DD(\cc_\b  \vv_\b^*), \vv_\b \rangle+...
  \equiv {\bf Q}_\b(\cc_\b)+... \, ,
 }
 \ee
 where ${\bf Q}_\b(\cc_\b)$
  is a  vector quadratic form in $\re^3$.
 Looking for a standard solution with a power decay yields a
 quadratic algebraic system on the coefficients,
  \be
  \label{dd1}
   \tex{
   \cc_\b(\t)= \frac {\AAA_\b} \t+... \LongA {\bf
   Q}_\b(\AAA_\b)=-\AAA_\b.
    }
    \ee
    Assuming that the matrix equation in \ef{dd1} has a solution
    $\AAA_\b \not = 0$, we obtain the following asymptotic pattern:
 \be
 \label{zzz2}
 \tex{
 \hat \YY_0(y,\t) \sim  \frac {\AAA_\b}{\t} \, \vv_\b^*(y)+... \, \asA \t \to +
 \infty \quad (\l_k=0).
 }
 \ee
 In this way, we observe a slow logarithmic drift, with $\t = -
 \ln(-t)$, relative to the centre subspace of $\HH'(\UU)$.

For $\l_k<0$, with $k=|\b|$, we arrive at more standard stable
subspace patterns
 \be
 \label{zzz3}
 \tex{
 \hat \YY_k(\t) \sim {\mathrm e}^{\l_k\, \t} \CC_\b\,  \vv_\b^*+... \quad (\l_k<0).
}
 \ee
 Both \ef{zzz11} and \ef{zzz3} describe for $\t \gg 1$
 asymptotically small perturbations in \ef{z6} of a steady singular
 solenoidal field $\UU$, so we have to choose solenoidal
 eigenfunctions.


 \ssk

\noi\underline{\em On outer matching: Outer Region}.  The main
idea of matching of these patterns with the {\em outer region},
which is most remote is explained in \cite{GalJMP},
 where we are assuming that $\vv_\b^*(y)$ have a polynomial
 decay, for some coefficients $\d_k>0$, $k=|\b|$,
  \be
  \label{pol1}
  |\vv_\b^*(y)| \sim b_k \, y^{-\d_k}+... \asA y \to +\infty.
 \ee
 Here we arrive at a simpler
first-order matrix operator, since in the new rescaled variables
of the  outer region, all the Laplacians form asymptotically small
singular perturbations simply meaning that Euler's equations are
dominant here. Of course, the resulting singularly perturbed
dynamical system is difficult to tackle to pass to the limit $\t
\to +\infty$.

Note the following crucial property of such an extension into the
next outer region II, where solution is much less: we are not now
obliged to keep there the self-similar power-like decay at
infinity.
   Moreover, we now can match this blow-up
behaviour with the exponentially decay
 associated with the heat semigroup ${\mathrm e}^{\D_2 \t}$ for
 the equations such as \ef{w1} and others. This matching is not
 a key
 difficulty of the blow-up twistor construction and
 deals with  not that singular
 orbits. Again,
 we take into account patterns \ef{Sol1}.

\ssk

\noi\underline{\em On principles of singular inner matching: Inner
Region II}. The last crucial step is to match the centre \ef{zzz2}
and stable \ef{zzz3} subspace patterns with the regular and
bounded solenoidal flow in the {\em singular inner region} close
to the origin $y=0$. This is the most involved matter, where the
divergence-free asymptotics  satisfying something like \ef{Sol1}
become key.
  {\em The singular inner matching in a countable
number of various cases is responsible for existence or
nonexistence of solenoidal blow-up patterns and is the main open
problem.}

 Assuming for definiteness that such matching exists, we predict some
 other evolution features of the resulting blow-up patterns.
 Here we formally suppose that the presented asymptotic expansions
 go  along the solenoidal classes like in \ef{Sol1},
   so that such
 spatial-temporal structures can be used for some preliminary
 estimates of blow-up evolution.

  Assume first that in the
stable subspace representation \ef{zzz3}, the eigenfunction is
singular at $y=0$ with the following behaviour ($\g_k>1$):
 \be
 \label{pol2}
 |\vv_\b^*(y)| \sim d_k \, y^{-\g_k}+... \asA y \to 0.
 \ee
Therefore, the pattern \ef{zzz3} has the asymptotic behaviour in
the singular inner region of the type (we omit here non-essential
constants and neglect other minor multipliers)
 \be
 \label{pol3}
  \tex{
\hat \uu_k(y,\t) \sim \frac \AAA y -{\mathrm e}^{\l_k \t} \CC_\b
y^{-\g_k}+... \,.
 }
  \ee
 Calculating as in \ef{WW9} (see  \cite{Dold1}) the maximum point in $y$ of the
 function in \ef{pol3} yields
  \be
  \label{pol4}
   \tex{
  \sup_{y>0} \, |\hat \uu_k(y,\t)| \sim {\mathrm e}^{ -\frac{\l_k}{\g_k-1} \,
  \t} \to + \infty \asA \t \to + \infty \quad (\l_k<0).
  }
  \ee
  Under the above assumptions,
this is about a right bound on regular hypothetical continuation
of the pattern $\hat \uu_k$ from an exponentially small
neighbourhood of the singular state $\UU$ to smooth bounded
solutions for $y \approx 0$. Note that the behaviour of the
pattern as $y \to 0$ can be very complicated including blow-up
swirl-like features.
  This gives a preliminary estimate of the rate of this non
  self-similar blow-up in the $\uu$-variable
   \be
   \label{pol5}
   \| \uu_k(\cdot, t)\|_\infty \sim (-t)^{-\frac 12 + \frac{\l_k}{\g_k-1}}
   \asA t \to 0^- \quad (k=|\b| \gg 1).
    \ee
Similarly, for the centre subspace pattern \ef{zzz2} (if any),
this very rough estimate imply the following perturbation of the
self-similar rate:
 \be
   \label{pol51}
   \tex{
   \|\uu_0(\cdot,t)\|_\infty \sim  \frac 1{\sqrt{-t}}\, |\ln(-t)|^{\frac
   1{\g_k-1}}
   \asA t \to 0^-.
   }
    \ee
 See Section \ref{S7.5} for further necessary properties of such a matching.

    \ssk

\noi{\bf Remark: back to reaction-diffusion.}
 For single parabolic equations such
as \ef{FK1}, \ef{SEHam1} (for $N >6$), or a more general
combustion model
 \be
 \label{Com1}
 u_t = \D u + |u|^{p-1}u \inA \quad (p>1),
  \ee
 such examples of creating non-similarity
asymptotics by matched asymptotic expansion techniques are known
and have been rigorously justified; see \cite{Dold1, Fil00, GKSob,
GV97, HVsup}, and references therein.
 Even for these
 reasonably simple parabolic PDEs,
 in the supercritical Sobolev range as in \ef{ST2},
 the set of singular stationary
 solutions can be rather complicated that affect the evolution
 properties of blow-up and other asymptotics; see
 \cite{Mat04, Miz04, Miz06} and Section \ref{S7.7}.
 For solutions of \ef{Com1} of changing sign (that are always taking place in our
 case), the blow-up patterns for $p=p_{\rm S}$ become more involved and
 their actual existence is still not fully justified rigorously, \cite{Fil00}.

 \ssk

  The reaction-diffusion equation
 \ef{Com1} does not assume any conservation of energy, but,
 however, there are examples of the so called {\em incomplete
 blow-up}, when the solution $u(x,t)$ becomes unbounded in $L^\iy(\ren)$ at a
 single moment $t=T$ only; see \cite{GV97} and earlier references
 therein concerning equation \ef{FK1}, as well as
 \cite{Miz06, Miz07} for more recent extensions and achievements.
 For the model \ef{NS1},  blow-up moments are
 {\em incomplete}, in view of the {\em a priori} bound \ef{ss1}.

 It seems this
quite a bounded mathematical experience cannot be directly
 translated to the Navier--Stokes rescaled equations \ef{z5}.

\ssk

We have discussed the second principle to generate another
countable set of blow-up twistors. No doubt that there are other
ways to get various countable families of such patterns, which we
are still not aware of and cannot even imagine.

\subsection{Discussion: on a very formal way  to create a vertex
and branching to localized smooth blow-up twistors;  a twistor
ring}
 \label{S3.9}

Here we discuss another important issue concerning the twistors on
the 2D subspace ${\rm Span}\,\{1,z\}$.

\ssk

(i) A vertex is  by truncation of the $W$-component in \ef{c5} via
the positive part,
 \be
 \label{pos1}
 W = (z)_+ \, \tilde W,
 \ee
 so that $W=0$ for $z<0$. This leads to a weak solution of
 (\ref{NS1}), so a certain effort is necessary to check whether this
 gives a Leray--Hopf  solution in the sense of the inequality in
 \ef{ss1} for all $t<0$, i.e.,  before the blow-up occurs.
 On the other hand, in view of the existence of a local semigroup
 of smooth bounded solutions, existence of such solutions with a
 vertex assumes performing a certain smoothing at the stagnation
 point at $x=0$. This cannot be done in the above self-similar
 manner, so the proof of such a ``smooth branching" (if any) from the given
 non-smooth blow-up twistor represents a typical difficult open problem
 of modern theory of nonlinear PDEs.

 \ssk

 (ii) Similarly, there occurs another, not simpler, open problem of
 ``smooth branching" from this unbounded twistor of
 other, partially and asymptotically self-similar blow-up twistors that might be sufficiently spatially localized
 as $t \to 0^-$. As we have mentioned, a possible way to check a
 possibility of such branching consists of using  asymptotic
 expansion theory to create solutions via series like \ef{em1}.
 Obviously, this leads to difficult nonlinear (first ones) and
 further linearized systems for the expansion coefficients and
 represents another face of this  open branching problem.

\ssk

(iii) There is another  formal geometric way to create blow-up
patterns by deforming the axis of their symmetry. For instance,
let us assume, in view of our pretty local analysis of such
blow-up behaviour for small $z \approx 0$, that the $z$-axis can
be bend into a finite ring such that at some points similar
blow-up structures oriented in the tangential $z$-direction occur
with a certain periodic behaviour along the
 ring generatrix\footnote{For the Euler equations \ef{ee1}, see \cite{Pom05}.}.
This leads to ``twistor rings" with a very difficult open
accompanying mathematics.


\subsection{On twistors in  Euler equations}
 \label{S3.14}

Similar  blow-up twistors, with analogous concepts of their
 swirl-like and other extensions,
  can be formally
detected for the Euler equations \ef{ee1}, where the calculations
are supposed to be simpler but actually they are not.
   Without the ``curse
of local smoothness" for the Navier--Stokes equations \ef{NS1},
for the Euler ones \ef{ee1} containing first-order differential
operators only without a typical interior regularity, the
necessary branching phenomena are easier to justify, so there is a
hope of doing these  more rigorously. However, one can see that,
neglecting Laplacians $\D_2$ in the systems \ef{c7} and \ef{c7Ns},
indeed simplifies the analysis but the remaining PDEs are still
very difficult to understand rigorously that demands new concepts
of solutions and entropy regularity. The stationary system
\ef{q1}--\ef{q3} with no Laplacian operators $\D_2$'s
 takes the form
 \be
 \label{q2E}
  \left\{
   \begin{matrix}
  - \frac 12\, y \hat w_y -\hat w -\sigma  \hat w_\mu + \hat w^2
  \qquad\qquad\qquad\qquad\qquad\qquad\qquad\qquad\qquad\qquad\,\,\,
 \ssk\\
 - \frac 1y
 \, \big(\int_0^y z \hat w\,{\mathrm d}z \big) \hat w_y +
 \frac 1y \, \big(\int_0^y \hat v_\mu \, {\mathrm d}z \big)
\hat w_y- \frac 1y \, \hat v \hat
w_\mu=0.\qquad\qquad\qquad\qquad\qquad \ssk\ssk\ssk\\
   \big[ \frac 12\, y \hat u_y + \frac 12\,\hat u +\sigma  \hat u_\mu + \hat u \hat u_y
    +\frac 1y \,  \hat v
 \hat u_\mu - \frac 1{y} \, \hat v^2  + \frac
 2{y^2} \, \hat v_\mu + \frac 1{y^2} \, \hat u\big]_\mu \qquad\qquad\qquad\quad \ssk\\
=
  \big[y\big( \frac 12\, y \hat v_y + \frac 12\,\hat v +\sigma  \hat v_\mu + \hat u \hat v_y + \frac 1y \,  \hat
  v
 \hat v_\mu + \frac 1{y} \, \hat u \hat v  - \frac
 2{y^2} \, \hat u_\mu + \frac 1{y^2} \, \hat v\big)\big]_y,
 \qquad\qquad
 \end{matrix}
 \right.
  \ee
with the same nonlocal constraint \ef{q1}. It is still a difficult
system  of a first and a second-order nonlocal evolution PDEs,

  i.e., actually, is a
system of higher-order equations.

\ssk

As we know, regardless its solvability, the linearization of this
system about constant, singular, or other equilibria  defines the
linear operators that may generate other countable families of
non-self-similar blow-up patterns with the evolution on centre
and/or stable eigenspaces, as in Sections \ref{S3.7} and
\ref{S3.8}.

\ssk

 It is important to mention that, following the lines of this
construction, we necessarily arrive at the non-stationary system
corresponding to \ef{q2E}.
 Since odd-order operators are dominated here, a proper
``entropy" setting for solutions with weak and strong
discontinuities will be necessary to check evolution consistency
 of the blow-up patterns constructed\footnote{Recall that
 there is still no a successful notion of weak solutions of Euler equations in 3D.}. This is expected to be a
  hard problem that, with  a clear inevitability, accompanies
 this pointwise blow-up analysis of the Euler equations.


\ssk

 Finally, we believe that general concepts of swirling
rotations,
 axis precessions, and vertex motion of periodic or
quasi-periodic nature developed in Section \ref{SLanFin} can be
applied to the Euler equations \ef{ee1}, but we will not  develop
these here and return to the Navier--Stokes ones \ef{NS1}; see
\cite{Pom05} for further refreshing and rather exotic ideas.
Recall that typical ``rolling-up mechanisms" for appeared bubble
caps and other swirling features of formation of blow-up
singularities were observed numerically even in the axisymmetric
setting; see \cite{Gra95, Pum92}, etc., though these scenarios of
blow-up seem remain still under scrutiny, while  some of them have
been ruled out in other works.




\section{On non-radial blow-up patterns on $W_2$: eigenfunctions
with swirl}
 \label{Snonrad}


In this section, we show how to extend some ideas coming from the
simpler model proposed in Ohkitani \cite{Oh07} with a countable
set of blow-up solutions constructed in \cite{GalJMP}.
 We return to the NSEs  restricted to the subspace $W_2={\rm
Span}\{1,z\}$ in Section \ref{S4.2}.

\subsection{The system on the subspace $W_2$ with swirl and angular dependence}
 \label{S6.1}

 For convenience, we  replace
 \be
 \label{WWW11}
 \tilde W \mapsto W,
  \ee
 so that the system of three equations from Section \ref{S4.2} takes the form
 \be
  \label{Sys11}
   \left\{
   \begin{matrix}
   W_t= \D_2  W - U  W_r - \frac 1r \, V  W_\var +  W^2,
   \qquad\qquad\qquad\qquad\qquad
  \ssk \\
 U=  \frac 1r \, \int_0^r z  W(z,\var,t) \, {\mathrm d}z
- \frac 1r \, \int_0^r  V_\var(z,\var,t) \, {\mathrm d}z,
 \qquad\quad\qquad\quad\,\,
\ssk \\
   \big[U_t+ UU_r+ \frac 1r \, V U_\var - \frac 1{r} \,
V^2 -\D_2 U + \frac 2{r^2}\, V_\var + \frac 1{r^2} \, U\big]_\var
\qquad\quad \ssk\ssk\\
 =\big[r\big(V_t + UV_r+ \frac 1r\, VV_\var + \frac 1r\, UV-
\D_2 V - \frac 2{r^2}\, U_\var + \frac 1{r^2} \, V)\big]_r.
 \end{matrix}
  \right.
 \ee
Indeed, this is a very difficult system. As we have mentioned, its
axi-symmetric version admits further, more rigorous study,
\cite{Oh07, GalJMP}. Here we demonstrate other aspects of such
solutions.

Thus, we apply the same scaling as in \ef{c8}, \ef{c9} to the
general system \ef{Sys11} to get
 \be
  \label{Sys11r}
   \left\{
   \begin{matrix}
   w_\t= \D_2  w - \frac 12\, y w_y-w - \sigma  w_\mu - u  w_y - \frac 1y \, v  w_\mu +  w^2,
   \qquad\qquad\qquad\qquad\,\,\quad\,\,\,
  \ssk \\
 u=  \frac 1y \, \int_0^y z  w(z,\mu,\t) \, {\mathrm d}z
- \frac 1y \, \int_0^y  v_\mu(z,\mu,\t) \, {\mathrm d}z,
 \qquad\quad\qquad\quad\,\,\qquad\qquad\qquad\quad\,\,\,
\ssk \\
   \big[u_\t+\sigma  u_\mu + u u_y+ \frac 1y \, v u_\mu - \frac 1{y} \,
v^2 -\D_2 u+ \frac 12 \, y u_y + \frac 12\, u + \frac 2{y^2}\,
v_\mu + \frac 1{y^2} \, u\big]_\mu \qquad\quad \ssk\ssk\\
 =\big[y\big(v_\t + \sigma  v_\mu+ u v_y+ \frac 1y\, v v_\mu + \frac 1y \, u v-
\D_2 v + \frac 12\, y v_y + \frac 12\, v  - \frac 2{y^2}\, u_\mu +
\frac 1{y^2} \, v)\big]_y.
 \end{matrix}
  \right.
 \ee
 The linearization $w=1+Y$ as in \ef{w3} yields the linearized
 equation (cf. \ef{w4})
  \be
  \label{gg1}
   \tex{
   Y_\t= \BB^*Y +Y - \sigma  Y_\mu - u Y_y- \frac 1y\, v Y_\mu + {\bf
   D}(Y),
    }
    \ee
    where $\DD$ is as in \ef{w4} and $\BB^*$ is the adjoint
    Hermite operator \ef{w8}. We next use the second equation for
    $u$ in \ef{Sys11r} to get
    \be
    \label{pppp1}
     \tex{
    u= \frac y2 + \frac 1y \, \int_0^y zY\,{\mathrm d}z - \frac 1y
    \, \int_0^y v_\mu\,{\mathrm d}z.
     }
     \ee

  Assume now that
  \be
  \label{kk2}
   Y(y,\mu,\t) \andA v(y,\mu,\t) \quad \mbox{are exponentially
   small for $\t \gg 1$.}
    \ee
    Then by \ef{pppp1} $u= \frac y2+...$ up to exponentially small
    perturbations, so that small solutions of \ef{gg1} are still
    governed by the linear operator \ef{w6}
 \be
 \label{kk3}
  \tex{
  \BB^*Y +Y - \sigma  Y_\mu
   \equiv \tilde
  \BB^* Y - \sigma  Y_\mu+... \asA \t \to +\iy.
  }
  \ee

Consider first the full operator with the rotational part with the
adjoint one in $L^2(\re^2)$
 \be
 \label{kk4}
  \tex{
 {\bf L}_\sigma ^*=
 \tilde \BB^* - \sigma  D_\mu \andA {\bf L}_\sigma =
 \tilde \BB + \sigma  D_\mu \quad \big(\tilde \BB=\D_2 +
 y D_y+3I\big).
 }
  \ee
  Similar to $\BB^*$ (similar properties are available for
  $2m$th-order counterparts \cite{Eg4}),
  it can be shown that ${\bf L}_\sigma ^*$ is a bounded operator from
  $H^2_{\rho^*}(\re^2)$ to $L^2_{\rho^*}(\re^2$, has compact
  resolvent, and a point spectrum only.
  The only, but a key difference, is that, it seems, for $\sigma  \not = 0$,
  \ef{kk4} does not admit a proper
self-adjoint setting, so that the whole point spectrum is not
necessarily real. Writing down in detail the spectral problem in
$L^2_{\rho^*}(\re^2)$,
 \be
 \label{kk41}
  \tex{
   {\bf L}_\sigma ^* \psi^*
  \equiv \psi^*_{yy}+ \frac 1 y \,
 \psi^*_y + \frac 1{y^2} \psi^*_{\mu\mu}- \frac 12\, y \psi^*_y+
 \psi^*
 - \sigma  \psi^*_\mu= \l \psi^*,
  }
  \ee
we see that the essential use of the angular operators makes it
rather difficult. Obviously, any radial Hermite polynomial remains
an eigenfunction for any $\sigma  \in \re$.
 Finding a suitable point spectrum and non-radial eigenfunctions
 with essential $\mu$-dependence for
\ef{kk41} can be  a difficult problem. One can see that a standard
separation of variables is  non-applicable for \ef{kk41}.

 Thus, a principal question is whether the operator \ef{kk4} has
 enough real eigenvalues for construction of linearized blow-up
 patterns. Let us show that in this sense \ef{kk4} is quite
 suitable.
 As an illustration of a  local approach to \ef{kk41}, we apply the classic
 techniques \cite{VainbergTr} to trace out  branching of eigenfunctions from $\psi^*_\b$ of $\tilde \BB^*$ at
 $\sigma =0$. 
 Let the kernel of $\BB^*-\l_\b I$ has the dimension $M \ge 2$.
 Then looking for the expansion
  \be
  \label{kk43}
   \tex{
   \psi^*(\sigma )= \sum_{k=1}^M c_k(\sigma )\psi^*_{\b,k}+ \sigma  \var+...
   \,\,(\var \bot \psi_\b^*) \andA \l(\sigma )=\l_\b+ \sigma  s+... \, ,
 }
    \ee
  we obtain the following problem ($\l_\b$ are as given in \ef{kk5}):
   \be
   \label{kk44}
    \tex{
   (\tilde \BB^*- \l_\b I)\var = \sum_{(k)} c_k(0)\big[
   s \psi_{\b,k}^* + (\psi_{\b,k}^*)'_\mu \big].
   }
 \ee
 This yields $M$
 orthogonality conditions of solvability
 \be
 \label{kk45}
  \tex{
\big\langle \sum_{(k)} c_k(0)\big[
   s \psi_{\b,k}^* + (\psi_{\b,k}^*)'_\mu \big], \psi_{\b,j}
\big\rangle=0 \quad \mbox{for} \quad j=1,2,\, ... \,,M,
 }
 \ee
 that define the eigenfunctions of $\tilde \BB^*$, from which branching
 is available.
  Then \ef{kk44}
 yields the unique solution $\var$ that shows the branching
 evolution in \ef{kk43}
 (as usual, the dimension of the kernel then can play a
 part).
  In other words, \ef{kk4} in the inverse
 integral form can be treated as a compact perturbation of a
 self-adjoint operator, and classic perturbation theory applies
 \cite{Kato} to calculate the deformation of the real point spectrum.

  This shows that there exists branching for small
 angular speeds $|\sigma |$ of
 eigenfunctions with swirl from standard Hermite polynomials for $\sigma =0$.
 For compact inverse integral operators involved, these are
 continuous $\sigma $-curves are indefinitely extensible and  can end up at other bifurcation
 point only or can be unbounded; see \cite{BufTol}.
  In other words, the swirl-dependent operators \ef{kk44} have
  enough real eigenfunctions for using in necessary construction
  of linearized blow-up patterns. Some of the related questions of
  linear operator theory for \ef{kk44} remain obscure; e.g.,
  eigenfunction closure for such non-symmetric cases (but indeed this
  looks doable though can be  rather technical).

\ssk

 It is curious that even in the case $\sigma =0$, we can find essentially non-radial
patterns, which can give an insight into a swirl structure of
blow-up patterns that have a clear $\var$-dependence. Namely, we
now take into account all the non-radial Hermite polynomial
eigenfunctions \ef{G33} of $\BB^*$ with the spectrum
 \be
 \label{kk5}
  \tex{
   \s(\tilde \BB^*)=\big\{\l_\b=1-|\b|, \,\, |\b|=0,1,2,...\big\}.
   }
   \ee
Therefore, there exist solutions of \ef{gg1} with the asymptotic
behaviour as $\t \to +\iy$,
 \be
 \label{kk6}
  \tex{
  w=1 -{\mathrm e}^{\l_\b \t} \psi_\b(y,\mu)+...\, \quad \mbox{for
  any}
  \quad |\b|>1,
   }
   \ee
   provided that \ef{kk2} holds on compact subsets in $y$. Obviously, the
   third equation in \ef{Sys11r} admits exponentially
   small solutions $v$ for $u$ being an exponentially small
   perturbation of $\frac y2$ according to \ef{pppp1}. Then
    patterns \ef{kk6} in the inner region make sense.

   The extensions of such patterns into the outer region is
   similar. We have :
    \be
    \label{kk7}
    \psi_\b(y,\mu) \sim y^\b = |y|^{|\b|} f_j(\mu) \asA y \to +\iy,
     \ee
     where the first entry $y$ means $(y_1,y_2)^T$ and $f_j$ are homogeneous
     harmonic polynomials being the
     eigenfunctions of the Laplace--Beltrami operator $\D_\mu
     =D_{\mu}^2$ on the circle ${\mathbb S}^1 \subset \re^2$,
      \be
      \label{kk8}
       \D_\mu f_j=-j^2 f_j \quad \mbox{on} \quad {\mathbb S}^1.
 \ee
 Therefore, \ef{kk6} yields the outer variable
  \be
  \label{kk9}
   \tex{
   w \sim 1- |\zeta|^{|\b|} f_j(\mu) \whereA \zeta=y\,
    {\mathrm e}^{\rho_\b \t}, \quad \rho_\b= \frac{\l_\b}{|\b|}<0.
    }
    \ee

Finally, rewriting the first equation in \ef{Sys11r} in terms of
the new variable $\zeta>0$ by the change
corresponding to the region governed by the new spatial rescaled
variable $\zeta$, we recast the original equation
 by
setting (this change is explained in detail in \cite{GalJMP})
 $ \th(\zeta,\t)= w(y, \t)$.
 This yields the
following perturbed Hamilton--Jacobi (Euler) equation:
 \be
 \label{kk10}
 \begin{matrix}
 \th_\t=\AAA_\b(\th) +  {\mathrm e}^{2\rho_\b \t} \D_2 \th
 +  {\mathrm e}^{\rho_\b \t}\big[ \frac 1 \zeta\, \big(\int_0^\zeta
 v_\mu\big) - v w_\mu \big],\ssk\ssk\ssk \\
 \mbox{where} \quad \AAA_\b(\th)=-\big(\rho_\b+ \frac 12\big)\zeta \th_\zeta
 - \frac 1 \zeta \, \big(\int_0^\zeta z \th \big)\, \th_\zeta
  - \th +
 \th^2.
  \end{matrix}
 \ee
 Thus, eventually, on passage to the limit $\t \to \iy$, we obtain the steady problem
 \be
 \label{kk12}
 \tex{
-\big(\rho_\b+ \frac 12\big)\zeta h_\zeta
 - \frac 1 \zeta \, \big(\int_0^\zeta z h \big)\, h_\zeta
 - h +
 h^2=0, \quad h(\zeta)\sim 1-\zeta^\b f_j(\mu),  \,\, \zeta \to 0.
 }
  \ee
For each fixed $\mu \in[0,2\pi)$, this is an ODE, which was the
object of a detailed study in \cite{GalJMP}.
  Therefore, for any $f_j(\mu)>0$, such a profile always exists and is
compactly supported with non-radial support.
 Unfortunately, for $f_j(\mu)<0$, the ODE \ef{kk12} gives a monotone
 increasing solution which moreover blow-up in finite $\zeta$.
 Therefore, such non-radial patterns can be constructed in
 any connected sector $S_j=\{\mu: f_j(\mu)>0\}$ with zero Dirichlet
 conditions at the boundary rays.
 There holds:
 \be
 \label{kk21}
  \zeta= x (-t)^{\frac{|\b|-2}{2|\b|}}
 \LongA \mbox{separable standing-wave blow-up for}\,\,\,|\b|=2.
  \ee
This case includes the radial pattern already studied and others,
which are no symmetric. For any $|\b|>2$, the patterns exhibit
global blow-up as usual in the fixed angular sectors.

We hope that using non-trivial rotation $\sigma  \not = 0$ will
supply us, {\em via} the eigenvalue problem \ef{kk41}, some extra
new eigenfunctions that allow us to get a non-radial blow-up
pattern with swirl in the whole space.

\section{On  blow-up patterns concentrated about  
Slezkin--Landau singular solutions of a ``submerged jet": first
example of linearization}
 \label{SLan}

We now discuss some  ways to construct necessary blow-up patterns
to be applied in greater detail in Section \ref{SLanFin} by using
the spherical geometry. Namely, we
 introduce implications of using the following classic singular
solutions of the NSEs:

\subsection{Singular homogeneous equilibria and their properties}

In 1934, Slezkin \cite{Slez34} (see also comments in
\cite{Slez54}) showed that
 \ef{NS1} admit special stationary solutions
 with the singularity and spatial
decay of the velocity field $\sim \frac 1r$ by reducing the
problem to a linear hypergeometric-type ODE; see more details in
Appendix C.
   In 1944\footnote{In Landau--Lifshitz's monograph \cite[p.~81]{LandM},
   this solution dated 1943, with no reference shown. Possibly,
   this reflects the  year, when the solution was actually constructed by Landau, but not published.}, Landau
 \cite{Lan44} found  a family of explicit
 solutions of that type
  describing  steady flows induced by
 a point source, which leads
  to the setting of a {\em submerge
jet} that is oriented along the positive part $Oz$ of the
$z$-axis.
This gives the following one-parameter family of the explicit {\em
Slezkin--Landau singular stationary  solutions} $\{\uu_{\rm
SL}({\bf x}),p_{\rm SL}({\bf x})\}$, of \ef{NS1}:
 \be
 \label{Lan1}
  \left\{
   \begin{matrix}
u({\bf x})=  \frac {2(c z- r)x}{(c r-z)^2 r},\quad
\,\,\,\,\,\,\,\,\ssk
   \\
    v({\bf x})=  \frac {2(c
     z -r)y}{(c r- z)^2 r},\quad\,\,\,\,\,\,\,\,
    \ssk
   \\
    w({\bf x})=  \frac {2(c r^2 - 2 z r + c
   z^2)}{(c r-z)^2r}, 
   \\
    p({\bf x})=  \frac {4(c z- r)}{(c r-z)^2 r},
    \quad\,\,\,\,\,\,\,\,
    \end{matrix}
 \right.
  \ee
where $r^2=x^2+y^2+z^2 \, (= |{\bf x}|^2)$, and  $c \in \re$ is a
constant such that $|c|>1$. In view of a  strong singularity at
$x=0$, the existence of such a solution demands an extra force at
the origin, so actually \ef{Lan1} is a {\em fundamental solution}
of the stationary operators of the NSEs satisfying, in the sense
of distributions (this computation was already done in Landau's
original work \cite[p.~300]{Lan44}; see also
\cite[pp.~2-9]{Bat70}, \cite[p.~182]{LandM}, or
\cite[p.~250]{Can04K}),
 \be
 \label{delta1}
  \begin{matrix}
 (\uu \cdot \n)\uu + \n p - \D u= b(c) \d(x) {\bf j}, \quad \n \cdot
 \uu =0 \whereA {\bf j}=(0,0,1)^T\,\,\mbox{and} \qquad
 \ssk\ssk\ssk\\
  b(c)= \frac{8 \pi c}{3(c^2-1)} \, \big[2+6
 c-3c(c^2-1) \ln \big( \frac{c+1}{c-1} \big) \big]=
 \frac {16 \pi}c- \frac {32 \pi}{3 c^3}+... \asA c \to +
 \iy.\qquad
  \end{matrix}
  \ee
 Thus, physically speaking, this means that
  \be
  \label{SLinj}
   \mbox{the steady S--L solution demands a permanent fluid
   injection at the origin $\xx=0$.}
    \ee
 Thus,  $\uu_{\rm SL}$ satisfies (this will be
used in the rescaled blow-up variables $y$):
  \be
  \label{ccc1}
    \begin{matrix}
   \uu_{\rm SL}({\bf x}) \sim \frac 1 {c|{\bf x}|} \asA x \to 0, \quad
   \mbox{and, more precisely,}\qquad\quad\ssk\\
     \uu_{\rm SL}({\bf x}) = \frac 1c \, \uu_0 + O\big(\frac 1{c^2}\big) \asA
  c \to + \iy \whereA
      \uu_0=\big(
     \frac{2zx}{r^3}, \frac {2zy}{r^3}, \frac 2
     r\big)^T,\qquad\quad
     \ssk\ssk\\
     |\uu_{\rm SL}({\bf x})| \to \iy \asA c \to 1^+ \,\,\,
     \mbox{on the semiaxis $Oz^+=\{x=y=0,\, z \ge
     0\}$}.\qquad\quad
      \end{matrix}
      \ee
 To illustrate the last property of unbounded steady profiles,
 consider the first component in \ef{Lan1} for a fixed $z>0$ (for $z=0$, the estimate is similar
 with $4 \mapsto 2$ at the end):
  \be
  \label{Last1}
   \tex{
    u({\bf x})\big|_{c=1}= \frac
    {2x}{(z-\sqrt{x^2+y^2+z^2})\sqrt{x^2+y^2+z^2}}
    =- \frac {4x(1+o(1))}{x^2+y^2} \to \iy \asA x^2+y^2 \to 0.
    }
    \ee
 This is rather impressive (and promising for blow-up around): it turns
 out that arbitrarily large vector fields (of special structure)
 can be locked in a singular {\em steady} pattern.

In the spherical coordinates, the flow \ef{Lan1}, being symmetric
about the polar ($Oz$-) axis,  reads (this is the actual Slezkin
\cite{Slez34} and Landau form \cite{Lan44}; see also
\cite[p.~82]{LandM})
 \be
 \label{Lan71}
  \left\{
   \begin{matrix}
     u(r,\th)= \frac{1+ \cos^2 \th-2 c \cos \th}
    {r(c-\cos \th)^2},\qquad\ssk\\
 v(r,\th)= \frac{2 \sin \th}
{r(c-\cos \th)}, \quad w=0,\ssk\\
 p(r,\th)= \frac{4(c \cos \th-1)} {r^2(c-\cos \th)^2}.\qquad\quad\,\,
 \end{matrix}
 \right.
 \ee
 Then, in accordance with the divergence \ef{Last1}, we have
  \be
  \label{Last2}
   \tex{
   v(r,\th)\big|_{c=1} = \frac 2r \, \frac{\cos( \frac \th 2)}{\sin( \frac \th 2)} \to
   \iy \asA \th \to 0.
   }
   \ee

The singular S--L solutions  are $L^2$ locally, but not globally:
 \be
 \label{loc11}
  \uu_{\rm SL} \in L^2_{\rm loc}(\re^3), \quad \mbox{but}
  \quad
\uu_{\rm SL} \not\in L^2(\re^3) \quad \mbox{for all} \,\,\,\,
|c|>1.
 \ee
  It is also worth noting  for further
 blow-up use of \ef{Lan71} that the total mass flux
 through any closed surface around the origin is equal to zero
 \cite[p.~83]{LandM}; cf. the ``vanishing oscillatory property"
 \ef{q5}.
 It turns out
 that  \ef{Lan1} are the only possible  stationary
 homogeneous of degree $-1$, regular except the
 origin
 $(0,0,0)$ solutions of \ef{NS1} \cite{Sv06}; see also earlier result \cite{Tsai98}
 proved under the $z$-axis symmetry assumption.
 One can see that \ef{Lan1} are invariant under the scaling group
 in \ef{ll2}.
 Note that the
 whole set of possible steady singularities of the NSEs is still
 not fully known; even in the class
 \be
 \label{Cl1}
  \tex{
  |\uu(x)| \le \frac {C_*}{|{\bf x}|}, \quad \mbox{where $C_*>0$ is
  a constant}.
 }
   \ee
 However, an essential first step in this direction
 was done in Miura--Tsai \cite{MiuTs08}, who proved: if $\uu$ is any {\em
 very weak} steady solution of the NSE in $\re^3\setminus\{0\}$,
 then
  \be
  \label{WWW1}
  \tex{
   \mbox{$\uu(x)$ satisfies (\ref{Cl1}) for some small $C_*>0$}
  \LongA \uu= \uu_{\rm SL} \,\,\, \mbox{for some $|c|>1$},
   }
   \ee
   where, by \ef{Last1}, \ef{Last2}, $c$ is assumed to be large enough.
In any case, the S--L solutions \ef{Lan1} should play a
 crucial role, since these are expected to be isolated from other singular ones (if any);
  see further details
  in \cite{Sv06, MiuTs08}.
 In addition, there are no smaller singular equilibria (see
 the result and earlier references in \cite{Kim06}): if $\uu, \,
 p \in C^\iy(B_R\setminus\{0\})$, then
  \be
  \label{Kim061}
   \tex{
  |\uu(x)| = o \big( \frac 1{|\xx|}\big) \asA \xx \to 0 \LongA
  \mbox{0 \,\,is removable and $\uu \in C^\iy$ there}.
  }
  \ee
  A similar removable singularity theorem for the nonstationary
NSEs \cite{Koz98}, where \ef{Kim061} is assumed to be valid for a
weak solution $\uu(\xx,t)$ in some set $B_r \times (0,t_0)$,
includes a  smallness condition on $\uu$ in $L^\iy(0,t_0;L^3)$.

Note that in addition to  the symmetry of these solutions about
the $z$-axis, \ef{Lan1} also implies that such a flow does not
exhibit any {\em swirl}, i.e., the angular component $w=u_\var
\equiv 0$ in \ef{Lan71}. This, in view of global existence results
for cylindrically axisymmetric (irrotational) flows  \cite{Ukh68,
Lad68A} (see \cite{Zaj07} for details),
 indeed,
reduces the chances to get a reasonable blow-up pattern moving
along the corresponding quasi-stationary manifolds (a
finite-dimensional centre, if any,  or a stable one of infinite
dimension). Nevertheless, the situation with the steady manifold
induced by \ef{Lan1} is not that hopeless. For instance, instead
of the ``standard" similarity swirl given by the log-law \ef{c11},
one can use a slower rotational mechanism by setting
 \be
 \label{sl1}
 \var= \mu + \sigma  \kappa(t) \whereA (-t) \kappa'(t) \to 0 \asA t \to
 0^-.
  \ee
 Then,  the stationary term $\sigma  \uu_\mu$ in \ef{c7N} is replaced
 by the asymptotically vanishing one,
  \be
  \label{sl2}
  \sigma  \uu_\mu \mapsto \sigma  (-t) \kappa'(
 t) \uu_\mu \to 0
   \quad \big(\mbox{e.g., \, $\kappa(t)= -\frac {\ln(-t)}{|\ln|\ln (-t)|\, |^\d}$ for $t
   \approx
   0^-$, with $\d>0$} \big),
   \ee
   on bounded smooth orbits. Therefore, for a suitable class of solutions, passing to the limit along a
   sequence $\{\t_k\} \to +\iy$ by setting $\t \mapsto \t_k + \t$
   will fix us the previous limit irrotational equations admitting the
   ``stationary"  singularity \ef{Lan1}. Of course, this leads to a
   delicate matching procedure of connecting the slow whirling
   flow characterized by \ef{sl1} with  Slezkin--Landau's  ``quasi-steady"
    solutions \ef{Lan1}, which actually should determine the
   function $\kappa(\t)$.
 We return to other extensions and use of \ef{Lan1} in
Section \ref{SLanFin}.


\subsection{On some hypothetical extensions}

On the other hand, it is principal for blow-up patterns to know
whether the rescaled NSEs admit other solutions with swirl.
 In general, classification of singular states for the NSEs is a difficult open
 problem (see \cite{Sv06}), though it is a necessary step for
 understanding of existence/nonexistence of finite energy blow-up
 patterns for such dynamical systems.

 For
instance, it is crucial to check existence of singular solutions
with the TW dependence in the angular direction, with
 \be
 \label{vv10}
 \varphi \mapsto \varphi + \sigma  \, t.
  \ee
 This gives the nonlinear eigenvalue problem (we assume that the
 singularity has not been essentially changed, so we may keep  the same its
  $\d$-interpretation, for simplicity)
 \be
 \label{delta2}
  \begin{matrix}
 \sigma  \uu_\varphi+(\uu \cdot \n)\uu + \n p - \D u= b(\sigma ) \d(x) {\bf j}, \quad \n \cdot
 \uu =0.
  \end{matrix}
  \ee
For $\sigma _0=0$ and $b(0)=b(c)$, this gives the S--L solutions
satisfying
 \ef{delta1}. Do other eigenvalues $\{\sigma _k \in \re\}$ exist?
 A partial negative answer is available (see \cite{MiuTs08} as a guide), but it does
 not cover the whole range.
 The corresponding singular states $\{\uu_{\rm SL}^k\}$ then can
 be used for constructing various blow-up patterns.

 In the rescaled NSEs \ef{ww1}, \ef{HHH21}, one needs to know existence of
 regular states with the blow-up angular dependence \ef{c9} (for
 $\s=0$, no solutions \cite{Nec96})
  \be
  \label{delta3}
  \sigma  \tilde \uu_\mu + \HH(\tilde \uu)=0 \inB \re^3, \quad \tilde \uu
  \in L^2_{\rm loc}(\re^3) \quad (\tilde \uu \not
  \in L^2(\re^3)\,\,\,\mbox{in general}.)
   \ee
Here $\tilde \uu(y)$ may be assumed to be bounded at $y=0$, so
that
 it is natural to suppose  that, along some subsequences or
 continuous branches (see Section \ref{S5.5} for Type II using),
 \be
 \label{delta4}
  \tilde \uu(0)= \tilde \CC \andA \tilde \uu \to \hat \uu_{\rm
  SL}(c)
  \asA |\tilde \CC| \to + \iy \quad(|c|>1)
   \ee
uniformly on compact subsets in $\re^3 \setminus \{0\}$. In other
words, the singular S--L solutions could serve as ``envelopes"
 of the regular ones (this is obscure and questionable). Then, it seems, at least a
 two-parameter family of such regular profiles might be expected.
 In any case,
 the extra nonlinear eigenvalues $\{\sigma _k\}$ (or other perturbation-like mechanisms)
  can essentially
 affect the structure of those hypothetical regular swirling
 states $\{\tilde \uu_k\}$ solving \ef{delta3}.
We again state that these speculations are made under the clear
absence of any clue on existence of those singular quasi-steady
patterns with various blow-up swirls. Both clear mathematically
existence or nonexistence conclusions are desperately needed here.

 \subsection{On possible blow-up  patterns: towards ``swirling tornado"
 about 
 S--L
 singular equilibria.
   Inner Region I}
 \label{S7.3}

We now perform first steps of the blow-up strategy according
to the blow-up scenario developed in Section \ref{S1.5} now
applied to the S--L singular stationary solution. So we are going
to check whether a ``blow-up swirling tornado" (not a twistor,
since we do not always apply directly the logTW mechanism, though
it is not excluded) can appear by approaching the
$L^\iy$-singularity at $\xx=0$  by drifting as $t \to T^-$ by
``screwing in" about this singular equilibrium structure. There is
a clear suspicion that this can hardly happen without, according
to  \ef{SLinj}, a permanent fluid injection, which is not
available for bounded $L^2$-solutions. Note that the hypothetical
blow-up tornado is going to occur very fast, during a miserable
time scale, when the total injection is negligible.
 Anyway,
 this physical negative motivation
 should find some clear mathematical issues in support or not.

 Recall that the
behaviour in {\em Inner Region I} is assumed to be characterized
by the linearized problem \ef{As12}.
Thus, following the same standard principles of matching,
 linearization of the
general equation \ef{ww1} about $\uu_{\rm SL}$ via \ef{z6}
 yields, in  {Inner Region I}, the equation \ef{z7} with quite a tricky linear
 pseudo-differential
 (integro-differential) operator obtained from \ef{HHH21}, which  in
 $C_0^\iy(\re^3\setminus\{0\})$ reads
 \be
 \label{Lan2}
  \begin{matrix}
\HH'(\hat \uu_{\rm SL})\YY =  \big(\BB^*- \frac 12 \, I \big) \YY
- \mathbb{P} \CC \YY, \,\,\, \mbox{where} \,\,\,
 \CC \YY= {\rm div}\,\big( \hat \uu_{\rm SL}\otimes \YY+ \YY
   \otimes\hat \uu_{\rm SL}\big)\qquad
\ssk\ssk\\
   \equiv
     (\hat
\uu_{\rm SL} \cdot \n )\YY + (\YY \cdot \n) \hat \uu_{\rm SL},
\,\,\,\, \mbox{so}\,\,\,\,\HH'(\hat \uu_{\rm SL})\YY = \big(\BB^*-
\frac 12 \, I \big) \YY -  \CC \YY \qquad
 \ssk\ssk\\
 + C_3
\int\limits_{\re^3} \frac {y-z}{|y-z|^3}\,\, \sum\limits_{(i,j)}
\,\big( \hat u_{{\rm SL}z_j}^i Y^j_{z_i} +  \hat u_{{\rm SL}z_i}^j
Y^i_{z_j}\big) \equiv {\bf J}_1 \YY + {\bf J}_2 \YY.\qquad
 \end{matrix}
 \ee
 Here, ${\bf J}_1$ includes all the local differential terms,
 while ${\bf J}_2$  the integral one.
 Note that, in general, according to \ef{delta3}, the above
 linearized operator must include the operator
  \be
  \label{kkk1}
  ...-\s D_\mu \YY +...\, , \quad \fbox{$\mbox{so (\ref{Lan2}) contains
  {\bf two} free parameters: $c$ and $\s$,}$}
  \ee
  with a nonlinear eigenvalue $\s$ coming from the stationary
  problem already containing the parameter $c>1$.
   In a most general setting, the linearization is performed
  relative to any unknown S--L-type singular steady states
  delivered by \ef{delta3}. For simplicity, we will continue to
  deal with the standard S--L solutions, actually meaning that,
 in addition to other branching ideas, one needs also (cf. Section \ref{S6.1}):
  \be
  \label{brsig}
   \fbox{$
  \mbox{to develop branching of eigenfunctions from the logTW speed $\s=0$,}
  $}
 \ee
 though possibly this is not an issue since  demanding too much from the operator \ef{Lan2}
 ($\l=0$ must be in the spectrum).
 On the other hand, even if \ef{brsig} makes sense,
 it is then also plausible that such a branching is
 available (or not) for $|\s|$ small only, and proper eigenfunctions may occur
 via
 ``saddle-node" bifurcations for large $|\s|$, which are very
 difficult to detect.
In any case, further detailing of these arguments seems excessive,
since, even
 without (\ref{brsig}) and/or others, the study gets very complicated with already too many
 open problems to appear.

 In all the cases,
  we observe in (\ref{Lan2})
  singular unbounded coefficients as $y \to \infty$ as well as $y \to
  0$.
  A proper treatment of $y=\iy$ is settled via the space $L^2_{\rho^*}(\re^3)$, as usual.
  The singular point $y=0$ is much worse.

\ssk

\noi\underline{\sc Local operator: discrete spectrum by improved
Hardy--Leray inequality}. Namely, using \ef{ccc1} yields that the
singularity of the potentials in the local terms,
 \be
 \label{loc1}
  \tex{
 \CC \YY \equiv
  (\hat
\uu_{\rm SL} \cdot \n )\YY +(\YY \cdot \n) \hat \uu_{\rm SL}
 \sim  (\frac y{c|y|^2} \cdot \n)\YY +(\YY \cdot \n) \frac 1{c|y|},
 }
 \ee
is well covered by Hardy's classic inequality  at least for all
large $c>0$, since both terms ``act" like the inverse square
potential. Thus,  according to \ef{WW6}, we need
 \be
 \label{loc2}
  \tex{
  {\bf J}_1 \YY \sim \D \YY +
    \frac \YY{c|y|^2} \asA y \to 0 \LongA \sim \frac 1c  \le
   \frac{(N-2)^2}4\big|_{N=3}= \frac 14.
   }
   \ee
   Note that, for axi-symmetric solenoidal flows, the Hardy
   optimal
   constant $c_{\rm H}=\frac{(N-2)^2}4$ becomes better, and, as shown in Costin--Maz'ya
   \cite{CosMaz08}, can be replaced by
 \be
 \label{Maz1}
 \tex{
  c_{\rm H,\,axi-sol.}=\frac{(N-2)^2}4\, \frac{N^2+2N +4}{N^2+2N-4}\big|_{N=3}= \frac{25}{68}=0.3676... > \frac 14.
 }
 \ee
  In other words, the operator $\D \YY + \frac {25}{68}
  \frac{\YY}{|x|^2}$ in the axi-symmetric solenoidal field admits
 a proper setting in $L^2$ for $x \approx 0$. Note that the inequality
 \be
 \label{Ler223}
  \tex{
  \frac 14 \, \int_{\re^3} \frac {|\YY|^2}{|x|^2} \le
  \int_{\re^3}|\n \YY|^2 \quad (N=3)
  }
  \ee
 already appeared in the same Leray's pioneering  paper in 1934 \cite{Ler34}.
 Thus, $\frac 14$ in \ef{Ler223} can be replaced by $\frac{25}{68}$ that  improves the admissible range of $c$'s
  according to \ef{loc2}. Hence, we obtain the first important
  and even principal conclusion:
  \be
  \label{CC11}
   \tex{
   \mbox{for any $c \gg 1$, the local linearized operator ${\bf
   J}_1$ is well-posed.}
   }
   \ee
Therefore, this local part of the linearized operator has a
compact resolvent, a discrete spectrum in the same weighted
$L^2$-space, etc., though several questions including completeness
and closure of the eigenfunction set and others remain open and
difficult in this non self-adjoint case. Moreover, it follows from
\ef{Lan2} and \ef{loc1} that  for $c \gg 1$, with a proper
functional setting, ${\bf J}_1$ has eigenvalues that are close to
those in \ef{bbb1} for $\BB^*- \frac 12\, I$. In addition, by
classic perturbation theory (see Kato \cite{Kato}), it can then be
shown that
 \be
 \label{Br1}
  \begin{matrix}
  \mbox{eigenfunctions $\{\hat \vv^*_\b\}$ of ${\bf J}_1$ can be obtained
  from generalized} \ssk\\
   \mbox{Hermite polynomials $\{\vv_\b^*\}$ via branching at $c=+\iy$}.
   \end{matrix}
   \ee
   This justifies existence of infinitely many real eigenvalues of
   ${\bf J}_1$, at least, for $c \gg 1$ (global extension of
   bifurcation $c$-branches to finite values of $c>1$ is a difficult open problem).

 Since in \ef{ww1}
$\mathbb{P}$ is a projector in $(L^2)^3$,  $\|\mathbb{P}\|=1$, it
seems reasonable to expect that some part of discrete spectrum may
be governed by the local differential operator ${\bf J}_1$.
 Consider the full linearized equation
  \be
  \label{YY12}
   \tex{
   \YY_\t = \big(\BB^*- \frac 12\, I\big) \YY - \mathbb{P}\CC \YY +
   \DD(\YY) \equiv \HH'(\hat \uu_{\rm SL})\YY + \DD(\YY),
   }
   \ee
   where $\DD$ is a quadratic perturbation. Assuming that there
   exists a proper eigenvalue $\l_\b \in \s_{\rm p}( \HH'(\uu_{\rm SL})) \ne
   \emptyset$, with ${\rm Re}\, \l_\b \le 0$, we expect the
   behaviour $\YY(\t) \sim \cc\,{\mathrm e}^{\l_\b \t} \hat\vv_\b^*$ for
   $\t \gg 1$, if ${\rm Re}\,\l_\b<0$, or (as in Section \ref{S3.8}, $\YY(\t) \sim \cc_\b(\t)\hat \vv_\b^*$
    close to a centre subspace, if ${\rm Re}\,\l_\b=0$. [Both to be matched with more regular flows
    around.]
Substituting this into \ef{YY12} yields the eigenvalue problem,
where the nonlocal term implies a tricky setting,
 \be
 \label{YY16}
 \tex{
 \big(\BB^*- \frac 12\, I\big)\hat \vv_\b^* - \mathbb{P}\CC
 \hat\vv_\b^*= \l_\b \hat\vv_\b^* \quad\big(\mbox{$\hat\vv_\b^* \in \hat H^2_{\rho^*}(\re^3)$, at
 least.} \big)
 }
 \ee
 Again, it should be noted that the angular operator \ef{kkk1}
 must enter this eigenvalue problem. We claim that
  this extra technicality can be tackled by separation in the
  angular variable, though this leads to more difficult calculus
  and a proper adaptation of further conclusions are not easy but
  doable.
 However, this shows
 a natural opportunity to
be retained within local spectral theory for ${\bf J}_1$.
Since, as we have seen, $\mathbb{P}\CC \to 0$ as $c \to +\iy$,
this demands solving the following problem for generalized Hermite
polynomials $\{\vv_\b^* \in \Phi_k^{*3}, \, k=|\b|\}$ as in
Section \ref{SPer}:
 \be
 \label{Br2}
 \fbox{$
\mbox{$c \gg 1$: \,\, for which  $\vv_\b^*$,  the linearized term
 $\CC \vv_\b^*$
is or ``almost" solenoidal?}
 $}
   \ee
 Obviously,   \ef{YY16} shows that if $\CC \vv_\b^*$ is fully divergence-free,
 then   $\mathbb{P}=I$ on it, so the local operator ${\bf J}_1$ takes the full
   power, and this settles the spectral problem for  $c \gg 1$. As
   we expect, this will allow us to start a series-like expansion
   in terms of the small parameter $\frac 1c$
   of the necessary div-free eigenfunction of the full problem
   \ef{YY16}.

 We did not check if such eigenfunctions are
actually existent, but, in a whole,  the problem \ef{Br2} looks
doable. Indeed, we do not expect to find such polynomials
$\vv_\b^*$ satisfying \ef{Br2} for sufficiently  small  $|\b|=0,1$
or 2, etc. The first such $l=|\b|$, for which \ef{Br2} holds with
necessary accuracy, would show a complicated geometry (with
remnants of swirl, axis precession, and so on), which would
eventually affect the global structure of a possible blow-up
pattern.
 Once we have found such solenoidal generalized Hermite
polynomials satisfying \ef{Br2}, this starts the not less simpler
procedure of constructing those ``twistors" by matching with more
regular flows close to the singular point $y=0$ and outside the
pattern blow-up region (the latter is supposed to be simpler). Of
course, it is too early discussing any seriously this nonlinear
matching matter (but we will, noting also that we have already
mentioned some formal principles of such procedures that are well
established in blow-up R--D theory). But indeed the linearized
spectral analysis about the rescaled S--L solution is assumed to
show some basics of formation of such mysterious blow-up patters
including their strong angular motion (the eigenfunctions are
supposed to be essentially non-radial), ``precession" features of
their axis of rotation (to be discussed in greater detail), and
even their possible merging into a closed path of several
individual twistors (hopefully, not infinite).

\ssk

\noi\underline{\sc Nonlocal operator: an obscure functional
setting}. Anyway, despite speculations concerning the more or less
 standard local  operator ${\bf J}_1$ (with basics related to
 Hardy--Leray inequalities), one cannot avoid
analyzing the nonlocal integral operator ${\bf J}_2$. The spectral
questions for the nonlocal operator in \ef{Lan2} become  more
involved.
Consider ${\bf J}_2$ in $C_0^\iy(\re^3 \setminus \{0\})$, where,
on integration by parts,
 \be
 \label{Lan2N}
  \begin{matrix}
 {\bf J}_2 \YY \equiv  C_3
\int\limits_{\re^3} \frac {y-z}{|y-z|^3}\,\, \sum\limits_{(i,j)}
\,\big( \hat u_{{\rm SL}z_j}^i Y^j_{z_i} +  \hat u_{{\rm SL}z_i}^j
Y^i_{z_j}\big)
  \ssk\\
  =
 - C_3
\int\limits_{\re^3} \, \sum\limits_{(i,j)} \,\big[\, \big(\frac
{y-z}{|y-z|^3}\,\, \hat u^i_{{\rm SL}z_j}\big)_{z_i} Y^j +
\big(\frac {y-z}{|y-z|^3}\,\, \hat u^j_{{\rm SL}z_i}\big)_{z_j}
Y^i\big].
 \end{matrix}
 \ee
We obtain singular integral operators with the kernels, formally
exhibiting the behaviour
 \be
 \label{KK11}
  \tex{
   K(y,z)=K_1(y,z)+K_2(y,z) \whereA |K_1| \sim \frac 1c\, \frac 1{|y-z|^2} \,\, \frac 1{|z|^3},
   \,\, |K_2| \sim \frac 1c\, \frac 1{|y-z|^3} \,\, \frac 1{|z|^2}.
   }
   \ee

   \noi\underline{\em First view of $K_1$: too singular, a
   condition at 0 needed}.
The kernel $K_1$  is standard polar relative to the first
multiplier having $\sim \frac 1{|y-z|^\mu}$, with $\mu=2<N=3$.
However,  the second multiplier, at first sight, exhibits  a
non-integrable singularity $\frac 1{|z|^3}$ at $z=0$ in $\re^3$
demands special conditions at the origin e.g., of the type
(otherwise, the integrals in \ef{Lan2N} ought to be understood in
a generalized  {\em v.p.} sense, see below)
 \be
 \label{spcond1}
 \tex{
 |\YY(z)| =O( \rho(|z|))\to 0 \asA z \to 0 \whereA \int_0 \frac{\rho(r)}r\, {\mathrm
 d}r < \iy.
 }
 \ee
 On the one hand, it seems clear that
 any such  condition \ef{spcond1} on the behaviour at the origin can reduce a hope to find a proper
eigenfunction of \ef{Lan2}. On the other hand, it is well known
that eigenvalues of pseudo-differential operators can be of
infinite multiplicity (e.g., for those with constant symbol),
i.e., the admitted behaviour at singularity points can be
extremely various (unlike the differential operators, for which
such a behaviour is restricted by their orders and asymptotics of
the coefficients). However, even for the local differential
operators, imposing such conditions is not that hopeless, as the
following example shows:

\ssk

\noi{\bf Example: differential operator with a locally
non-integrable potential.} Consider for simplicity  a scalar
differential operator associated with ${\bf J}_1$ in \ef{Lan2},
with an extra inverse cubic  potential ($\not \in L^2_{\rm
loc}(\re^3)$) as in \ef{KK11}:
 \be
 \label{Pot1}
  \tex{
  {\bf E} = \BB^*- \frac 12\, I -\CC - \frac A{|z|^3}\, I \quad (A \not = 0).
   }
   \ee
 Obviously, looking for possible eigenfunctions, the behaviour of
 solutions close to the singular origin 0 is defined by the
 principal part, so that, as $z \to 0$,
  \be
  \label{Pot2}
  \tex{
  {\bf E} \psi^*=\l \psi^* \LongA
   \D \psi^* -\frac A{|z|^3}\,\psi^*+...=0
   \LongA \psi^*(z) \sim \o\big( \frac z{|z|}\big)\,|y|^{-\d} {\mathrm
   e}^{a|z|^\a},
    }
    \ee
    where  $\o$ on the unit sphere ${\mathbb S}^2$ in
    $\re^3$ makes the angular separation.
  Substituting the last WKBJ-type expansion term into the operator
 $(\psi^*)'' + \frac 2{|z|} (\psi^*)'- \frac A{|z|^3}\,\psi^* = 0$ yields
 \be
 \label{Pot3}
  \tex{
   \a=- \frac 12, \quad \d=\frac 34, \andA a^2 = 4A \LongA a_\pm= \pm 2 \sqrt A \forA
   A>0
   }
   \ee
   ($A>0$ corresponds to  ``stable" potential, i.e.,  monotone
  principal operator).
Thus, close to 0, there exist two asymptotic bundles of solutions
with essentially different behaviours:
 \be
 \label{Pot4}
  \tex{
 \psi^*_-(z) \sim |z|^{-\frac 34}\, {\mathrm
   e}^{-  2 \sqrt {\frac A {|z|}}} \to 0 \andA
    \psi^*_+(z) \sim  |z|^{-\frac 34}\, {\mathrm
   e}^{ 2 \sqrt {\frac A {|z|}}} \to \iy.
   }
   \ee
 Therefore, posing the condition (cf. \ef{spcond1})
  \be
  \label{Pot5}
  Y(0)=0
   \ee
  does not spoil at all the eigenvalue problem for ${\bf E}$ and
  just eliminates those singular behaviours at 0 (half of all of
  them) making the spectrum discrete. In linear operator theory, the same
  is usually and naturally
  done by assuming $Y \in L^2(B_1)$, which
  nevertheless has come true after certain asymptotic analysis
  presented above.
  For the principal radial ordinary differential operator, this
  means that the singular point $|z|=0$ is in the {\em limit point
  case} (an index characterization is also available for the elliptic
  setting).

   On the contrary (and this case seems more correctly describes the nature of such an
   ``unstable" potential in \ef{KK11}, though not that obviously), if $A<0$, then \ef{Pot3} yields
  $a_\pm= \pm 2 {\rm i}\, \sqrt {|A|}$, so that both bundles are
  singular and
  oscillatory: as $z \to 0$,
   \be
   \label{Pot6}
    \tex{
    \psi^*_-(z) \sim  |z|^{-\frac 34}\, \cos\big(2 \sqrt {\frac A {|z|}}\big),
    \quad \psi^*_+(z) \sim  |z|^{-\frac 34}\, \sin\big(  2 \sqrt {\frac A {|z|}}\big), \quad
  \psi^*_\pm \in L^2_{\rm loc};
 }
 \ee
the {\em limit circle case} for the operator $\D - \frac
A{|z|^3}\,I$. Then the discrete spectrum of ${\bf E}$ is also
achieved by posing special conditions at the singularity,
\cite{Nai1}, and then any eigenfunctions exhibit the singular
behaviour \ef{Pot6}. It is another hard problem to check whether
\ef{Pot6} allows to match such an oscillatory bundle with smoother
$L^\iy$-flow near the origin $z=0$ for $\t \gg 1$; see below.
 Therefore, in the case of total
oscillatory bundle \ef{Pot6}, the condition such as \ef{Pot5}
makes no sense and, evidently, destroys the eigenvalue problem
 giving $\s({\bf E})= \emptyset$.

 Note finally, that the oscillatory-type behaviour similar to
 \ef{Pot6}, in general, is not an absolute obstacle for getting
 from such ``eigenfunctions" a proper blow-up pattern. We discuss
 this for some R--D equations in Section \ref{S7.7}. Then, for the behaviour like
 \ef{Pot6}, the integral operator as in \ef{Lan2N} should be
 understood in a (canonical) regularized sense, see below.


\ssk

\noi\underline{\em Singular kernel $K_2$ is of Calder\'on--Zygmund
type}.
 The integral operator with the kernel $K_2$ in \ef{KK11} is singular.
However, in view of the divergence form of the operator $\CC$ in
\ef{Lan2},
 there is a hope (not yet fully justified) that
 this kernel falls into the scope of
 Calder\'on--Zygmund's classic result (1952) \cite{Cald52} saying that such an operator can
 be bounded in $L^p(\re^3)$ (in the local sense, meaning that, as usual,
  we cut-off the infinity by the appropriate weight $\rho$) for some $p>1$;
 see \cite[Ch.~5, \S~5]{Tay} for modern overview and references. Most of the text-books on pseudo-differential
 operator theory quote  such
 fundamental results as being its origin; see \cite[p.~278]{Mikhlin}. Note that the conditions of
  the fundamental
 Calder\'on--Zygmund result in the standard form \cite[p.~16]{Tay}
 do not directly cover  the singularities of $K_2$, so that an
 extra hard work is essential. In addition, one surely needs
 extensions of such boundedness results to $L^p$-spaces of
 vector-valued functions $\YY$ (such a study was already initiated by Calder\'on
 himself in 1962); these questions being  well
 understood; see references and a survey on further operator-valued issues in \cite{Hyt07}.

\ssk

  \noi\underline{\em Back to $K_1$: the operator at 0 in a {c.r.}
  sense}.
Indeed, such a possibility to ``overrun" (if possible) the very
restrictive condition \ef{spcond1} is to take into account the
specific divergence part of the $ \frac 1{|z|^3}$-terms in
\ef{Lan2N}. To explain this, let us fix $i=j=1$ in the first term,
where for $z \approx 0$ and $y \not =0$, for $Y^1 \in C_0^\iy$, we
 have that the singular part in the regularized value
 sense vanishes:
 \be
 \label{Int11}
  \tex{
\sim C_3 \frac {y}{|y|^3} Y^1(0)  \int_{B_\e}
     \big(\hat u^1_{\rm SL}\big)_{z_1z_1}\, {\mathrm d}z=0,
 }
 \ee
 since the second-order $z_1$-derivative, according to \ef{Lan1},
 is an odd function in $z_1$. Namely, the corresponding
 part of the integral operator in \ef{Lan2N} acts like the
 standard distribution ${\mathcal P} \frac 1{z_1^3}$ with the
  regularization
  \be
  \label{Int12}
   \tex{
   \langle {\mathcal P} \frac 1{z_1^3}, \varphi \rangle=
   \int_0^\iy \frac{\varphi(s)- \varphi(-s)-2 s
   \varphi'(0)}{s^3}\, {\mathrm d}s \forA \varphi \in
   C_0^\iy(\re),
   }
   \ee
   which is {\em canonical}, {\em c.r.} (i.e.,
   this regularization keeps
   the linear properties of the functionals, as well as the
   differentiation).
A similar {\em c.r.}-property is observed for the whole $Y^1$-term
appeared in \ef{Lan2N}   that consists of three members: $2[(\hat
u^1_{\rm SL})_{z_1z_1} + (\hat u^2_{\rm SL})_{z_1z_2}+ (\hat
u^3_{\rm SL})_{z_1z_3}]$. Note that all of them are odd relative
$z_1$, so the indefinite integral is even. It seems that direct
$L^p$-theory is not applicable to such {\em c.r.}-integral
operators that thus deserve further study. Elsewhere, in other,
less singular terms, the standard {\em v.p.}-sense {\em c.r.}
 $$
  \tex{ \langle {\mathcal P} \frac 1{z_1},
\varphi \rangle=\int_{-\iy}^{\iy} \frac{\varphi(s)}s \, {\mathrm
d}s \equiv \int_{0}^{\iy} \frac{\varphi(s)-\varphi(-s)}s \,
{\mathrm d}s
 }
  $$
   in the Cauchy sense can occurs, leading to more
standard integral operators with typical singularities  of the
{\em Hilbert transform} in $z_1$ (a Calder\'on--Zygmund operator
in $\re$), etc.

\ssk

Continuing the fruitful idea of eigenfunction branching at
$c=+\iy$ as in the local case above, we obtain an open problem for
the whole pseudo-differential operator \ef{Lan2}:
 \be
 \label{Pr55}
  \fbox{$
 \mbox{In which operator topology, does \,\,$\HH'(\hat \uu_{\rm SL}) \to
 \BB^*- \frac 12\, I$ as $c \to +\iy$?}
 $}
  \ee
Or the limit always contains a certain ``singular" part? The
answer is principal.


  Note that any extra condition  such as \ef{spcond1} can be
 inconsistent with admitted ``weakly singular" behaviour
 corresponding to the local operator ${\bf J}_1$
 composed from $\BB^*$ and $\CC$ in
 \ef{loc1}
 provided that the latter ones   are ``dominant".
 Then this case becomes rather standard and, as in blow-up  R--D theory (see Section
 \ref{S1.5}),
  is governed by improved Hardy's inequalities discussed above.

 However, if the nonlocal
  operator ${\bf J}_2$ is ``dominant" close to the
 origin\footnote{Actually, as seen from $\mathbb{P}$ in  \ef{ww1}, both local and nonlocal operator parts are of
 a  ``similar power".}, which
 seems unavoidable and  represents the main difficulty of blow-up in the NSEs (cf.
 the uncertain polynomial micro-structure of multiple zeros detected in Section
 \ref{SHerm} by the same reason), then the conditions such as
 \ef{spcond1} can be natural, giving rise to the necessity of different spectral
 theory, which is entirely unknown and represents an open problem.
We do not know whether the linearized pseudo-differential operator
\ef{Lan2} with  such a strongly singular kernel (and hence with a
nonregular full symbol; see an overview in  \cite[p.~45]{Tay}) has
a nontrivial discrete spectrum in a proper functional setting
involving conditions at the singularity and Lorentz--Marcinkiewicz
or Zygmund-type  spaces (see \cite[p.~37]{Tay} with applications
to nonlinear PDEs in this Taylor's volume).

 Note that, plausibly, by some mysterious
reason, such a pseudo-differential operator
 may admit just a single proper eigenfunction (for further use in
 blow-up matching, with a positive or negative conclusion, it does
 not matter), which cannot be detected in principle by any general advanced
spectral theory.
 This
is a difficult problem, where numerics for such involved nonlocal
operators, can be key, though a definite justified answer: ``yes"
or ``no" to existence of proper eigenfunctions and hence
centre-stable eigenspaces, can be very questionable.

Overall, spectral properties of the pseudo-differential operator
\ef{Lan2N} are of great demand (its symbol is not still well
understood) and can be key for existence/nonexistence of proper
blow-up patterns on centre-stable manifolds created by  the S--L
exact steady singular solutions. Actually, even the negative
result on nonexistence would play a role for concentrating on
other more involved scenarios of blow-up to be focused on.
 Of course, as usual, since the S--L
solutions \ef{Lan71} are axi-symmetric with no swirl, our analysis
is assumed to include blow-up swirling mechanism as in \ef{c11},
so that the resulting operators  contain typical angular terms
like $-\s D_\mu$ as in \ef{kkk1}, with, possibly, sufficiently
large angular speeds $\s$ (this assumes complicated separation
angular techniques).
Thus,  {Questions (i)--(vi)} from Section \ref{S3.8} do deserve
further study in this case, and, with a certain luck, will give a
first insight into  blow-up singularities for the NSEs.

Let us   note that, after matching of the linearized patterns on
the stable (centre) manifold for ${\bf J}_{1,2}$ with more regular
flow close to the origin (this is even more difficult open
problem, which makes no sense if spectral theory is still
unavailable; Sections \ref{S1.5} and future \ref{S7.5} may be
consulted for an idea), a typical shape of the resulting patterns
(according to \ef{spcond1} or other restrictions) will have a form
of {\em swirling ``tornado" about a S--L singular steady
solution},
which is self-focused onto the origin $x=0$ as $t \to T^-=0^-$.

\ssk

\noi\underline{\em Final conclusion: spectral results are expected
by branching at $c=+\iy$}.
 Despite various difficulties and open
problems already detected, our final conclusion is not fully
negative.
 Though the extra condition such as \ef{spcond1} looks rather
 restrictive and even frightening, the branching idea at $c=+\iy$
correlated with this  well, since:
 \be
 \label{YY65}
 \fbox{$
\mbox{there exist infinitely many  Hermite
 polynomials $\vv^*_\b$ satisfying (\ref{spcond1}).}
 $}
 \ee
For instance, \ef{HP1} implies that those are $\vv^*_{11}$,
$\vv^*_{12}$, $\vv^*_{13}$, $\vv^*_{24}$, $\vv^*_{25}$,
$\vv^*_{26}$, $\vv^*_{27}$, $\vv^*_{28}$, etc. However, obviously
conditions such as \ef{spcond1} can mean that $\hat \uu(y,\t)$ is
always singular at $y =0$, so \ef{Br2} is indeed more preferable.

 Thus, within the positive reflection of
the hard demand \ef{Pr55} and in view of the desired feature
\ef{YY65},
 we expect
 that
\be
\label{YY13}
 \mbox{for $c \gg 1$, $\HH'(\hat \uu_{\rm SL})$ may have several real
 eigenvalues close to $\l_\b=- \frac{|\b|}2$,}
  \ee
 and the total number of those gets infinite as $c \to +\iy$.
 Indeed, in this limit, we observe a ``convergence" to the Hermite
 operator $\BB^*- \frac 12 \, I$. In other words, the above
 analysis makes it possible  to start a real procedure of checking
 whether at least one from this countable set of linearized
 structures
 admits a proper matching to get a finite energy blow-up pattern
 for the NSEs (or all of them are hopeless, which nevertheless would not prove
  nonexistence of blow-up since there are other blow-up scenarios).
  In other words,
   \be
   \label{YY14}
   \mbox{S--L solutions with $c \gg 1$ are not still
   forbidden
   for blow-up evolution around,}
 \ee
 while  extensions of the branches to finite $c>1$ is even more promising, but
 indeed extremely difficult.
Then the still mysterious matching procedure gets principal, which
itself can cross out all the previous ``spectral" speculations and
illusive  achievements.

\subsection{First application of spectral theory: towards periodic
blow-up patterns}
 \label{S.7.3AH}

As a first elementary application of the above spectral
discussion, returning to simpler periodic blow-up orbits inducing
\ef{ww5}, we pose a straightforward question on the Andronov--Hopf
(A--H) classic scenario of bifurcation of periodic orbits from the
singular equilibrium $\hat \uu_{\rm SL}$. Namely, we state:
 \be
 \label{AH1}
  \begin{matrix}
 \mbox{to check if an A--H bifurcation
 can occur at some $\s=\s_{\rm AH}$, $c=c_{\rm AH}$} \qquad\quad
  \ssk\\
 \mbox{for the operator (\ref{Lan2}), (\ref{kkk1}), i.e.,
  $\exists \,\, {\rm i}\, \o \in \s(\HH'(\hat \uu_{\rm
 SL})+\s_{\rm AH}D_\mu)$, $\o \ne 0$}.\qquad\quad
 \end{matrix}
  \ee
 The asymptotic behaviour of the corresponding eigenfunction
 $\vv^*_\b(y)$ as $y \to 0$ is also of crucial importance to get
 a periodic pattern by matching with the regular bounded
 flow for $y \approx 0$. In case of both positive answers, this
 would lead to a periodic twistor blow-up pattern for $\s \approx
 \s_{\rm AH}$ that gives rise to the $\o$-limit set \ef{ww5}.
 Respectively, existence of multiple eigenvalues ${\rm i}\,
 \o_1$,..., ${\rm i}\,
 \o_n$ may lead to quasi-periodic orbits obeying \ef{cou1}.
These are hard open hypothetical  questions.

Let us finish this brief discussion with a negative result
concerning bifurcation from 0:

\begin{proposition}
 \label{Pr.Bif1}
For the rescaled equation $(\ref{ww1})$ with the swirl operator
$(\ref{kkk1})$,
 an A--H bifurcation from 0 is impossible.
 \end{proposition}

 \noi{\em Proof.} Consider the eigenvalue problem for the
 linearized operator about 0:
  \be
  \label{jh1}
   \tex{
  \big(\BB^*- \frac 12\, I\big) \vv^* - \s D_\mu \vv^*= \l \vv^*
  \inB L^2_{\rho^*}.
  }
  \ee
Using the symmetry \ef{mm1} and divergence of $D_\mu$ on ${\mathbb
S}^1=(0,2\pi)$, we get
 \be
 \label{jh2}
  \tex{
  \frac{\l + \bar \l}2\,
  \|\vv^*\|_{L^2_{\rho^*}}=-\|\n \vv^*\|_{L^2_{\rho^*}}- \frac 12\,
  \|\vv^*\|_{L^2_{\rho^*}} \LongA {\rm Re} \, \l<0 \,\,\, \mbox{for any $\s \in \re$}.
  \qed
   }
   \ee


\subsection{Inner Region II: matching with smoother solenoidal flow
near the origin}

Assume that the previously posed spectral problem has been solved
successfully, so we have found the actual Inner Region I, where
\ef{As12} holds. For simplicity, we assume a stable subspace
behaviour:
 \be
 \label{PPP1}
  \tex{
  \mbox{Inner Region I}:
  \quad \hat \uu(y,\t) \sim \hat \uu_{\rm SL}(y) + {\mathrm e}^{\l_\g \t}
  \hat \vv^*_\g(y)
  \forA \t \gg 1,
   }
   \ee
   on a certain ``maximal" set $y \in \Upsilon(\t)$,
    where such an
   expansion is applicable.
 Roughly speaking, it is given by
  \be
  \label{HH66}
\Upsilon(\t) \sim \{r(\t) < |y| < R(\t)\},
 \ee
 where $R(\t)$ characterizes the outer matching with smooth and
 almost regular flow, while much smaller $r(\t) \to 0$ is the
 crucial sphere, on which the inner matching with the flow bounded
 at the origin is supposed to occur.
 Thus, in particular,
    \be
    \label{PPP0}
 \mbox{set $\Upsilon(\t)$ essentially depends on the eigenfunction
$\hat \vv^*_\g(y)$
 behaviour as $y \to 0$.}
  \ee
 Indeed, as a proper illustration, $\Upsilon(\t)$ is different for
 the first monotone patterns as in \ef{Pot4} and for the singular
 oscillatory ones in \ef{Pot6}. For other types of possible
 eigenfunctions, the actual matching as in \ef{pol2}--\ef{pol4}
 defines the inner boundary of $\Upsilon(\t)$.

 Concerning the outer boundary of the radius $R(\t)$, which is not
 that essential, it is estimated as follows: since $\hat
 \vv^*_\g(y)$  is assumed to be close to a Hermite polynomial for
 $y\gg 1$ of order $k$ large and hence $\l_\g(c)= \frac {k+1}2$, with
 $|\g|=k$, this radios is characterized by matching with zero: for
 $\t \gg 1$,
 \be
 \label{Mat11}
  \tex{
 \uu_{\rm SL}(y) + {\mathrm e}^{\l_\g \t}
  \hat \vv^*_\g(y) \sim 0 \LongA \frac 1{|y|} - {\mathrm e}^{-
  \frac {k+1}2 \t}|y|^k \sim 0 \LongA R(\t) \sim {\mathrm
  e}^{\frac \t 2},
   }
   \ee
   as should be via the basic kinetic energy estimate
   \ef{ss1} in the rescaled variables \ef{ll2}. However, as we
   will see, this outer region cannot be important for the inner
   matching.

   In Inner Region II, we have to return to the full original
   rescaled
   equation \ef{ww1}, where for convenience we set ${\mathbb P}= I
   -{\mathbb Q}$,
    \be
    \label{PPP2}
     \tex{
   \hat \uu_\t=\big(\BB^*- \frac
    12\,I\big)\hat\uu - (\hat\uu \cdot \n)\hat\uu +{\mathbb Q}(\hat\uu \cdot
    \n)\hat \uu,
     }
     \ee
  where the convection term $(\hat\uu \cdot \n)\hat \uu$ is negligible in
  comparison with the linear $\BB^*$-term: this is the actual definition of
  Region II. The non-local term keeps be influential and even
  leading therein (we have seen a similar nonlocal phenomenon in
  the study of polynomial micro-structure of multiple zeros in Section
  \ref{SHerm}, where the integral term was shown to essentially
  deform the Hermitian solenoidal fields).

 Thus, in Region II, the following asymptotic equation occurs:
 \be
    \label{PPP3}
     \tex{
    \hat\uu_\t=\big(\BB^*- \frac
    12\,I\big)\hat \uu  + {\bf f}^*(\hat\uu) + ...\, \whereA
    {\bf f}^*(\hat\uu)=
    {\mathbb Q}(\hat\uu \cdot
    \n)\hat\uu \asA \t \to +\iy.
     }
     \ee
   It follows (or this can be a key assumption to be checked)
   that, in view of the expansion \ef{PPP1},
   the forcing term in \ef{PPP3} can be estimated as follows:
 \be
 \label{PPP4}
  \tex{
  {\bf f}^*(\hat\uu) \sim {\bf f}^*(\hat\uu_{\rm SL}) \sim {\mathbb
  Q}_{\Upsilon(\t)}(\hat\uu_{\rm SL} \cdot \n)\hat\uu_{\rm SL},
  }
  \ee
  where in the last term we assume integration in \ef{HHH21} over $\Upsilon(\t)$
  only, which represents the leading-order expansion term of the
  nonlocal term. If the latter is not true, one needs to include
  also the nonlocal portion that depends on the still unknown
  solution expansion in Region I, that makes the analysis more
  complicated (but formally doable). The {\em Outer Region}, where
  the solution is much smaller, is assumed to produce no
  essential influence on the integral.

Thus, the possibility of matching of Regions I and II depends on
existence of a smooth bounded solution of the following {\em limit
problem}:
\be
    \label{PPP5}
     \tex{
    \hat\uu_\t=\big(\BB^*- \frac
    12\,I\big)\hat\uu  + {\bf f}^*_{\Upsilon(\t)}(\t) \whereA {\bf f}^*_{\Upsilon(\t)}(\t)=
{\mathbb
  Q}_{\Upsilon(\t)}(\hat\uu_{\rm SL}) \forA \t \gg 1.
     }
     \ee
It is convenient to use the generalized Hermite polynomials to
describe the solution:
 \be
 \label{PPP6}
  \tex{
   \hat\uu(\t) = \sum \cc_\b \vv^*_\b \whereA \dot \cc_\b= \big(\l_\b-
   \frac 12\big)
   \cc_\b + \langle {\bf f}^*_{\Upsilon(\t)}(\t), \vv_\b \rangle.
   }
   \ee

The first condition (by no means, a necessary and/or sufficient))
of proper matching reads then as follows:
 \be
 \label{PPP7}
  \fbox{$
  \mbox{the auxiliary problem (\ref{PPP6}) has an $L^\iy$-solution
  for all $\t \gg 1$}
  $}
  \ee
  (otherwise, $\uu_0 \not \in L^\iy$).
  As a second one, for a rough checking of matching,
   consider the
  equation for the leading Fourier coefficient for $\b=0$, with $\l_0=0$:
   \be
   \label{Mat20}
    \tex{
\dot \cc_0= -
   \frac 12\,
   \cc_0 + \langle {\bf f}^*_{\Upsilon(\t)}(\t), \vv_0 \rangle,
   }
    \ee
    where $\vv_0$ is a constant vector (a polynomial of degree
    zero). Integrating this and assuming a slow growth divergence,
    one can expect a ``quasi-stationary" behaviour given by
    \ef{Mat20}:
    \be
    \label{Mat21}
    \cc_0 \sim 2 \langle {\bf f}^*_{\Upsilon(\t)}(\t), \vv_0
    \rangle\forA \t \gg 1.
    \ee
Assuming that $r(\t) \to 0$ (otherwise, no blow-up as $ \t \to
+\iy$),  matching the stable manifold behaviour  about the S-L
profile $\hat u_{\rm SL}$ assumes, at least, that, for $\t \gg 1$,
 \be
 \label{Mat22}
  \tex{
  \hat u(0,\t) \sim \cc_0(\t) \sim \hat u_{\rm SL}(\xx)\big|_{|\xx|
  \sim r(\t)} \LongA |\cc_0(\t)| \sim \frac 1{r(\t)}.
 }
   \ee
Involving other expansion coefficients $\{\cc_\b(\t)\}$ will lead
to a similar but more complicated relation with $r(\t)$.
Altogether, \ef{Mat21} and \ef{Mat22} define a complicated
nonlinear integral equation for the expansion coefficient
$\cc_0(\t)$ and the matching radius $r(\t)$. On integration in
 the non-local term given via \ef{HHH21}, i.e.,
  \be
  \label{Mat23}
  \tex{
{\bf f}^*_{\Upsilon(\t)}(\t) \sim \int_{\Upsilon(\t)}
\frac{y-z}{|y-z|^3} \sum_{(i,j)} \hat u^i_{{\rm SL},z_j} \hat
u^j_{{\rm SL},z_j} \whereA |\hat u^i_{{\rm SL},z_j}| \sim \frac
1{|z|^2}, }
    \ee
it follows that the RHS in \ef{Mat21} behaves as $\sim \frac 1r$
as $r \to 0$, which is satisfactory for a possible matching
purpose (at least, an obvious contradiction for matching is not an
immediate option). Surely, these are just rough matching
estimates, and further more difficult matching study along the
lines of \ef{PPP6}, \ef{Mat21} is necessary.

 Thus, due to \ef{PPP0}, the ``quality" of the eigenfunction $\hat
\vv^*_\g(y)$ as $y \to 0$ is essentially involved into \ef{PPP7}
via the set of integration $\Upsilon(\t)$ for $\t \gg 1$. This
somehow  involves verifying  a certain ``integrability"
(sufficient integral convergence) over $\Upsilon(\t)$, where the
asymptotics of the eigenfunction $\hat \vv^*_\g(y)$ as $y \to 0$
is  key. In other words, this eigenfunction must belong to a
``proper functional class" to make the above computations
meaningful. In general, similar to the construction in Section
\ref{S1.5}, the actual blow-up structure in Region II will indeed
 depend on the behaviour of $\hat \vv^*_\g(y)$ as $y \to 0$ (unknown), and hopefully will look like a
smoother``quasi-stationary"
 (driven by $\D$ only)
 evolution.
 On the other hand, even for very bad, singular and/or
highly oscillatory eigenfunctions $\hat \vv_\g^*$, there is still
some plausible hope that the integrals in \ef{PPP6} still properly
converge
  in such a manner
that \ef{PPP7} holds (this can be checked), and the matching with
a bounded flow for $y \approx 0$ can be purely ``accidental". This
means that \ef{PPP7} can be valid for some special and possibly
very singular eigenfunctions $\hat \vv_\g^*$, so that the
functional setting for the eigenvalue problem for the occurred
pseudo-differential operator can be treated in a wide sense.

Further matching conditions, involving also the already obtained
Region I expansions (and the Outer Region if necessary), may be
important, which eventually give the matched solutions in terms of
a converging bounded functional series.
 Finally,
in view of the precaution \ef{SLinj}, another fruitful idea is to
diminish the fluid injection by assuming that
 \be
 \label{ccinf}
  c=c(\t) \to +\iy \asA \t \to +\iy,
  \ee
which is an individual subject of Section \ref{S5.4}.

\subsection{A scaling view to formation of Type II blow-up
solutions: heteroclinic orbits are necessary}
  \label{S7.5}

We then need to consider the full rescaled equation taking into
account the angular operator \ef{kkk1}:
 \be
 \label{kkk2}
  \tex{
  \hat \uu_\t + \s \hat \uu_\mu + {\mathbb P}(\hat \uu \cdot \n)\hat
  \uu=\D \hat \uu - \frac 12 \, y \cdot \n \hat \uu - \frac 12\, \hat \uu.
   }
   \ee
 Assuming that the blow-up is faster than the self-similar one, i.e., is of Type II and
  \be
  \label{kkk3}
   \tex{
   \sup_y|\hat \uu (y,\t_k)|=C_k \to \iy \asA \{\t_k\} \to \iy,
 }
    \ee
we apply the $C_k$-scaling technique as in Section \ref{SEE2} by
setting
 \be
 \label{kkk4}
  \tex{
  \hat \uu = C_k \ww, \,\, y = y_k+ a_k z, \,\, \t = \t_k+a_k^2 s, \,\, \mu \mapsto a_k^2 \mu
  \whereA a_k = \frac 1{C_k} \to 0.
  }
\ee
  The resulting equation for $\ww_k$ takes the
 following perturbed form:
  \be
  \label{kkk5}
   \tex{
  \ww_s+  \s  \ww_\mu + {\mathbb P}( \ww \cdot \n)
  \ww=\D  \ww - \frac 1{C_k^2}\, \big( \frac 12 \, z \cdot \n  \ww + \frac 12\,  \ww\big).
 }
 \ee
In view of uniform boundedness and further regularity of the
sequence $\{\ww_k\}$ on compact subsets in $\re^3 \times \re$, we
can pass to the limit $k \to \iy$ in the weak sense in \ef{kkk5}
to conclude that $\ww_k(s)$ must approach a regular  solution
$\WW$ satisfying the NSEs
 \be
 \label{kkk6}
\WW_s+\s  \WW_\mu + {\mathbb P}( \WW \cdot \n)
  \WW=\D  \WW \inB \re^3\times \re, \quad \|\WW(0)\|_\iy=1.
   \ee
 Note that $\WW(z,s) \not \equiv 0$ is an {\em ancient solution}, which is defined for
 all $s \le 0$. At the same time, by construction, it is also a
 {\em future solution}, which must defined for all $s>0$. Indeed,
 one can see that if $\WW(s)$ blows up at some finite $s=S^->0$,
 this would contradict the Type II solution $\ww(y,\t)$ is globally
 defined for all $\t > 0$.
  Thus, by scaling of the Type II blow-up orbit (\ref{kkk4}), we arrive at the problem (\ref{kkk6}),
  which defines:
 \be
 \label{nnN1}
 \mbox{$\{$\underline{heteroclinic} solution $\WW(z,s) \not =0\}$ = $\{$\underline{ancient} for
 $s<0\}\,\,\cup\,\,\{$\underline{future} for $s>0\}$}.
 \ee

 Note that constant stationary solutions $\WW_0$ of (\ref{kkk6}) are the
 simplest possibilities, which being rather suspicious (too
 elementary) have not still been ruled out.
  We do
   not know whether or not non-constant or non-steady global bounded regular solutions of \ef{kkk6}
exist    for some $\s \not =0$.
 Note that in the stabilization problem to a constant solution
 $\WW_0$, which appear by linearization:
  \be
  \label{bi21}
  \WW(s)=\WW_0+\vv(s) \LongA
  \vv_s+\s  \vv_\mu + {\mathbb P}( \WW_0 \cdot \n)
  \vv=\D  \vv,
   \ee
 the linear operator does not have a clear discrete spectrum,
 unlike the one appeared in Region I. Therefore, the actual
 exponential stabilization rate in Region II is inherited from the
 spectral problem in Region I, assumed to be properly
 solved beforehand.
    In any case, this analysis shows a
   quasi-stationary nature of formation of Type II blow-up
   patterns on shrinking spatial subsets around blow-up points. Of
   course, the actual details of such a behaviour can be much more
   complicated involving various extra singularity mechanisms.

Regardless a clear lack of rigorous arguments here on existence of
Type II blow-up patters,  a {\em negative} conclusion occurs
(surely, not surprisingly):
 \be
 \label{bp11}
  \fbox{$
   \begin{matrix}
  \mbox{if the hypothetical heteroclinic patterns $W(z,s)$ from
  (\ref{kkk6}) {\em do not match}}
 \\
 \mbox{for $|z|  \gg 1$  typical spatial
 structures in (\ref{Lan1}) or those in (\ref{delta3}),}
 \\
  \mbox{tornado-type blow-up  around the S--L singularities is
 not possible.}
 \end{matrix}
 $}
 \ee


\subsection{Blow-up drift on the ${\mathbf c}$-manifold of S--L solutions}
 \label{S5.4}

As we know from R--D theory, the above blow-up scenario does not
exhaust all the types of possible singular patterns.
  As a
new evolution possibility, we again consider the S--L solutions,
but currently we assume a centre-subspace-like evolution on their
whole manifold.

We begin with the first linearized procedure to
 calculate a first approximation of such a behaviour.
Namely, following \cite{GKSob}, we assume that slow blow-up
evolution occurs ``along" an unknown functional dependence
$c=c(\t)$ in \ef{Lan2}, i.e., via \ef{ccinf},
 \be
 \label{ttt1}
  \hat \uu(\t)= \hat \uu_{\rm SL}(c(\t))+ \YY(\t) \in
  L^2_{\rho^*}(\re^3)
 \whereA c(\t) \to +\iy \asA \t \to +\iy.
  \ee
Then, instead of \ef{z7}, using the second asymptotics in
\ef{ccc1}, we substitute \ef{ttt1} into the NSEs with the
 corresponding extra force on the right-hand side as in \ef{delta1} of the form
 $= \frac{16}{c(\t)}\, \d(y) {\bf j} +...\,.$
Omitting higher-order terms, we obtain the following perturbed
inhomogeneous PDE:
 \be
 \label{ttt2}
  \begin{matrix}
 \YY_\t - \frac {\dot c(\t)}{c^2(\t)} \, \hat \uu_0=
  \HH'(\hat \uu_{\rm SL}(c(\t)) \YY + \DD(\YY)+...
  \qquad\qquad
   \ssk\ssk\\
  \equiv \big(\BB^*- \frac 12\,  I\big)\YY - \frac 1c\,[(\hat \uu_0
  \cdot \n)\YY+(\YY \cdot \n)\hat \uu_0] - \frac 1c \,C_3
  \int_{\re^3}(\cdot) (\hat \uu_0,\YY) +... \, .\qquad\qquad
  \end{matrix}
   \ee
In the last term, we mean the nonlocal operator in \ef{Lan2}
defines (still formally)  at the constant vector field $\hat
\uu_0$, with a special truncation at the origin.

We next looking for a solution governed by the leading
perturbation in \ef{ttt2},
 \be
 \label{ttt3}
  \tex{
   \YY(\t)= \var(\t) \Psi^*(y)+ \ww(\t)
   \whereA \var(\t)= - \frac {\dot c(\t)}{c^2(\t)}  \to 0
   \asA \t \to +\iy,
    }
    \ee
    and we assume that $\ww(\t)\bot \Psi^*$ in the dual metric,
    i.e., with a proper definition of the adjoint element (a
    linear functional),
     \be
     \label{ttt4}
     \langle \ww(\t),\Psi \rangle \equiv 0 \quad (\Psi \in
     L^2_\rho(\re^3)).
 \ee
 Substituting \ef{ttt3} into \ef{ttt2} yields
  \be
  \label{tt10}
   \tex{
   \dot \var \Psi^* + \var \hat \uu_0= \var \big(\BB^*- \frac 12\,
   I\big)\Psi^* - \frac{ \var}c[\cdot](\hat \uu_0,\Psi^*)+...\, ,
   }
   \ee
   where $[\cdot](\hat \uu_0,\Psi^*)$ denotes the linear operators in \ef{ttt2} at $\hat \uu_0$ (with
    possibly
  a necessary truncation at $0$).

\ssk

\noi\underline{\em Countable set of exponential patterns}.
Formally, it then follows that \ef{tt10} admits a countable set of
linear solenoidal blow-up patterns (here $\l_k=-\frac k2$)
 \be
  \label{ex11}
   \tex{
  \YY_k(y,\t)=a_k(\t) \vv_k^*(y)+...\,
  \whereA \dot a_k - \frac{\dot c}{c^2}\, \langle \uu_0, \vv_k
  \rangle= \big(\l_k- \frac 12)a_k.
  }
  \ee
In particular, we have the following special family of functions
$\{c_k(\t)\}$:
 \be
 \label{ex12}
  \tex{
  c_k(\t)={\mathrm e}^{(\frac 12-\l_k)\t} \LongA
  a_k(\t)= \big(\frac 12-\l_k\big)\, \t \,{\mathrm e}^{(\l_k-\frac 12)\t}
  \forA \t \gg 1.
  }
  \ee

\ssk

 \noi\underline{\em Power decay pattern}.
  Assuming now that in \ef{ttt3},
 \be
 \label{tt11}
  |\dot \var(\t)| \ll |\var(\t)| \forA \t \gg 1,
   \ee
  we mean that $c(\t)$ is a slow  growing
 function  not of an exponential form.
 Then
 the only possible way to balance the terms in \ef{tt10} is
  the vector $\Psi^*$ to satisfy:
  \be
  \label{tt12}
   \tex{
\big(\BB^*- \frac 12\,
   I\big)\Psi^*= \hat \uu_0.
   }
    \ee
Since the operator is invertible, there exists a unique solution
$\Psi^*$ constructed by an eigenfunction expansion via solenoidal
Hermite polynomials, so $\Psi^*$ is also solenoidal. It follows
from \ef{ccc1} that this linear procedure gives  a similar
singularity at the origin (cf. below with a nonlinear one),
 \be
 \label{tt13}
  \tex{
 | \Psi^*(y)| \sim \frac 1{|y|} \asA y \to 0.
   }
   \ee

The further balance in \ef{tt10} is then obtained, as usual, by
multiplication by the adjoint orthonormal element $\Psi$ that
yields an asymptotic ODE for $\var(\t)$ for $\t \gg 1$:
 \be
 \label{tt14}
  \tex{
  \dot \var =  \g_0 \frac {\var}{c}+...\, \whereA \g_0=\langle
  [\cdot](\hat \uu_0,\Psi^*),\Psi \rangle.
  }
   \ee
 We must admit that actually, $\g_0$ is not supposed to be a
 constant, since in this dual product, a further matching of
 the pattern with a bounded flow for $y \approx 0$ should be assumed.
 This matching (a cut-off for $|y| \ll 1$) procedure makes
 the integral in \ef{tt14}  finite for all $\t \gg 1$, but then we
 conclude that
 \be
 \label{tt15}
  \g_0= \g_0(\t) \to \infty \asA \t \to +\iy \quad (\mbox{e.g.,
  \,$\g_0 \sim \ln \t$}),
   \ee
 though the divergence turns out to be
  slower than that for $\var(\t)$, so this does not essentially affect
 the leading term in \ef{tt14}.
Thus, we arrive at the equation for $\var(\t)$:
 \be
 \label{tt16}
  \tex{
  \big( \frac {c'}{c^2}\big)'= \g_0 \, \frac {c'}{c^3}+...
  \LongA c(\t) \sim -\frac {\g_0}2\, \t+... \asA \t \to +\iy
   }
   \ee
 (recall that we do not exclude the case, e.g., $c(\t) \sim \t \ln
 \t$ for $\t \gg 1$).

Hence, up to lower-order multipliers, we get the following
expansion of such solenoidal  patterns, which is convenient to
write in terms of a series of the form (again, lower-order,
logarithmic-like factors are not taken into account):
 \be
 \label{MM1}
  \tex{
   \hat \uu(y,\t) \sim \frac 1 \t \, \hat \uu_0 + \frac 1{\t^2}\, \hat \uu_1
+ \frac 1{\t^3}\, \hat \uu_2+... \, ,
 }
 \ee
 where for convenience we put $\hat \uu_1 \sim \Psi^*$.
 To get the coefficients $\{\hat \uu_k\}$ of this nonlinear
 expansion, one needs to use the full rescaled equation \ef{ww1}
 with the operator \ef{HHH21}. It then follows that
 the equation for the third coefficient $\hat \uu_2$ takes
 the form
  \be
  \label{MM2}
   \tex{
   - \hat \uu_1 = \big(\BB^*- \frac 12\, I\big) \hat \uu_2- (\hat
   \uu_0,\n)\hat \uu_1 - (\hat
   \uu_1,\n)\hat \uu_0 - C_3 \int_{\re^3}(\cdot)(\hat \uu_0, \hat \uu_1),
    }
    \ee
 where in the last term we again assume a certain evolution cut-off
 procedure near the origin.  Then we obtain the  singularity for $\hat
 \uu_2(y)$:
 \be
 \label{MM3}
 \tex{
 - (\hat
   \uu_0,\n)\hat \uu_1 \sim \frac 1{|y|^3} \LongA \hat \uu_2(y) \sim
    \frac {\ln |y|}{|y|} \asA y \to 0,
    }
    \ee
    i.e, we have found that the third term is more singular at the
    origin than two previous ones.
 Therefore, similar to the procedure described in our first Type II blow-up
structure in \ef{WW8} for the Frank--Kamenetskii equation, we
observe that on small compact subsets in $y$, a model one-sided
(say from above) behaviour of the rescaled vector field can be
roughly estimated as
 \be
 \label{RRR1}
  \tex{
  |\hat \uu(y,\t)| \sim U(y,\t) \sim \frac  1{\t}\, \frac 1{|y|}+ \frac  1{\t^2}\, \frac 1{|y|} + \frac 1{\t^3}
  \, \frac{\ln |y|}{|y|}+... \, .
   }
   \ee
   Taking into account the first and the third terms,
  the absolute positive maximum of the scalar function $U(y,\t)$ is  attained at
   \be
   \label{yy1}
    \tex{
   |y| \sim {\mathrm e}^{-\t^2}
  \to 0, \quad \mbox{and hence}\,\,\,
    \sup_y U(y,\t) \sim
   \frac {{\mathrm e}^{\t^2}}\t \forA \t \gg
   1.
 }
 \ee
 Once we have known the radius (\ref{yy1}) of truncation of the singularity
 structure of $U$ at the origin, all computations can be redone
 to see its actual  influence on the final asymptotics.
 In view of scaling \ef{ll2}, overall, this yields
 a formal  expansion \ef{MM1} of the solenoidal vector field
 for a quadratic NSEs system, which is singular at the origin,
  but can admit some kind of critical surface inflection-like
 behaviour at the level set
 ($T=1$)
  \be
   \label{RRR3}
    \tex{
   |\uu(x,t)| \sim
    \frac {{\mathrm e}^{|\ln(T-t)|^2}}{\sqrt{T-t}\,|\ln(T-t)|}
   \asA t \to T^-.
 }
 \ee
 In other words, values \ef{RRR3} for $t \approx
 T^-$, the singular vector field $\uu(x,t)$ is expected to exhibit certain
 special transitional behaviour, which possibly can be used
 for further necessary matching and branching to create more
 reasonable blow-up patterns.

 The pattern with special level sets as in \ef{yy1} is obtained
 in the simplest situation, where no ``resonance" between
 different terms of the expansion \ef{RRR1}, which itself
can generate logarithmic factors, is assumed. In general, for
quadratic dynamical systems as \ef{ww1}, there can be several
types of asymptotic centre-manifold-type expansions, which are
very sensitive and depend on resonance conditions of various terms
involved; see invariant manifold theory in \cite{Lun}. The
validity of such conditions are very difficult to check,
especially in the presence of extra outer matchings involved.
However, in our opinion, \ef{yy1} correctly characterizes some
features of the blow-up behaviour of such patterns up to some
extra exponents and/or slower growing factors.
 A  justification of such a behaviour requires a quite involved matching analysis
of fully using eigenfunction expansions on the vector solenoidal
Hermitian polynomials and related delicate and  technical matching
procedures involved.


\subsection{On a possibility of blow-up on regular equilibria}
 \label{S5.5}

We now follow the discussion in  Section \ref{SRD1}. These
structures are more exotic but also reasonably well-known in
reaction-diffusion theory; cf. \cite{Fil00} and \cite{GKSob}.
 For simplicity, we take $\sigma =0$ in \ef{delta3} (or \ef{delta2}, then a new perturbation occurs)
  and  study a
possibility of blow-up behaviour
 \be
 \label{bb1}
 \hat \uu(y,\t)= \var(\t) \tilde \uu(\var(\t)y)+\YY \whereA
 \var(\t) \to \iy \asA \t \to +\iy.
  \ee
To this end, we first introduce the blow-up variables
 \be
 \label{bb2}
 \hat \uu(y,\t)= \var(\t) \hat \vv(z,s), \quad z= \var(\t)y, \quad
 \var^2(\t)\,{\mathrm d} \t= {\mathrm d}s
  \ee
 where $\hat \vv$ solves the following rescaled non-autonomous
 PDE:
  \be
  \label{bb3}
   \tex{
  \hat \vv_s =\HH(\hat \vv) - \rho(s)[(z \cdot \n) \hat \vv + \hat
  \vv] \whereA \rho(s)= \big(\frac{\dot \var}{\var^3}+ \frac 12\,
  \frac 1{\var^2}\big)(\t).
  }
  \ee
 If $\sigma  \not =0$,
  we introduce the angular TW dependence
  $\mu=\hat \mu +s$, so that
 the extra term
 $ - \sigma  \hat \vv_{\hat \mu}$
  should be put into the right-hand side (cf. \ef{Stat1}), which does not change
  the concepts of matching.
 It then follows that we may estimate $\hat \vv(s)$ as follows:
  \be
  \label{bb4}
   \tex{
  \hat \vv(s)= \tilde \uu + \rho(s) W, \,\,\, \mbox{provided that}
  \,\, \frac{\rho'}\rho \to 0, \,\,\, s \to + \iy
 }
   \ee
   (i.e., $\rho(s)$ is not exponentially decaying and is
   algebraic),
 where $W$ is a proper solutions of the linearized inhomogeneous
 equation
  \be
  \label{bb5}
  \HH'(\tilde \uu) W -[(z \cdot \n) \tilde \uu+ \tilde \uu]=0.
   \ee

In order to get possible acceptable families of functions
$\{\var_k(\t)\}$, we return to the original rescaled variable $y$
 and use \ef{bb1} in \ef{ww1}, \ef{HHH21} to get
 for $Y$ the equation
  \be
  \label{bb6}
   \tex{
  \YY_\t= \big(\BB^*- \frac 12\, I \big) \YY + \HH(\var(\t) \tilde
  \uu(\var(\t)y))+...\, ,
   }
    \ee
   where we omit higher-order terms. By the assumption
   \ef{delta4}, we have from \ef{delta1} that the main
   inhomogeneous term in \ef{bb6} can be estimates as follows:
 \be
 \label{bb7}
  \tex{
 \HH\big(\var(\t) \tilde
  \uu(\var(\t)y)\big) = -\frac {16 \pi}{c(\t)}\, \d(y)+...\,
  \LongA
  \YY_\t= \big(\BB^*- \frac 12\, I\big) \YY
 -\frac {16 \pi}{c(\t)}\, \d(y)
+...\, .
   }
   \ee
 We then balance the terms similar to \ef{ex12}:
 \be
 \label{bb9}
  \tex{
   \YY_k(\t)= {\mathrm e}^{(\l_k-\frac 12)\t} \vv_k^*+... \andA
  c_k(\t)={\mathrm e}^{(\frac 12-\l_k)\t}
  \forA \t \gg 1.
  }
  \ee

 Comparing \ef{bb1}, \ef{bb9} and \ef{bb2}, \ef{bb4} yields the
 matching condition, which we take in the simplest form to catch
 the exponential factors only:
 \be
 \label{bb10}
  \tex{
 \rho(s) \equiv \frac{\dot \var}{\var^3}+ \frac 12\,
  \frac 1{\var^2} \sim \frac 1 \var \, {\mathrm e}^{(\l_k-\frac
  12)\t} \LongA
 \var_k(\t) \sim \big(\frac 12- \l_k\big) {\mathrm e}^{(\frac
 12-\l_k)\t} \to \iy.
  }
  \ee
Let us first check the consistency of the expansion \ef{bb9} for
\ef{bb2}:
 \be
 \label{bb91}
  \tex{
  s \sim {\mathrm e}^{(k+1)\t} \LongA \rho_k(s) \sim \frac 1s \forA s
  \gg 1,
  }
\ee
 so that each $\rho_k(s)$ has the desired algebraic (and non-integrable) decay at infinity.

Recall that these are blow-up patterns constructed via a
hypothetical evolution on the manifolds of ``steady" (or
``TW-swirl" for $\sigma  \not = 0$) solutions, which includes both
families of the singular S--L profiles and the regular ones as in
\ef{delta3}. Including the angular eigenvalue $\sigma  \not =0$
will make the formal analysis more complicated to say nothing
about a rigorous justification (a finite energy interpretation of
such blow-up patterns is also hard).

\subsection{Towards Burnett equations}

Finally, we note that, for the Burnett equations \ef{NS1m}, a
similar formal analysis can be performed relative to the {\em
singular solutions} with a different behaviour near the origin,
 \be
 \label{Dif1}
  \tex{
  \uu_{\rm SS}(x) \sim \frac 1{|x|^5} \asA x \to 0,
 }
 \ee
 though proving existence of such steady structures is a
 difficult problem, as well as a rigorous mathematical
 justification of the expansion and matching procedures for both
 families of singular and bounded
 states. Concerning the linearization approach as in Section
 \ref{S5.5}, it can be performed in similar lines by using
 spectral theory and the generalized Hermite polynomials for the
operator \ef{NBN1}, \cite{Eg4}. The possible ways of
``interaction" with the spectral characteristics of the nonlocal
operator ${\bf J}_2$ are out of question here.

  On the
 other hand, it is not excluded that \ef{NS1m} can admit a purely
 self-similar blow-up with various multiple patterns
 (a finite or countable set? -- seems should be finite as higher-order parabolic flows
 suggest \cite{BGW1}; however, the solenoidal restriction can indeed easily spoil
 Leray's-type ``similarity party" even here). Of course, the nonexistence proof based on the MP ideas from
 \cite{Nec96} is not applicable here.
  Thus, \ef{NS1m} can admit much more
  complicated  families of
 blow-up patterns (it seems that numerics can help here being however very difficult),
  though the classification problem
 \ef{Turb1} remains a principle difficult issue, which, as expected, will never
 be completely solved (too difficult, and not that essential?).

 \subsection{Complicated oscillatory singular equilibria for an
 R--D equation}
  \label{S7.7}

Here, as a key illustration to some of our speculations, we
briefly consider \ef{ST2}. In radial geometry,
 such Emden (1907)--Fowler (1914) equations
 have been most carefully studied in ODE
 theory. First detailed classification of solutions were obtained  in Gel'fand \cite{Gel59}
 (this ODE section is known be written by
  Barenblatt) and
  Joseph--Lundgren \cite{JL73}; see further
  references and applications for blow-up in
 \cite[\S~6.5]{GalGeom}. We briefly comment on still difficult and
 seems not completely well-understood
 oscillatory properties of singular solutions.

First of all, \ef{ST2} admits the {\em homogeneous singular
stationary solutions} (SSS):
 \be
 \label{EF1}
  \tex{
  U(r)= \pm \,C_* r^{-\mu} \whereA \mu= \frac 2{p-1}, \quad
  C_*=\big[\mu(N-2-\mu)\big]^{\frac 1{p-1}} \,\,\,\big(p> \frac N{N-2}\big).
 }
   \ee
To describe others of changing sign, we introduce the {\em
oscillatory component} $\var$ by
 \be
 \label{EF2}
  \tex{
u(r)=r^{-\mu} \var(- \ln r): \quad
 \var''+(2\mu+2-N) \var'+\mu(\mu+2-N)\var+|\var|^{p-1} \var=0.
 }
 \ee

\begin{figure}
\centering
\subfigure[shooting of $u_0(r)$]{
\includegraphics[scale=0.52]{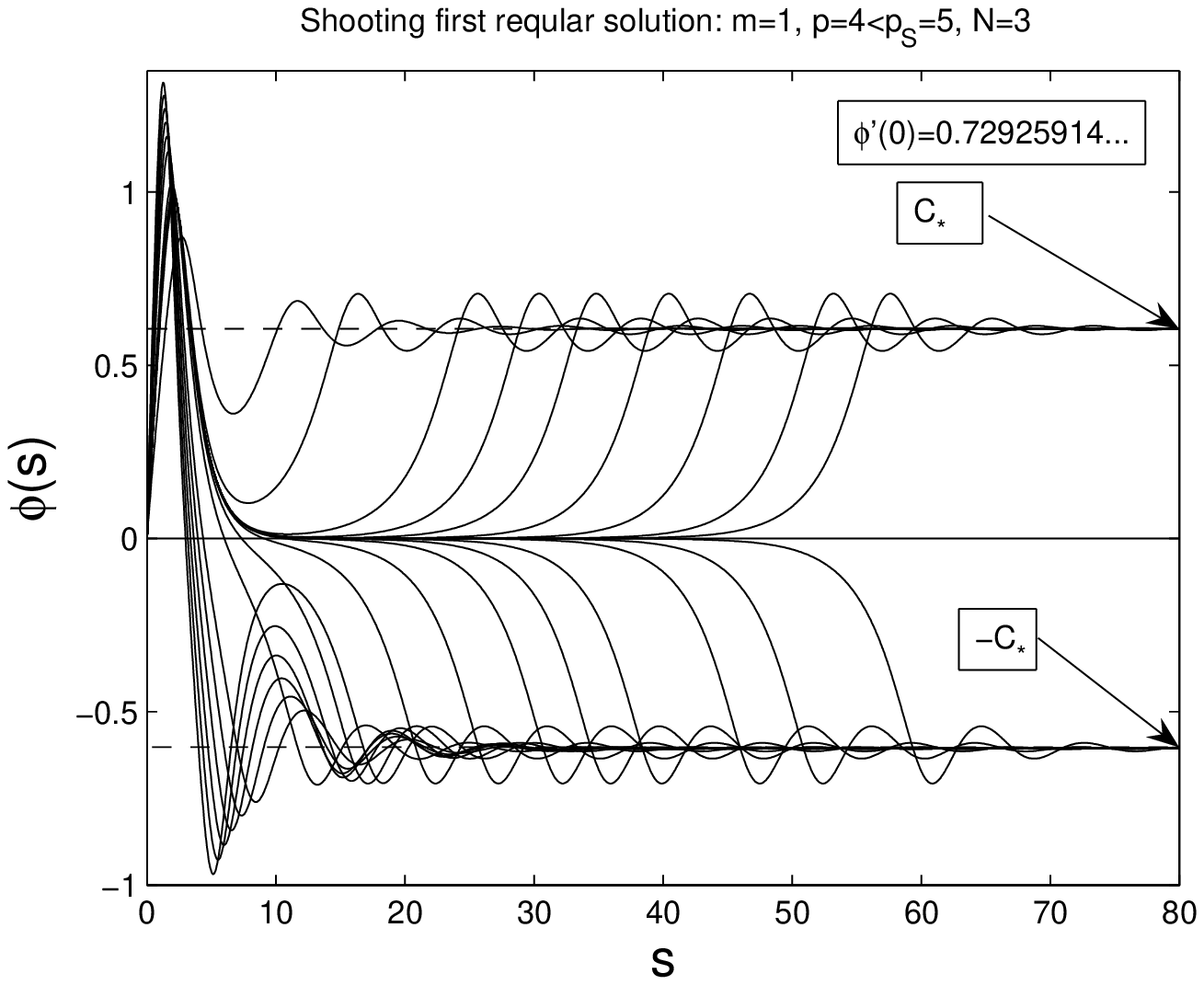} 
}
\subfigure[shooting of $u_1(r)$]{
\includegraphics[scale=0.52]{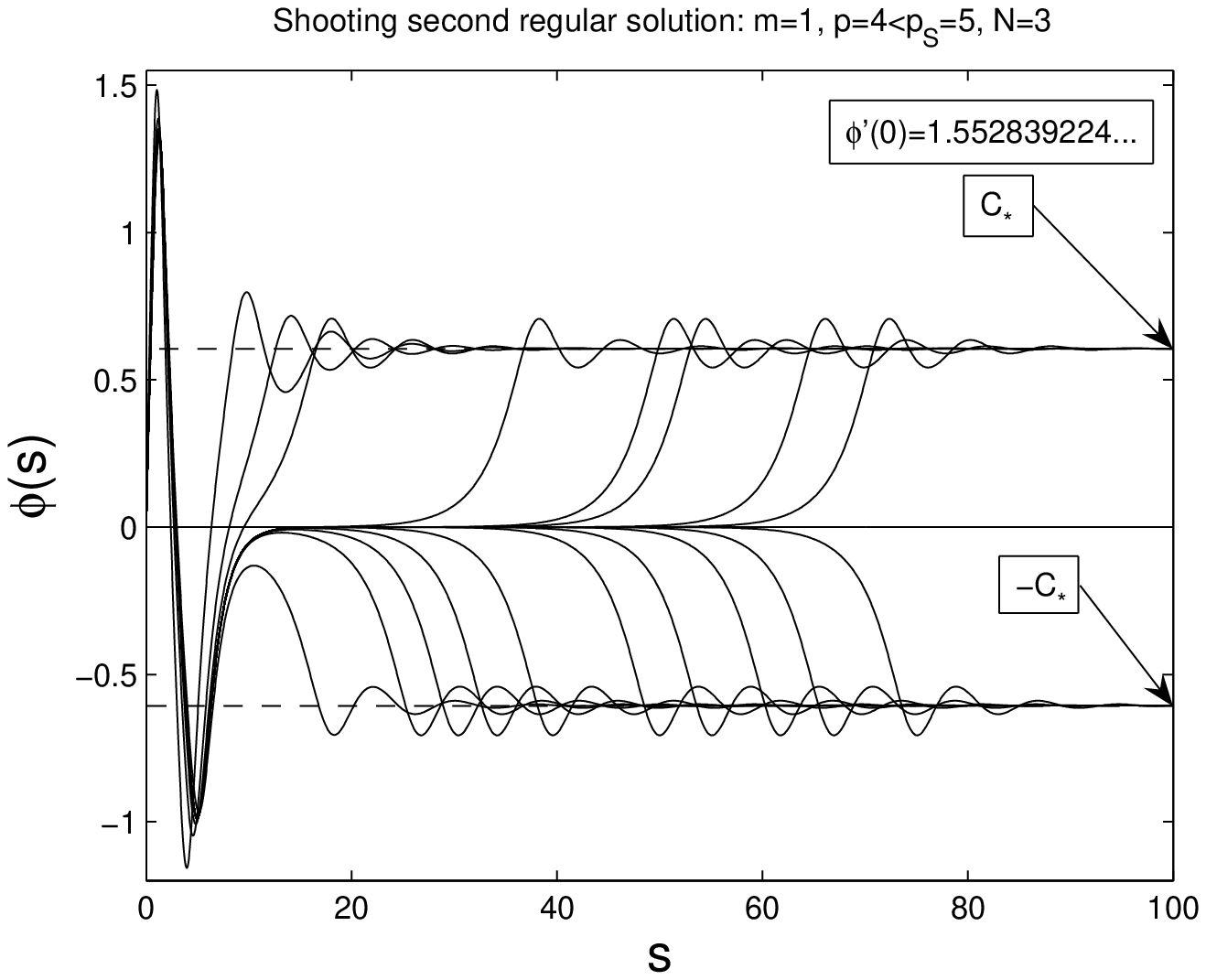}
}
   \vskip -.3cm
    \caption{\small Shooting the first (a) and the second (b) regular equilibria
   via ODE in (\ref{EF2})  for $p=4<p_{\rm S}=5$, $N=3$.}
   \vskip -.3cm
 \label{Fg2}
\end{figure}

Figure \ref{Fg2} shows the non-oscillatory character of singular
solutions in the subcritical Sobolev range $p<p_{\rm S}$ by
shooting from $s=0$ ($r=1$). Both Figures \ref{Fg2}(a) and (b)
explain that as $r \to 0$, i.e., $s \to + \iy$, the oscillatory
component stabilizes to constants $\pm C_*$ as in \ef{EF1},
excluding countable set of regular patterns $\{u_k, \, k \ge 0\}$,
which are
 bounded at $r=0$ and have a fixed number of sign changes.
Here (a) shows shooting the first profile $u_0(r)>0$ with no
zeros, while (b) corresponds to shooting of $u_1(r)$ with a single
zero on $r \in (0,1)$ ($s \in (0,+\iy)$). Surrounding those
regular $u_k(r)$ are continuous families of singular solutions
exhibiting the behaviour \ef{EF1} as $r \to 0$.

\begin{figure}
\centering
\subfigure[$p_{\rm S}<p<p^*$]{
\includegraphics[scale=0.52]{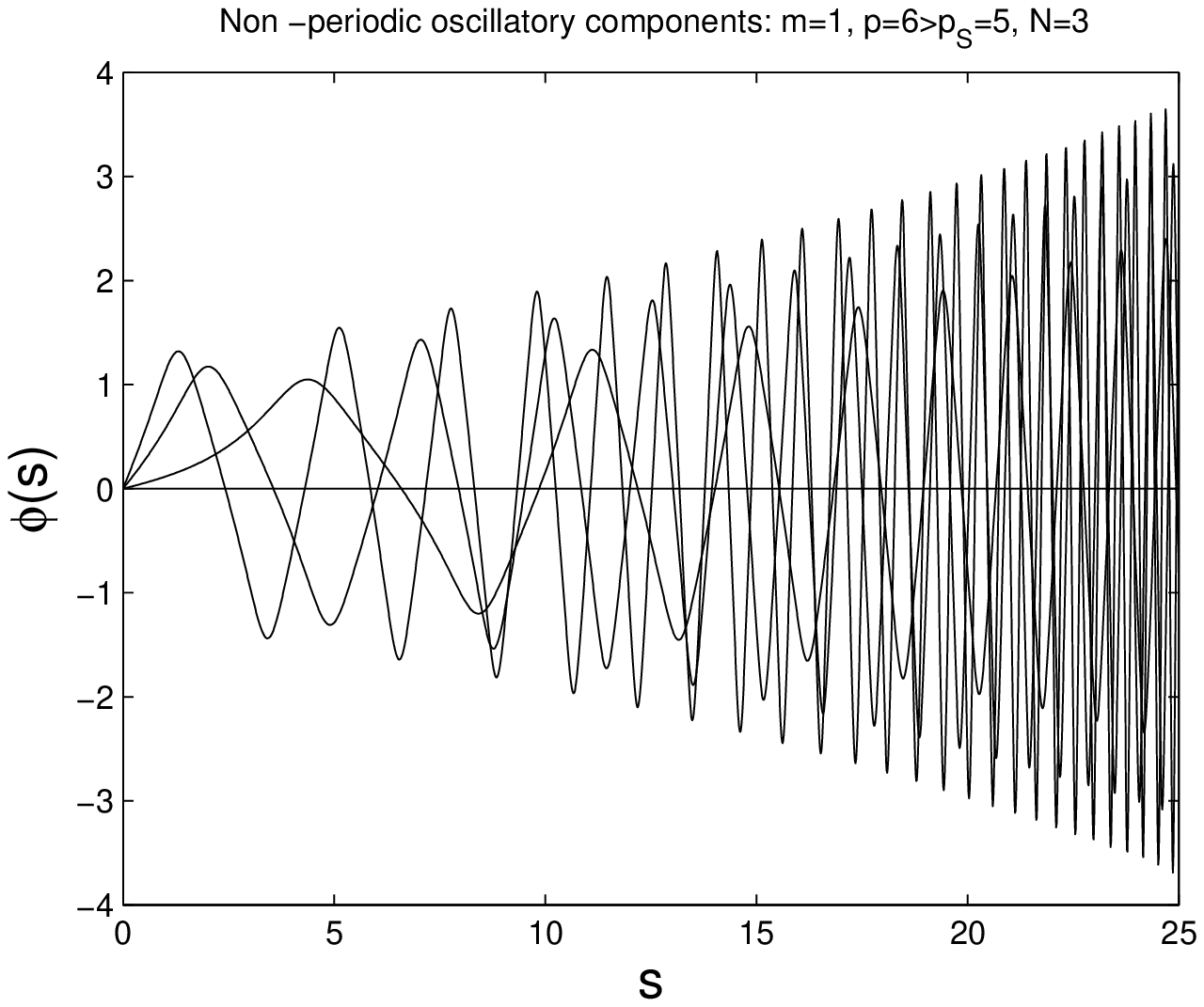} 
}
\subfigure[$p>p^*$]{
\includegraphics[scale=0.52]{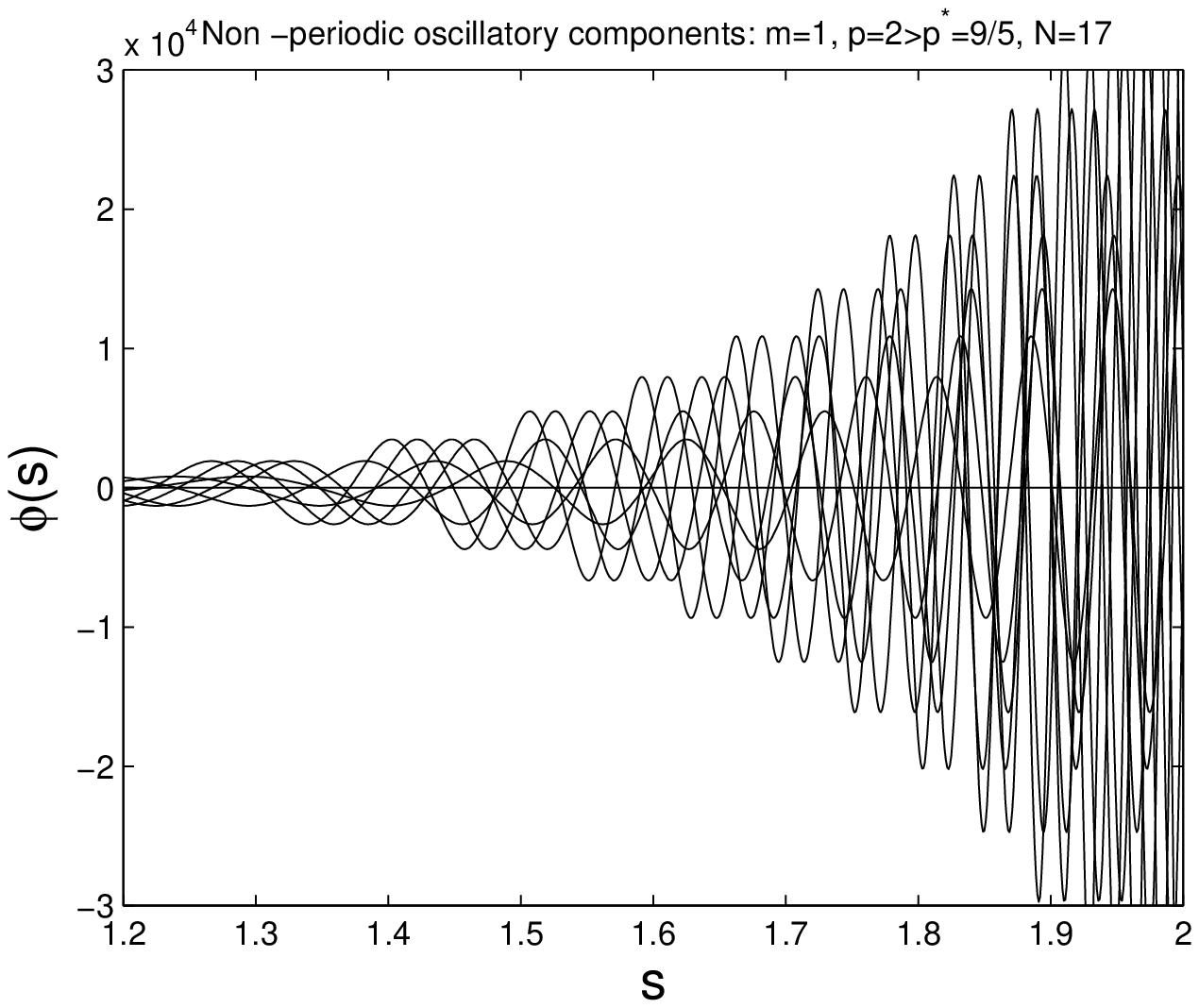}
}
   \vskip -.3cm
    \caption{\small Oscillatory singularities of
   the ODE in (\ref{EF2})  for
    $p=6>p_{\rm S}=5$, $N=3$ (a) and for
    $p=2>p^*= \frac 95$, $N=17$ (b).}
   \vskip -.3cm
 \label{Fg4}
\end{figure}

Next Figure \ref{Fg3} shows the crucial appearance of infinitely
many {\em oscillatory periodic orbits} $\var(s)$ at the critical
value $p=p_{\rm S}$, where \ef{EF2} takes the autonomous
 variational form
 \be
 \label{VAR1}
  \tex{
  \var''- \frac{(N-2)^2}8 \, \var + |\var|^{p-1}\var=0.
  }
  \ee
  Finally, in
Figure \ref{Fg4}, we show how the oscillatory character of
singular solutions dramatically changes for $p>p_{\rm S}$, where
nonlinear spiral out behaviour replaces the bounded periodic one.
In (a), we take $p=6>p_{\rm S}=5$ for $N=3$. In (b),  the highly
oscillatory behaviour is shown for $p=2$ and $N=17$, i.e., above
the uniqueness \cite{GV97} critical exponent
 \be
 \label{uni1}
  \tex{
  p^*= 1+ \frac 4{N-4-2 \sqrt{N-1}} \forA N \ge 11
   \quad \big(p^*= \frac 95<2 \forA N=17 \big),
  }
  \ee
  which plays a role in a number of blow-up problems for this
  R--D equation; see \cite{GV97} and \cite{Mat04} for further
  details. Thus, the spiral-type oscillations for large $p>p_{\rm
  S}$ get arbitrarily large, as $r \to 0$, so singular equilibria stay
  arbitrarily far from the homogeneous one \ef{EF1} (is this possible for the steady NSEs \ef{ST1}
  relative to the homogeneous S--L solutions \ef{Lan1}?).

\begin{figure}
\centering
\includegraphics[scale=0.55]{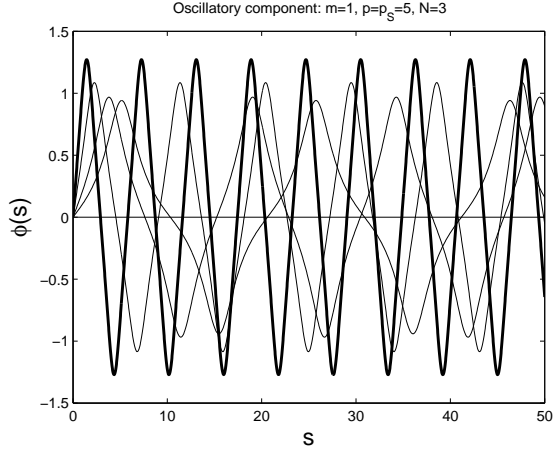}
 \vskip -.3cm
\caption{\small Periodic oscillations in the ODE in
(\ref{EF2}) for $p=5=p_{\rm S}$, $N=3$.}
   \vskip -.3cm
\label{Fg3}
\end{figure}

Thus, the singular equilibrium manifold for \ef{ST2} is rather
complicated even in the radial case. Of course, in the radial
geometry, the blow-up analysis of the corresponding parabolic
equation \ef{Sob1} is  essentially simplified by using the
Sturmian argument of intersection comparison (the number of
intersections of different solutions cannot increase with time;
C.~Sturm, 1836), and this can prohibit some scenarios of blow-up;
see various applications in \cite{GalGeom, GV97, Mat04}. In the
non-radial geometry, even for \ef{Sob1}, this advantage is almost
nonexistent, the singular equilibria can get much more complicated
 and are  unknown.
 On the other hand, the MP remains then in place, but its application
 to control, simplify, or prohibit possible
 blow-up structures is unclear.

  Thus, the non-radial geometry in such parabolic problems  can result in existence of new
blow-up patterns with evolution ``close" to such singular
manifolds, where specific ``swirl-torsion" blow-up phenomena may
occur. Of course, the matching conditions such as \ef{PPP7} for
such singular equilibria may prohibit, via obvious divergence,
almost all types of such blow-up patterns (the sum in \ef{PPP6} is
not from $L^\iy$), but we cannot always rule out a possibility of
``accidental convergence" of the corresponding integral over
$\Upsilon(\t)$.

 We suggest that this R--D singular equilibria experience can be
 taken into account for the NSEs, and then the open problem \ef{ST1}
 turns out to be  key (does it admit local oscillatory
 singular solutions that have nothing to do with
 the homogeneous S--L one \ef{Lan1}?--a principal question).
  Recall that regardless the presence of the standard Laplacian $\D$
  in both \ef{ST1} and \ef{NS1}, the MP arguments
  cannot be somehow essential for such
   nonlocal problems (cf. \ef{ww1}). [All speculations are {\em
   modulo} \ef{PPP7}.]

\ssk

\noi{\bf Remark:  singular equilibria in the bi-harmonic case and
open problems.} The analogous elliptic problem for the Burnett
equations \ef{NS1m} is the bi-harmonic one,
 \be
 \label{bi1}
 -\D^2 u+|u|^{p-1}u=0,
 \ee
 which admits a similar to \ef{EF2} substitution with $\mu= \frac
 4{p-1}$ and the ODE:
  \be
  \label{bi2}
 \begin{matrix}
  -\var^{(4)}- A \var'''-B \var''-C \var' - D \var + |\var|^{p-1}
  \var=0 \whereA\qquad\qquad\qquad\quad \ssk \\
  A=2(2\mu+4-N), \quad B= 6 \mu^2+18\mu+11+(N-1)(N-9-6\mu),\quad\,\,\qquad\ssk\\
  C=2\big[2\mu^3+9\mu^2+11\mu+3+(N-1)[(N-3)(\mu+1)-3\mu^2-6\mu-2]\big],\quad\,\,\, \ssk\\
  D=\mu(\mu+2)[(\mu+1)(\mu+3)+(N-1)(N-5-2\mu)].\qquad\qquad\quad\qquad\qquad\,\,
   \end{matrix}
 \ee
This is a harder equation than \ef{EF2}, and its 
complexity shows how difficult proofs of global or blow-up bounds
on solutions of the corresponding parabolic PDE
 \ef{ho1}, $m=2$,
 can be for
$p \ge p_{\rm S}=\frac{N+4}{N-4}$, with $N>4$. The homogeneous SSS
of \ef{bi1} is
  \be
  \label{bi5}
  \begin{matrix}
 U(r)= \pm C_*\, r^{-\mu} \whereA \mu= \frac 4{p-1}, \,\,
  C_*=D^{\frac 1{p-1}}>0,
   \ssk\\
 \mbox{existing for}\,\,\,
  p>
 \frac N{N-4}, \,\,\,N>4, \quad \mbox{or} \quad p < \frac {N+2}{N-2},\,\,\,N>2.
 \end{matrix}
 \ee
 Global radial solutions of \ef{bi1} for $p> p_{\rm S}$ are
 obtained in \cite{Gaz06}, where further references can be found.

 Again, blow-up
patterns can be created by the manifold of singular equilibria
with completely unknown non-radial structure (with no traces of
the MP); the analogy of the matching \ef{PPP7} is supposed to be
taken into account; see \cite{Gal5Bl} for a discussion. Any
general estimates of Type I or II blow-up for \ef{ho1} for $p \ge
p_{\rm S}$ are absent. Moreover, we claim
 that a full and complete description of all the (non-radial)
 blow-up patterns for \ef{ho1} will never be achieved.
 Similar difficulties concerning existence and nonexistence
of various blow-up patterns
 occur for
supercritical nonlinear Sch\"odinger equations such as \ef{Sr1}
(see key references
  on this subject
  in  Merle--Raphael \cite{Mer05} and Visan \cite{Vis07})
 and for many other important
higher-order PDEs (see the list around \ef{Sr1}) and systems of
the twenty-first century PDE mathematics/application.



\section{On complicated blow-up patterns with swirl and precessions}
 \label{SLanFin}

 Using the previous ideas of possible concepts of blow-up in the
 equations \ef{NS1}, we continue the construction  of more refined
structures of a full complexity. We now use spherical coordinates
that allow us to fix more complicated singular stationary
solenoidal fields, in a neighbourhood of which some non-steady
blow-up phenomena may occur.

\subsection{Basic solutions with swirl in spherical coordinates}
 \label{S5.1}

Consider the equations \ef{NS1} in the spherical polar coordinate
system, with $\uu=(u_r,u_\th,u_\var) \equiv (u,v,w)$:
 \be
 \label{A1}
  \left\{
  \begin{matrix}
 u_t+ uu_r+ \frac 1r \, vu_\th + \frac 1{r \sin \th} \, w u_\var -
 \frac 1r\, v^2 - \frac 1{r} \, w^2\qquad\qquad\,\,\,\,\,\,\,\,\ssk\\
 =-p_r + \D_3 u - \frac 2{r^2} \, u - \frac 2{r^2} \, v_\th -
 \frac {2 \cot \th}{r^2} \, v - \frac {2}{r^2 \sin \th} \, w_\var,
 \qquad\,\,
 \ssk\\
v_t+ uv_r+ \frac 1r \, u v + \frac 1{r} \, v v_\th +
 \frac 1{r \sin \th}\, w v_\var  - \frac {\cot \th}{r} \, w^2\qquad\qquad\,\,\,\,\ssk\\
 =- \frac 1r \,p_\th + \D_3 v - \frac 1{r^2 \sin^2 \th} \, v + \frac 2{r^2} \, u_\th -
 \frac {2 \cot \th}{r^2 \sin \th} \, w_\var, \qquad\qquad
 \ssk\\
 w_t+ u w_r+ \frac 1r \, u w + \frac 1{r} \, v w_\th +
 \frac {\cot \th}r\, v w + \frac 1{r \sin\th} \, w w_\var \qquad\,\,\,\,\,\,\,\ssk\\
 =- \frac 1{r \sin\th}\, p_\var + \D_3 w - \frac 1{r^2 \sin^2 \th} \, w + \frac 2{r^2 \sin \th}\, u_\var
 +
 \frac {2 \cot \th}{r^2 \sin \th} \, v_\var,
 \ssk\\
  u_r+ \frac 2r \, u + \frac 1r \, v_\th + \frac{\cot \th}r \, v +
  \frac 1{r \sin \th} \, w_\var=0. \qquad \qquad\qquad\quad
 \end{matrix}
  \right.
  \ee
 Here $\D_3$ denotes the spherical $\re^3$-Laplacian:
 \be
 \label{A10}
  \tex{
  \D_3 h= \frac 1{r^2}\, ( r^2 h_r)_r + \frac 1{r^2 \sin \th} \,
  (\sin\th \, h_\th)_\th + \frac 1{r^2 \sin^2 \th} \,
  h_{\var\var}.
  }
  \ee


\noi\underline{\em  Homogeneous singular stationary states}.
Similar to \ef{z1}, we begin with simpler {\em singular stationary
solutions} of \ef{A1} of the homogenuity $-1$,
  \be
  \label{vv3}
   \tex{
  (u,v,w)= \frac 1r \, ( \hat u, \hat v, \hat w)(\var,\th), \quad p= \frac
  1{r^2} \, \hat p(\var, \th).
 }
   \ee
 This solution representation means that
the operators of the stationary equations \ef{A1}, in view of
their homogenuity, perform the following mappings of linear
subspaces:
 \be
 \label{vv4}
  \tex{
  W_1= {\rm Span} \, \big\{\frac 1r\big\} \to
  \hat W_2= {\rm Span} \, \big\{ \frac 1{r^2}, \, \frac 1{r^3}\big\}.
   }
   \ee
   Substituting yields the following PDE system for the four unknowns $(\hat
   u, \hat v, \hat w, \hat p)$:
 \be
 \label{A10n}
  \left\{
  \begin{matrix}
 -\hat u^2 +  \hat v \hat u_\th + \frac 1{\sin \th} \, \hat w  \hat u_\var -
 \hat  v^2 - \hat w^2\qquad\qquad\,\,\,\,\,\,\,\,\qquad\qquad\qquad\qquad\qquad\quad\ssk\\
 =2\hat p + \hat u_{\th\th} +  \cot \th \,\hat u_\th + \frac 1{\sin^2 \th} \, \hat u_{\var\var}
  -  2 \hat u -  2\hat v_\th -
 2 \cot \th \hat v - \frac {2}{\sin \th} \, \hat w_\var,
 \qquad\,\,
 \ssk\\
 \hat  v \hat v_\th+\frac 1{r \sin \th}\, \hat w \hat v_\var  - \cot \th \hat w^2\qquad\,\,
 \qquad\qquad\qquad\qquad\qquad\quad\qquad\qquad\quad\quad\ssk\\
 =- \hat p_\th + \hat v_{\th\th} + \cot \th \hat v_\th  + \frac 1{\sin^2\th} \, \hat v_{\var\var}
 - \frac 1{\sin^2 \th} \, \hat v +  2\hat u_\th -
 \frac {2 \cot \th}{ \sin \th} \, \hat w_\var, \qquad\qquad
 \ssk\\
  \hat v \hat w_\th+
  {\cot \th}\, \hat v \hat w + \frac 1{\sin\th} \, \hat w \hat  w_\var \qquad\,\,\,\,\,\,\,
  \qquad\qquad\qquad\qquad\qquad\quad\qquad\qquad\,\,\,\ssk\\
 =- \frac 1{\sin\th}\, \hat p_\var + \hat w_{\th\th} + \cot \th \, \hat w_\th  +
 \frac 1{\sin^2 \th} \, \hat w_{\var\var} - \frac 1{\sin^2 \th} \, \hat w +
 \frac 2{\sin \th}\, \hat u_\var
 +
 \frac {2 \cot \th}{\sin \th} \, \hat v_\var,
 \ssk\\
 \hat  u+ \hat v_\th + \cot \th \, \hat v +
  \frac 1{\sin \th} \, \hat w_\var=0. \qquad \qquad\qquad\quad
  \qquad\qquad\qquad\qquad\qquad
 \end{matrix}
  \right.
  \ee
 Indeed,  the system looks rather frightening for studying in
 general. As we have mentioned, fortunately, it was proved in \cite{Sv06}
 that, up to an isometry, the only non-trivial $(-1)$-homogenuity
 stationary
solutions in $\re^3 \setminus \{0\}$  are the Slezkin--Landau ones
\ef{Lan1} (it was also conjectured there that branching from
$\uu_{\rm SL}$ is impossible). Further extensions, showing the
exceptional role of the S--L solutions, are obtained in
\cite{MiuTs08}. In any case, \ef{A10n} shows the range of typical
difficulties concerning systems that inevitably occur while
searching for other types of stationary or non-stationary
singularity manifolds for linearization.

 Nevertheless,
 according to our ``linearized
strategy" of further construction of blow-up patterns (and for
answering questions \ef{Pos} or \ef{Neg}),
 \be
 \label{mmmm1}
  \fbox{$
 \mbox{all properties of all solutions of (\ref{A10n}) should be known
 in detail.}
  $}
  \ee
 Indeed, singular stationary solutions \ef{vv3},
 \ef{A10n} can be key for a possible successful matching to create
 a blow-up pattern.
The second step is then to answer:
\be
 \label{mmmm2}
  \fbox{$
 \mbox{{Questions (i)--(vi)}, \S\, \ref{S3.8}, for
 operators linearized about all solutions of
  (\ref{A10n}).}
  $}
  \ee
  To underline the complexity of the problems
  \ef{mmmm1} and \ef{mmmm2},
 we will stress below  the attention to some known
  examples of particular solutions of that type.

  \ssk

  \noi\underline{\em  On general singular stationary states: elliptic evolution}.
 Of course, the system \ef{A10n} does not contain all the necessary singular
stationary states for \ef{A1}, which might be important for
blow-up constructions. The general representation of such
solutions takes the form
 \be
  \label{vv3s}
   \tex{
  (u,v,w)= \frac 1r \, ( \hat u, \hat v, \hat w)(\var,\th,s), \quad p= \frac
  1{r^2} \, \hat p(\var, \th,s),
 }
   \ee
where $s= - \ln r \to + \iy$ as $ r \to 0^+$
plays the role of the time-variable. Then
  \be
  \label{tau2}
  \tex{
 u_r=- \frac 1{r^2}\,(\hat u_s+ \hat u), \quad u_{rr}= \frac
 1{r^3}\,(\hat u_{ss} + 3 \hat u_s +2 \hat u), \quad \mbox{etc.,}
 }
  \ee
  so that the stationary system \ef{A1}, according to the variable
  $s$, takes the form of  ``elliptic evolution" equations for $s
  \gg 1$,
   \be
   \label{ev11}
    \left\{
u_{ss}=...\,,\quad
    v_{ss}=...\,, \quad
    w_{ss}=...\,, \quad
    u_s=...\, .
     \right.
   \ee
In general regularity linear PDEs theory, such ``blow-up" scalings
lead to complicated spectral theory of pencils of linear
operators, whose spectrum and root functions classify all possible
types of singularities occurred. We refer to seminal Kondrat'ev's
papers in the 1960s \cite{Kond66, Kond67} and monographs by Maz'ya
with collaborators \cite{KMR1, KozMaz01} (further extensions via
{\tt MathSciNet}).

  As is well-known (since Hadamard's classic  example), such an elliptic-like evolution
 is ill-posed and almost all of the orbits
  are destroyed before reaching the singularity point $s=+\infty$
   ($r=0$), but anyway  the rest of the  orbits that are defined for all $s \gg
   1$ describe all possible singular steady states for \ef{A1}.
   Such an analysis has been effectively implemented for a number
   of single semilinear elliptic problems, even in the case of
   non-Lipschitz nonlinearities; see e.g., a machinery and a full list of references in
   \cite{BidVer99}. Indeed, for the stationary system of four PDEs as
   in \ef{A1}, the problem of identifying {\em all} possible
   global evolution trajectories reaching $s=+\iy$ is
    extremely difficult.
        Recall
    again that each such a non-trivial singular stationary state
    with possible inclinations and precessions of the swirl axis
    (see further developments below) can be key for construction
    of a blow-up pattern by a some kind of a linearization and
    matching procedures.
Of course, eventually we cannot escape  the problem \ef{mmmm2} for
each of those singular stationary solutions \ef{vv3s}.

In other words, the singular stationary problem \ef{ST1} plays a
first key role for understanding of the blow-up mechanism of the
NSEs, to say nothing about a similar problem with rotations to be
represented later on.

\subsection{Scaling and introducing blowing up angular singularity with precessions}

We return to the general system \ef{A1} in the standard rescaled
blow-up variables
 \be
 \label{A2}
 \tex{
 \uu=\frac 1{\sqrt{-t}}\, \hat \uu, \quad p=\frac 1{(-t)}\, \hat p, \quad y=
 \frac r{\sqrt{-t}}, \quad \t=- \ln (-t),
  }
  \ee
 where, as one of the possibilities, we first assume the similarity swirling mechanism
 of  accelerating  rotation
 in the angle
 $\var$,
 \be
 \label{A3}
 \var= \mu - \sigma \,  \ln(-t) \equiv \mu + \sigma  \, \t \quad (\sigma  \not = 0).
  \ee
 As we have seen, using the variable \ef{A3} introduces into the system an extra parameter $\sigma
 \in \re$ being a nonlinear eigenvalue that increases
the probability of successful matching of flows in  various
regions.


Substituting \ef{A2}, \ef{A3} into \ef{A1} yields the following
system:
 \be
 \label{A20}
  \left\{
  \begin{matrix}
\hat u_\t + \frac 12 \, y \hat u_y + \frac 12 \, \hat u + \sigma
\hat u_\mu + \hat u \hat u_y+ \frac 1y \, \hat v \hat u_\th +
\frac 1{y \sin \th} \, \hat w \hat u_\mu
-
 \frac 1y\, \hat v^2 - \frac 1{y} \, \hat w^2\qquad\qquad\,\,\,\,\,\ssk\\
 =-\hat p_y + \D_3 \hat u - \frac 2{y^2} \, \hat u - \frac 2{y^2} \, \hat v_\th -
 \frac {2 \cot \th}{y^2} \, \hat v - \frac {2}{y^2 \sin \th} \, \hat w_\mu,
 \qquad \qquad\qquad\qquad\qquad\,\,\,\,\,
 \ssk\ssk\\
\hat v_\t+ \frac 12 \, y \hat v_y  + \frac 12 \, \hat v + \sigma
\hat v_\mu+ \hat u \hat v_y+ \frac 1y \, \hat u \hat v + \frac
1{y} \, \hat v \hat v_\th +
 \frac 1{y \sin \th}\,\hat  w \hat v_\mu  - \frac {\cot \th}{y} \, \hat w^2\qquad\qquad
 \,\,\,
 \ssk\\
 =- \frac 1y \,\hat p_\th + \D_3 \hat v - \frac 1{y^2 \sin^2 \th} \, \hat v + \frac 2{y^2} \,
 \hat u_\th -
 \frac {2 \cot \th}{y^2 \sin \th} \, \hat w_\mu, \qquad\qquad
 \qquad\qquad\qquad\,\qquad\,
 \ssk\ssk\\
 \hat w_\t+ \frac 12\, y \hat w_y  + \frac 12 \, \hat w+ \sigma  \hat w_\mu+ \hat u \hat w_y+ \frac 1y \, \hat u \hat
 w
  + \frac 1{y} \, \hat v \hat w_\th +
 \frac {\cot \th}y\, \hat v \hat w + \frac 1{y \sin\th} \, \hat w \hat w_\mu \qquad\,\,\,
 \,\,\,
 \ssk\\
 =- \frac 1{y \sin\th}\, \hat p_\mu + \D_3 \hat w - \frac 1{y^2 \sin^2 \th} \, \hat w
 + \frac 2{y^2 \sin \th}\, \hat u_\mu
 +
 \frac {2 \cot \th}{y^2 \sin \th} \,\hat v_\mu, \qquad\,\qquad \qquad\qquad\,\,
 \ssk\ssk\\
 \hat u_y+ \frac 2y \, \hat u + \frac 1y \, \hat v_\th + \frac{\cot \th}y \, \hat v +
  \frac 1{y \sin \th} \, \hat w_\mu =0. \qquad \qquad\qquad\quad
  \qquad\qquad\qquad\,\,\qquad\,
 \end{matrix}
  \right.
  \ee
 The rescaled Laplacian is now
 \be
 \label{A101}
  \tex{
  \D_3 h= \frac 1{y^2}\, ( y^2 h_y)_y + \frac 1{y^2 \sin \th} \,
  (\sin\th \, h_\th)_\th + \frac 1{y^2 \sin^2 \th} \,
  h_{\mu\mu}.
  }
  \ee

Recall that, in view of the ban \ef{Sim},
for $\sigma =0$ non-trivial self-similar, i.e., stationary,
solutions of \ef{A20} are in fact non-existent. More generally, we
are supposed to perform a matching asymptotic expansion
construction of blow-up patterns for \ef{A20} using some already
known quasi-stationary manifolds. It should be mentioned  the
possibility of taking $\sigma =0$ for the $\var$-independence
solutions
 \be
 \label{prr1}
 \hat \uu=\uu(r,\th,\t), \quad p=p(r,\th,\t).
  \ee
  The long history of difficulties in proving global existence
  of solutions \ef{prr1} \cite{SerZaj07} suggests that a blow-up pattern
  can be reveled even in this restricted geometry.


However, it seems that a most reliable approach to blow-up
patterns should comprise:


(i) either the blow-up swirl mechanism \ef{A3} with $\sigma  \not
=0$, or


 (ii) the asymptotically slowing down mechanism such as in
\ef{sl1}, \ef{sl2}, formally corresponding to the case
 \be
 \label{prr2}
 \sigma =\sigma (\t) \to 0 \asA \t \to  +\iy  \inB (\ref{A3}).
  \ee


Both cases seem can
suit the linearized construction about the Slezkin--Landau
solutions \ef{Lan71} or others more regular,  which assumes a
difficult spectral analysis of the linear operator as in
\ef{Lan2}, where the linearization is now performed relative to
the nonlinear operators in \ef{A20}.
Note that for the asymptotically $\sigma =0$ in any similarity
annulus $\big\{\,\frac 1C \sqrt{-t} \le r \le C\sqrt{-t}\,\big\}$,
with arbitrary $C>1$, the construction remains the same and
includes general eigenfunctions of \ef{Lan2} having a
$\var$-dependence. Again, {Questions (i)--(vi)} from Section
\ref{S3.8} appear. Assume that we can answer these questions for
some particular setting of singular quasi-stationary manifold. For
instance,  then as a by-product,
 this would actually mean
that we may look for a point spectrum of $\BB^*$ with the
eigenfunction behaviour of the third $w$-component of the
linearization \ef{z6} about the singular equilibrium, say,
$\UU=\uu_{\rm SL}$,
 \be
 \label{prr3}
  Y_3(\t) \sim {\mathrm e}^{\l_k \t} \psi_\b^*(r,\th,\var) + ...
  \to 0
 \asA \t \to + \iy,
  \ee
 where $\l_k \in \s_{\rm p}(\BB^*)$, $k=|\b|$, is such that ${\rm Re} \,\l_k <
 0$ (or simply $\l_k<0$ for a real eigenvalue that also can be
 expected). In other words, the swirling blow-up behaviour then will occur
 on smaller shrinking subsets $\big\{r = o\big(\sqrt{-t}\,\big)\big\}$ as $t \to 0^-$,
 where one can expect both scenarios associated with  \ef{A3} or
 \ef{sl1}, \ef{sl2}, or with a constant, independent  of $\var$
 rotation as for solutions \ef{prr1}.

Let us mention the following aspect from Section \ref{S4.8}: can
the possible zero mean oscillation property similar to \ef{q5} (in
the original $(x,t)$-variables) affect the asymptotic behaviour of
solutions as $x \to \iy$ in such a way that the patterns will
attain a finite energy? Actually, this is related to the hard
problem of the asymptotic behaviour for the NSEs in a compliment
of a bounded smooth domain, on the boundary of which blow-up swirl
rotations with the zero mean as in \ef{q5} is prescribed.
 Possibly, this could affect the asymptotics of the solutions
 which become better localized in the $L^2$-sense.
 The
difficulty  of this question dramatically increases if a possible
axis precession (see below) is taken into account.

 Thus, it seems that, according to the above scenarios,
  there is no chance to study possible admissible solutions
 of the resulting  systems such as \ef{A20} in a reasonable and reliable
 generality and mathematical strictness. So we  will   continue to develop
 some necessary formal  arguments in an attempt to give  further hints for
 understanding such potential complicated blow-up patterns.
 As in Sections
\ref{S3.8} and \ref{SLan}, we will assume that in some
intermediate region the blow-up behaviour goes along the
quasi-stationary manifolds of singular steady solutions of
\ef{A20}, which, hopefully, have asymptotically the form \ef{vv3},
with $r \mapsto y$, i.e., do not exhibit somehow essential
dependence on the $\var$-torsion  (which may get crucial for
smaller $y$). Then, as we have seen, a complicated linearized
operator such as in \ef{z7} and \ef{Lan2} occurs. The main
difficulty is not studying its spectral properties in suitable
weighted topology, but the matching procedures with bounded orbits
for $y \approx 0$. For large $y \gg 1$, we always assume that
there exists a possibility of matching the inner region with a
properly deformed rescaled kernel of the fundamental solution of
the operator $D_\t-\D_3$, since all the evolution equations in
\ef{A20} contain the necessary counterparts. Observe that the
resulting patterns
are not supposed to be of a simple self-similar form, so they do
not exhibit any uniform homogenuity as in \ef{vv3} for steady
profiles and/or $L^\infty$-boundedness as non-stationary
similarity solutions after scaling \ef{A2}.


\noi\underline{\em An axis precession mechanism}. In order to
further increase the probability of such a matching of various
manifolds for systems like \ef{A20}, in order to suit \ef{cou1},
it is necessary and natural for such swirling/vortex flows to
introduce an additional  {\em axis precession} mechanism for these
blow-up twistors. This cannot be described explicitly or by a
system on a lower-dimensional subspace. However, there is a
standard asymptotic approach.

Thus, we first need to plug  into system another parameter in
order to change the axisymmetric (with swirl) geometry of the
solutions involved. As it has been mentioned, any strong symmetry
(or quasi-symmetry) constraint in the Navier--Stokes equations
reduces the dimension of the space, so approach them closer to the
regular $\re^2$-case.
 Therefore, thinking about creating even more complicated geometry and
  not hesitating  to perform rather weird transformations
 with such unknown  solutions, we will show a way how to perform
 a necessary perturbation of those solution structures.
 As a first simple formal illustration, let us  take into
account an extra slow motion in the $\th$-angle by setting
 \be
 \label{pr1}
  \th= \rho + \e \, \t \whereA |\e|\ll 1.
   \ee
This will mean introducing into the system \ef{A20} extra
$O(\e)$-operators, i.e.,
 \be
 \label{pr2}
 ...+ \e \hat u_\rho+... \, ,\quad ...+ \e \hat v_\rho+... \, ,\quad...+ \e \hat w_\rho+... \quad
 \mbox{into the LHSs of equations}.
  \ee
Using in all the equations the corresponding $\e$-expansions
 such as
  \be
  \label{pr3}
   \tex{
  \sin\th=\sin(\rho+\e \t)= \sin\rho + \e \t \cos \rho+...\,,
  \,\,\,
   \frac 1{\sin \th}= \frac 1{\sin \rho}- \e \t \, \frac{\cot
   \rho}{\sin \rho} +...\, ,
 }
    \ee
   etc.,
gives extra entries of $\e$ into the system. According to
 asymptotic expansion theory (see Il'in \cite{Il92}
for typical difficult methods and further references) and not
taking into account at this moment singularities introduced by
\ef{pr3} into the PDE system, this makes it possible to look for
solutions in the form of formal expansions\footnote{We discuss
this simplified version of an $\e$-expansion for convenience. In
general, depending on the kernel of the linearized operator and
other factors, the expansion can depend on  the small parameter
$\e^{1/l}$, where integer $l \ge 1$ is defined from the
solvability of the corresponding nonlinear algebraic systems on
the coefficients; see general bifurcation-branching theory, e.g.,
\cite{VainbergTr}.}
 \be
 \label{pr4}
 (\hat u, \hat v, \hat w, \hat p)= (\hat u_0, \hat v_0, \hat w_0, \hat p_0) + \e
 (\hat u_1, \hat v_1, \hat w_1, \hat p_1) +... \, .
  \ee
By a standard procedure,  being fully developed, the series
\ef{pr4} gives a unique formal representation of a certain formal
solution. The main step is the first one, for $\e=0$ that leads to
the above nonlinear systems (which were assumed can be solved).
The next terms $(\hat u_k, \hat v_k, \hat w_k, \hat p_k )$  are
obtained by iterating some linear systems, that is easier, and
existence of a solution can be checked by standard (but not always
simple) arguments.
The convergence of such series, as usual, is not required, and
is often extremely difficult to guarantee  even for lower-order
PDEs, where the rate of convergence or asymptotics are also hardly
understandable.
 Typically, such series serve as asymptotic ones, i.e., each next
 term correctly describes the behaviour of the actual solution
  as $\e \to 0$.

  Thus, \ef{pr4} is now assumed to describe {\em branching} of
  suitable solutions from the unperturbed  one  at $\e=0$ with no fast
  swirl. Here, we refer to classic bifurcation-branching theory
 for equations with  compact nonlinear integral operators,
   \cite{Berger,  Deim, KrasZ, VainbergTr}, etc. Proving the
   actual branching is a deadly difficult problem, especially,
   since we are interested in quite special solutions only.
   Therefore, we do not concentrate on branching phenomena and
   continue to reveal possible  features of such blow-up
    twistors.


 Note that
\ef{pr1} indicates a tendency of  the desired  precession with
potentially unbounded deviation of the swirl axis. To avoid such
an unpleasant (and seems non-realistic?) pattern, we
assume that the actual precession is governed by (cf.  \ef{sl1})
 \be
 \label{pr5}
 \th=\rho+ \e \kappa(\t,\e) \whereA \kappa(\t,\e) \,\,\, \mbox{is
 bounded},
  \ee
  so that \ef{pr1} contains the first term of the $\t$-expansion,
  and extra difficult asymptotic matching theory occurs.


\noi\underline{\em Periodic and quasi-periodic axis precession}.
 However, even the axis precession according to the simplified
 dependence \ef{pr1} can give insight into the actual behaviour of
 such matching blow-up patterns. In particular, as a formal illustration,
 ignoring at the moment a periodic-like
 behaviour in $\t$,  consider both
 angular dependencies of the swirl axis on the unit sphere ${\mathbb S}^2 \subset
 \re^3$,
  \be
  \label{pr6}
   \left\{
    \begin{matrix}
     \var=  \sigma  \, \tau, \\
\th=  \e \, \tau.\,
 \end{matrix}
  \right.
   \ee
Since, as in the classic representation of periodic and
quasi-periodic motion on a torus in $\re^3$, the angle behaviour
is understood {\em modulo} $2\pi$, so we  can have both
 scenarios
 on the sphere ${\mathbb S}^2$.
Namely,  assume for a moment that a matching procedure has turned
out to be successful for a given pair of eigenvalues $(\sigma
,\e)$. Then \ef{pr6} indicates that the {\em periodic
 scenario} of the axis evolution associated with \ef{ww5}
 (but not entirely; see below)
  takes place provided that
   \be
   \label{pr7}
   \tex{
    \frac \sigma  \e \in \mathbb Q\,\,\,\mbox{is rational},
    }
    \ee
and, {\em vice versa}, we have a quasi-periodic precession if
 \be
   \label{pr8}
   \tex{
    \frac \sigma  \e  \quad \mbox{is an irrational number}.
    }
    \ee
 The periodic scenario \ef{pr7} looks being   the simplest one
to create a blow-up singularity in the dynamical system \ef{ww1}
provided stabilization to a point is forbidden by \ef{Sim}
 (if  to follow
 the principle of Occam's
Rasor\footnote{W.~Ockham's {\sc lex parsimoniae}: ``entia non sunt
multiplicanda praeter necessitatem".}).
 But of course, this does not
  rule  out
 other patterns, which actually, do not look more complicated
and admit entirely similar characterization \ef{pr8}
 (though we do not exclude the
non-precision case $\e=0$). For \ef{pr5}, we can observe a more
realistic scenario of very small precession exhibiting periodic or
quasi-periodic motion of the axis.



\noi\underline{\em Briefly on precession of the vertex}. It is
quite natural that  including both swirl and axis precession
mechanisms will also require an extra ``precession of the vertex"
(as an acceptable analogy with a rotating top on a sufficiently
smooth surface suggests); otherwise such a complicated
non-symmetric vortex evolution with swirl and precession
 having  a fixed stagnation
point
 would not be possible.
 This assumes the slow-variable change of the original rescaled
coordinate system,
 \be
 \label{xx1}
  \hat \xx \mapsto \hat \xx + \aaa(\t),
   \ee
   where hopefully $\aaa(\t)$ can be a bounded function with also periodic or
   quasi-periodic behaviour as $\t \to \iy$\footnote{For reaction-diffusion equations such as \ef{FK1} or
   \ef{SEHam1}, the {\em blow-up set} is introduced
    $
    B[u_0]=\{x_0\in \ren: \,\,\, \exists \,\,\{x_n\}\to x_0 \,\,
   \mbox{and \, $\{t_n\}\to T^-$ \,such as \, $|u(x_n,t_n)| \to
   \iy$}\}.
    $
     It is then proved that $B$ is closed (e.g., is a
   point; see Friedman--McLeod \cite{FM} for a pioneering approach), and next blow-up scaling is performed relative to an
   $x_0 \in B[u_0]$ by setting $y=\frac{x-x_0}{\sqrt{T-t}}$, etc.}.
   This will involve into the systems such as \ef{A2} extra
   perturbed operators according to the change
    \be
    \label{ch1}
     \uu_\t \mapsto \uu_\t + (\aaa'\cdot \n)\uu,
      \ee
 so that these extra terms can be responsible for additional evolution blow-up phenomena
 (even in the asymptotically vanishing case $\aaa'(\t) \to 0$,
  since this also might support a quasi-periodic-like
 motion; the integrable case $\int^\iy |\aaa'(\t)|\, {\mathrm d \t}
 <\iy$ makes \ef{xx1} non-essential)
 similar to those induced by swirling \ef{sl1} and precession
 \ef{pr5}. All three mechanisms taking altogether lead to a
 possibility to attempt to construct a blow-up twistor pattern originated in
 the outer region by the Slezkin--Landau singular steady solutions \ef{Lan71} or others, which
  possibly are still unknown.


In conclusion, let us mention that we do not stress attention
to a possible more physical-mechanical interpretation of the  blow-up
patterns under speculations. Namely, we do not know and cannot
imagine how many actual twistors should be involved in such a
``hypothetical fluid configuration" to produce such a  pattern
under the fixed blow-up frame and zoom. We just recall that, as a consistent part of
our ``swirling-like philosophy", we mean that complicated blow-up patterns can be locally trapped as $t \to T^-$ by the quasi-stationary singularity manifold described by
the problems \ef{ST1}, \ef{Stat1}, and others. Such a blow-up drift along those singular manifolds will itself define the type of generalized swirling that is necessary to support the evolution.

Mechanical
interpretations
 of various solutions and patterns have  always been
very effective, difficult, and involved techniques of modern
applied and mathematical fluid dynamics, which can be also
efficient after a blow-up scaling, but possibly could fail in view
of the fact that the blow-up micro-structure of the NSEs might
have nothing to do with their classic macro-coherent structures
studied during almost two centuries. A possibly acceptable example
is as follows: the local atomic (micro-scaled) structure of a
desk, where this paper is about to be finished, has nothing to do
with its global (macro) properties.
 This is up to the
obvious fact that here the micro- and macro-mechanics are
different: the quantum and Newton's ones. A similar phenomenon
occurs for the NSEs:
the operators of the original PDEs (\ref{NS1}) and of the blow-up
rescaled ones (\ref{ww1}) live in completely different spaces,
$L^2(\re^3)$ and $L^2_{\rho^*}(\re^3)$.



\noi\underline{\em
Again on $L^\iy$-bounded or unbounded rescaled orbits $\{\hat
\uu(\t)\}$}. Finally, we add to this formal description of
possible blow-up patterns of twistor type the following necessary
observation based on our previous analysis. Namely, according to
matching concepts revealed in Sections \ref{S3.7}, \ref{S3.8}, and
\ref{SLan}, there exist two cases (sub-scenarios relative to the
above):



(i) {\sc Type I: the rescaled orbit $\{\hat \uu(\t)\}$ is bounded
in $L^\infty$ as $\t \to + \iy$}. Then this corresponds to
scenarios \ef{ww5} and \ef{cou1}, as before, and



(ii) {\sc Type II:  $\|\hat \uu(\t)\|_\iy \to \iy$
 as $\t \to + \iy$}. This can happen according to the matching
 analysis and can be driven by an exponential \ef{pol4} or power-like
 as in  \ef{pol51} divergence (or others) depending on the
 manifold that was used for matching purposes. Of course, this
 does not ruled out  the quasi-periodic $L^\iy$ behaviour
 \ef{S3.8} of bounded orbits in the necessary new rescaled
 variables.





 In both cases,
  the blow-up patterns remain exhibiting similar
properties of swirl and precession but, possibly, on different
spatio-temporal subsets.


\section{Final remarks}
  \label{SFin}








\subsection{Micro-structure of turbulence}

 As the attentive Reader has noticed, the main goal of
the present essay is not about achieving or even essentially
approaching a definite answer to the fundamental open problem:
\ef{Pos} or \ef{Neg} for the NSEs \ef{NS1}.
 Actually, this is more about approaching better understanding the
 {\bf Goal} \ef{Turb1}.
In other words, this is about presenting some ideas on a
description of
  \be
 \label{ab2}
  \fbox{$
 \mbox{existing micro-scaled fluid configurations appearing
 from smooth data.}
 $}
 \ee
 We recall that even for bounded smooth solutions, revealing such
 micro-structure of multiple zeros at regular points (Section \ref{SHerm}) led us to
 some difficult problems, though their Hermite polynomial solenoidal structure was
 partially justified.


Concerning singular points,
 in \ef{ab2}, the term ``micro-scaled" is key meaning
to look for fluid configurations that are seen via microscopic
``blow-up rescaling zoom" as in \ef{NS2}\footnote{In general, this
reflects a small part of fluid dynamic problems of fundamental
importance; e.g., cf. A.M.~Lyapunov Master's Thesis ``On Stability
of Spheroidal Equilibrium Forms of a Rotating Fluid"
(S.-Petersburg University, 27th January, 1885; Supervisor:
P.L.~Chebyshov), which was a forerunner of Lyapunov's stability
theory and a starting point  for his correspondence since 1885
with H.~Poincar\'e, who also in 1885
obtained a linearized system of Euler equations about a rotating
fluid pattern as a rigid body around the $z$-axis; see
\cite{Poin10} for an extra account.}. This means a description and
a classification of the admitted ``local turbulent
micro-structures" of the Navier--Stokes equations \ef{NS1}. The
problem of ``micro-structure" can be posed for any linear or
nonlinear evolution PDE of parabolic, hyperbolic, nonlinear
dispersion, etc., types, and gets very complicated even for simple
models (cf. \cite[\S~9]{2mSturm} as a short introduction to this
involved subject).


As a by-product, solving \ef{ab2} would also mean the negative
answer \ef{Neg} provided that the family of those configurations
would include some blow-up patterns. Not pretending at all to
giving any comprehensive insight into the problem \ef{ab2}, we
just have shown that including a standard similarity
``log-torsion" mechanism produces  a number of very complicated
singular stationary (cf. \ef{mmmm1}) or evolutionary dynamical
systems as some lower-dimensional reduction of the NSEs. This also
implies that {\em a detailed study of such reduced, but still very
complicated, dynamical systems and corresponding equilibria are
necessary and unavoidable steps that should have been passed
before even thinking about attacking the Millennium Problem}
(along the proposed lines).


\subsection{A final pessimistic expectation in general PDE theory}

 As a consequence of all the above blow-up speculations, we first
 state the simplest claim: for the NSEs,
  \be
 \label{kk00}
  \fbox{$
   \mbox{from the side of (\ref{Neg}), checking all blow-up
   configurations is impossible}.
   $}
    \ee
 More precisely, involving the positive part \ef{Pos},
it can be
 emphasized that,
 for a sufficiently wide class of complicated dynamical systems ``$M \times (N+1)$" such
  as \ef{ww1} with
  $M=4$, $N=3$ (the numbers of dependent and independent variables involved),
 with  a similar mathematics including a global energy control of solutions (not enough
  to guarantee $L^\iy$-bounds by embedding and/or interpolation), divergence of operators, scaling laws, and
 other related and necessary properties\footnote{Just in case, {\em modulo} (\ref{meth2});
 recent almost purely ``parabolic" and MP-type results in \cite{Koch07, Chen07}
 (I am sure that further, even stronger papers in this directions will follow
 soon)
for nonexistence of Type I blow-up for axi-symmetric flows
inspired some optimism,
 though, as we have tried
to show in Sections \ref{SAxi}, \ref{SNAxi}, and \ref{SLan}, this
can be just the beginning of a very long road to success.},
 \be
 \label{kk1}
  \fbox{$
\mbox{a definite answer to claims like (\ref{Neg}) and
(\ref{Pos}), in general, is
 impossible.}
 $}
  \ee

 As we have seen, for the NSEs \ef{NS1}, some ideas of blow-up focusing
 can be understood via standard asymptotic language, though
 a full justification could also take years or do not admit such at all.
 For the fourth-order bi-harmonic operator as in (\ref{NS1m}),
 being also parabolic but with no order-preserving and nonlocal properties, a
 similar proof often can be characterized as being  completely illusive.
  Moreover, here  a self-similar blow-up is then rather
 plausible as for \ef{ho1}, $m=2$, \cite{BGW1}, though will be extremely difficult to
 prove.
 In other words,
 there is a huge probability that
  \be
  \label{PP11}
  \fbox{$
 \mbox{Problem: (\ref{Pos}) or (\ref{Neg}),  can be  {\em analytically
 non-solvable}.}
  $}
  \ee
 Rephrasing the above,
  problems of such complexity from the side (\ref{Pos})  can be
  {\em ``accidentally solvable"}\footnote{Or, at least, it then could be not an exaggeration to comment
   that, posing
  this Millennium Problem, it would be also useful to specify, which millennium it is supposed
   to be.},
 when special tricks associated with these equations and operators {\sc only}
 (not {\em robust} techniques admitting perturbations, i.e., ``structurally stable") can rule out
  some essential part of the core difficulties; all {\em modulo} (\ref{meth2}).
  But this ``accidental" feature is rather unlikely: both \ef{Pos} and \ef{Neg} are too much related
  to each other, and proofs should pass through a lot of similar extremely difficult stages.
  We should accept that, in modern nonlinear higher-order PDE theory, many standard results,
  which had been perfectly solved for lower-order counterparts by
  inventing great mathematical methods by great mathematicians, do
  not admit rigorous setting in principle. In this case, looking and searching for results admitting
  rigorous proofs would be a wrong idea, also contradicting Kolmogorov's thoughts; see a
  footnote in Section \ref{S1.4}.


   To this end,
we recall that a definite, positive or negative, answer in
particular assumes a rigorous checking whether \ef{cou1} is true
or not for any $n \ge 2$. Let us also remind that, unlike
\ef{pr6}, in general, we are talking about a quasi-periodic motion
in an infinite-dimensional functional space for $\hat \uu(\t)$, so
that possible localization of that is a very difficult problem.
 Recall that such  quasiperiodic motions can be trapped in a vicinity of
 singular  ``equilibria"  \ef{Stat1} or \ef{ST1} of very complicated unknown
  structure.
In
its turn, via a matching approach, this means checking whether the
corresponding asymptotic bundles do or do not overlap for {\bf
any} $n=1,2,3,... \, $.
 This creates a restriction that cannot be checked analytically
 in general because it depends on unknown  and unpredictable
 conditions of matching of various pair of infinite-dimensional local vector bundles.
 In view of the energy control, those matching look like checking
 if, for a suitable 2D restriction, there exists a heteroclinic path connecting two saddles
 for a given dynamical system on a 2D manifold. For general DSs,
 existence/nonexistence of such a path cannot be predicted,
  and moreover, DSs with such heteroclinic
 orbits are known to be structurally unstable
 (Andronov--Pontriagin--Peixoto's  theorem, 1937--57). It then seems that
 such a matching is not possible for almost all DSs like that.
 This is true, but since, according to \ef{cou1}, at least a
 countable (or more than that)
 number of such possibilities is supposed to be checked, this
 changes the dimension of the parameter space and, eventually,  makes the
 problem to be analytically non-solvable.


If, in reality, the justification of the negative claim \ef{Neg} would
have been proved of being of a geometric ``configuration" as a
matching of the type {\em saddle--saddle} on some manifold for
establishing existence of a blow-up pattern with finite energy, it
would clearly suggest that any enhancements of functional space
and corresponding facilities for proving \ef{Pos}
 by more and more refined interpolation-embedding techniques
 or similar would be entirely hopeless.
Unless there is a parameter ``gap" prohibiting such matchings {\em
uniformly} in $n \ge 1$ in \ef{cou1} and for other types of
patterns. Hence, the existence problem takes a
principally other background and becomes the question of
blow-up scaling.

In other words, then the NSEs regularity problem should be
classified as being
 of ``pointwise sense", i.e., its global solvability is not controlled by any
 possible {\em a priori} bounds in $L^p$ or related Sobolev spaces, so it is hopeless to try to derive those.
 Nowadays, more and more nonlinear evolution PDEs (for instance, supercritical
 nonlinear Schr\"odinger-type or sufficiently multi-dimensional Burnett equations)
  penetrate into this class, which demands
 principally new mathematics of the truly twenty first century.
 One should be ready to recognize that, in this PDE class, most of typical
  mathematical problems will be analytically non-solvable. This by no means
   diminishes the role of the pure mathematics. On the contrary, this implies that the
    whole mathematical culture will be needed to built a necessary well-organized understanding
    of the problem, under the pressure that no even a hope for any definite  rigorous answers exists.

\ssk

Eventually, the author quite bravely expresses his personal
opinion, which in view of the possible feature \ef{PP11} might
have some sense:
 \be
 \label{FF11}
 \mbox{for (\ref{NS1}), global existence is more plausible than
 blow-up.}
  \ee
Actually, this is a pure probability based upon above  long
speculations and discussions. However, global existence issues,
though expected to be somehow related to blowing up ones, were not
the subject of the paper.



\ssk


\noindent{\bf Acknowledgments.} The  author   thanks C.~Budd for
 a short but stimulating discussion. The author also thanks
 A.V.~Aksenov and S.R.~Svirshchevskii for
 getting some original Russian papers by N.A.~Slezkin and others from the
 1930s-50s.




\begin{small}
\begin{appendix}
\section*{Appendix A. English Translation of N.A.~Slezkin's paper \cite{Slez34}}
 \setcounter{section}{1}
\setcounter{equation}{0}

\begin{center}
{\bf ON AN INTEGRABILITY CASE OF FULL
  DIFFEREN-}\\
  {\bf TIAL EQUATIONS OF THE MOTION OF A VISCOUS FLUID}
 \end{center}

 \ssk

\begin{center}
 N. A. Slezkin
 \end{center}

\ssk

\begin{center}
[Scientific-Research Institute of Mathematics and
Mechanics]
 \end{center}

\ssk\ssk

If the motion of viscous fluid is stationary and the axisymmetry
of the flow takes place, then, as is known, the stream function
satisfies the following differential equation in the cylindrical
coordinates:
 \be
 \label{1}
  \tex{
   - \frac 1r \, \frac{\pa \psi}{\pa z}\, \frac{\pa D \psi}{\pa r}+
 \frac 1r \, \frac{\pa \psi}{\pa r}\, \frac{\pa D \psi}{\pa z}
 + \frac 2{r^2} \, \frac{\pa \psi}{\pa z}\, D \psi= \nu  D D \psi,
 }
 \ee
 where
 $$
 \tex{
 D= \frac{\pa^2}{\pa r^2} - \frac 1r \, \frac{\pa}{\pa r} +
 \frac{\pa^2}{\pa z^2},
 }
 $$
 or, on introduction of the conical coordinate $\rho$ and $\cos
 \th=\t$, the equation \ef{1} takes the form:
  \be
  \label{2}
  \tex{
  \frac 1{\rho^2} \, \big[\frac{\pa \psi}{\pa \rho}\, \frac{\pa D \psi}{\pa \t}
  -\frac{\pa \psi}{\pa \t}\, \frac{\pa D \psi}{\pa \rho} + 2\big(
   \frac \t{1-\t^2}\, \frac{\pa \psi}{\pa \rho} + \frac 1 \rho \,
   \frac{\pa \psi}{\pa \t} \big) D \psi \big]= \nu  D D \psi,
   }
   \ee
   where
$$
 \tex{
 D= \frac{\pa^2}{\pa \rho^2} + \frac {1-\t^2}{\rho^2} \,
 \frac{\pa^2}{\pa \t^2}.
 }
 $$
Under a certain assumption on the form of the function $\psi$,
equation \ef{2} can be reduced to an ordinary differential
equation of fourth order, and the latter to a Riccati equation.
Indeed, set
 \be
 \label{3}
 \psi= \rho f(\t).
 \ee
 Then we have:
  $$
  \tex{
  D \psi= \frac{1-\t^2}\rho \, f'',
  }
  $$
  and the equation \ef{2} takes the form:
   \be
   \label{4}
   f f'''+3 f'f''= \nu[(1-\t^2) f^{IV}- 4 \t f'''].
   \ee
The left-hand side of this equation can be represented as
 $$
  \tex{
 ff'''+3 f'f''= \frac 12\, (f^2)''',
 }
 $$
while the right-hand one as
 $$
 \nu[(1-\t^2) f' + 2 \t f]''',
 $$
 and then \ef{4} is rewritten as follows:
  \be
  \label{5}
  \tex{
  \big(\frac12\, f^2 \big)'''= \nu[(1-\t^2)f'+ 2 \t f]'''.
  }
  \ee
 Integrating it three times, we obtain:
  \be
  \label{6}
  \tex{
  \frac 12\,f^2- \nu[(1-\t^2)f'+2 \t f]= C_0+C_1 \t + C_2 \t^2,
  }
  \ee
  where $C_0$, $C_1$, $C_2$ are constants of integration.

Equation \ef{6} can be rewritten as:
 \be
 \label{7}
 \tex{
 f'= \frac 1{2\nu(1-\t^2)}\, f^2- \frac {2 \t}{1-\t^2} \, f +
 \frac{C_0+C_1 \t + C_2 \t^2}{1-\t^2} \, .
 }
 \ee

 Thus, we have obtained a Riccati differential equation of the
 form that is not explicitly integrated. On substitution
  $$
  \tex{
  f= - 2 \nu(1-\t^2) \frac{d \ln y}{d \t}
  }
  $$
  this reduces to the linear differential equation of the 2nd
  order:
   \be
   \label{8}
    \tex{
 \frac{d^2 y}{d \t^2} +
 \frac{C_0+C_1 \t + C_2 \t^2}{1-\t^2} \,y=0.
}
 \ee

The solution of the latter one can be studied by using methods of
analytic theory of differential equations.

\ssk

\begin{center}
--------------
\end{center}

\begin{center}
[Uchenije Zapiski MGU, No. II, 1934.]
\end{center}

\ssk\ssk

\begin{center}
--------------
\end{center}

\begin{center}
{\bf \"UBER EINEN INTEGRIERBAREN FALL DER VOLLST\"ANDIGEN BEWE-}
 \end{center}
 \begin{center}
  {\bf GUNGSGLEICHUNGEN EINER Z\"AHEN FL\"USSIGKEIT}
 \end{center}

 \ssk

\begin{center}
 N. Slioskin
 \end{center}

\ssk

\begin{center}
[Institut f\"ur  Mathematik und Mechanik]
 \end{center}

\ssk

\begin{center}
(Z u s a m m e n f a s s u n g)
\end{center}

\ssk\ssk

Es wird der Fall untersucht, in dem die Gleichung f\"ur die
Stromfunk
 tion sich auf eine gew\"ohnliche Gleichung, und diese
letztere auf eine
 Riccatische Gleichung zur\"uckf\"uhren
l\"asst.

\ssk

\begin{center}
--------------
\end{center}

\begin{center}
[Wissenschaftliche Berichte der Moskauer Staatsuniversit\"at, H.
II, 1934.]
\end{center}

\ssk\ssk

\begin{center}
--------------
\end{center}

\end{appendix}
\end{small}

\begin{small}
\begin{appendix}
\section*{Appendix B. English Translation of N.A.~Slezkin's paper \cite{Slez54}}
 \setcounter{section}{2}
\setcounter{equation}{0}

\begin{center}
{\bf Remark on the notes of Yu. V.~Rumer, ``The
  problem of a submerged jet"$^{[1]}$}
 \end{center}
 \begin{center}
  {\bf and of L.G.~Loitianskii,
  ``Propagation of a whirling jet into an infinite spaces filled
  with the same fluid"$^{[2]}$}
 \end{center}

 \ssk

\begin{center}
 N. A.\, S l e z k i n \,(Moscow)
 \end{center}

\ssk\ssk

In both notes, it was pointed out that the problem on a laminar
submerged jet was first considered by L.D.~Landau (Mechanics of
Continuum Media$^{[3]})$ and that the solution of this problem
 represents a new case of explicit integration of equations of
 motion of a viscous fluid. The note of such kind in the second
 part is not correct.

 In the note by N. A.~Slezkin
 ``On an integrability case of full
  differential equations of the motion of a viscous fluid"$^{[4]}$
  published in Uchenije Zapiski MGU, No. II, 1934, it was shown
  that the full differential equation for the stream function of a
  steady axisymmetric motion in the spherical coordinates $r$ and
  $\t=\cos \th$
 $$
  \begin{matrix}
 \frac 1{r^2} \, \big[\frac{\pa \psi}{\pa r}\, \frac{\pa D \psi}{\pa \t}
  -\frac{\pa \psi}{\pa \t}\, \frac{\pa D \psi}{\pa r} + 2\big(
   \frac \t{1-\t^2}\, \frac{\pa \psi}{\pa r} + \frac 1 r \,
   \frac{\pa \psi}{\pa \t} \big) D \psi \big]= \nu D D \psi
\ssk\ssk
  \\
D= \frac{\pa^2}{\pa r^2} + \frac {1-\t^2}{r^2} \,
 \frac{\pa^2}{\pa \t^2}
  \end{matrix}
 $$
 under the assumption $\psi=r f(\t)$ reduces to the equation
 $$ 
 \tex{
 f'= \frac 1{2\nu(1-\t^2)}\, f^2- \frac {2 \t}{1-\t^2} \, f +
 \frac{C_0+C_1 \t + C_2 \t^2}{1-\t^2} \,
 }
 $$ 
which by the substitution\footnote{Clearly, a misprint: should be
 $(1-\t^2)$.--{\em VAG}.}
 $
  f= - 2 \nu(1-\t)^2 \frac{d \ln y}{d \t} $
  was
   reduced to the linear differential equation
 $$ 
    \tex{
 \frac{d^2 y}{d \t^2} +
 \frac{C_0+C_1 \t + C_2 \t^2}{1-\t^2} \,y=0
 }
 $$ 

Equating  the constants $C_0$, $C_1$, $C_2$ to zero, we arrive at
the solution
 $$
 \tex{
 f=-2 \nu(1-\t^2) \, \frac1{A+\t}
 }
 $$
 which was obtained by L. D. Landau within the study of the
 problem of a submerged jet. This case of integration of full
 differential equations of the motion of a viscous fluid was
 referred to in L. I. Sedov's book$^{[5]}$ (p.~104) and in
 Rosenblatt's paper$^{[6]}$.

 \ssk

 Received 3 I 1954

 \ssk

 \begin{center}
 LITERATURE
 \end{center}

 \ssk

 \noi 1. R u m e r\, Yu. B.
  {The problem of a submerged jet},
 {Prikl. Mat. Meh.,} Vol. {XVI}, No. 2, pp. 255--256, 1952.

\noi 2.
 L o \u{i} c y a n s k i i\, L. G.  {Propagation of a whirling jet into an
 infinite space filled with the same fluid}, Prikl. Mat. Meh.,
 {Vol. XVII}, No. 1, pp. 3--16, 1953.

\noi 3. L a n d a u\, L. and L i f s h i t z \,E. \, Mechanics of
continuum media, GTTI, Ed. 1, \S\, 19, 1944, Gostechizdat, Ed 2,
\S\, 23, 1953.

\noi 4.
  S l e z k i n\, N. A.  {On an integrability case of full
  differential equations of the motion of a viscous fluid},
 Uchen.
Zapiski MGU, No. II,  1934.

\noi 5. S e d o v \,L. I. Methods of similarity and dimension in
mechanics. Gostechteorizdat, Ed. 2, Moscow-Leningrad. 1951, p.
104.

\noi 6. R o s e n b l a t t\, A. Sur certaines classes de
mouvements permanents d'un liquide visqueux sym\'etrique par
rapport \'a un axe. M\'em. Soc. Roy. Sci. Li\'ege, IV S. 1--1936
(Zentralblatt f. Mech., Bd. 6, H. 4).

\end{appendix}

\begin{appendix}
\section*{Appendix C. On some simpler exact singular steady solutions: more history}
 \setcounter{section}{3}
\setcounter{equation}{0}


For further application to blow-up patterns, we will need to
revise some stationary solutions of \ef{NS1} that are {\em
singular} (unbounded) at the origin $x=0$. As above, such singular
equilibria can be used for attempting to construct  different
blow-up patterns that are not bounded in the rescaled similarity
variables as in \ef{NS2}. Constructing various exact or explicit
particular solutions of the NSEs is a classic, effective, diverse,
and, at the same time, difficult way of understanding fluid flows.
Some basic ideas go back to Oseen and von K\'arm\'an.
 There are various other approaches to constructing
 other, exact or
explicit, less and non-singular solutions of the Navier--Stokes
equations; see monographs \cite{LandM, PaiB} and surveys in
\cite{Gib99, Ohk04}, with a number of references therein.

 Singular solutions \ef{vv3} of the homogenuity
$-1$ have a long history. Actually, Slezkin--Landau's explicit
solutions \ef{Lan1}
 belong to the same class.
 As we have mentioned, the first ODE reduction was due to Slezkin
 \cite{Slez34} in 1934. For convenience of the Reader,
 in Appendix A, we present the English translation of
 Slezkin's short note \cite{Slez34}, while Appendix B contains the
 translation of his further note \cite{Slez54} (1954), so both
 rather convincingly show that Slezkin's ODE analysis includes the
 crucial explicit solution.
  Later, Landau \cite{Lan44} in 1944
 derived a particular explicit solution using a deep physical
 motivation and
 interpretation associated with submerged jet (the derivation in
 Landau--Lifshitz \cite[pp.~81-83]{LandM} well corresponds to the
 that in the original paper \cite{Lan44}).
 Some other aspects of derivation of Slezkin--Landau solutions can
 be found in Sedov \cite[pp.~99-104]{Sedov} and Batchelor
 \cite[p.~206]{Bat70}.


It seems that, before Slezkin, the first attempt to construct
exact stationary solutions was performed by Strakchowitsch (1931)
\cite{Str31} by looking for stream functions
 in the forms
 $$
 \tex{
  \psi=\chi(r) + \a z
   \andA
  \psi= \chi(r) + \frac {a z^2}2.
 }
  $$
  Based on these  Strakchowitsch's and Slezkin's \cite{Slez34}
  representations,
  Rosenblatt (1936) \cite{Ros36} constructed solutions with
  \be
  \label{RR1}
   \psi = f(r) z + f_1(r) \LongA
     \left\{
     \begin{matrix}
     \big(\frac {f'}r\big)''' + \big[ \frac 1r \, \big(\frac
     {f'}r\big)'\big]' + f \big[ \frac 1r \, \big(\frac
     {f'}r\big)'\big]' - \frac {f'}r\,  \big(\frac {f'}r\big)'=0,
     \ssk\ssk\\
\big(\frac {f_1'}r\big)''' + \big[ \frac 1r \, \big(\frac
     {f_1'}r\big)'\big]' + f \big[ \frac 1r \, \big(\frac
     {f_1'}r\big)'\big]' - \frac {f_1'}r\,  \big(\frac
     {f'}r\big)'=0.
      \end{matrix}
       \right.
       \ee
       Both equations admit straightforward integration ones.
It was shown that the equations admit non-singular solutions given
by infinite power series. Therefore, it seems that Slezkin
\cite{Slez34} was the first who looked for exact  stationary fluid
flows that are {\em singular} at the origin $x=0$.

 In 1950,  Yaceev
\cite{Yac50}
 performed a detailed construction
 (in fact, quite similar to Slezkin's one; cf. Appendix A)
  of such
 exact  solutions without torsion and
$w=0$, where
 \be
 \label{Y1}
  \tex{
  u= \frac {\hat u(\th)}r,\quad  v= \frac {\hat v(\th)}r,
  \quad w=0,\quad p= \frac {\hat p(\th)}{r^2}.
 }
  \ee
These assumptions essentially simplify the system \ef{A10n} that
now takes the form
 \be
 \label{A11}
  \left\{
  \begin{matrix}
 -\hat u^2 +  \hat v \hat u'  -
 \hat  v^2
 =2\hat p + \hat u'' +\cot \th \, \hat u'
  -  2 \hat u -  2\hat v' -
 2 \cot \th \hat v,
 \\
\hat v \hat v' =- \hat p' + \hat v'' + \cot \th\, \hat v'
 - \frac 1{\sin^2 \th} \, \hat v +  2\hat u',
 \qquad\qquad\quad\qquad\,\,\,\,
 \\
 \hat  u+ \hat v' + \cot \th \, \hat v=0. \qquad \qquad\qquad\quad
  \qquad\qquad\qquad\qquad\,\,\,\,\,
 \end{matrix}
  \right.
  \ee
The last equation gives the first velocity component
 \be
 \label{vv4nn}
 \hat u= - \hat v'-\cot \th \, \hat v,
  \ee
and substituting into the second one, on integration once, yields
the pressure,
 \be
 \label{vv5}
  \tex{
  \hat  p=- \hat v'- \cot \th \, \hat v - \frac 12 \, \hat v^2+C_0 \quad (C_0
 \in \re).
 }
  \ee

Finally, substituting into the first ODE yields the following
solutions of \ef{A11} \cite{Yac50}:
 \be
 \label{vv6}
  \tex{
 \hat u= - \hat v'-\cot \th \,\hat v, \quad
 \hat p= - 2 \hat v' + \frac{2(b \cos \th- a)}{\sin^2 \th},\quad
 \hat v= - \frac {2 \chi'(\th)}{\chi(\th)},
 }
  \ee
where $\chi(\th)$ is given by the hyper-geometric functions
 \be
 \label{vv7}
  \begin{matrix}
 \chi(\th)= \big( \cos \frac \th 2 \big)^\g \big(\sin \frac \th 2
 \big)^{1+\a+\b-\g}\qquad\qquad \ssk\ssk\\
 \times \, \big[
 c_1 F(\a,\b,\g, \cos^2 \frac \th 2) +
 c_2 F(\a+1-\g,\b+1-\g,2-\g, \cos^2 \frac \th 2)\big],\qquad\qquad
  \end{matrix}
  \ee
  where $c_{1,2} \in \re$.
 The relations between constants $a,b,c$ and $\a,\b,\g$ are as
 follows:
  \be
  \label{vv8}
   \left\{
  \begin{matrix}
 a= \g^2-(1+\a+\b)\g+ \frac 12\,(\a+\b)^2- \frac 12,\qquad\qquad\qquad\,\, \ssk\\
 b=(\a+\b-1)\g-\frac 12\, (\a+\b)+ \frac 12,\,\,
 c= \frac 12\,[(\a-\b)^2-1].
 \end{matrix}
 \right.
  \ee
The hyper-geometric function $F=F(z)$ in \ef{vv7} satisfying the
ODE
 \be
 \label{HG1}
 z(1-z)F''+[c-(a+b+1)z]F'-ab F=0,
   \ee
 is regular at
the origin and is given by the Kummer power series,
 \be
 \label{vv9}
  \tex{
 F(a,b,c,z)= 1+ \frac{ab}c \, z + \frac 1{2!} \,
 \frac{a(a+1)b(b+1)}{c(c+1)} \, z^2 +... \, ,
 }
  \ee
where $c \not =0$ and  $c \not = -l$, with $l \in \mathbb N$. This
converges uniformly in $\{|z|<1\}$ and also on the unit circle
$\{|z|=1\}$ if $a+b-c<0$. Note that the reduction to a linear ODE
related to \ef{HG1} was already available in Slezkin
\cite{Slez34}.  It follows from the expression for $\hat p$ in
\ef{vv6} that
 sufficiently regular solutions demand $a=b=0$.
 The
solution \ef{Lan1} in these coordinates yields (\ref{Lan71}) and
is a particular case of \ef{vv6}; see \cite{Lan44} and
\cite[p.~82]{LandM}.


 Similar
   exact singular solutions were obtained by
Squire (1951) \cite{Squ51} (see Pai \cite[p.~72]{PaiB} for extra
comments), where (cf. \ef{Y1})
 \be
 \label{Pa1}
 \tex{
 u= \frac {f'(\th)}{r \sin \th}, \quad v= - \frac {f(\th)}{r \sin
 \th},
 \quad w=0.
 }
 \ee
 Then, similarly to \ef{vv4nn}, we have from the second equation in
 \ef{A1} that
  \be
  \label{pa2}
   \tex{
   p_\th= - v v_\th + \frac 1r \, u_\th \LongA p=- \frac 12 \,
   v^2 + \frac 1r \, u + \frac {c_1}r
   \quad (c_1 \in \re),
    }
    \ee
   where the constant of integration is taken in the form $\frac{c_1}r$
   for further convenience.
  The first equation in \ef{A10n} then
  yields for $f=f(\xi)$, with $\xi=\cos\th$, the following equation:
   \be
   \label{pa3}
   (f')^2+f f''=2 f' +[(1-\xi^2)f'']'-2c_1.
    \ee
    Integrating leads to the first-order quadratic Bernoulli-type equation
  $f^2= 4 \xi f+2(1-\xi^2)f'-2(c_1 \xi^2+c_2 \xi+c_3)$,  so that
   ($\a$, $\b$, $b$ depend on
  $c_{1,2,3,}$)
 \be
 \label{pa4}
 \tex{
 f(\xi)= \a(1+\xi)+\b(1+\xi)+
  \frac{2(1-\xi^2)(1+\xi)^\b}{(1-\xi)^\a} \,
  \big[b- \int_1^\xi \frac{(1+\eta)^\b}{(1-\eta)^\a}\, {\mathrm
  d}\eta\big]^{-1}.
  }
  \ee
 The regular case
  $\a=\b=0$ can be interpreted as a jet issuing from a nozzle.


A further asymptotic extension of the exact solutions \ef{Y1} was
performed two years later by Rumer \cite{Rum52}, who
looked for solutions on the subspace $W_3={\rm Span} \,\{\frac 1r,
\, \frac 1{r^2},\, \frac 1{r^3}\}$,
\be
 \label{Y10}
  \tex{
  u= \frac {\hat u(\th)}r + \frac {\tilde u(\th)}{r^2},\quad  v= \frac {\hat v(\th)}r
  + \frac {\tilde v(\th)}{r^2},
  \quad w=0,\quad p= \frac {\hat p(\th)}{r^2}+ \frac {\tilde p(\th)}{r^3}.
 }
  \ee
It turns out  that a certain extension of  solutions \ef{vv6},
\ef{vv7} exists in the case $a=b=c=0$ in \ef{vv8}, but these are
not exact solutions so that \ef{Y10} gives the next asymptotic
term of expansion in terms of $\frac 1r$.
Further, Lo\u{i}cyanskii \cite{Loi53} in 1953 proposed to look for
solutions with torsion in the cylindrical coordinates
$\{r,\var,z\}$ in terms of power series for $u$, $w \not =0$
small,
 \be
 \label{vv51}
  \tex{
  (u,v,w): \quad \sum
  _{n=1}^\iy \, \frac { a_n(\eta)}{z^n} \whereA
  \eta= \frac r z,
   }
   \ee
   which are not converging and are singular at the origin.
The formula \ef{vv51} may supply us with some other singular
solutions, that, though  not admitting  explicit representations,
can be used in blow-up analysis.
 In general,
    it
   seems that the method of formal asymptotic expansions is an
   appropriate way to treat our more general and difficult system \ef{A10n} and
   further related ones to be introduced.
 We refer to Kurdyumov \cite{Kurd05} for more recent
 developments associated with such exact solutions of \ef{NS1}.


\end{appendix}

\end{small}

\end{document}